\documentclass[11pt]{article}

\usepackage[margin=1.3in]{geometry}

\usepackage{graphicx}
\usepackage{amsmath}
\usepackage{amssymb}
\usepackage{amsthm}
\usepackage{textcomp}
\usepackage{esint}
\usepackage{xcolor}
\usepackage{color} 
\usepackage{mathrsfs}
\usepackage{tikz}
\usetikzlibrary{arrows.meta,
	positioning,
	quotes
}
\usepackage{subfigure}
\usepackage[colorlinks=true]{hyperref}
\hypersetup{
	colorlinks=blue,%
	citecolor=red,%
	filecolor=black,%
	linkcolor=blue,%
	urlcolor=blue
}
\usepackage[colorinlistoftodos]{todonotes}
\newtheorem{remark}{Remark}[section]
\newtheorem{theorem}{Theorem}[section]
\newtheorem{lemma}[theorem]{Lemma}
\newtheorem{proposition}[theorem]{Proposition}

\newtheorem*{claim}{Claim}
\newtheorem*{convention}{\textbf{Convention}}
\newtheorem{corollary}{Corollary}[section]
\newtheorem{definition}{Definition}[section]
\def \no{\nonumber}
\def\ve{\varepsilon}
 
\newcommand{\R}{\mathbb{R}}

\newcommand{\F}{\mathcal{F}}
\newcommand{\Rn}{\mathbb{R}^n}
\newcommand{\ud}{\mathrm{d}}
\newcommand{\bcw}{\mathbin{\bigcirc\mkern-15mu\wedge}}

\newcommand{\Sn}{\mathbb{S}^n}
\newcommand{\Sp}{\mathbb{S}}
\newcommand{\N}{\mathbb{N}}

\newcommand{\pa}{\partial}

\numberwithin{equation}{section}

\makeatletter

\newdimen\bibspace
\renewenvironment{thebibliography}[1]{%
	\section*{\refname 
		\@mkboth{\MakeUppercase\refname}{\MakeUppercase\refname}}%
	\list{\@biblabel{\@arabic\c@enumiv}}%
	{\settowidth\labelwidth{\@biblabel{#1}}%
		\leftmargin\labelwidth
		\advance\leftmargin\labelsep
		\itemsep\bibspace
		\parsep\z@skip     %
		\@openbib@code
		\usecounter{enumiv}%
		\let\p@enumiv\@empty
		\renewcommand\theenumiv{\@arabic\c@enumiv}}%
	\sloppy\clubpenalty4000\widowpenalty4000%
	\sfcode`\.\@m}
{\def\@noitemerr
	{\@latex@warning{Empty `thebibliography' environment}}%
	\endlist}

\makeatother

\makeatletter

\allowdisplaybreaks

\begin{document}
	\title{\textbf{\Large{The $\sigma_{2}$-curvature equation on a compact manifold with boundary}}}
	\author{Xuezhang Chen\thanks{X. Chen:xuezhangchen@nju.edu.cn, is partially supported by NSFC (No.12271244). }~
		and Wei Wei\thanks{W. Wei:wei\_wei@nju.edu.cn, is partially supported by NSFC (No.12201288, No.12271244) and BK20220755.}\\
		{\small $^{\ast}$$^{\dag}$Department of Mathematics \& IMS, Nanjing University, Nanjing 210093, P. R. China }}
        \date{}
	\maketitle
	
	\begin{abstract}
		\medskip
		We first establish local $C^2$ estimates of solutions to the $\sigma_2$-curvature equation with nonlinear Neumann boundary condition. Then,  under assumption that the  mean curvature of a background metric is nonnegative on totally non-umbilic boundary, for dimensions three and four there exists a conformal metric  having a prescribed positive $\sigma_2$-curvature and a prescribed nonnegative  boundary mean curvature.  The local estimates play an important role in the blow up analysis for the latter existence result.
				\medskip
		
		\medskip{}
		
	\end{abstract}
	{\footnotesize \tableofcontents}
	\section{Introduction}
	
	The $\sigma_k$-Yamabe problem on closed manifolds has been studied extensively.  See \cite{CGY2, Chen3, Ge-Wang1,Ge-Lin-Wang,Ge-Wang2,Guan-Wang0,Guan-Wang,GVW, Guan-Wang2,GV0,Gursky-Viaclovsky,Li-Li1,Li-Luc2,Li-Luc1,Li-Luc3,Li-Luc4,STW,TW,Via1,Via2,Wang} etc. This problem is not only an interesting and difficult existence
	problem in geometric PDEs, especially in the field of fully nonlinear PDEs, but also provides many applications in differential geometry. One among its important applications is  four dimensional Conformal Sphere Theorem in the seminal work of  Chang, Gursky and Yang \cite{CGY1} with use of $\sigma_2$-curvature. For compact manifolds with boundary, a classical problem in conformal geometry is the boundary Yamabe problem, which traced back to Escobar \cite{escobar1,escobar4} in 1992. A more general version than Escobar's was first proposed in 1999 and investigated by Z. C. Han and Y. Y. Li \cite{han-li2,han-li1}. For recent developments, readers are referred to Chen, Ruan and Sun \cite{Chen-Ruan-Sun}, Chen and Sun \cite{Chen-Sun} and references therein.
	
	
	Similar to the aforementioned two closely related problems, the prescribing $\sigma_k$-curvature problem with boundary, including boundary $\sigma_k$-Yamabe problem, is also one of fundamental and active subjects in conformal geometry. 
	
	Let $(M,g)$ be a smooth compact Riemannian manifold of dimension
	$n\geq3$ with boundary $\pa M$. Denote    the Ricci curvature by $\mathrm{Ric}_g$ and scalar curvature by $R_g$. For any symmetric two-tensor $W$ and $W_j^i=g^{ik}W_{kj}$, we define
	\begin{equation}\label{def:sigma2}
		\sigma_2(W)=\frac{1}{2}
		\delta\left(\begin{matrix}
			i_1&i_2       \\
			j_1 &j_2
		\end{matrix}
		\right)
		W_{i_1}^{j_1}W_{i_2}^{j_2} ,
	\end{equation}
	where the Kronecker symbol 
	$\delta \left(\begin{matrix}
		i_1&i_2       \\
		j_1 &j_2
	\end{matrix}
	\right)$
	has value  $1~ (-1)$ for an even (odd) permutation in the index set $\{1,2,\cdots,n\}$ when $i_1,i_2$ are distinct; otherwise it has value $0$. Moreover, it is straightforward to show
	\begin{equation}\label{def:F_j^i}
		(T_1)_j^i:=\frac{\pa \sigma_2(W)}{\pa W_i^j}=\sigma_1(W)\delta_j^i-W_j^i.
	\end{equation}

	Define
	$$\Gamma_1^+:=\{\lambda=(\lambda_1,\cdots,\lambda_n)\in \Rn; \sum_{i} \lambda_i>0\}$$
	and
	$$\Gamma_2^+:=\{\lambda\in \Gamma_1^+; \sum_{i<j} \lambda_i\lambda_j >0\}.$$
	
	\begin{convention} 
		The Latin letters like $1\leq i,j,k \leq n$ stand for the full indices; the Greek letters like $1 \leq \alpha,\beta,\gamma\leq n-1$ stand for tangential indices. 
		
	\end{convention}

	 We define the second fundamental form at some $P \in \pa M$ by 
	\[
	L(X,Y)=-\langle\nabla_{X}\vec{n},Y\rangle
	\]
	and its trace-free part by 
	\[
	\mathring{L}(X,Y)=L(X,Y)-h_{g}\langle X,Y\rangle
	\]
	for all $X,Y \in T_P(\pa M)$, where $\vec{n}$ is the inward unit normal at $P$, $h_{g}=\sum_{\alpha}L(e_\alpha,e_\alpha)/(n-1)$   is the mean curvature with $\{e_1,\cdots,e_{n-1}\}$  being an orthonormal basis of $T_P(\pa M)$.
	
	For a conformal metric $g_{u}:=e^{-2u}g$,
	a conformal deformation of the Schouten tensor 
	\[
	A_{g}=\frac{1}{n-2}(\mathrm{Ric}_{g}-\frac{R_{g}}{2(n-1)}g)
	\]
	gives 
	\begin{equation}\label{conf_change_Schouten}
		A_{u}:=A_{g_{u}}=\nabla^{2}u+\ud u\otimes\ud u-\frac{1}{2}|\nabla u|^{2}g+A_{g}.
	\end{equation}

	Under the above conformal change of metrics there hold
	\begin{equation}
		\mathring{L}_{g_{u}}=e^{-u}\mathring{L}_{g}\label{trace_free-2nd_ff}
	\end{equation}
	and 
	\[
	h_{g_{u}}=(\frac{\pa u}{\pa\vec{n}}+h_{g})e^{u}.
	\]

	We are interested in the following $\sigma_2$-curvature equation with nonlinear Neumann boundary condition:
	\begin{align}\label{PDE:sigma_2_with_bdry}
		\begin{cases}
			\sigma_{2}^{1/2}(A_{u})=f(x,u), \quad \lambda(A_u)\in \Gamma_2^+ & \qquad\mathrm{in~~}M,\\
			\frac{\pa u}{\pa\vec{n}}=ce^{-u}-h_{g} & \qquad\mathrm{on~~}\partial M,
		\end{cases}
	\end{align} 	
	where  $c=h_{g_{u}}$ and $\lambda(A_u)$ is the eigenvalues of $g^{-1}A_u$. Throughout the paper, we don't distinguish $\sigma_{2}(A_{u})$ and $\sigma_2(\lambda(A_{u}))$, the second elementary symmetric function of $\lambda(A_u)$. For convenience, we call $\sigma_2^{1/2}(A_g)$ as the $\sigma_2$-curvature of metric $g$.
	
	In the particular case of
	\[
	f(x,u)=f(x)e^{-2u(x)}, 
	\]
	the above problem \eqref{PDE:sigma_2_with_bdry}
	is exactly the prescribing $\sigma_{2}$-curvature and boundary mean curvature problem:
	\begin{align}
		\begin{cases}
			\sigma_{2}^{1/2}(A_{u})=f e^{-2u}, \quad \lambda(A_u)\in \Gamma_2^+& \qquad\mathrm{in~~}M,\\
			\frac{\pa u}{\pa\vec{n}}=ce^{-u}-h_{g} & \qquad\mathrm{on~~}\partial M.
		\end{cases}\label{PDE:prescribing_sigma_2}
	\end{align}
As first pointed out by S. Chen in  \cite[Lemma 5]{Chen}, under the assumptions that $0<f \in C^\infty(\overline M)$ and $c=0$, on compact manifolds with umbilic boundary the above equation \eqref{PDE:prescribing_sigma_2} serves as the Euler-Lagrange equation of a geometric functional $\mathcal{F}_2[g]$ (see \eqref{def:F_2} for the definition). We are going back to the variational characterization in Section \ref{Sect:degree_theory}.

	In comparison of the Yamabe equation, the $\sigma_2$-curvature equation is a fully nonlinear elliptic  equation, to solve which, the local  $C^1$ and $C^2$ estimates are necessary. 
	A more general fully nonlinear PDE with Neumann boundary condition including \eqref{PDE:sigma_2_with_bdry} had been studied by Q. Jin, A. Li and Y. Y. Li \cite{Jin-Li-Li}, where the authors \cite[Theorem 1.4]{Jin-Li-Li} established local $C^{1}$ estimates assuming a local $C^0$  lower bound. 
	
	Once the local $C^2$ estimates were achieved,  the compactness and existence of solutions to (\ref{PDE:prescribing_sigma_2}) come out with proper $C^0$ estimates. There exist several existence theorems  for \eqref{PDE:prescribing_sigma_2} on manifolds with umbilic boundary. When $(M,g)$  is  locally conformally flat  with umbilic boundary having multiply connected components,  Li and Nguyen \cite{Li-Luc3}  established the compactness of solutions to $\sigma_k$-curvature equation   for nonnegative prescribed mean curvature.  He and Sheng \cite{HS} obtained the compactness of solutions to $\sigma_k$-curvature equation on manifolds with totally geodesic boundary for $c=0$.  On a smooth compact four-manifold with umbilic boundary, which is not conformally equivalent to the standard hemi-sphere $(\Sp^4_+,g_{\Sp^4})$, S. Chen \cite[Theorem  2]{Chen} proved the existence of solutions to \eqref{PDE:prescribing_sigma_2} with $c=0$   provided that the boundary Yamabe constant $Y(M,\pa M,[g])$ and $\mathcal{F}_2[g]$ are both positive. Here
	\begin{equation}\label{def:bdry_Yamabe_constant}
		Y(M,\partial M,[g]):=\inf_{\tilde{g}\in[g]}\frac{\int_{M}R_{\tilde{g}}\ud\mu_{\tilde{g}}+2(n-1)\int_{\partial M}h_{\tilde{g}}\ud\sigma_{\tilde{g}}}{(\int_{M}\ud\mu_{\tilde{g}})^{\frac{n-2}{n}}}.
	\end{equation}
	With the local estimates in Theorem \ref{Thm:local estimates} and the argument of S. Chen in \cite[Theorem  2]{Chen}, we extend to prove the existence result of \eqref{PDE:prescribing_sigma_2} on manifolds with umbilic boundary when  $c\geq 0$.   However, the existence of  \eqref{PDE:prescribing_sigma_2} on manifolds with non-umbilic boundary is largely open. 
	
	Besides local $C^1$ estimates, boundary local $C^2$ estimates  are difficult obstructions to the existence for non-umbilic boundary. Our first purpose is to derive local $C^2$ estimates  of solutions to \eqref{PDE:sigma_2_with_bdry} on compact manifolds with \emph{generic} boundary.
	
	\begin{theorem}\label{Thm:local estimates}
		Let $(M,g)$ be a smooth compact Riemannian manifold of dimension
		$n\geq3$ with boundary $\pa M$ and assume $0\leq c \in C^\infty(\pa M)$. Suppose $f(x,z): M \times \R \to \R_+$ is a smooth function. Let  $\mathcal{O}$ and $\mathcal{O}'$ be two neighborhoods near $\pa M$ with $\mathcal{O}' \Subset \mathcal{O}$. Assume $u \in C^4(\overline M)$ is a solution to \eqref{PDE:sigma_2_with_bdry}. Then there exists a positive constant $C$ depending on $n, g, \mathrm{dist}_g(\mathcal{O}', \partial\mathcal{O}\backslash \pa M)$, $\sup_{\mathcal O}(c+|\tilde \nabla c|+|\tilde \nabla^2c|)e^{-u}$,  $\sup_{\mathcal{O}} (f+|f_x|+|f_{xz}|+|f_{xx}|+|f_{zz}|)$ and $\sup_{\mathcal{O}}\frac{1}{f(x,u)}$, such that
		$$|\nabla^2 u| \leq C \qquad \mathrm{~~in~~} \quad \mathcal{O}'.$$
		Here $\tilde{\nabla}$ denotes the
	Levi-Civita connection with respect to the induced metric  on $\pa M$.
	\end{theorem}
	The boundary local  $C^2$ estimates for \eqref{PDE:sigma_2_with_bdry} are even \emph{new} for both Euclidean domains and compact manifolds with umbilic boundary for $c>0$. Jin, Li and Li \cite{Jin-Li-Li} obtained local  $C^2$ estimates when $\partial M$ is umbilic and $(M,g)$ is locally conformally flat near $\pa M$. Meanwhile, with different technique, S. Chen \cite{Chen3} obtained the boundary local  $C^2$ estimates under the same assumption except for $c=0$. Jin \cite{Jin} obtained the boundary local  $C^2$ estimates for manifolds with totally geodesic boundary.   Jiang and Trudinger \cite{JT} introduced the so-called $(\Gamma_k, A,G)$-convex boundary condition, which is equivalent to the condition $c e^{-u}>0$ in $\Gamma_2^+$  (cf. \cite[p.401]{JT}), and obtained the global $C^2$ estimates for a Euclidean domain $\Omega$ with non-umbilic boundary $\pa \Omega$.  Specifically, $\inf_{\partial \Omega} (c e^{-u})>0$ plays a crucial role in their proof.

	 The nonnegativity of  prescribing mean curvature for umbilic boundary is \emph{necessary} to boundary local $C^2$ estimates, which can be  seen from the counter-examples constructed by Li and Nguyen \cite{Li-Luc4}.   It is shown  in \cite{Li-Luc4} that on an annulus $B_{R_0}\backslash \overline{B_1}(0)$,  a sequence of radially symmetric solutions $\{u_j\}$ to positive constant $\sigma_k$-curvature  ($2 \leq k \leq n$) equation  enjoys the property:  Given any  constant  mean curvature $c<0$ on $\pa B_1$, $\pa_r^2 u_j$   goes to infinity  when approaching $\pa B_1$. For non-umbilic boundary, this necessity is ensured by our counter-examples. One of our important ingredients is to replace the inner boundary by an appropriate ellipsoid touching  $\pa B_1(0)$ at exactly two non-umbilic points. See Proposition \ref{prop:counterexample_new}.  It is no doubt that the sign of $c$ is essential in the estimates of double normal derivative on the boundary.   	 

For non-umbilic boundary case, one of missing elements in the existence for \eqref{PDE:sigma_2_with_bdry} is exactly the local  $C^2$ estimates in Theorem \ref{Thm:local estimates}, although  the Neumann boundary problem for fully nonlinear PDEs has been developed a lot, see  Lions, Trudinger and Urbas \cite{LTU} and recent work of X. N. Ma and G. Qiu \cite{MQ} etc. Generally, the boundary curvatures play an essential role on the boundary $C^2$ estimates for the most fully nonlinear equations. For example, in \cite{CMW, CW, LTU, MQ} etc., the boundary $C^2$ estimates only work for strictly convex domains. Without the convexity of the boundary, Dong and Wei \cite{DW} proved  $C^2$ estimates for certain type of fully nonlinear equations, which highly depend on the strictly elliptic property of these equations. For most of fully nonlinear elliptic equations we would not expect  $C^2$ estimates without any condition on equations or curvatures on manifolds.  The classical existing method consists of two steps: The first step is to reduce second  tangential derivative estimates to double normal derivative estimates; the second step is to obtain double normal derivative estimates.  However, this reduction argument does not work well to obtain the local $C^2$ estimates for $\sigma_2$-curvature equation on manifolds with non-umbilic boundary. Our strategy to derive the boundary local $C^2$ estimates is quite different from these works just described.
	
	For boundary local $C^2$ estimates, we split our proof into two steps.  The first step is to derive the trace of second tangential derivative estimates. We first realize that an inaccuracy occurs in the computations of second covariant derivatives of the barrier function by S. Chen \cite{Chen}.  This brings us some unavoidable `\emph{hard}' terms involving third covariant derivate terms. Our method combines some overlooked ingredients from the study of $\sigma_2$-Yamabe problem: As the traditional way for $\sigma_2$-Yamabe equation, only the concavity of $\sigma_2^{1/2}(W)$ with $W=A_u$ is not enough to ours, instead we use a complete expansion of the term $\frac{\pa^2 \sigma_2(W)}{\pa W_i^j \pa W_r^s}$, which enables us to take back  some `\emph{positive}' terms (originally thrown away) to get a good control of these `\emph{hard}' terms. Readers are referred to \eqref{key-est:3rd_derivatives} for details. In fact, the $\sigma_2$-curvature equation has interior $C^2$ estimates, which are rare in the fully nonlinear elliptic equations and make second tangential  derivative estimates reasonable. In the second step, we construct a delicate barrier function to obtain double normal derivative estimates and then obtain the trace of second derivative estimates.  To derive the local $C^2$ estimates, we take advantage of a beautiful structure hidden behind two geometric PDEs: The interior $\sigma_{2}$-curvature equation and boundary mean curvature equation. We emphasize that the dependence of the above constant $C$ on $\sup_{\mathcal{O}}\frac{1}{f(x,u)}$ only originates in double normal derivative estimates. This is a different phenomenon from the local $C^1$ estimates in \cite{Jin-Li-Li}. Moreover, the nonnegativity of prescribing boundary mean curvature is crucial in the proof.  This in turn can be confirmed by Li-Nguyen's and our counterexamples.

	Based on the local $C^1, C^2$ estimates and the blow-up analysis, we are able to obtain $C^0$ estimates, as well as the following existence result.
	\begin{theorem}\label{Thm:existence}
		Let $(M,g)$ be a compact Riemannian manifold of dimension $n\geq 3$ with totally non-umbilic boundary $\pa M$. Assume that the Schouten tensor satisfies $\lambda(A_g) \in \Gamma_2^+$ and the mean curvature $h_g\geq 0$  on $\pa M$. Suppose $n=3,4$ and $0<f \in C^\infty(\overline M), 0\leq c \in C^\infty(\pa M)$. Then the boundary $\sigma_2$-curvature equation \eqref{PDE:prescribing_sigma_2} admits a smooth solution.
	\end{theorem}
	
	Here `\emph{totally non-umbilic boundary}' means every boundary point is non-umbilic, which is essentially used to exclude the occurrence of boundary isolated blow-up points as in Definition \ref{def:bdry_blowup_pts}. Guan, Viaclovsky and Wang \cite{GVW} established that  for $k\ge n/2$, $\lambda(A_g)\in \Gamma_k^+$ implies the nonnegativity of  Ricci curvature, which is crucial in the study of $\sigma_k$-Yamabe problem on closed manifolds; see \cite{Gursky-Viaclovsky,Li-Luc2,TW} etc. For this reason, we limit ourselves to dimensions three and four.  Moreover, the assumption $\lambda(A_g) \in \Gamma_2^+$ and $h_g\geq 0$  directly implies $Y(M,\partial M,[g])>0$.

	Schoen \cite{Schoen} initiated the blow up analysis for the scalar curvature equation on closed manifolds. Such  blow-up theory of boundary Yamabe problem has been delicately explored  in Han and Li \cite{han-li2}.  
	For $n\geq 3$, some basic results for the boundary $\sigma_2$-curvature equation \eqref{PDE:prescribing_sigma_2} are  developed for  the analogue of Schoen and Zhang's \cite[Lemma 3.1]{Schoen-Zhang}.  However, our strategy is technically quite different from that of Schoen  on prescribing scalar curvature  problem. On the basis of and Liouville-type theorems by Y. Y. Li and collaborators \cite{Li-Li1,Li-Li2,Han-Li-Tei} etc,  for $n=3,4$  we are able to complete the blow up analysis together with our local estimates in Theorem \ref{Thm:local estimates}, and thus the existence result follows. Another advantage of $n=3,4$ is that for any sequence of solutions $u_i$ to \eqref{PDE:prescribing_sigma_2} with $\min_{\overline M} u_i \to -\infty$ (such a $\{u_i\}$ is called a blow-up sequence), we define a rescaled sequence of  functions $v_ i=u_i-u_i(y_0)$ for some interior `\emph{regular}' point $y_0$ (cf.  Section \ref{Subsect:compactness}). The sequence of conformal metrics $g_{v_i}:=v_i^{4/(n-2)}g$ has nonnegative Ricci curvature by virtue of Guan, Viaclovsky and Wang \cite{GVW}. This plays a central role in excluding interior isolated blow-up points as in Definitions \ref{def:isolated blow up} and \ref{def:interior_isolated_simple} for the above blow-up sequence $\{u_i\}$.

	The following is the organization  of the paper, as well as a sketch of our strategy. 
	
	Section \ref{Sect:local_est} is devoted to the a priori $C^2$ estimates for  \eqref{PDE:sigma_2_with_bdry}, which is accomplished by the combination of estimates of trace of  second tangential derivatives and double normal derivative on $\pa M$. In Subsection \ref{subsect:counter-examples}, counter-examples to local boundary $C^2$ estimates on \emph{generic} boundary are presented. In Section \ref{Sect:first_bdry_eigenvalue}, we study first boundary $\sigma_2$-eigenvalue problem, which parallels  first boundary Yamabe/$\sigma_1$-curvature eigenvalue problem. As a byproduct, this reduces the assumption $h_g\geq 0$ to $h_g=0$ in the study of existence result of \eqref{PDE:prescribing_sigma_2} in the next two sections. In Section \ref{Sect:degree_theory}, we adopt the degree theory to prove the existence of solutions to \eqref{PDE:prescribing_sigma_2} and encounter great difficulty obtaining the uniform $C^0$ estimates  to a $1$-parameter family of boundary $\sigma_2$-curvature equations. The uniform $C^0$ bounds are obtained via the blow-up analysis in Section \ref{Sect:blowup_analysis}. 
	
	 In Section \ref{Sect:blowup_analysis}, we build the blow-up analysis for a slightly different boundary $\sigma_2$-curvature equation \eqref{PDE:sigma_2_new}, which enjoys a similar structure of each equation in the $1$-parameter family \eqref{eq:path}.  In Subsection \ref{Subsect:Basic facts in all dimensions}, we introduce relevant definitions, such as (interior and boundary)  isolated and isolated simple blow-up points etc, and establish some common facts in the blow-up procedure for all $n\geq 3$.  Starting from Subsection \ref{Subsect:isolated_singularity}, we restrict consideration to $n=3,4$. Given any blow-up sequence $\{u_i\}$ for \eqref{PDE:sigma_2_new} (if exists), we describe the isolated singularity behavior near each interior and boundary isolated blow-up point. We next  establish that every isolated blow-up point is indeed simple and the singular set $\mathcal S:=\mathcal{S}[\{u_i\}]$, consisting of only isolated blow-up points, is a finite set. Keep in mind that isolated blow-up points do not occur on the boundary due to the assumption that the boundary is totally non-umbilic. We are able to adapt some insights from Gursky and Viaclovsky's argument in \cite{Gursky-Viaclovsky}  to show that $\mathcal{S}$ has exactly one interior isolated blow-up point.  Using the volume estimate by Perales \cite{Perales} for manifolds with boundary under nonnegative Ricci curvature and a lower  mean curvature bound, with some extra efforts we  demonstrate that the unique interior isolated blow-up point could not happen either. In summary, such a blow-up sequence actually does not exist.
	
	Appendix \ref{Appendix:A} is of elementary nature in differential geometry, and we give proofs of Propositions \ref{prop:ellipsoid} and \ref{prop:umbilic_pts_ellipsoid} involving description of \emph{possible} umiblic points for higher dimensional ellipsoids. To make our proof more transparent, we postpone the proof of  Lemmas \ref{one critical point} and \ref{lem:four-dim lower bound-boundary} to Appendix \ref{Appendix:B}.
	For readers' convenience,  the proof of two results involving degenerate $\sigma_2$-curvature equations with Neumann data is left to Appendix \ref{Appendix:C}: One concerns the limit of a  blow-up sequence of solutions to the above family of boundary $\sigma_2$-curvature equations; the other is about an even reflection of a viscosity solution to a degenerate $\sigma_2$-curvature equation in a Euclidean half-ball with zero Neumann boundary condition.
	
	\bigskip
	
	\noindent{\bf Acknowledgments.} The second named author would like to thank Professor Hao Fang and Professor Xinan Ma for enlightening discussions and constant support. The authors are grateful to Professors Yan Yan Li and Xinan Ma for bringing their attention to Jiang-Trudinger's paper \cite{JT}. They greatly appreciate Yan Yan Li's encouragement to construct the counter-examples on generic boundary for local $C^2$ estimates. The authors thank their colleague at Nanjing University, Professor Yalong Shi for fruitful discussions about possible umbilic points of higher dimensional ellipsoids. The authors are grateful to Professor Neil S. Trudinger's interest in our work and encouragement, and thank Professor Feida Jiang for private communications about their results in \cite{JT} relative to those in this paper.
	\bigskip
	
	\section{Local estimates}\label{Sect:local_est}
	In this section, we establish the a priori $C^2$ estimates for solutions to \eqref{PDE:sigma_2_with_bdry} and thus complete the proof of Theorem \ref{Thm:local estimates}.

	Denote by
	$\nabla$ the covariant derivative with respect to  the Levi-Civita connection of $g$ in $M$  with Christoffel symbols $\Gamma_{ij}^{k}$, and by $\tilde{\nabla}$ the covariant derivative with respect to the
	Levi-Civita connection of the induced metric of $g$ on $\pa M$, respectively. From now on, given $u \in C^\infty(\overline{M})$, we simplify $\nabla_{j}\nabla_{i}u$
	as $u_{ij}$ or $u_{,ij}$ and $\tilde{\nabla}_{\beta}\tilde{\nabla}_{\alpha}u$
	as $u_{;\alpha\beta}$, so are higher order covariant derivatives
	of $u$; for tensors with $A_{ij}$ and $L_{\alpha\beta}$ for example,
	we shall use $A_{ij,k}=\nabla_{k}A_{ij},L_{\alpha\beta;\gamma}=\tilde{\nabla}_{\gamma}L_{\alpha\beta}$.
	Denote by $\pa_{i}u$ the first order Euclidean derivatives of $u$,
	so are higher order Euclidean derivatives of $u$.
	
	Let us fix $x_{0}\in\pa M$, choose geodesic normal coordinates
	$z'=(z_{1},\cdots,z_{n-1})$ on $\pa M$ around $x_{0}$, corresponding to the
	point $x'\in\pa M$ near $x_0$. The local coordinates $z=(z',z_{n})$ with small
	$z_{n}\geq0$ are called Fermi coordinates if $x=\exp_{x'}\{z_{n}\vec{n}(x')\}$
	in a neighborhood of $x_{0}$. For simplicity, we denote such a map by $G_{x_{0}}(z)=\exp_{x'}\{z_{n}\vec{n}(x')\}:B_{\rho}^{+}:=B_{\rho}(0)\cap\Rn_{+}\to\overline{M}$. Under Fermi
	coordinates around $x_{0}$, the metric has the form  
	\[
	g=\ud z_{n}^{2}+g_{\alpha\beta}\ud z_{\alpha}\ud z_{\beta}\qquad\mathrm{~~for~~}z=(z',z_{n})\in B_{\rho}^{+}.
	\]
	We sometimes implicitly  identify the
	coordinates $x$ with $z$ via use $G_{x_{0}}$. Readers are referred to \cite[Section 3]{escobar1}
	for further properties of Fermi coordinates. The second fundamental form on $\pa M$ and its trace-free part are given by
	\[
	L_{\alpha\beta}=-\langle\nabla_{\pa_{x_{\alpha}}}\vec{n},\pa_{x_{\beta}}\rangle \qquad \mathrm{~~and~~} \qquad  \mathring{L}_{\alpha\beta}=L_{\alpha\beta}-h_{g}g_{\alpha \beta},
	\]
	respectively, and we simplify $H:=H_g=(n-1)h_{g}$.
	The Codazzi equation states that
	\[
	R_{\alpha\beta\gamma n}=L_{\alpha\gamma;\beta}-L_{\beta\gamma;\alpha} \qquad \mathrm{~~on~~} \partial M.
	\]
	
	To simplify notations, we denote $W=A_u$ as in \eqref{conf_change_Schouten} and by $W^\top$ the tangential part of $W$, and define
	$$F_j^i:=\frac{\pa \sigma_2(W)}{\pa W_i^j},\quad F^{ij}=g^{ik}F_{k}^j \quad \mathrm{~~and~~} \quad F_{js}^{ir}:=\frac{\pa^2 \sigma_2(W)}{\pa W_i^j \pa W_r^s}.$$
It sounds more natural and convenient to adopt these notations above in the study of fully nonlinear partial differential equations on manifolds. Throughout this section, we adopt the Einstein summation.

	We first present the following elementary fact.
	
	\begin{lemma}\label{lem:3rd_derivative} 
	Let $u \in C^4(\overline{M})$ satisfy $\frac{\pa u}{\pa \vec n}=ce^{-u}-h_{g}$
		on $\pa M$, then under Fermi coordinates around $P \in \pa M$, at $P$ there holds 
		\begin{align*}
			u_{\alpha\beta n}
			= & L_{\beta}^{\gamma}u_{\gamma\alpha}+L_{\alpha}^{\gamma}u_{\gamma\beta}-L_{\alpha\beta}u_{nn}-ce^{-u}u_{\alpha\beta}\\
			& +L_{\alpha;\beta}^{\gamma}u_{\gamma}+L_{\alpha}^{\gamma}L_{\beta\gamma}u_{n}-e^{-u}(c_{\alpha}u_{\beta}+c_{\beta}u_{\alpha}-cu_\beta u_\alpha+cL_{\alpha\beta}u_{n}-c_{;\alpha\beta})\\
			&-(h_{g})_{;\alpha\beta}+R_{\alpha\beta n}^{i}u_{i}.
		\end{align*}
	\end{lemma} 
	\begin{proof}
		Under Fermi coordinates around $P$  the metric can be expressed
		as 
		\[
		g=\ud x_{n}^{2}+g_{\alpha\beta}\ud x_{\alpha}\ud x_{\beta}\qquad\mathrm{in~~}B_{\rho}^{+},
		\]
		where $(x_{1},\cdots,x_{n-1})$ are the geodesic normal coordinates
		on $\pa M$ and $\vec n=\pa_{x_{n}}$ is the inward unit normal on $\pa M$.
		In other words, at $P$ we have 
		\[
		g_{\alpha\beta}=\delta_{\alpha\beta},\quad\Gamma_{\alpha\beta}^{\gamma}=0\qquad\mathrm{for~~}1\leq\alpha,\beta,\gamma\leq n-1.
		\]
		Moreover, in a boundary neighborhood  near $P$ there hold 
		\[
		\Gamma_{\alpha\beta}^{n}=L_{\alpha\beta},\quad\Gamma_{\alpha n}^{\beta}=-L_{\alpha\gamma}g^{\gamma\beta}=-L_{\alpha}^{\beta},\quad\Gamma_{\alpha n}^{n}=0.
		\]
		
		In the following calculations, all involved quantities are evaluated
		at $P$. Notice that 
		\begin{align*}
			u_{\gamma\beta}= & \pa_{\beta}\pa_{\gamma}u-\Gamma_{\beta\gamma}^{i}u_{i}=\pa_{\beta}\pa_{\gamma}u-\Gamma_{\beta\gamma}^{n}u_{n}=u_{;\gamma\beta}-L_{\beta\gamma}u_{n},\\
			u_{\alpha n}= & u_{n\alpha}=\pa_{\alpha}\pa_{n}u-\Gamma_{\alpha n}^{\gamma}u_{\gamma}=\pa_{\alpha}u_{n}+L_{\alpha\gamma}u_{\gamma}
		\end{align*}
		and 
		\begin{align*}
			u_{\alpha\beta n}=u_{\alpha n\beta}+R_{\alpha\beta n}^{i}u_{i}=u_{n\alpha\beta}+R_{i\alpha\beta n}u_{i}.
		\end{align*}
		Then we have 
		\begin{align*}
			u_{n\alpha\beta}= & \pa_{\beta}(u_{n\alpha})-\Gamma_{n\beta}^{i}u_{i\alpha}-\Gamma_{\alpha\beta}^{i}u_{ni}\\
			= & \pa_{\beta}\pa_{\alpha}(u_{n})+L_{\alpha\gamma;\beta}u_{\gamma}+L_{\alpha\gamma}\pa_{\beta}\pa_{\gamma}u-\Gamma_{n\beta}^{\gamma}u_{\gamma\alpha}-\Gamma_{\alpha\beta}^{n}u_{nn}\\
			= & \pa_{\beta}\pa_{\alpha}(u_{n})+L_{\alpha\gamma;\beta}u_{\gamma}+L_{\alpha\gamma}(u_{\gamma\beta}+L_{\beta\gamma}u_{n})+L_{\beta\gamma}u_{\gamma\alpha}-L_{\alpha\beta}u_{nn}
		\end{align*}
		and 
		\begin{align*}
			\pa_{\beta}\pa_{\alpha}(u_{n})= & \pa_{\beta}\pa_{\alpha}(ce^{-u}-h_{g})\\
			= & \pa_{\beta}(e^{-u}\pa_{\alpha}c-ce^{-u}\pa_{\alpha}u)-(h_{g})_{;\alpha\beta}\\
			= & c_{;\alpha\beta}e^{-u}-c_{\alpha}e^{-u}u_{\beta}-c_{\beta}e^{-u}u_{\alpha}+ce^{-u}u_\beta u_\alpha-ce^{-u}\pa_{\beta}\pa_{\alpha}u-(h_{g})_{;\alpha\beta}\\
			= & c_{;\alpha\beta}e^{-u}-c_{\alpha}e^{-u}u_{\beta}-c_{\beta}e^{-u}u_{\alpha}+ce^{-u}u_\beta u_\alpha-ce^{-u}(u_{\alpha\beta}+L_{\alpha\beta}u_{n})-(h_{g})_{;\alpha\beta}.
		\end{align*}
		
		Therefore, putting these facts together the desired assertion follows. 
	\end{proof}

	The proof of Theorem \ref{Thm:local estimates} is accomplished through a series of technical lemmas, mainly Lemmas \ref{lem:tangential trace estimate} and \ref{lem:double_normal_derivatives}.

	\subsection{Trace of second tangential derivatives}
	
	\begin{lemma}\label{lem:tangential trace estimate} 
		Suppose $u \in C^4(\overline M)$ is a solution to \eqref{PDE:sigma_2_with_bdry}.
		Let $\mathcal{O}$ and $\mathcal{O}_1$ be two  neighborhoods 
		near $\pa M$ with $\mathcal{O}_1 \Subset \mathcal{O}$, then  there exist two positive constants  $C _1$ depending on $n, \|g\|_{C^2(\overline{\mathcal{O}})}$, and $C$ depending on $n, \|g\|_{C^4(\overline{\mathcal{O}})}$, $\mathrm{dist}_g(\mathcal{O}_1, \partial\mathcal{O}\backslash \pa M)$, $\sup_{\mathcal O}(c+|\tilde{\nabla}c|+|\tilde{\nabla}^2c|)e^{-u},\sup_{\mathcal{O}}(f+|f_{x}|+|f_{z}|+|f_{xz}|+|f_{zz}|)$ such that
		under Fermi coordinates around some boundary point in $\overline{\mathcal{O}_1}$ there holds
		\[
		-C_1 \leq\sum_{\alpha=1}^{n-1}u_{~~\alpha}^{\alpha}\le C \qquad \mathrm{~~in~~}\quad \overline{\mathcal{O}_1}.
		\]
	\end{lemma}
	\begin{proof}
For the lower bound, since $\lambda(W)\in\Gamma_{2}^{+}$, we have
\[
0<\frac{\partial\sigma_{2}(W)}{\partial W_{nn}}=\sum_{\alpha}W_{\alpha}^{\alpha}=\sum_{\alpha}(u_{~~\alpha}^{\alpha}+u_{\alpha}u^{\alpha}+A_{~~\alpha}^{\alpha})-\frac{n-1}{2}|\nabla u|^{2}.
\]
This implies 
\[
\sum_{\alpha}(u_{~~\alpha}^{\alpha}+A_{~~\alpha}^{\alpha})\geq\frac{n-3}{2}|\nabla u|^{2}
\]
and thereby the lower bound follows.

For the upper bound,  under Fermi coordinates around some boundary point in $\overline{\mathcal{O}}$, the metric is expressed as 
\[
g=\ud x_{n}^{2}+g_{\alpha\beta}\ud x_{\alpha}\ud x_{\beta}\qquad\mathrm{~~in~~}\quad B_{\rho}^{+}.
\]
The boundary functions $c,h_{g}$ are extended to the interior by
letting $c(x)=c(x'),h_{g}(x)=h_{g}(x')$ for $x=(x',x_{n})\in B_{\rho}^{+}$.
In the following proof, we use $\alpha,\beta, \gamma,\zeta,\xi$ as the tangential coordinates from $1$ to $n-1$.

Let $\eta$ be a smooth cut-off function in $\mathcal{O}$ such that
$\frac{\partial\eta}{\partial x_{n}}=0$ on $\partial M\cap\mathcal{O}$,
$\eta=0$ on $\partial\mathcal{O}\cap M$, $\eta=1$ on $\mathcal{O}_{1}\Subset\mathcal{O}$,
and $|\nabla\eta|\leq\frac{C\eta^{1/2}}{r}$, $|\nabla^{2}\eta|\leq\frac{C}{r^{2}}$,
where $r=\mathrm{dist}(\mathcal{O}_{1},\pa\mathcal{O})$.

In $\mathcal{O}$ we define 
\begin{align*}
G & =\eta e^{ax_{n}}E :=\eta e^{ax_{n}}(\sum_{\alpha}u_{~~\alpha}^{\alpha}+\sum_{\alpha} u_\alpha u^\alpha+au_{n}),
\end{align*}
where $a\in\R_{+}$ is  to be determined later. We remark that the above barrier function $G$ is local, depending on Fermi coordinates.

Suppose 
\[
\max_{\overline{\mathcal{O}}}G=G(x_{0})>0
\]
and further assume $\left(\eta e^{ax_{n}}E\right)(x_{0})\gg 1$. For simplicity, we use $u_{~~\alpha}^{\alpha}$
instead of $\sum_{\alpha}u_{~~\alpha}^{\alpha}$.

\vskip 4pt 
\emph{Case 1:} $x_{0}\in\partial M\cap \overline{\mathcal{O}}$. 
\vskip 4pt

At the maximum point $x_{0}$, by Lemma \ref{lem:3rd_derivative} we have 
\begin{align*}
0\ge & \frac{\partial G}{\partial x_{n}}\\
= & a\eta E+\eta\big(u_{~~\alpha n}^{\alpha}+2 u^\alpha u_{\alpha n}+au_{nn}\big)\\
\ge & \eta\big(au_{~~\alpha}^{\alpha}+au_{nn}+2L^{\alpha\gamma}u_{\alpha\gamma}-Hu_{nn}-ce^{-u}u_{~~\alpha}^{\alpha}\big)\\
 & -C\eta\big[|\nabla u|^{2}(1+ce^{-u})+((c+|\tilde \nabla c|)e^{-u}+|\tilde \nabla^2 c|)|\nabla u|\big]\\
\ge & \eta\bigg[a\Delta u-\sup_{\mathcal{O}\cap\partial M}(ce^{-u}+2|L|)\sum_{\alpha,\gamma}|u_{\alpha\gamma}|-(\sup_{\mathcal{O}\cap\partial M}|H|)|u_{nn}|\bigg]\\
 &-C\eta\big[|\nabla u|^{2}(1+ce^{-u})+(ce^{-u}+|\tilde \nabla c| e^{-u}+|\tilde \nabla^2 c|)|\nabla u|\big],
\end{align*}
where $C$ depends on  $\|L\|_{C^{1}(\overline{\mathcal{O}}\cap\partial M)}$ and $g$.

Since $\lambda(W)\in\Gamma_{2}^{+}$, we know that for all $1\leq i,j \leq n$,
\[
|W_{ij}|\le\sigma_{1}(W)\le\Delta u+C
\]
and also
$$|u_{nn}| \leq \Delta u+C.$$
We may assume $\Delta u(x_{0})\geq1$. Then there exists a positive
constant $C_1$ such that 
\[
\sum_{\alpha,\gamma}|u_{\alpha\gamma}|\le C_1\Delta u \qquad \mathrm{~~at~~} x_0.
\]

Take
\[
a=1+2C_{1}\sup_{\mathcal{O}\cap\partial M}(|ce^{-u}|+2|L|)+\sup_{\mathcal{O}\cap\partial M}|H|
\]
to make
\[
\Delta u(x_{0})\leq C\sup_{\mathcal{O}\cap\partial M}\big(|\nabla u|^{2}+(|c|+|\tilde \nabla c|+|\tilde \nabla^2 c|)^{2}e^{-2u}+1\big).
\]
This yields 
\[
G(x)\leq C\sup_{\mathcal{O}\cap\partial M}\big(|\nabla u|^{2}+(|c|+|\tilde \nabla c|+|\tilde \nabla^2 c|)^{2}e^{-2u}+1\big)\qquad\mathrm{~~in~~}\mathcal{O}_{1}\cap\overline{M}.
\]
This in turn implies the desired estimate.

\vskip 4pt 
\emph{Case 2:} $x_{0}\in M\cap\mathcal{O}$. 
\vskip 4pt

At the maximum point $x_{0}$ we have
\begin{align*}
0=\partial_{i}(\log G)=  \frac{\eta_{i}}{\eta}+a\delta_{ni}+\frac{\partial_{i}E}{E} \qquad \mathrm{~~for~~} \quad 1 \leq i \leq n
\end{align*}
and $(\pa_{j}\pa_{i}\log G)$ is negative semi-definite.

 In the following, calculations are evaluated at $x_{0}$ unless otherwise stated, and denote by $C$
 a positive constant depending on $\|\nabla u\|_{C^{1}(\overline{\mathcal{O}})},n,g$ and $\sup_{\mathcal{O}}(ce^{-u})$.

Notice that 
\begin{align*}
\pa_i E=&\partial_{i}u_{~~\alpha}^{\alpha}+\pa_i(|\nabla u|^2-u_n^2)+a\partial_{i}u_{n}\\
=&u_{~~\alpha i}^{\alpha}+\Gamma_{\alpha i}^{k}u^{\alpha}_{~~k}-\Gamma_{ki}^\alpha u^k_{~~\alpha}+2u^{\alpha}u_{\alpha i}-2\Gamma_{ni}^\alpha u_n u_\alpha+a(u_{ni}+\Gamma_{ni}^{\beta}u_{\beta})\\
=&u_{~~\alpha i}^{\alpha}+\Gamma_{\alpha i}^{n}u^{\alpha}_{~~n}-\Gamma_{ni}^\alpha u^n_{~~\alpha}+2u^{\alpha}u_{\alpha i}-2\Gamma_{ni}^\alpha u_n u_\alpha+a(u_{ni}+\Gamma_{ni}^{\beta}u_{\beta}).
\end{align*}
It follows that 
\begin{equation}\label{3rd_derivative_u_nth}
u_{~~\alpha n}^{\alpha}=\pa_n (u_{~~\alpha}^{\alpha})=-(\frac{\eta_n}{\eta}+a)E-2u^{\alpha}u_{\alpha n}-au_{nn}
\end{equation}
and for $1\le\gamma\le n-1$,
\begin{align}\label{3rd_derivative_u}
u_{~~\alpha \gamma}^{\alpha}  =&-\frac{\eta_{\gamma}}{\eta}E-a(u_{n\gamma}+\Gamma_{n\gamma}^{\beta}u_{\beta})\no\\
&-\Gamma_{\alpha \gamma}^{n}u^{\alpha}_{~~n}+\Gamma_{n\gamma}^\alpha u^n_{~~\alpha}-2u^{\alpha}u_{\alpha \gamma}-\Gamma_{\gamma\alpha}^{k}u_{k}u^\alpha+\Gamma_{k\gamma}^\alpha u^ku_\alpha\no\\
=& -\frac{\eta_{\gamma}}{\eta}E-\Gamma_{\alpha\gamma}^{n}u^{\alpha}_{~~n}+\Gamma_{n\gamma}^\alpha u^n_{~~\alpha}-2u^{\alpha}u_{\alpha\gamma}-au_{n\gamma}+O(1),
\end{align}
where we use $O(1)$ to denote the constant depending on $\|\nabla u\|_{C^0(\mathcal{\overline O})}$ and $g$.

A direct computation yields 
\begin{align*}
 &\pa_{j}\pa_{i}\log G=(\log G)_{ij}\\
 =& \frac{\eta_{ij}}{\eta}-\frac{\eta_{i}\eta_{j}}{\eta^{2}}+\frac{E_{ij}}{E}-\frac{E_{i}E_{j}}{E^{2}}+a(x_{n})_{,ij}\\
 =& \frac{\eta_{ij}}{\eta}-\frac{\eta_{i}\eta_{j}}{\eta^{2}}+\frac{E_{ij}}{E}-(\frac{\eta_{i}}{\eta}+a\delta_{ni})(\frac{\eta_{j}}{\eta}+a\delta_{nj})+a(x_{n})_{,ij}\\
=&\frac{\eta_{ij}}{\eta}-2\frac{\eta_{i}\eta_{j}}{\eta^{2}}-a\frac{\eta_{i}}{\eta}\delta_{nj}-a\frac{\eta_{j}}{\eta}\delta_{ni}-a^{2}\delta_{ni}\delta_{nj}+\frac{E_{ij}}{E}+a(x_{n})_{,ij}.
\end{align*}
This together with Cauchy inequality follows that
\begin{align}\label{est:F_log G}
0\ge & F^{ij}(\log G)_{ij}\no\\
\geq & F^{ij}(\frac{\eta_{ij}}{\eta}-3\frac{\eta_{i}\eta_{j}}{\eta^{2}})-2a^{2}F^{nn}+aF^{ij}(x_{n})_{,ij}+F^{ij}\frac{E_{ij}}{E}.
\end{align}

An inaccuracy occurs in the previous work \cite[p.1057]{Chen} involving
calculations to derive that $(\pa_{j}\pa_{i}G)$ is negative semi-definite at some local maximum point of $G$, where $G$ is the barrier function therein. For this reason, we give a correction of such an inaccuracy. 

To this end, notice that
\begin{align*}
(u_n)_{,ij} =&(\langle \nabla u,\nabla x_n\rangle)_{,ij}\\
 =&\big(u^k_{~~i}(x_n)_{,k}+u^k (x_n)_{,ki}\big)_{,j}\\
 =&u^k_{~~ij}(x_n)_{,k}+u_{i}^k(x_n)_{,kj}+u_{j}^k (x_n)_{,ki}+u^k (x_n)_{,kij}\\
 =&u_{nij}-\Gamma_{kj}^n u^k_{~~i}-\Gamma_{ki}^n u^k_{~~j}+u^k (x_n)_{,kij}
\end{align*}
and
\begin{align*}
(u^\alpha u_\alpha)_{,ij}=&(|\nabla u|^2- u_n^2)_{,ij}\\
=&2 u^k u_{kij}+u^k_{~~i}u_{kj}+u^k_{~~j}u_{ki}-2 \pa_i u_n \pa_j u_n-2 u_n (u_n)_{,ij}\\
=&2 u^k u_{kij}+2 u^k_{~~i}u_{kj}-2 (u_{ni}+\Gamma_{ni}^k u_k) (u_{nj}+\Gamma_{nj}^l u_l)\\
&-2 u_n\big[u_{nij}-\Gamma_{kj}^n u^k_{~~i}-\Gamma_{ki}^n u^k_{~~j}+u^k (x_n)_{,kij}\big]\\
=& 2  u^\alpha u_{\alpha ij}+2 u^\alpha_{~~i}u_{\alpha j}\\
&-2 \Gamma_{ni}^k u_k u_{nj} -2 \Gamma_{nj}^k u_k u_{ni}+2 u_n \Gamma_{kj}^n u^k_{~~i}+2 u_n \Gamma_{ki}^n u^k_{~~j}\\
&-2 \Gamma_{ni}^k \Gamma_{nj}^l u_k u_l-2 u_n u^k (x_n)_{,kij}.
\end{align*}

If we write 
\[
u_{~~\alpha}^{\alpha}=\Delta u-u_{nn},
\]
then 
\begin{align*}
(u_{~~\alpha}^{\alpha})_{,ij}=(\Delta u)_{,ij}-(u_{nn})_{,ij}
\end{align*}

Observe that 
\[
u_{nni}=\pa_{i}(u_{nn})-2\Gamma_{in}^{k}u_{kn}
\]
and 
\begin{align*}
u_{nnij}= & \pa_{j}(u_{nni})-\Gamma_{jn}^{k}u_{kni}-\Gamma_{jn}^{k}u_{nki}-\Gamma_{ji}^{k}u_{nnk}\\
= & \pa_{j}\pa_{i}(u_{nn})-2\pa_{j}(\Gamma_{in}^{k}u_{kn})-2\Gamma_{jn}^{k}u_{nki}-\Gamma_{ji}^{k}u_{nnk}\\
= & \pa_{j}\pa_{i}(u_{nn})-2\pa_{j}\Gamma_{in}^{\beta}u_{\beta n}-2\Gamma_{in}^{\beta}(u_{\beta nj}+\Gamma_{j\beta}^{l}u_{ln}+\Gamma_{jn}^{l}u_{\beta l})\\
 & -2\Gamma_{jn}^{\beta}u_{n\beta i}-\Gamma_{ji}^{k}u_{nnk}.
\end{align*}
On the other hand, we have 
\begin{align*}
(u_{nn})_{,ij}= & \pa_{j}\pa_{i}(u_{nn})-\Gamma_{ij}^{k}\pa_{k}(u_{nn})\\
= & \pa_{j}\pa_{i}(u_{nn})-\Gamma_{ij}^{k}(u_{nnk}+2\Gamma_{kn}^{l}u_{ln}).
\end{align*}
Combining  above terms we obtain 
\begin{align*}
(u_{~~\alpha}^{\alpha})_{,ij}= & u_{~~\alpha ij}^{\alpha}-2\pa_{j}\Gamma_{in}^{\beta}u_{\beta n}-2\Gamma_{in}^{\beta}(u_{\beta nj}+\Gamma_{j\beta}^{l}u_{ln}+\Gamma_{jn}^{l}u_{\beta l})\\
 & -2\Gamma_{jn}^{\beta}u_{n\beta i}-\Gamma_{ji}^{k}u_{nnk}+\Gamma_{ij}^{k}(u_{nnk}+2\Gamma_{kn}^{l}u_{ln})\\
= & u_{~~\alpha ij}^{\alpha}-2(\Gamma_{jn}^{\beta}u_{n\beta i}+\Gamma_{in}^{\beta}u_{\beta nj})\\
 & -2\pa_{j}\Gamma_{in}^{\beta}u_{\beta n}-2\Gamma_{in}^{\beta}(\Gamma_{j\beta}^{l}u_{ln}+\Gamma_{jn}^{l}u_{\beta l})+2\Gamma_{ij}^{k}\Gamma_{kn}^{l}u_{ln}.
\end{align*}

Consequently, putting these terms together we obtain
\begin{align*}
E_{ij}=&u_{~~\alpha ij}^{\alpha}-2(\Gamma_{jn}^{\beta}u_{n\beta i}+\Gamma_{in}^{\beta}u_{\beta nj})+2  u^\alpha u_{\alpha ij}+au_{nij}+2u^\alpha_{~~i}u_{\alpha j}\\
 & -2\pa_{j}\Gamma_{in}^{\beta}u_{\beta n}-2\Gamma_{in}^{\beta}(\Gamma_{j\beta}^{l}u_{ln}+\Gamma_{jn}^{l}u_{\beta l})+2\Gamma_{ij}^{k}\Gamma_{kn}^{l}u_{ln}\\
&-2 \Gamma_{ni}^k u_k u_{nj} -2 \Gamma_{nj}^k u_k u_{ni}+2 u_n \Gamma_{kj}^n u^k_{~~i}+2 u_n \Gamma_{ki}^n u^k_{~~j}\\
 &-a\Gamma_{kj}^n u^k_{~~i}-a\Gamma_{ki}^n u^k_{~~j}\\
&-2 \Gamma_{ni}^k \Gamma_{nj}^l u_k u_l-2 u_n u^k (x_n)_{,kij}+au^k (x_n)_{,kij}.
\end{align*}

Taking covariant derivatives of $E$ under Fermi coordinates around a boundary point,
we could not regard these terms $u_{~~\alpha}^{\alpha}$ and $u^\alpha u_\alpha$, $u_{n}$
in $E$ as $(1,1)$-tensor and $(0,1)$-tensor, but just functions.
This is the remarkable difference between ours and those of \cite{Chen}.

Keep in mind that as in \eqref{def:F_j^i}:
   $$F_j^i:=\frac{\pa \sigma_2(W)}{\pa W_i^j}=(T_1)_j^i\quad \mathrm{~~and~~} \quad F_j^i=g^{ik}F_{kj},$$
then it is not hard to see that
\begin{equation}\label{def:2nd_D_sigma_2}
F_{js}^{ir}:=\frac{\pa^2 \sigma_2(W)}{\pa W_i^j \pa W_r^s}=\frac{\pa}{\pa W_r^s}(\sigma_1(W)\delta_{j}^i-W_j^i)=\delta_{s}^r\delta_j^i-\delta^{i}_{s}\delta^{r}_j.
\end{equation}

Consequently, by \eqref{est:F_log G} we obtain 
\begin{align}\label{est:F*E}
F^{ij}E_{ij}\ge & F^{ij}u_{~~\alpha ij}^{\alpha}-4F^{ij}\Gamma_{jn}^{\beta}u_{n\beta i}+2F^{ij}u_{\alpha j}u^\alpha_{~~i}+aF^{ij}u_{nij}+2F^{ij}u^{\alpha}u_{\alpha ij}\no\\
 & -C(\sum_i F_i^{i})\sum_{k,\beta}|u_{\beta k}|-C(\sum_i F_i^{i})|u_{nn}|-C\sum_{i}F_i^{i}.
\end{align}

By definition of $\sigma_2(W)$ we know that
\begin{equation}\label{def:sigma_2}
 \sigma_{2}(W^{\top})+\sigma_{1}(W^{\top})W_{nn}-\sum_{\alpha}W_{n\alpha}W_{~~n}^{\alpha}=f^2.
 \end{equation}
This together with the fact that $u_{~~\alpha}^{\alpha}\gg 1$ yields
\begin{equation}
\sigma_{1}(W^{\top})|W_{nn}|=|f^2+\sum_{\alpha=1}^{n}W_{n\alpha}^{2}-\sigma_{2}(W^{\top})|\le2\sum_{k=1}^{n}\sum_{\alpha=1}^{n-1}|u_{\alpha k}|^{2}.\label{eq:contorl 1}
\end{equation}
 Keep in mind that $\sigma_{1}(W^{\top}) \gg 1$ as $u_{~~\alpha}^{\alpha} \gg 1$. 

We combine (\ref{eq:contorl 1}) and 
\begin{equation}\label{eq:upper bound of F^ii}
\sum_i F_i^{i}=(n-1)\sigma_{1}(W)\ge(n-1)|W_{ij}|\qquad\mathrm{~~for~~}\quad1\le i,j\le n
\end{equation}
to show
\begin{equation}\label{eq:good bound of W_nn^2}
  (\sum_i F_i^{i})\sum_{k,\alpha} |u_{\alpha k}|^{2}
\ge \frac{1}{2} (\sum_i F_i^{i})\sigma_{1}(W^{\top})|W_{nn}|
\ge \frac{n-1}{2}\sigma_{1}(W^{\top})|W_{nn}|^{2}.
\end{equation}

We first deal with 
\begin{align*}
 & F^{ij}u^\alpha_{~~\alpha ij}\\
= & F^{ij}(u_{ij~~\alpha}^{~~\alpha}+R_{~~~~\alpha i}^{m\alpha}u_{mj}+2R^m_{~~i\alpha j}u_{m}^{~~\alpha}+R^{m\alpha}_{~~~~\alpha j}u_{mi}+R^{m\alpha}_{~~~~\alpha i,j}u_{m}+R^{m~\alpha}_{~~i~~j,\alpha}u_{m})\\
\geq & F^{ij}u_{ij~~\alpha}^{~~\alpha}+2F^{ij}R_{~~~~\alpha i}^{m\alpha}u_{mj}+2F^{ij}R^m_{~~i\alpha j}u_{m}^{~~\alpha}-C\sum_{i}F_i^{i}.
\end{align*}
Notice that
\begin{align*}
|2F^{ij}R_{~~~~\alpha i}^{m\alpha}u_{mj} |=&|2F^{in}R_{~~~~\alpha i}^{m\alpha}u_{mn}+2F^{i\beta}R_{~~~~\alpha i}^{m\alpha}u_{m\beta}|\\
 \le& C(\sum_i F_i^{i})\sum_{m,\alpha}|u_{m\alpha}|+C(\sum_i F_i^{i})|u_{nn}|
\end{align*}
and 
\[
2\left|F^{ij}R^m_{~~i\alpha j}u_{m}^{~~\alpha}\right|\le C(\sum_i F_i^{i})\sum_{m,\alpha}|u_{m\alpha}|	.
\]

We are ready to handle the first three terms on the right hand side
of (\ref{est:F*E}). Denote 
\[
I=F^{ij}u_{ij~~\alpha}^{~~\alpha}.
\]

A direct computation yields 
\[
W_{ij,~~\alpha}^{~~\alpha}=u_{ij~~\alpha}^{~~\alpha}+2u_{i}^{~\alpha}u_{j\alpha}+u_{i}u_{j~~\alpha}^{~\alpha}+u_{j}u_{i~~\alpha}^{~\alpha}-\left(u^{k}u_{k~~\alpha}^{~\alpha}+u^{k\alpha}u_{k\alpha}\right)g_{ij}+A_{ij,~~\alpha}^{~~~\alpha}.
\]
We employ second covariant derivatives to the equation in (\ref{PDE:sigma_2_with_bdry}) to show
\[
(f^2)^{\alpha}_{~~\alpha}-F_{js}^{ir}W_{i,}^{j~\alpha}W^s_{r,\alpha}=F^{ij}W_{ij,~~\alpha}^{~~\alpha},
\]
where $(f^2)^{\alpha}_{~~\alpha}=\nabla_\alpha \nabla^\alpha \left(f^2(x,u)\right)$. Then we have 
\begin{align*}
I\geq & F^{ij}\left(W_{ij,~~\alpha}^{~~~\alpha}-2u_{i}^{~\alpha}u_{j\alpha}-2u_{i}u_{j~~\alpha}^{~\alpha}+\left(u^{k}u_{k~~\alpha}^{~\alpha}+u^{k\alpha}u_{k\alpha}\right)g_{ij}\right)-C\sum_{i}F_i^{i}\\
\geq & F^{ij}W_{ij,~~\alpha}^{~~~\alpha}+F^{ij}(-2u^\alpha_{~~\alpha i}u_{j}+u^{k}u^{\alpha}_{~~\alpha k}g_{ij})+F^{ij}(-2u_{i}^{~\alpha}u_{j\alpha}+u^{k\alpha}u_{k\alpha}g_{ij})-C\sum_{i}F_i^{i}\\
= &-F_{js}^{ir}W_{i,}^{j~\alpha}W^s_{r,\alpha}+ (f^2)^{\alpha}_{~~\alpha}+F^{ij}(-2u^\alpha_{~~\alpha i}u_{j}+u^{k}u^{\alpha}_{~~\alpha k}g_{ij})\\
&+F^{ij}(-2u_{i}^{~\alpha}u_{j\alpha}+u^{k\alpha}u_{k\alpha}g_{ij})-C\sum_{i}F_i^{i}\\
:=& -F_{js}^{ir}W_{i,}^{j~\alpha}W^s_{r,\alpha}+I_{1}.
\end{align*}
 
We apply \eqref{def:2nd_D_sigma_2} to show
\[
-F_{js}^{ir}W_{i,}^{j~\alpha}W^s_{r,\alpha}=W_{i,}^{j~\alpha}W^i_{j,\alpha}-W_{i,}^{i~\alpha}W^j_{j,\alpha}.
\]

A key ingredient is to control the following terms involving third covariant derivatives.
\vskip 8pt
\noindent\fbox{\begin{minipage}[t]{1\columnwidth - 2\fboxsep - 2\fboxrule}%
\begin{claim}
There holds
\begin{align}\label{key-est:3rd_derivatives}
 & -4F^{ij}\Gamma_{jn}^{\beta}u_{n\beta i}+W_{i,}^{j~\alpha}W^i_{j,\alpha}-W_{i,}^{i~\alpha}W^j_{j,\alpha}\no\\
\ge & -C(\sum_i F_i^{i})(|u_{nn}|+\sum_{k,\alpha}|u_{k\alpha}|+\frac{\sum_{k,\alpha}|u_{k\alpha}|^{2}}{\sigma_{1}(W^{\top})})-C\frac{1}{\eta}(E^{2}+W_{nn}^{2}+\sum_{k,\beta}|u_{k\beta}|^{2}).
\end{align}
\end{claim}
\end{minipage}}
\vskip 8pt

To this end, we write
\begin{align*}
 & -4F^{ij}\Gamma_{jn}^{\beta}u_{n\beta i}+W_{i,}^{j~\alpha}W^i_{j,\alpha}-W_{i,}^{i~\alpha}W^j_{j,\alpha}\\
= & (W_{i,}^{j~\alpha}W^i_{j,\alpha}-W_{n,}^{n~\alpha}W^n_{n,\alpha})-W_{\beta,}^{\beta~\alpha}W^\gamma_{\gamma,\alpha}-2W_{\beta,}^{\beta~\alpha}W^n_{n,\alpha}-4F^{ij}\Gamma_{jn}^{\beta}u_{n\beta i}
\end{align*}
and simplify
$$W_{i,}^{j~\alpha}W^i_{j,\alpha}-W_{n,}^{n~\alpha}W^n_{n,\alpha}=|W^{n}_{\gamma,\alpha}|^{2}+|W^{\gamma}_{\zeta,\alpha}|^{2}=|W^{k}_{\gamma,\alpha}|^{2}.$$
Our purpose is to show that the first term $|W^{k}_{\gamma,\alpha}|^{2}$ can absorb all other three terms.

By \eqref{3rd_derivative_u} we have
\begin{align}\label{eq:1st_derivative_W}
W^{\beta}_{\beta,\alpha}=&u^{\beta}_{~~\beta\alpha}+2u^{\beta}u_{\beta\alpha}-(n-1)u^{k}u_{k\alpha}+A^{\beta}_{\beta,\alpha},\no\\
=&-\frac{\eta_{\alpha}}{\eta}E-\Gamma_{\beta\alpha}^{n}u^{\beta}_{~~n}+\Gamma_{n\alpha}^\beta u^n_{~~\beta}-au_{n\alpha}-(n-1)u^{k}u_{k\alpha}+O(1).
\end{align}
Then, it is not hard to see that
\begin{align*}
 & W_{\beta,}^{\beta~\alpha}W^\gamma_{\gamma,\alpha}\leq C\frac{|\nabla\eta|^{2}}{\eta^{2}}E^{2}+C\sum_{k,\beta}|u_{k\beta}|^{2}.
\end{align*}

Differentiating the equation in \eqref{PDE:sigma_2_with_bdry} in the variable $x_\alpha$ to show
\[
\sigma_{1}(W^{\top})W_{nn,\alpha}-2W_{n}^{\gamma}W^{n}_{\gamma,\alpha}+F^{\gamma\zeta}W_{\gamma\zeta,\alpha}=F^{ij}W_{ij,\alpha}=(f^2)_{\alpha},
\]
we have
\[
W_{nn,\alpha}=\frac{(f^2)_{\alpha}-F^{\gamma\zeta}W_{\gamma\zeta,\alpha}+2W_{n}^{\gamma}W^{n}_{\gamma,\alpha}}{\sigma_{1}(W^{\top})}.
\]

We split the following difficult term as follows:
\begin{align*}
 & 2W_{\beta,}^{\beta~\alpha}W^n_{n,\alpha}\\
= & 2\left(-\frac{\eta_{\alpha}}{\eta}E-\Gamma_{\beta\alpha}^{n}u^{\beta}_{~~n}+\Gamma_{n\alpha}^\beta u^n_{~~\beta}-au_{n\alpha}-(n-1)u^{k}u_{k\alpha}+O(1)\right)\\
&~~\cdot\frac{(f^2)^{\alpha}-F^{\gamma\zeta}W_{\gamma\zeta,}^{~~~\alpha}+2W_{n}^{\gamma}W^{n~\alpha}_{\gamma,}}{\sigma_{1}(W^{\top})}\\
:= & J_{1}+J_{2}+J_{3}+J_{4}+J_{5},
\end{align*}
where 
\begin{align*}
J_{1}=&-2\frac{\eta_{\alpha}}{\eta}E\frac{(f^2)^{\alpha}-F^{\gamma\zeta}W_{\gamma\zeta,}^{~~~\alpha}+2W_{n}^{\gamma}W^{n~\alpha}_{\gamma,}}{\sigma_{1}(W^{\top})},\\
J_{2}=&2(-\Gamma_{\beta\alpha}^{n}u^{\beta}_{~~n}+\Gamma_{n\alpha}^\beta u^n_{~~\beta})\frac{(f^2)^{\alpha}-F^{\gamma\zeta}W_{\gamma\zeta,}^{~~~\alpha}+2W_{n}^{\gamma}W^{n~\alpha}_{\gamma,}}{\sigma_{1}(W^{\top})},\\
J_{3}=&-2au_{n\alpha}\frac{(f^2)^{\alpha}-F^{\gamma\zeta}W_{\gamma\zeta,}^{~~~\alpha}+2W_{n}^{\gamma}W^{n~\alpha}_{\gamma,}}{\sigma_{1}(W^{\top})},\\
J_{4}=&-2(n-1)u^{k}u_{k\alpha}\frac{(f^2)^{\alpha}-F^{\gamma\zeta}W_{\gamma\zeta,}^{~~~\alpha}+2W_{n}^{\gamma}W^{n~\alpha}_{\gamma,}}{\sigma_{1}(W^{\top})},\\
J_{5}=&O(1)\frac{(f^2)^{\alpha}-F^{\gamma\zeta}W_{\gamma\zeta,}^{~~~\alpha}+2W_{n}^{\gamma}W^{n~\alpha}_{\gamma,}}{\sigma_{1}(W^{\top})}.
\end{align*}

We next deal with $J_{1},\cdots,J_{5}$ one by one. Without loss of generality, we assume $\frac{1}{2}\sigma_{1}(W^{\top})\le E\le2\sigma_{1}(W^{\top})$ recalling that $\sigma_{1}(W^{\top})\gg 1$.

We estimate
\begin{align}\label{eq:J_1}
|J_{1}| \le& C\frac{|\nabla\eta|}{\eta}|(f^2)_{\alpha}-F^{\gamma\zeta}W_{\gamma\zeta,\alpha}+2W_{n}^{\gamma}W^{n}_{\gamma,\alpha}|\no\\
 \le&\frac{1}{16}|W^{\gamma}_{\zeta,\alpha}|^{2}+\frac{1}{16}|W^{n}_{\zeta,\alpha}|^{2}+CW_{nn}^{2}\frac{|\nabla\eta|^{2}}{\eta^{2}}+C\frac{|\nabla\eta|^{2}}{\eta^{2}}\sum_{k,\beta}|u_{k\beta}|^{2}.
\end{align}

Since $J_3$ has the same order as $J_2$, it suffices to estimate $J_2$.

Notice that
\begin{align*}
\frac{1}{2}|J_{2}|\le& \left|(\Gamma_{n\alpha}^\beta u^n_{~~\beta}-\Gamma_{\beta\alpha}^{n}u^{\beta}_{~~n})\frac{(f^2)^{\alpha}-F^{\gamma\zeta}W_{\gamma\zeta,}^{~~~\alpha}+2W_{n}^{\gamma}W^{n~\alpha}_{\gamma,}}{\sigma_{1}(W^{\top})}\right|\\
\le & C\frac{\sum_{\beta}|u_{n\beta}|}{\sigma_{1}(W^{\top})}+\frac{|F^{\gamma\zeta}W_{\gamma\zeta,}^{~~~\alpha}(\Gamma_{\beta\alpha}^{n}u^{\beta}_{~~n}-\Gamma_{n\alpha}^\beta u^n_{~~\beta})|}{\sigma_{1}(W^{\top})}\\
&+\frac{2|W_{n}^{\gamma}W^{n~\alpha}_{\gamma,}(\Gamma_{n\alpha}^\beta u^n_{~~\beta}-\Gamma_{\beta\alpha}^{n}u^{\beta}_{~~n})|}{\sigma_{1}(W^{\top})}\\
:=& C\frac{\sum_{\beta}|u_{n\beta}|}{\sigma_{1}(W^{\top})}+|J_{21}|+|J_{22}|,
\end{align*}
where 
\[
J_{21}=\frac{F^{\gamma\zeta}W_{\gamma\zeta,}^{~~~\alpha}(\Gamma_{\beta\alpha}^{n}u^{\beta}_{~~n}-\Gamma_{n\alpha}^\beta u^n_{~~\beta})}{\sigma_{1}(W^{\top})},\qquad J_{22}=\frac{2W_{n}^{\gamma}W^{n~\alpha}_{\gamma,}(\Gamma_{n\alpha}^\beta u^n_{~~\beta}-\Gamma_{\beta\alpha}^{n}u^{\beta}_{~~n})}{\sigma_{1}(W^{\top})}.
\]

By \eqref{def:sigma_2} and Young inequality we  get the control of
\begin{align}
\frac{1}{2}\sum_{\beta}|u_{n\beta}|^{2} \le &\sigma_{1}(W^{\top})W_{nn}+\sigma_{2}(W^{\top})+C\le\sigma_{1}(W^{\top})(W_{nn}+\frac{n-2}{2n-2}\sigma_{1}(W^{\top}))+C\no\\
  \le&\sigma_{1}(W^{\top})\sum_i F_i^{i}+C.\label{eq:upper bound of u_nbeta}
\end{align}

By \eqref{eq:1st_derivative_W} we have
\begin{align}\label{sum:3rd_derivative_W}
&F^{\gamma\zeta}W_{\gamma\zeta,\alpha}\no\\
  =&\sigma_{1}(W)W^{\beta}_{\beta,\alpha}-W_{\gamma}^{\zeta}W^{\gamma}_{\zeta,\alpha}\no\\
  =&\sigma_{1}(W)[-\frac{\eta_{\alpha}}{\eta}E-\Gamma_{\beta\alpha}^{n}u^{\beta}_{~~n}+\Gamma_{n\alpha}^\beta u^n_{~~\beta}-au_{n\alpha}-(n-1)u^{k}u_{k\alpha}+O(1)]-W_{\gamma}^{\zeta}W^{\gamma}_{\zeta,\alpha}.
\end{align}
It follows that
\begin{align*}
& J_{21}\\
=&\frac{\sigma_{1}(W)}{\sigma_{1}(W^{\top})}[-\frac{\eta^{\alpha}}{\eta}E-g^{\alpha\gamma}\Gamma_{\beta\gamma}^{n}u^{\beta}_{~~n}+g^{\alpha\gamma}\Gamma_{n\gamma}^\beta u^n_{~~\beta}-au_n^{\alpha}-(n-1)u^{k}u_k^{\alpha}+O(1)](\Gamma_{\beta\alpha}^{n}u^{\beta}_{~~n}-\Gamma_{n\alpha}^\beta u^n_{~~\beta})\\
 & +\frac{W_{\gamma}^{\zeta}W^{\gamma~\alpha}_{\zeta,}(\Gamma_{n\alpha}^\beta u^n_{~~\beta}-\Gamma_{\beta\alpha}^{n}u^{\beta}_{~~n})}{\sigma_{1}(W^{\top})}.
 \end{align*}
This together with \eqref{eq:upper bound of u_nbeta} yields
 \begin{align*}
  |J_{21}|\le& C\frac{|\nabla\eta|}{\eta}\sigma_{1}(W)\sum_{\beta}|u_{n\beta}|+C\frac{\sigma_{1}(W)}{\sigma_{1}(W^{\top})}(\sum_{k,\beta}|u_{k\beta}|^{2})+\frac{1}{16}|W_{\gamma,\alpha}^{\zeta}|^{2}+C\frac{|W^{\top}|^{2}(\sum_{\beta}|u_{n\beta}|^{2})}{\sigma_{1}(W^{\top})^{2}}\\
\le& \frac{1}{16}|W_{\gamma,\alpha}^{\zeta}|^{2}+C\frac{|\nabla\eta|}{\eta}\sigma_{1}(W)\sum_\beta |u_{n\beta}|+C\frac{\sigma_{1}(W)}{\sigma_{1}(W^{\top})}(\sum_{k,\beta}|u_{k\beta}|^{2})\\
&+C\frac{(\sum_i F_i^{i})(\sum_{\gamma,\zeta}|u_{\gamma\zeta}|^{2})}{\sigma_{1}(W^{\top})}+C\frac{\sum_{\gamma,\zeta}|u_{\gamma\zeta}|^{2}}{\sigma_{1}^{2}(W^{\top})}.
\end{align*}

Again by (\ref{eq:upper bound of u_nbeta}) we have
\begin{align*}
|J_{22}|= & \frac{|2W_{n}^{\gamma}W^{n~\alpha}_{\gamma,}(\Gamma_{n\alpha}^\beta u^n_{~~\beta}-\Gamma_{\beta\alpha}^{n}u^{\beta}_{~~n})|}{\sigma_{1}(W^{\top})}\\
\le & \frac{1}{16}|W^{n}_{\gamma,\alpha}|^{2}+C\frac{(\sum_\beta |u_{n\beta}|^{2})(\sum_\gamma|W_{n\gamma}|^{2})}{\sigma_{1}^{2}(W^{\top})}\\
\le & \frac{1}{16}|W^{n}_{\gamma,\alpha}|^{2}+C\frac{(\sum_\gamma|W_{n\gamma}|^{2})(\sigma_{1}(W^{\top})\sum_i F_i^{i}+1)}{\sigma_{1}^{2}(W^{\top})}\\
\le & \frac{1}{16}|W^{n}_{\gamma,\alpha}|^{2}+C\frac{(\sum_\gamma|W_{n\gamma}|^{2})(\sum_i F_i^{i})}{\sigma_{1}(W^{\top})}+C\frac{\sigma_{1}(W^{\top})\sum_i F_i^{i}+1}{\sigma_{1}^{2}(W^{\top})}.
\end{align*}

Hence, putting these estimates together and recalling that $J_3$ and $J_2$ have the same order, we are ready to show that
\begin{align}\label{eq:J_2&J_3}
&\frac{1}{2}|J_{2}|+|J_3|\no\\
\le & \frac{1}{8}|W^{k}_{\gamma,\alpha}|^{2}+C\frac{\sum_{\beta}|u_{k\beta}|}{\sigma_{1}(W^{\top})}+C\frac{|\nabla\eta|}{\eta}\sigma_{1}(W)\sum_{\beta}|u_{n\beta}|\no\\
&+C\frac{\sigma_{1}(W)}{\sigma_{1}(W^{\top})}(\sum_{k,\beta}|u_{k\beta}|^{2})+\frac{(\sum_i F_i^{i})(\sum_{\gamma,\zeta}|u_{\gamma\zeta}|^{2})}{\sigma_{1}(W^{\top})}+C\frac{\sum_{\gamma,\zeta}|u_{\gamma\zeta}|^{2}}{\sigma_{1}^{2}(W^{\top})}\no\\
 & +C\frac{(\sum_\gamma|W_{n\gamma}|^{2})(\sum_i F_i^{i})}{\sigma_{1}(W^{\top})}+C\frac{\sigma_{1}(W^{\top})\sum_i F_i^{i}+1}{\sigma_{1}^{2}(W^{\top})}.
\end{align}

By \eqref{sum:3rd_derivative_W} we have
\begin{align*}
&-\frac{1}{2(n-1)}J_{4}\\
 =&u^{k}u_k^{\alpha}\frac{(f^2)_{\alpha}-F^{\gamma\zeta}W_{\gamma\zeta,\alpha}+2W_{n}^{\gamma}W^{n}_{\gamma,\alpha}}{\sigma_{1}(W^{\top})}\\
 =&u^{k}u_k^{\alpha}\frac{(f^2)_{\alpha}+2W_{n}^{\gamma}W^{n}_{\gamma,\alpha}}{\sigma_{1}(W^{\top})}-(u^{\xi}u_{\xi}^\alpha+u_{n}u_n^{\alpha})\frac{F^{\gamma\zeta}W_{\gamma\zeta,\alpha}}{\sigma_{1}(W^{\top})}\\
 =&\frac{u^{k}u_{k}^{\alpha}(f^2)_{\alpha}}{\sigma_{1}(W^{\top})}+\frac{2u^{k}u_k^{\alpha}W_{n}^{\gamma}W^{n}_{\gamma,\alpha}}{\sigma_{1}(W^{\top})}\\
 &-u^{\xi}u_\xi^{\alpha}\frac{\sigma_{1}(W)[-\frac{\eta_{\alpha}}{\eta}E-\Gamma_{\beta\alpha}^{n}u^{\beta}_{~~n}+\Gamma_{n\alpha}^\beta u^n_{~~\beta}-au_{n\alpha}-(n-1)u^{k}u_{k\alpha}+O(1)]-W_{\gamma}^{\zeta}W^{\gamma}_{\zeta,\alpha}}{\sigma_{1}(W^{\top})}\\
 &-u_{n}u_{n}^{\alpha}\frac{\sigma_{1}(W)[-\frac{\eta_{\alpha}}{\eta}E-\Gamma_{\beta\alpha}^{n}u^{\beta}_{~~n}+\Gamma_{n\alpha}^\beta u^n_{~~\beta}-au_{n\alpha}-(n-1)u^{k}u_{k\alpha}+O(1)]-W_{\gamma}^{\zeta}W^{\gamma}_{\zeta,\alpha}}{\sigma_1(W^{\top})}.
\end{align*}
This together with \eqref{eq:upper bound of u_nbeta} gives
\begin{align*}
|J_{4}|\le & \frac{1}{16}|W^{k}_{\gamma,\alpha}|^{2}+C\frac{\sum_{k,\beta}|u_{k\beta}|}{\sigma_{1}(W^{\top})}+C\frac{|\nabla\eta|}{\eta}\sigma_{1}(W)\sum_{\beta}|u_{n\beta}|\\
&+C\frac{\sigma_{1}(W)}{\sigma_{1}(W^{\top})}(\sum_{k,\beta}|u_{k\beta}|^{2})+C\frac{(\sum_i F_i^{i})(\sum_{\gamma,\zeta}|u_{\gamma\zeta}|^{2})}{\sigma_{1}(W^{\top})}+C\frac{\sum_{\gamma,\zeta}|u_{\gamma\zeta}|^{2}}{\sigma_{1}^{2}(W^{\top})}\\
 & +C\frac{(\sum_{\gamma}|W_{n\gamma}|^{2})(\sum_i F_i^{i})}{\sigma_{1}(W^{\top})}+C\frac{\sigma_{1}(W^{\top})\sum_i F_i^{i}+1}{\sigma_{1}^{2}(W^{\top})}+\frac{|u^{\beta}u_{\beta\alpha}W_{\gamma}^{\zeta}W^{\gamma}_{\zeta,\alpha}|}{\sigma_{1}(W^{\top})}.
\end{align*}
Note that 
\begin{align*}
 \frac{|u^{\beta}u_{\beta}^{\alpha}W_{\gamma}^{\zeta}W^{\gamma}_{\zeta,\alpha}|}{\sigma_{1}(W^{\top})}
\le C\frac{|W^\top|^2\sum_{\beta,\alpha}|u^{\beta}u_{\beta\alpha}|^{2}}{\sigma_{1}(W^{\top})^{2}}+\frac{1}{16}|W^{\gamma}_{\zeta,\alpha}|^{2}.
\end{align*}

By \eqref{def:sigma_2} we have
\[
|W^{\top}|^{2}=\sigma_{1}(W^{\top})^{2}-2\sigma_{2}(W^{\top})=-2(f^2+W_{n\alpha}W_{~~n}^{\alpha}-\sigma_{1}(W^{\top})W_{nn})+\sigma_{1}(W^{\top})^{2}
\]
and then
\begin{align*}
\frac{|W^\top|^2}{\sigma_1(W^{\top})^{2}}\leq 1+2\frac{W_{nn}}{\sigma_1(W^{\top})}.
\end{align*}
Thus, we obtain 
\[
\frac{|W^{\top}|^{2}\sum_{\beta,\alpha} |u^{\beta}u_{\beta\alpha}|^{2}}{\sigma_{1}(W^{\top})^{2}}\le C(\sum_{\beta,\alpha}|u_{\beta\alpha}|^{2})\left(\frac{\sum_i F_i^{i}}{\sigma_{1}(W^{\top})}+1\right)
\]
and then 
\begin{align}\label{eq:J_4}
|J_{4}|\leq&\frac{1}{8}|W^{k}_{\gamma,\alpha}|^{2}+C\frac{\sum_{\beta}|u_{k\beta}|}{\sigma_{1}(W^{\top})}+C\frac{|\nabla\eta|}{\eta}\sigma_{1}(W)|\sum_{\beta}|u_{n\beta}|+C\frac{\sigma_{1}(W)}{\sigma_{1}(W^{\top})}(\sum_{k,\beta}|u_{k\beta}|^{2})\no\\
&+C\frac{(\sum_i F_i^{i})(\sum_{\gamma,\zeta}|u_{\gamma\zeta}|^{2})}{\sigma_{1}(W^{\top})}+C\frac{\sum_{\gamma,\zeta}|u_{\gamma\zeta}|^{2}}{\sigma_{1}^{2}(W^{\top})}+C\frac{(\sum_{\gamma}|W_{n\gamma}|^{2})(\sum_i F_i^{i})}{\sigma_{1}(W^{\top})}\no\\
&+C\frac{\sigma_{1}(W^{\top})\sum_i F_i^{i}+1}{\sigma_{1}^{2}(W^{\top})}+C(\sum_{\beta,\alpha}|u_{\beta\alpha}|^{2})\left(\frac{\sum_i F_i^{i}}{\sigma_{1}(W^{\top})}+1\right).
\end{align}

It is relatively easier to estimate
\begin{align}\label{eq:J_5}
|J_{5}| \le &\frac{1}{16}|W^{k}_{\gamma,\alpha}|^{2}+ \frac{C}{\sigma_{1}(W^{\top})}+C\frac{|\nabla\eta|}{\eta}\sigma_{1}(W)+C\frac{\sigma_{1}(W)}{\sigma_{1}(W^{\top})}(\sum_{k,\beta}|u_{k\beta}|)\no\\
&+C\frac{\sum_{\gamma,\zeta}|u_{\gamma\zeta}|^{2}}{\sigma_{1}^{2}(W^{\top})} +C\frac{\sum_{\gamma}|W_{n\gamma}|^{2}}{\sigma_{1}(W^{\top})}.
\end{align}

Notice that
\[
u_{ni\beta}=W_{ni,\beta}-(u_{n\beta}u_{i}+u_{n}u_{i\beta}-u^{k}u_{k\beta}g_{ni}+A_{ni,\beta}).
\] 
Then we have
\begin{align*}
 & -4F^{ij}\Gamma_{jn}^{\beta}u_{n\beta i}\\
\ge & -4F^{i\gamma}\Gamma_{\gamma n}^{\beta}u_{ni\beta}-C\sum_i F_i^{i}\\
= & -4F^{n\gamma}\Gamma_{\gamma n}^{\beta}u_{nn\beta}-4F^{\zeta\gamma}\Gamma_{\gamma n}^{\beta}u_{n\zeta\beta}-C\sum_i F_i^{i}\\
= & -4F^{n\gamma}\Gamma_{\gamma n}^{\beta}u_{nn\beta}-4F^{\zeta\gamma}\Gamma_{\gamma n}^{\beta}W_{n\zeta,\beta}-C(\sum_i F_i^{i})(1+\sum_{m,\alpha}|u_{m\alpha}|)\\
\ge & -4F^{n\gamma}\Gamma_{\gamma n}^{\beta}W_{nn,\beta}-\frac{1}{16} |W^{n}_{\zeta,\beta}|^{2}-C|W|^2-C(\sum_i F_i^{i})(1+\sum_{m,\alpha} |u_{m\alpha}|).
\end{align*}

A similar argument to the estimate of $J_2$ also yields 
\begin{align}\label{eq:boundary bad third derivative}
&-4F^{n\gamma}\Gamma_{\gamma n}^{\beta}W_{nn,\beta} \no\\
=&4W_{n}^{\gamma}\Gamma_{\gamma n}^{\beta}W_{nn,\beta}\no\\
 =&4W_{n}^{\gamma}\Gamma_{\gamma n}^{\beta}\frac{(f^2)_{\beta}-F^{\alpha\zeta}W_{\alpha\zeta,\beta}+2W_{n}^{\alpha}W^{n}_{\alpha,\beta}}{\sigma_{1}(W^{\top})}\nonumber \\
  =&4\frac{W_{n}^{\gamma}\Gamma_{\gamma n}^{\beta}(f^2)_{\beta}}{\sigma_{1}(W^{\top})}-\frac{4W_{n}^{\gamma}\Gamma_{\gamma n}^{\beta}}{\sigma_{1}(W^{\top})}F^{\alpha\zeta}W_{\alpha\zeta,\beta}+\frac{8W_{n}^{\gamma}\Gamma_{\gamma n}^{\beta}W_{n}^{\alpha}W^{n}_{\alpha,\beta}}{\sigma_{1}(W^{\top})}\nonumber \\
  \ge&-\frac{1}{16}|W^{k}_{\gamma,\alpha}|^{2}-C\frac{\sum_{\beta}|u_{k\beta}|}{\sigma_{1}(W^{\top})}-C\frac{|\nabla\eta|}{\eta}\sigma_{1}(W)\sum_{\beta}|u_{n\beta}|-C\frac{\sigma_{1}(W)}{\sigma_{1}(W^{\top})}(\sum_{k,\beta}|u_{k\beta}|^{2})\no\\
 &-C\frac{(\sum_i F_i^{i})(\sum_{\gamma,\zeta}|u_{\gamma\zeta}|^{2})}{\sigma_{1}(W^{\top})}-C\frac{\sum_{\gamma,\zeta}|u_{\gamma\zeta}|^{2}}{\sigma_{1}^{2}(W^{\top})}\nonumber \\
 & -C\frac{(\sum_{\gamma}|W_{n\gamma}|^{2})(\sum_i F_i^{i})}{\sigma_{1}(W^{\top})}-C\frac{\sigma_{1}(W^{\top})\sum_i F_i^{i}+1}{\sigma_{1}^{2}(W^{\top})}.
\end{align}

Let us combine these third derivative terms: (\ref{eq:J_1}),(\ref{eq:J_2&J_3}),(\ref{eq:J_4}),(\ref{eq:J_5}),(\ref{eq:boundary bad third derivative}) to show
\begin{align*}
 &  -4F^{ij}\Gamma_{jn}^{\beta}u_{n\beta i}+W_{i,}^{j~\alpha}W^i_{j,\alpha}-W_{i,}^{i~\alpha}W^j_{j,\alpha}\\
\ge &\frac{1}{2}|W^{k}_{\gamma,\alpha}|^{2}-CW_{nn}^{2}\frac{|\nabla\eta|^{2}}{\eta^{2}}-C\frac{\sigma_{1}(W)}{\sigma_{1}(W^{\top})}(\sum_{k,\beta}|u_{k\beta}|^{2})-C\frac{|\nabla\eta|^{2}}{\eta^{2}}\sum_{k,\beta}|u_{k\beta}|^{2}\\
&-C\frac{|\nabla\eta|}{\eta}\sigma_{1}(W)\sum_{\beta}|u_{n\beta}|
 -C\frac{\sum_{\beta}|u_{k\beta}|}{\sigma_{1}(W^{\top})}-C\frac{(\sum_i F_i^{i})(\sum_{\gamma,\zeta}|u_{\gamma\zeta}|^{2})}{\sigma_{1}(W^{\top})}-C\frac{\sum_{\gamma,\zeta}|u_{\gamma\zeta}|^{2}}{\sigma_{1}^{2}(W^{\top})}\no\\
 & -C\frac{(\sum_{\gamma}|W_{n\gamma}|^{2})(\sum_i F_i^{i})}{\sigma_{1}(W^{\top})}-C\frac{\sigma_{1}(W^{\top})\sum_i F_i^{i}+1}{\sigma_{1}^{2}(W^{\top})}-C(\sum_{\beta,\alpha}|u_{\beta\alpha}|^{2})\no\\
  \ge & -C(\sum_i F_i^{i})(|u_{nn}|+\sum_{k,\alpha}|u_{k\alpha}|+\frac{\sum_{k,\alpha}|u_{k\alpha}|^{2}}{\sigma_{1}(W^{\top})})-C\frac{1}{\eta}(E^{2}+W_{nn}^{2}+\sum_{k,\beta}|u_{k\beta}|^{2}).
\end{align*}
This completes the proof of the Claim.

By (\ref{3rd_derivative_u}) we have 
\begin{align*}
 & F^{ij}(-2u^{\alpha}_{~~\alpha i}u_{j}+u_{k}u^{\alpha}_{~~\alpha k}g_{ij})\\
= & -2u_{j}F^{ij}\big[-E(\frac{\eta_{i}}{\eta}+a\delta_{ni})-\Gamma_{\alpha i}^{n}u^{\alpha}_{~~n}+\Gamma_{ni}^\alpha u^n_{~~\alpha}-2u^{\alpha}u_{\alpha i}+2\Gamma_{ni}^\alpha u_n u_\alpha-a(u_{ni}+\Gamma_{ni}^{\beta}u_{\beta})\big]\\
 & +(\sum_{i}F_i^{i})u^{k}\big[-E(\frac{\eta_{k}}{\eta}+a\delta_{nk})-\Gamma_{\alpha k}^{n}u^{\alpha}_{~~n}+\Gamma_{nk}^\alpha u^n_{~~\alpha}-2u^{\alpha}u_{\alpha k}+2\Gamma_{nk}^\alpha u_n u_\alpha-a(u_{nk}+\Gamma_{nk}^{\beta}u_{\beta})\big]\\
= & 2u_{j}F^{ij}\Gamma_{\alpha i}^{n}u_{n\alpha}-2u_{j}F^{ij}\Gamma_{ni}^\alpha u^n_{~~\alpha}+2au_{j}F^{ij}(u_{ni}+\Gamma_{ni}^{\beta}u_{\beta})+4u_{j}F^{ij}u_{\alpha}(u^{\alpha}_{~~i}-\Gamma_{ni}^\alpha u_n)\\
 & +(\sum_{i}F_i^{i})u^{k}(-\Gamma_{\alpha k}^{n}u^{\alpha}_{~~n}+\Gamma_{nk}^\alpha u^n_{~~\alpha})-a(\sum_{i}F_i^{i})u^{k}(u_{nk}+\Gamma_{nk}^{\beta}u_{\beta})\\
 &-2u^{k}(\sum_{i}F_i^{i})u_{\alpha}(u^{\alpha}_{~~k}-\Gamma_{nk}^\alpha u_n)+2Eu_{j}F^{ij}(\frac{\eta_{i}}{\eta}+a\delta_{ni})-E(\sum_{i}F_i^{i})u^{k}(\frac{\eta_{k}}{\eta}+a\delta_{nk})\\
\ge & 2Eu_{j}F^{ij}(\frac{\eta_{i}}{\eta}+a\delta_{ni})-E(\sum_{i}F_i^{i})u^{k}(\frac{\eta_{k}}{\eta}+a\delta_{nk})\\
 & +4u_{j}F^{ij}u^{\alpha}u_{\alpha i}-2(\sum_{i}F_i^{i})u^{k}u^{\alpha}u_{\alpha k}+2au_{j}F^{ij}u_{ni}+a(\sum_{i}F_i^{i})u^{k} u_{nk}\\
 & -C(\sum_i F_i^{i})\sum_{k,\beta}|u_{\beta k}|-C(\sum_i F_i^{i}).
\end{align*}

Using 
\begin{align*}
F^{ij}u_{\alpha ij} & =F^{ij}u_{ij\alpha}+F^{ij}R^k_{\alpha ij}u_{k} \ge F^{ij}u_{ij\alpha}-C\sum_{i}F_i^{i}
\end{align*}
and 
\[
W_{ij,l}=u_{ijl}+u_{i}u_{jl}+u_{j}u_{il}-u^{k}u_{kl}g_{ij}+A_{ij,l},
\]
we estimate 
\begin{align*}
II:= & 2F^{ij}u^{\alpha}u_{\alpha ij}\\
\geq & 2u^{\alpha}F^{ij}\bigg(W_{ij,\alpha}-(u_{i}u_{j\alpha}+u_{j}u_{i\alpha}-u^{k}u_{k\alpha}g_{ij}+A_{ij,\alpha})\bigg)-C\sum_{i}F_i^{i}\\
= & 2u^{\alpha}(f^2)_{\alpha}-2u^{\alpha}F^{ij}(u_{i}u_{j\alpha}+u_{j}u_{i\alpha}-u^{k}u_{k\alpha}g_{ij}+A_{ij,\alpha})-C\sum_{i}F_i^{i}
\end{align*}
and 
\begin{align*}
III:=&a  F^{ij}u_{nij}\\
\geq& a(f^2)_{n}-C\sum_{i}F_i^{i} -aF^{ij}(u_{i}u_{jn}+u_{j}u_{in}-u^{k}u_{kn}g_{ij}+A_{ij,n}).
\end{align*}

Consequently, we obtain 
\begin{align*}
 & I_{1}+II+III\\
\geq & (f^2)^{\alpha}_{~~\alpha}-2F^{ij}u_{\alpha i}u_{\alpha j}+(\sum_{i}F_i^{i})(\sum_{k,\alpha}|u_{k\alpha}|^{2})\\
 & +2Eu_{j}F^{ij}(\frac{\eta_{i}}{\eta}+a\delta_{ni})-E(\sum_{i}F_i^{i})u^{k}(\frac{\eta_{k}}{\eta}+a\delta_{nk})\\
 & +2u^{\alpha}(f^2)_{\alpha}+a(f^2)_{n}-C(\sum_i F_i^{i})\sum_{k,\beta}|u_{\beta k}|-C(\sum_i F_i^{i}).
\end{align*}

Since $\eta\sigma_{1}(W^{\top})\gg 1$ by assumption,  it follows from (\ref{eq:good bound of W_nn^2}) that
\begin{align*}
&\frac{1}{4}(\sum_{i}F_i^{i})(\sum_{m}u_{m\alpha}^{2})\\
\ge& C\frac{1}{\eta}(E^{2}+W_{nn}^{2}+\sum_{k,\beta}|u_{k\beta}|^{2})+C(\sum_iF_i^{i})(|u_{nn}|+\sum_{k,\alpha}|u_{k\alpha}|+\frac{\sum_{k,\alpha}|u_{k\alpha}|^{2}}{\sigma_{1}(W^{\top})}).
\end{align*}

Therefore, putting these facts together, we conclude from (\ref{est:F*E})
that 
\begin{align*}
 & F^{ij}E_{ij}\\
\geq & (f^2)^{\alpha}_{~~\alpha}+(\sum_{i}F_i^{i})(\sum_{k,\alpha}|u_{k\alpha}|^{2})\\
 & +2Eu_{j}F^{ij}(\frac{\eta_{i}}{\eta}+a\delta_{ni})-E(\sum_{i}F_i^{i})u_{k}(\frac{\eta_{k}}{\eta}+a\delta_{nk})-C\|f\|_{C^{1}(\overline{\mathcal{O}})}\\
 & -C(\sum_i F_i^{i})(|u_{nn}|+|u_{k\alpha}|+\frac{|u_{k\alpha}|^{2}}{\sigma_{1}(W^{\top})})-C\frac{1}{\eta}(E^{2}+W_{nn}^{2}+|u_{k\beta}|^{2})-C\sum_{i}F_i^{i}\\
\ge & (f^2)^{\alpha}_{~~\alpha}+\frac{3}{4}(\sum_{i}F_i^{i})(\sum_{k,\alpha}|u_{k\alpha}|^{2})-C(\sum_{i}F_i^{i})\frac{|\nabla\eta|}{\eta}E-CE\sum_{i}F_i^{i}\\
&-C\sum_{i}F_i^{i}-C\|f\|_{C^{1}(\overline{\mathcal{O}})}.
\end{align*}

Eventually, going back to (\ref{est:F_log G}) and collecting estimates above we have
\begin{align}
0\ge & EF^{ij}(\log G)_{ij}\no\label{eq:final tangent derivative}\\
\ge & EF^{ij}(\frac{\eta_{ij}}{\eta}-3\frac{\eta_{i}\eta_{j}}{\eta^{2}})-2a^{2}F^{nn}E\no\\
 & +(f^2)^{\alpha}_{~~\alpha}+\frac{3}{4}(\sum_{i}F_i^{i})(\sum_{k,\alpha}|u_{k\alpha}|^{2})\no\\
 & -C_{1}(\sum_{i}F_i^{i})\frac{|\nabla\eta|}{\eta}E-C_{2}E\sum_{i}F_i^{i}-C\sum_{i}F_i^{i}-C\|f\|_{C^{1}(\overline{\mathcal{O}})}\no\\
\ge & EF^{ij}(\frac{\eta_{ij}}{\eta}-3\frac{\eta_{i}\eta_{j}}{\eta^{2}})+\frac{1}{4}(\sum_{i}F_i^{i})(\sum_{k}|u_{k\alpha}|^{2})-C(\sum_{i}F_i^{i})\frac{|\nabla\eta|}{\eta}E-C,
\end{align}
where $C$ also depends on $\sup_{\mathcal{O}}(f+|f_{x}|+|f_{z}|+|f_{xz}|+|f_{zz}|)$.
Here we emphasize that $(\sum_{i}F_i^{i})(\sum_{k,\alpha}|u_{k\alpha}|^{2})$
is the dominated term, which can absorb most of other terms on the right hand side of (\ref{eq:final tangent derivative}).
Therefore, we can apply (\ref{eq:final tangent derivative}) to conclude
that 
\[
(\eta\sum_{k,\alpha}|u_{k\alpha}|^{2})(x_{0})\leq C\quad\Longrightarrow\quad(\eta e^{ax_{n}}E)(x_{0})\le C.
\]

In conclusion, we establish that in $\mathcal{O}_{1}$, there holds
\[
e^{ax_{n}}E\leq C\quad\Longrightarrow\quad\sum_{\alpha=1}^{n-1}u_{~~\alpha}^{\alpha}\leq C.
\]
This completes the proof. 
\end{proof}

	\subsection{Double normal derivative}
	
	For brevity, we let $u_{n}=ce^{-u}-h_{g}:=\varphi(x,u)$ on $\partial M$.
	
	\begin{lemma}\label{lem:double_normal_derivatives} 
		Assume the assumptions as in Theorem \ref{Thm:local estimates}. Then for any $x_{0}\in\partial M$ and a sufficiently small $r>0$, there exists a positive constant $C$ depending on $n,r,\sup_{M_{r}}|\nabla u|$,
		$\sup_{M_r}(f+|f_{x}|+|f_{z}|)$ and $(\inf_{M_{r}}f)^{-1}$, such that
		\[
		u_{nn}(x_{0})\le C,
		\]
		where $M_{r}=\{x\in\overline{M};d_{g}(x,x_{0})<r\}$. 
	\end{lemma}
	\begin{proof}
Fix $x_{0}\in\partial M$. In a tubular neighborhood of $x_0$ near $\partial M$
we define 
\begin{align*}
G(x)= & \langle\nabla u,\nabla d_{g}(x,\partial M)\rangle-ce^{-u}+h_{g}+\left(\langle\nabla u,\nabla d_{g}(x,\partial M)\rangle-ce^{-u}+h_{g}\right)^{7/5}\\
 & -\frac{1}{2}Nd_{g}(x,\partial M)-N_{1}d_{g}^{2}(x_{0},x),
\end{align*}
where $N,N_{1}\in\mathbb{R}_{+}$ are to be determined
later and $N>N_{1}$. 

For our purpose, it suffices to consider $G$ in $M_{r}:=\{x\in\overline{M};d_{g}(x,x_{0})<r\}$ for small $r>0$. We claim that with appropriate choice of $N$ and $N_1$,
$x_{0}$ is a maximum point of $G$ in $\overline{M_{r}}$. 

If we admit the above claim temporarily, then under Fermi coordinates around $x_0$  with $\vec n=\pa_{x_{n}}$ we have 
\begin{align*}
0  \ge G_{n}(x_{0})=&(u_{nn}-\varphi_{n})(x_{0})[1+\frac{7}{5}(u_{n}-\varphi)^{2/5}(x_{0})]-\frac{1}{2}N\\
 =& u_{nn}(x_{0})-\varphi_{n}(x_{0})-\frac{1}{2}N.
\end{align*}
This gives $u_{nn}(x_{0})\le\frac{1}{2}N+\|\varphi\|_{C^{1}(\overline{M_{r}})}$. 

We first choose $N_1$ such that
\begin{align}
N_{1}r^{2}> & \sup_{M_{r}}|\langle\nabla u,\nabla d_{g}(x,\partial M)\rangle-ce^{-u}+h_{g}+(\langle\nabla u,\nabla d_{g}(x,\partial M)\rangle-ce^{-u}+h_{g})^{7/5}|.\label{eq:condtion 0-1}
\end{align}
It is not hard to see that $G<0$ on $\{d_{g}(x,x_{0})=r\}\cap M$;
$G(x_{0})=0$ and $G<0$ on $(\partial M\cap M_{r})\setminus\{x_{0}\}$.
Let $\eta(x)=d_{g}^{2}(x,x_{0})$, then there exists a positive constant
$\delta=\delta(r)$ such that 
\[
\|\eta\|_{C^{2}(\overline{M_{r}})}\leq\delta.
\]

Until now the constant $N_{1}$ has been fixed, it remains to choose a sufficiently large $N$ depending on $r,N_{1},g$, $\sup_{M_{r}}(|c|+|\nabla u|+e^{-u})$, $\sup_{M_{r}}(f+|f_{x}|+|f_{z}|)$
and $\sup_{M_{r}}\frac{1}{f(x,u)}$ such that the maximum point of $G$ is
$x_{0}$. Without confusion, the following constant $C$ may be different line
by line and only depends on the aforementioned data.

By contradiction, we assume 
\[
\max_{\overline{M_{r}}}G=G(x_{1})>0\quad\mathrm{~~for~~}\quad x_{1}\in\mathring{M}_{r}.
\]

We may choose Fermi coordinates around $x_{0}$ such that 
\[
g=\ud x_{n}^{2}+g_{\alpha\beta}\ud x_{\alpha}\ud x_{\beta}\qquad\mathrm{~~in~~}B_{\rho}^{+}
\]
with $\nu=\pa_{x_{n}}$, meanwhile $W^{\top}=\big(W_{\beta}^{\alpha}\big)$
is diagonal at $x_{1}$.  For $x$ near $\partial M$, we let $d_{g}(x,x')=d_{g}(x,\partial M)$
for some $x'\in\partial M$ and extend $c(x)=c(x'),h_{g}(x)=h_{g}(x')$
for $x=(x',x_{n})\in B_{\rho}^{+}$.

Notice that $G(x_1)>0$ implies
\begin{align*}
 & \big(u_{n}(x_{1})-c(x_{1}')e^{-u(x_{1})}+h_{g}(x_{1}')\big)(1+(u_{n}(x_{1})-c(x_{1}')e^{-u(x_{1})}+h_{g}(x_{1}'))^{2/5})\\
> & \frac{1}{2}Nd_{g}(x_{1},\partial M)+N_{1}d_{g}^{2}(x_{0},x_{1})>0
\end{align*}
implies $A:=u_{n}(x_{1})-c(x_{1}')e^{-u(x_{1})}+h_{g}(x_{1}')>0$. More precisely, we have 
\begin{equation}
A\ge\frac{\frac{1}{2}Nx_{n}+N_{1}d_{g}^{2}(x_{0},x_{1})}{1+A_{1}^{2/5}},\label{eq:lower bound of A}
\end{equation}
where $A_{1}=\sup_{M_{r}}(|\nabla u|+|h_{g}|+ce^{-u})$.

Near $x_{1}$ there holds
\[
G(x)=u_{n}-ce^{-u}+h_{g}+\big(u_{n}-ce^{-u}+h_{g}\big)^{7/5}-\frac{1}{2}Nx_{n}-N_{1}d^{2}(x_{0},x).
\]
Then at the maximum point $x_1$, for all $1\le i\le n$ we have    
\begin{align*}
0=G_{i}= & \big(u_{ni}+u^{\alpha}(x_{n})_{\alpha i}-(ce^{-u}-h_{g})_{i}\big)(1+\frac{7}{5}\big(u_{n}-ce^{-u}+h_{g}\big)^{\frac{2}{5}})\\
&-\frac{N}{2}\delta_{ni}-2N_{1}d_{g}(x_{0},x)d_{g}(x_{0},x)_{i}.
\end{align*}

In particular, 
\begin{equation}
u_{nn}=\frac{\frac{1}{2}N+2N_{1}d_{g}(x_{0},x)d_{g}(x_{0},x)_{n}}{1+\frac{7}{5}A^{\frac{2}{5}}}+(ce^{-u}-h_{g})_{n}\ge\frac{\frac{1}{4}N}{1+\frac{7}{5}A^{\frac{2}{5}}},\label{lbd:u_nn}
\end{equation}
which requires 
\begin{equation*}
N\geq4N_{1}\delta+4(1+\frac{7}{5}A^{\frac{2}{5}})\|\nabla\varphi\|_{C^{0}(\overline{M_{r}})}\label{res_cod_N1}
\end{equation*}
and 
\begin{equation}
u_{n\alpha}=-u^{\beta}(x_{n})_{\beta\alpha}+\frac{2N_{1}d_{g}(x_{0},x)d_{g}(x_{0},x)_{\alpha}}{1+\frac{7}{5}A^{\frac{2}{5}}}+\varphi_{\alpha}.\label{u_n alpha}
\end{equation}

Since $\lambda(W)\in\Gamma_{2}^{+}$, by Lemma \ref{lem:tangential trace estimate} we have for all $1 \leq i,j \leq n$,
\begin{equation}
|W_{i}^{j}|\le \sigma_1(W)=\sigma_1(W^{\top})+W_{nn}\le2u_{nn}\le3W_{n}^{n} \qquad \mathrm{~~at~~}x_{1}.\label{eq:double tangential initial bound}
\end{equation}
At $x_{1}$, $(G_{ij})$ is negative semi-definite and a direct computation yields 
\begin{align*}
G_{ij}(x)= & \left(u_{nij}+u_{j}^{k}(x_{n})_{ki}+u^{k}(x_{n})_{kij}+u_{i}^{k}(x_{n})_{kj}-\varphi_{ij}\right)(1+\frac{7}{5}A^{\frac{2}{5}})\\
 & +\frac{14}{25}A^{-\frac{3}{5}}(u_{ni}+u^{k}(x_{n})_{ki}-\varphi_{i})(u_{nj}+u^{k}(x_{n})_{kj}-\varphi_{j})\\
 & -\frac{1}{2}N(x_{n})_{ij}-N_{1}\eta_{ij}\\
= & (u_{j}^{\beta}(x_{n})_{\beta i}+u^{k}(x_{n})_{kij}+u_{i}^{\beta}(x_{n})_{\beta j})(1+\frac{7}{5}A^{\frac{2}{5}})+(u_{nij}-\varphi_{ij})(1+\frac{7}{5}A^{\frac{2}{5}})\\
 & +\frac{14}{25}A^{-\frac{3}{5}}(u_{ni}+u^{k}(x_{n})_{ki}-\varphi_{i})(u_{nj}+u^{k}(x_{n})_{kj}-\varphi_{j})\\
 & -\frac{1}{2}N(x_{n})_{ij}-N_{1}\eta_{ij}.
\end{align*}
Consequently, we obtain
\begin{align}\label{est:crucial}
0\ge & F_{j}^{i}G_{i}^{j}\no\\
= & \frac{14}{25}A^{-\frac{3}{5}}F_{j}^{i}(u_{ni}+u^{k}(x_{n})_{ki}-\varphi_{i})(u_{n}^{j}++u^{l}(x_{n})_{l}^{j}-\varphi^{j})\no\\
 & +(1+\frac{7}{5}A^{\frac{2}{5}})F_{j}^{i}\left(u^{\beta j}(x_{n})_{\beta i}+u^{k}(x_{n})_{ki}^{j}+u_{i}^{\beta}(x_{n})_{\beta}^{j}\right)\no\\
 & +F_{j}^{i}(u_{ni}^{\quad j}-\varphi_{i}^{j})(1+\frac{7}{5}A^{\frac{2}{5}})-\frac{N}{2}F_{j}^{i}(x_{n})_{i}^{j}-\delta N_{1}\sum_{i}F_{i}^{i}.
\end{align}

By definition of $\sigma_{2}(W)$ we have 
\begin{equation}
\sigma_{1}(W^{\top})W_{n}^{n}+F_{\beta}^{\alpha}W_{\alpha}^{\beta}-2\sum_{\alpha}W_{n}^{\alpha}W_{\alpha}^{n}=F_{j}^{i}W_{i}^{j}=2\sigma_{2}(W)=2f^{2}.\label{def1:sigma_2}
\end{equation}

It also follows from that 
\begin{align*}
F_{\alpha}^{n}=-W_{\alpha}^{n}=-(u_{n\alpha}+u_{n}u_{\alpha}+A_{n\alpha})
\end{align*}
and 
\[
F_{n}^{n}=\sigma_{1}(W^{\top})=\sum_{\alpha}W_{\alpha}^{\alpha}=\frac{f^{2}-\sigma_{2}(W^{\top})+\sum_{\alpha}W_{n}^{\alpha}W_{\alpha}^{n}}{W_{n}^{n}}.
\]

Observe that 
\begin{align*}
W_{n}^{n}=W_{nn}= & u_{nn}+u_{n}^{2}-\frac{1}{2}|\nabla u|^{2}+A_{nn}\\
= & \frac{\frac{1}{2}N+2N_{1}d_{g}(x_{0},x)d_{g}(x_{0},x)_{n}}{1+\frac{7}{5}A^{\frac{2}{5}}}+O(1),
\end{align*}
where $O(1)$ depends on $|\nabla u|$.

Observe that
\begin{align}\label{bdd_trace_F-1-1}
\sum_{i}F_{i}^{i}=(n-1)\sigma_{1}(W)> & (n-1)W_{nn}\no\\
\geq & (n-1)u_{nn}-C\ge\frac{\frac{1}{4}N(n-2)}{1+\frac{7}{5}A^{\frac{2}{5}}}. 
\end{align}

For later use, by (\ref{lbd:u_nn}), (\ref{u_n alpha}) and (\ref{bdd_trace_F-1-1})
we collet a series of estimates here: 
\begin{equation}
\begin{split}F_{n}^{n}u_{n}^{n}= & \sigma_{1}(W^{\top})u_{nn},\\
|F_{\beta}^{\alpha}u_{n}^{\beta}|\le & C(\sum_{i}F_{i}^{i})\bigg(\frac{N_{1}d_{g}(x_{0},x)}{1+\frac{7}{5}A^{\frac{2}{5}}}+1\bigg),\\
|F_{n}^{\alpha}u_{n}^{n}|\le & C\left(1+\frac{N_{1}d_{g}(x_{0},x)}{1+\frac{7}{5}A^{\frac{2}{5}}}\right)\frac{N}{1+\frac{7}{5}A^{\frac{2}{5}}}\le C(\sum_{i}{F}_{i}^{i})\left(1+\frac{N_{1}d_{g}(x_{0},x)}{1+\frac{7}{5}A^{\frac{2}{5}}}\right),\\
|F_{\alpha}^{n}u_{n}^{\alpha}|\le & C\left(1+\frac{N_{1}d_{g}(x_{0},x)}{1+\frac{7}{5}A^{\frac{2}{5}}}\right)\frac{N_{1}d(x_{0},x)}{1+\frac{7}{5}A^{\frac{2}{5}}}\le C(\sum_{i}{F}_{i}^{i})\left(1+\frac{N_{1}d_{g}(x_{0},x)}{1+\frac{7}{5}A^{\frac{2}{5}}}\right),
\end{split}
\label{est:important-1}
\end{equation}
where the above constant $C$ only depends on $g$
and $\|\varphi\|_{C^{1}(\overline{M_{r}})}$.

Notice that 
\[
F_{j}^{i}u_{ik}^{j}=F_{j}^{i}\bigg(W_{i,k}^{j}-(u_{i}u_{k}^{j}+u^{j}u_{ik}-u_{p}u_{k}^{p}\delta_{i}^{j}+A_{i,k}^{j})\bigg)
\]
and 
\[
\varphi_{i}^{j}=c_{i}^{j}e^{-u}-c_{i}u^{j}e^{-u}-c^{j}u_{i}e^{-u}-cu_{i}^{j}e^{-u}+cu_{i}u^{j}e^{-u}-(h_{g})_{i}^{j}.
\]
By (\ref{def1:sigma_2}), (\ref{bdd_trace_F-1-1})
and (\ref{est:important-1}) we obtain 
\begin{align*}
 & F_{j}^{i}(u_{ni}^{\quad j}-\varphi_{i}^{j})\\
= & F_{j}^{i}u_{in}^{j}+F_{j}^{i}R_{nil}^{j}u^{l}-F_{j}^{i}\varphi_{i}^{j}\\
= & 2ff_{n}-2F_{j}^{i}u_{i}u_{n}^{j}+(\sum_{i}{F}_{i}^{i})u^{p}u_{pn}-F_{j}^{i}A_{i,n}^{j}+F_{j}^{i}R_{nil}^{j}u^{l}\\
 & -F_{j}^{i}(c_{i}^{j}e^{-u}-2u_{i}c^{j}e^{-u}+cu_{i}u^{j}e^{-u}-(h_{g})_{i}^{j})\\
 & +cF_{j}^{i}(W_{i}^{j}-u_{i}u^{j}+\frac{1}{2}|\nabla u|^{2}\delta_{i}^{j}-A_{i}^{j})e^{-u}\\
= & 2ff_{n}-2F_{n}^{n}u_{n}u_{n}^{n}-2F_{n}^{\alpha}u_{\alpha}u_{n}^{n}-2F_{\alpha}^{n}u_{n}u_{n}^{\alpha}-2F_{\beta}^{\alpha}u_{\alpha}u_{n}^{\beta}\\
 & +(\sum_{i}{F}_{i}^{i})u^{n}u_{nn}+(\sum_{i}{F}_{i}^{i})u^{\alpha}u_{\alpha n}-F_{j}^{i}A_{i,n}^{j}+F_{j}^{i}R_{nil}^{j}u^{l}\\
 & -F_{j}^{i}(c_{i}^{j}e^{-u}-2u_{i}c^{j}e^{-u}+cu_{i}u^{j}e^{-u}-(h_{g})_{i}^{j})\\
 & +cF_{j}^{i}(W_{i}^{j}-u_{i}u^{j}+\frac{1}{2}|\nabla u|^{2}\delta_{i}^{j}-A_{i}^{j})e^{-u}\\
\ge & [(n-1)W_{n}^{n}+(n-3)\sigma_{1}({W}^{\top})]u_{n}u_{nn}-C(\sum_{i}{F}_{i}^{i})-C(\sum_{i}{F}_{i}^{i})\frac{N_{1}d_{g}(x_{0},x)}{1+\frac{7}{5}A^{\frac{2}{5}}},
\end{align*}
where $C$ depends on $\sup_{M_{r}}(|\nabla u|+c+f+|\nabla f|)$
and $g$.

By (\ref{lbd:u_nn}) we have 
\begin{align}
 & (n-1)W_{n}^{n}u_{n}u_{nn}(1+\frac{7}{5}A^{\frac{2}{5}})\no\\
= & (n-1)W_{n}^{n}u_{n}\bigg(\frac{1}{2}N+2N_{1}d_{g}(x_{0},x)d_{g}(x_{0},x)_{n}+(1+\frac{7}{5}A^{2/5})\varphi_{n}\bigg)\no\\
\ge & \frac{n-1}{2}W_{n}^{n}Nu_{n}-C\delta N_{1}W_{n}^{n}-CW_{n}^{n}.\label{key5-1}
\end{align}

 Let $\widetilde{\varphi}_{i}=\varphi_{i}-u^{k}(x_{n})_{ki}$
and write 
\begin{align*}
 & F_{j}^{i}(u_{ni}+u^{k}(x_{n})_{ki}-\varphi_{i})(u_{n}^{j}++u^{l}(x_{n})_{l}^{j}-\varphi^{j})\\
= & F_{j}^{i}\widetilde{\varphi}_{i}\widetilde{\varphi}^{j}-2F_{j}^{i}u_{ni}\widetilde{\varphi}^{j}+F_{j}^{i}u_{ni}u_{n}^{j}.
\end{align*}

Notice that
\begin{align*}
F_{j}^{i}\widetilde{\varphi}_{i}\widetilde{\varphi}^{j}= & \sigma_{1}({W}^{\top}){\widetilde{\varphi}}_{n}^{2}+\sigma_{1}({W}){\widetilde{\varphi}}_{\alpha}\widetilde{\varphi}^{\alpha}-{W}_{\alpha}^{\beta}{\widetilde{\varphi}}^{\alpha}{\widetilde{\varphi}}_{\beta}-2{W}_{n}^{\alpha}{\widetilde{\varphi}}^{n}{\widetilde{\varphi}}_{\alpha}\nonumber \\
\ge & u_{nn}\sum_{\alpha}{\widetilde{\varphi}}_{\alpha}\widetilde{\varphi}^{\alpha}-C\sum_{\alpha}|u_{n\alpha}|+\sigma_{1}({W}^{\top}){\widetilde{\varphi}}^{i}{\widetilde{\varphi}}_{i}-{W}_{\alpha}^{\beta}{\widetilde{\varphi}}^{\alpha}{\widetilde{\varphi}}_{\beta}
\end{align*}
and
\begin{align*}
 & -2F_{j}^{i}u_{ni}\widetilde{\varphi}^{j}\nonumber \\
= & -2\sigma_{1}({W}^{\top})u_{nn}\widetilde{\varphi}^{n}+2{W}_{\alpha}^{n}u_{nn}\widetilde{\varphi}^{\alpha}+2{W}_{n}^{\alpha}u_{n\alpha}\widetilde{\varphi}^{n}-2(\sigma_{1}({W})\delta_{\alpha}^{\beta}-{W}_{\alpha}^{\beta})u_{n\beta}\widetilde{\varphi}^{\alpha}\nonumber \\
\geq & -2\sigma_{1}({W}^{\top})u_{nn}\widetilde{\varphi}_{n}+2[{W}_{\alpha}^{n}u_{nn}-\sigma_{1}({W})u_{n\alpha}]\widetilde{\varphi}^{\alpha}+2W_{\alpha}^{\beta}u_{n\beta}\widetilde{\varphi}^{\alpha}-C\sum_{\alpha}|u_{n\alpha}|^{2}\nonumber \\
\ge & -2\sigma_{1}({W}^{\top})u_{nn}\widetilde{\varphi}_{n}+2(u_{n}u_{\alpha}+A_{n\alpha})u_{nn}\widetilde{\varphi}^{\alpha}-2\sigma_{1}(W^{\top})u_{n\alpha}\widetilde{\varphi}^{\alpha}\\
&+2W_{\alpha}^{\beta}u_{n\beta}\widetilde{\varphi}^{\alpha}-C\sum_{\alpha}|u_{n\alpha}|^{2}\nonumber \\
\ge & -2\sigma_{1}({W}^{\top})\widetilde{\varphi}_{n}u_{nn}-u_{nn}\sum_{\alpha}\widetilde{\varphi}_{\alpha}\widetilde{\varphi}^{\alpha}-\sum_{\alpha}(u_{n}u_{\alpha}+A_{n\alpha})(u_{n}u^{\alpha}+A_{n}^{\alpha})u_{nn}\\
 & -2\sigma_{1}(W^{\top})u_{n\alpha}\widetilde{\varphi}^{\alpha}+2W_{\alpha}^{\beta}u_{n\beta}\widetilde{\varphi}^{\alpha}-C\sum_{\alpha}|u_{n\alpha}|^{2}.
\end{align*}
We turn to deal with the most difficult term: 
\begin{align}
 & F_{j}^{i}u_{ni}u_{n}^{j}\no\\
= & F_{n}^{n}u_{nn}u_{n}^{n}+F_{\beta}^{\alpha}u_{n\alpha}u_{n}^{\beta}+F_{\beta}^{n}u_{nn}u_{n}^{\beta}+F_{n}^{\beta}u_{n\beta}u_{n}^{n}\nonumber \\
= & \sigma_{1}({W}^{\top})u_{nn}^{2}+(\sigma_{1}({W})\delta_{\beta}^{\alpha}-W_{\beta}^{\alpha})u_{n\alpha}u_{n}^{\beta}+2(-W_{\beta}^{n})u_{nn}u_{n}^{\beta}\no\\
= & \sigma_{1}({W}^{\top})u_{nn}^{2}+W_{nn}\sum_{\alpha}u_{n\alpha}u_{n}^{\alpha}+\sigma_{1}({W}^{\top})\sum_{\alpha}u_{n\alpha}u_{n}^{\alpha}-W_{\beta}^{\alpha}u_{n\alpha}u_{n}^{\beta}-2W_{\beta}^{n}u_{nn}u_{n}^{\beta}\nonumber \\
= & \sigma_{1}({W}^{\top})u_{nn}^{2}+u_{nn}\sum_{\alpha}u_{n\alpha}u_{n}^{\alpha}-2W_{\beta}^{n}u_{nn}u_{n}^{\beta}\nonumber \\
 & +(u_{n}^{2}-\frac{1}{2}|\nabla u|^{2}+A_{nn})\sum_{\alpha}u_{n\alpha}u_{n}^{\alpha}+\sigma_{1}({W}^{\top})\sum_{\alpha}u_{n\alpha}u_{n}^{\alpha}-W_{\beta}^{\alpha}u_{n\alpha}u_{n}^{\beta}\nonumber \\
= & u_{nn}(\sigma_{1}({W}^{\top})u_{nn}+\sum_{\alpha}u_{n\alpha}u_{n}^{\alpha}-2W_{\beta}^{n}u_{n}^{\beta})\nonumber \\
 & +(u_{n}^{2}-\frac{1}{2}|\nabla u|^{2}+A_{nn})\sum_{\alpha}u_{n\alpha}u_{n}^{\alpha}+\sigma_{1}({W}^{\top})\sum_{\alpha}u_{n\alpha}u_{n}^{\alpha}-W_{\beta}^{\alpha}u_{n\alpha}u_{n}^{\beta}.\label{eq:key2-1}
\end{align}
We need the lower bound of the following
term 
\begin{align}
 & \sigma_{1}({W}^{\top})u_{nn}+u_{n\alpha}u_{n}^{\alpha}-2W_{\beta}^{n}u_{n}^{\beta}\nonumber \\
= & \sigma_{1}({W}^{\top})(W_{nn}-u_{n}^{2}+\frac{1}{2}|\nabla u|^{2}-A_{nn})+u_{n}^{\alpha}(u_{n\alpha}-2{W}_{n\alpha})\nonumber \\
= & \sigma_{1}({W}^{\top})(W_{nn}-u_{n}^{2}+\frac{1}{2}|\nabla u|^{2}-A_{nn})\no\\
 & +(W_{n}^{\alpha}-u_{n}u^{\alpha}-A_{n}^{\alpha})(-W_{n\alpha}-u_{n}u_{\alpha}-A_{n\alpha})\no\\
=& \sigma_{1}({W}^{\top})W_{nn}-\sum_{\alpha}{W}_{n\alpha}W_{n}^{\alpha}+\sigma_{1}({W}^{\top})(-u_{n}^{2}+\frac{1}{2}|\nabla u|^{2}-A_{nn})\no\\
 & +\sum_{\alpha}(u_{n}u_{\alpha}+A_{n\alpha})(u_{n}u^{\alpha}+A_{n}^{\alpha})\nonumber \\
\overset{\eqref{def:sigma_2}}{=}   & f^{2}-\sigma_{2}({W}^{\top})+\sigma_{1}({W}^{\top})(-u_{n}^{2}+\frac{1}{2}|\nabla u|^{2}-A_{nn})+\sum_{\alpha}(u_{n}u_{\alpha}+A_{n\alpha})(u_{n}u^{\alpha}+A_{n}^{\alpha})\no\\
\geq & f^{2}-\sigma_{2}({W}^{\top})-C\sigma_{1}({W}^{\top})+\sum_{\alpha}(u_{n}u_{\alpha}+A_{n\alpha})(u_{n}u^{\alpha}+A_{n}^{\alpha}).\label{est:key1-1}
\end{align}

Putting these facts together, we are able to show that 
\begin{align}
 & F_{j}^{i}(u_{ni}+u^{k}(x_{n})_{ki}-\varphi_{i})(u_{n}^{j}+u^{l}(x_{n})_{l}^{j}-\varphi^{j})\no\label{eq:extra term 1}\\
\ge & u_{nn}\left(f^{2}-\sigma_{2}({W}^{\top})-C_4\sigma_{1}({W}^{\top})\right)-C_4\sum_{\alpha=1}^{n-1}|W_{\alpha}^{\alpha}|.
\end{align}

Also 
\begin{align*}
F_{j}^{\gamma}u^{\beta j}L_{\beta\gamma} & =F_{n}^{\gamma}u^{\beta n}L_{\beta\gamma}+F_{\zeta}^{\gamma}u^{\beta\zeta}L_{\beta\gamma}\\
 & =F_{n}^{\gamma}u^{\beta n}L_{\beta\gamma}+(\sigma_{1}(W)\delta_{\zeta}^{\gamma}-W_{\zeta}^{\gamma})(W^{\beta\zeta}+O(1))L_{\beta\gamma}\\
 & =O(N_{1}^{2})+\sigma_{1}(W|W_{\alpha}^{\alpha})W_{\alpha}^{\alpha}L_{\alpha}^{\alpha}-C\sum_{\alpha=1}^{n-1}|W_{\alpha}^{\alpha}|-CW_{n}^{n}\\
 & =O(N_{1}^{2})+\sigma_{1}(W|W_{\alpha}^{\alpha})W_{\alpha}^{\alpha}L_{\alpha}^{\alpha}-C\sum_{\alpha=1}^{n-1}|W_{\alpha}^{\alpha}|-CW_{n}^{n}.
\end{align*}
and 
\begin{align}\label{eq:extra term2}
 & (1+\frac{7}{5}A^{\frac{2}{5}})F_{j}^{i}\left(u^{\beta j}(x_{n})_{\beta i}+u^{k}(x_{n})_{ki}^{j}+u_{i}^{\beta}(x_{n})_{\beta}^{j}\right)\no\\
\ge & -2(1+\frac{7}{5}A^{\frac{2}{5}})\sigma_{1}(W|W_{\alpha}^{\alpha})W_{\alpha}^{\alpha}L_{\alpha}^{\alpha}-C\sum_{\alpha=1}^{n-1}|W_{\alpha}^{\alpha}|-CW_{n}^{n}-CN_{1}^{2}.
\end{align}

For simplicity, we define
\begin{align*}
\mathcal{I}:= &\frac{14}{25}A^{-\frac{3}{5}}F_{j}^{i}(u_{ni}+u^{k}(x_{n})_{ki}-\varphi_{i})(u_{n}^{j}+u^{l}(x_{n})_{l}^{j}-\varphi^{j})\\
 & +(1+\frac{7}{5}A^{\frac{2}{5}})F_{j}^{i}\left(u^{\beta j}(x_{n})_{\beta i}+u^{k}(x_{n})_{ki}^{j}+u_{i}^{\beta}(x_{n})_{\beta}^{j}\right).
\end{align*}
Then, it follows from (\ref{eq:extra term 1}) and (\ref{eq:extra term2}) that
\begin{align*}
\mathcal{I}\ge & \frac{14}{25}A^{-\frac{3}{5}}u_{nn}\left(f^{2}-\sigma_{2}({W}^{\top})-C_{4}\sigma_{1}({W}^{\top})-\frac{C_4\sum_{\alpha=1}^{n-1}|W_{\alpha}^{\alpha}|}{u_{nn}}\right)\\
 & -2(1+\frac{7}{5}A^{\frac{2}{5}})\sigma_{1}(W|W_{\alpha}^{\alpha})W_{\alpha}^{\alpha}L_{\alpha}^{\alpha}-C\sum_{\alpha=1}^{n-1}|W_{\alpha}^{\alpha}|-CW_{n}^{n}-CN_{1}^{2}.
\end{align*}

We are ready to estimate the remaining term
\begin{align}\label{eq:key4-1}
-\frac{N}{2}F_{j}^{i}(x_{n})_{i}^{j}= & -\frac{N}{2}F_{\beta}^{\alpha}(x_{n})_{\alpha}^{\beta}=\frac{N}{2}F_{\beta}^{\alpha}(L_{\alpha}^{\beta}+O(x_{n}))\no\\
= & \frac{N}{2}(\sigma_{1}({W})\cdot H-{W}_{\beta}^{\alpha}L_{\alpha}^{\beta})+O(x_{n})N\sum_{\alpha,\beta}|F_{\alpha}^{\beta}|\no\\
\ge &  \frac{N}{2}(W_{nn}+\sigma_{1}({W}^{\top}))\cdot(n-1)h_{g}-\frac{N}{2}{W}_{\beta}^{\alpha}L_{\alpha}^{\beta}-C(\sum_{i}{F}_{i}^{i})Nx_{n}.
\end{align}

Therefore, going back to \eqref{est:crucial} and putting these facts together, we conclude that 
\begin{align}\label{eq:final 1}
0\ge & F_{j}^{i}G_{i}^{j}\no\\
\ge & \frac{14}{25}A^{-\frac{3}{5}}\big[u_{nn}\left(f^{2}-\sigma_{2}({W}^{\top})-C_{4}\sigma_{1}({W}^{\top})\right)-C_4\sum_{\alpha=1}^{n-1}|W_{\alpha}^{\alpha}|\big]\nonumber \\
 & -2(1+\frac{7}{5}A^{\frac{2}{5}})\sigma_{1}(W|W_{\alpha}^{\alpha})W_{\alpha}^{\alpha}L_{\alpha}^{\alpha}\nonumber \\
 & +(1+\frac{7}{5}A^{\frac{2}{5}})\left(-C(\sum_{i}{F}_{i}^{i})-C(\sum_{i}{F}_{i}^{i})\frac{N_{1}d_{g}(x_{0},x)}{1+\frac{7}{5}A^{\frac{2}{5}}}\right)\nonumber \\
 & +\frac{n-1}{2}W_{n}^{n}Nu_{n}+\frac{n-1}{2}NW_{n}^{n}h_{g}+\frac{n-1}{2}N\sigma_{1}(W^{\top})h_{g}\no\\
 &+(n-3)(1+\frac{7}{5}A^{\frac{2}{5}})\sigma_{1}({W}^{\top})u_{n}u_{nn}\nonumber \\
 & -\frac{N}{2}{W}_{\beta}^{\alpha}L_{\alpha}^{\beta}-C(\sum_{i}{F}_{i}^{i})Nx_{n}-C\delta N_{1}\sum_{i}F_{i}^{i}-C\sum_{\alpha=1}^{n-1}|W_{\alpha}^{\alpha}|-CW_{n}^{n}-CN_{1}^{2}\nonumber \\
\ge & \frac{14}{25}A^{-\frac{3}{5}}u_{nn}\left(f^{2}-\sigma_{2}({W}^{\top})-C_{4}\sigma_{1}({W}^{\top})-\frac{C_{4}\sum_{\alpha=1}^{n-1}|W_{\alpha}^{\alpha}|}{u_{nn}}\right)\no\\
&+\frac{n-1}{2}NW_{n}^{n}(ce^{u}+A)-C_{1}(\sum_{i}{F}_{i}^{i})Nx_{n}\nonumber \\
 & +\frac{n-1}{2}N\sigma_{1}(W^{\top})h_{g}+(n-3)(1+\frac{7}{5}A^{\frac{2}{5}})\sigma_{1}({W}^{\top})u_{n}u_{nn}\no\\
 & -2(1+\frac{7}{5}A^{\frac{2}{5}})\sigma_{1}(W|W_{\alpha}^{\alpha})W_{\alpha}^{\alpha}L_{\alpha}^{\alpha}-\frac{N}{2}{W}_{\beta}^{\alpha}L_{\alpha}^{\beta}\nonumber \\
 & -C(\sum_{i}{F}_{i}^{i})-C(\sum_{i}{F}_{i}^{i})N_{1}\delta-C\sum_{\alpha=1}^{n-1}|W_{\alpha}^{\alpha}|-CW_{n}^{n}\no\\
 \ge & \frac{14}{25}A^{-\frac{3}{5}}u_{nn}\left(f^{2}-\sigma_{2}({W}^{\top})-C_{4}\sigma_{1}({W}^{\top})-\frac{C_{4}\sum_{\alpha=1}^{n-1}|W_{\alpha}^{\alpha}|}{u_{nn}}\right)+\frac{n-2}{2}NW_{n}^{n}A\nonumber \\
 & +\frac{n-1}{2}N\sigma_{1}(W^{\top})h_{g}+(n-3)(1+\frac{7}{5}A^{\frac{2}{5}})\sigma_{1}({W}^{\top})u_{n}u_{nn}\no\\
 & -2(1+\frac{7}{5}A^{\frac{2}{5}})\sigma_{1}(W|W_{\alpha}^{\alpha})W_{\alpha}^{\alpha}L_{\alpha}^{\alpha}-\frac{N}{2}{W}_{\beta}^{\alpha}L_{\alpha}^{\beta}-C(\sum_{i}{F}_{i}^{i}),
\end{align}
where the last inequality follows from \eqref{eq:lower bound of A},
non-negativity of $c$, \eqref{eq:double tangential initial bound} and a larger $N$.

Let 
\begin{align}
\mathbb{P}:= & \frac{14}{25}A^{-\frac{3}{5}}u_{nn}\left(f^{2}-\sigma_{2}({W}^{\top})-C_{4}\sigma_{1}({W}^{\top})-\frac{C_{4}\sum_{\alpha=1}^{n-1}|W_{\alpha}^{\alpha}|}{u_{nn}}\right)\no\label{expression of P}\\
 & +\frac{n-1}{2}N\sigma_{1}(W^{\top})h_{g}+(n-3)(1+\frac{7}{5}A^{\frac{2}{5}})\sigma_{1}({W}^{\top})u_{n}u_{nn}\nonumber \\
 & -2(1+\frac{7}{5}A^{\frac{2}{5}})\sigma_{1}(W|W_{\alpha}^{\alpha})W_{\alpha}^{\alpha}L_{\alpha}^{\alpha}-\frac{N}{2}{W}_{\beta}^{\alpha}L_{\alpha}^{\beta}-C(\sum_{i}{F}_{i}^{i}),
\end{align}
then \eqref{eq:final 1} becomes 
\begin{equation}
0\ge F_{j}^{i}G_{i}^{j}\ge\mathbb{P}+\frac{n-2}{2}NW_{n}^{n}A.\label{final2}
\end{equation}
A basic fact is needed later:
\begin{equation}
\frac{13}{25}A^{-\frac{3}{5}}W_{n}^{n}f^{2}+\frac{1}{2}NW_{n}^{n}A\ge \frac{1}{2}W_{n}^{n}N^{\frac{3}{8}}f^{\frac{5}{4}}\ge\frac{1}{2}W_{n}^{n}N^{\frac{3}{8}}\inf_{M_{r}}f^{\frac{5}{4}},\label{eq:lower bound of big term}
\end{equation}
where the first inequality follows by regarding the left hand side
as a function of $A$ and seeking its minimum value.

\vskip 4pt
\emph{Case 1:} $\sigma_{2}(W^{\top})(x_{1})>0$.
\vskip 4pt
 This together with \eqref{def:sigma_2} yields 
\[
0<\sigma_{1}(W^{\top})\leq C\frac{N_{1}^{2}}{W_{n}^{n}}\quad\mathrm{~~and~~}\quad\sigma_{2}(W^{\top})\le\frac{(n-2)\sigma_{1}(W^{\top})^{2}}{2(n-1)}\le C\frac{N_{1}^{4}}{(W_{n}^{n})^{2}}.
\]
Actually, now $\lambda(W^{\top})\in\Gamma_{2}^{+}$ implies 
\begin{equation}
|W_{\beta}^{\alpha}|\le\sigma_{1}(W^{\top})\le C\frac{N_{1}^{2}}{W_{n}^{n}}.\label{CaseAtangentialbound}
\end{equation}
With a large $N$, by \eqref{CaseAtangentialbound} and
\eqref{lbd:u_nn} we have 
\begin{equation*}
u_{nn}\big(f^{2}-\sigma_{2}({W}^{\top})-C_{4}\sigma_{1}({W}^{\top})-\frac{C_{4}\sum_{\alpha=1}^{n-1}|W_{\alpha}^{\alpha}|}{u_{nn}}\big)\ge\frac{13}{14}f^{2}W_{nn}.
\end{equation*}
and 
\begin{equation}
\mathbb{P}\ge\frac{13}{25}f^{2}W_{nn}A^{-\frac{3}{5}}-C-C(\sum_{i}{F}_{i}^{i}).\label{CaseA-P}
\end{equation}

By (\ref{final2}), \eqref{CaseA-P} and (\ref{eq:lower bound of big term}),
we obtain 
\begin{align*}
0\ge & F_{j}^{i}G_{i}^{j}\\
\ge & \frac{1}{2} W_{n}^{n}N^{\frac{3}{8}}\inf_{M_r} f^{\frac{5}{4}}-C-C(\sum_{i}{F}_{i}^{i})>0,
\end{align*}
which yields a contradiction by choosing a larger $N$. We emphasize
that such an $N$ depends on $\sup_{M_{r}}\frac{1}{f(x,u)}$.

\vskip 4pt
\emph{Case 2:} $\sigma_{2}(W^{\top})(x_{1})\leq 0$.
\vskip 4pt

 Without loss of generality,
we assume that $W_{1}^{1}\ge W_{2}^{2}\ge\cdots \ge W_{n-1}^{n-1}$. Then
we have $W_{1}^{1}>0$ and $W_{n-1}^{n-1}<0$ due to $\sigma_{1}(W^{\top})(x_{1})>0$
and $\sigma_{2}(W^{\top})(x_{1})\leq 0$.

It follows from $0<\sigma_{1}(W^{\top})=W_{1}^{1}+\cdots+ W_{n-1}^{n-1}$ that
\[
0<-W_{n-1}^{n-1}\le W_{1}^{1}+\cdots+W_{n-2}^{n-2}\le(n-2)W_{1}^{1}
\]
and thus 
\begin{equation}
|W_{\alpha}^{\alpha}|\le(n-2)W_{1}^{1}\qquad\mathrm{~~for~~}1\le\alpha\le n-1.\label{the upper bound of W^T by W11}
\end{equation}

By Lemma \ref{lem:tangential trace estimate}, there exists a positive constant $C^\ast$ such that  $0<\sigma_{1}(W^{\top})\le C^{*}$.

\emph{Case 2.1:} $W_{1}^{1}>C^{**}:=2C^{\ast}+2\sqrt{(C^{\ast})^{2}+C_{4}C^{\ast}+C_{5}}+5\sqrt{C_{6}}A_{1}^{\frac{3}{10}}+25C_{6}A_{1}^{\frac{3}{5}}$,
where $C_{5},C_{6}$ are positive constants below depending on $g$ and
$n$.

By 
\[
0<W_{1}^{1}+\sigma_{1}(W^{\top}|W_{1}^{1})\le C^{*},
\]
we have $-W_{1}^{1}\le\sigma_{1}(W^{\top}|W_{1}^{1})\leq C^{*}-W_{1}^{1}<0$
by virtue of the choice of $C^{\ast\ast}$.

By 
\[
-C^{*}W_{n}^{n}\le f^{2}-\sigma_{1}(W^{\top})W_{n}^{n}+\sum_{\alpha}W_{n}^{\alpha}W_{\alpha}^{n}=\sigma_{2}(W^{\top})<0,
\]
we obtain 
\begin{align}
 & C^{*}W_{1}^{1}-\frac{n-1}{2(n-2)}(W_{1}^{1})^{2}\no\label{polynomialofW11}\\
 & \ge W_{1}^{1}(C^{*}-W_{1}^{1})+\frac{n-3}{2(n-2)}\sigma_{1}(W^{\top}|W_{1}^{1})^{2}\no\\
 & \ge W_{1}^{1}\sigma_{1}(W^{\top}|W_{1}^{1})+\sigma_{2}(W^{\top}|W_{1}^{1})=\sigma_{2}(W^{\top})\ge-C^{*}W_{n}^{n}.
\end{align}
This yields 
\begin{equation}
\sqrt{\frac{2(n-2)}{n-1}}\sqrt{C^{*}W_{n}^{n}+\frac{(C^{*})^{2}(n-2)}{2(n-1)}}+\frac{C^{*}(n-2)}{n-1}\ge W_{1}^{1}>0.\label{eq:subcase B.2 the bound of W^11}
\end{equation}
By \eqref{the upper bound of W^T by W11} we have 
\begin{equation}
|W_{\alpha}^{\alpha}|\le(n-2)\sqrt{\frac{2(n-2)}{n-1}}\sqrt{C^{*}W_{n}^{n}+\frac{(C^{*})^{2}(n-2)}{2(n-1)}}+\frac{C^{*}(n-2)^{2}}{n-1}.\label{tangential W_alphaalpha}
\end{equation}

By \eqref{polynomialofW11} and \eqref{the upper bound of W^T by W11}
we have 
\begin{align}
 & f^{2}-\sigma_{2}({W}^{\top})-C_{4}\sigma_{1}({W}^{\top})-\frac{C_{4}\sum_{\alpha=1}^{n-1}|W_{\alpha}^{\alpha}|}{u_{nn}}\no\label{lower-positive}\\
\ge & f^{2}+\frac{n-1}{2(n-2)}(W_{1}^{1})^{2}-C^{*}W_{1}^{1}-C_{4}C^{\ast}-C_{5}\no\\
\ge & f^{2}+\frac{n-1}{4(n-2)}(W_{1}^{1})^{2},
\end{align}
where the last inequality follows from the choice of $C^{\ast\ast}$.

By \eqref{lower-positive} and \eqref{the upper bound of W^T by W11},
we have 
\begin{align*}
\mathbb{P\ge} & \frac{12}{25}A^{-\frac{3}{5}}W_{nn}\left(f^{2}+\frac{n-1}{4(n-2)}(W_{1}^{1})^{2}\right) -C_{6}W_{1}^{1}W_{nn}-C_{6}(W_{1}^{1})^{2}-C_{6}W_{nn}.
\end{align*}

With a sufficiently large $N$ and the choice of $C^{\ast\ast}$ we arrive
at a contradiction via 
\begin{align}
0\ge & F_{j}^{i}G_{i}^{j}\no\label{eq:final subcase B.2}\\
\ge & \frac{12}{25}A^{-\frac{3}{5}}W_{nn}\left(f^{2}+\frac{n-1}{4(n-2)}(W_{1}^{1})^{2}\right)\nonumber \\
 & -C_{6}W_{1}^{1}W_{nn}-C_{6}(W_{1}^{1})^{2}-C_{6}W_{nn}+\frac{n-2}{2}NW_{n}^{n}A>0.
\end{align}

\emph{Case 2.2:} $W_{1}^{1}\le C^{**}=2C^{\ast}+2\sqrt{(C^{\ast})^{2}+C_{4}C^{\ast}+C_{5}}+5\sqrt{C_{6}}A_{1}^{\frac{3}{10}}+25C_{6}A_{1}^{\frac{3}{5}}$.

It follows from \eqref{the upper bound of W^T by W11} that 
\begin{equation}
|W_{\alpha}^{\alpha}|\le C\quad\text{for}\quad1\le\alpha\le n-1,\label{upper bound of double in the last case}
\end{equation}
and then 
\[
\sigma_{1}(W^{\top})=\frac{f^{2}-\sigma_{2}(W^{\top})+\sum_{\alpha}W_{n}^{\alpha}W_{\alpha}^{n}}{W_{n}^{n}}\le\frac{C}{W_{n}^{n}}.
\]
By choosing $N$ large enough we have
\begin{align*}
 & u_{nn}\big(f^{2}-\sigma_{2}({W}^{\top})-C_{4}\sigma_{1}({W}^{\top})-\frac{C_{4}\sum_{\alpha=1}^{n-1}|W_{\alpha}^{\alpha}|}{u_{nn}}\big)\\
\geq & u_{nn}\big(f^{2}-\frac{C_{7}}{W_{n}^{n}}\big)\geq\frac{13}{14}W_{n}^{n}\inf_{M_r} f^{2}
\end{align*}
and then 
\begin{equation}
\mathbb{P}\ge\frac{13}{25}f^{2}W_{nn}A^{-\frac{3}{5}}-C-C(\sum_{i}{F}_{i}^{i}).\label{CaseB-P lastsubcase}
\end{equation}
By (\ref{final2}), (\ref{CaseB-P lastsubcase}) and (\ref{eq:lower bound of big term})
we know that 
\begin{align*}
0 & \ge F_{j}^{i}G_{i}^{j}\\
 & \ge\frac{13}{25}W_{n}^{n}\inf_{M_{r}}f^{2} A^{-\frac{3}{5}}+\frac{n-2}{2}NW_{n}^{n}A-C(\sum_{i}{F}_{i}^{i})\\
 &\ge \frac{1}{2}W_{n}^{n}N^{\frac{3}{8}}\inf_{M_{r}}f^{\frac{5}{4}}-C(\sum_{i}{F}_{i}^{i})>0.
\end{align*}
This is a contradiction.

Therefore, we establish  that $x_0 \in \pa M$ is the maximum point of $G$ in $\overline{M_{r}}$. 
\end{proof}
	
	Since $\lambda(W)\in\Gamma_{2}^{+}$, by Lemma
	\ref{lem:tangential trace estimate} we have 
	\begin{equation}
		|u_{ij}|\le C(\Delta u+1)\le C(u_{nn}+1).\label{upper_bound_delta_u}
	\end{equation}
	and 
	\[
	-C\leq \sigma_1(W)-C\leq\Delta u=\sum_{\alpha=1}^{n-1}u_{~~\alpha}^{\alpha}+u_{nn}.
	\]
	Again by Lemma \ref{lem:tangential trace estimate} we see that $u_{nn}$ is bounded below. 

	The following corollary is a direct consequence of Lemma \ref{lem:double_normal_derivatives} and \eqref{upper_bound_delta_u}.
	
	\begin{corollary}\label{Delta boundary estimates}
		Suppose $c\ge 0$ on $\partial M$. Let $u$ be a solution to \eqref{PDE:sigma_2_with_bdry}.
		Then for any open sets $\mathcal{O'}\Subset \mathcal{O} \subset  \partial M$,  $$\sup_{\mathcal{O'}}\Delta u\le C, $$
		where $C$ depends on $n,g$, $\sup_{\mathcal{O}}((c+|\tilde \nabla c|+|\tilde \nabla^2 c|)e^{-u}),\sup_{\mathcal{O}} (f+|f_x|+|f_{xz}|+|f_{xx}|+|f_{zz}|)$,$\mathcal{O},\mathcal{O}'$, $\|c\|_{C^2(\overline{\mathcal{O}})}$ and $\sup_{\mathcal{O}} \frac{1}{f(x,u)}$.
	\end{corollary}
	
	We are now in a position to complete the proof of Theorem \ref{Thm:local estimates}.	
\begin{proof}[Proof of Theorem \protect\ref{Thm:local estimates}]
		Fix $P \in \pa M$. Under Fermi coordinates around $P$ we consider the function 
		\[F(x)=\eta(\Delta u+|\nabla u|^2), \qquad \mathrm{~~in~~} \quad B_r^+,\]
		where $\eta$ is the same cut-off
		function in $B_r^+$ used in the proof of Lemma \ref{lem:tangential trace estimate}. The maximum point $x_0$ of $F$ happens either on the boundary or in the interior.
		
		When $x_0 \in M$, the interior $C^2$-estimates for the $\sigma_k$-curvature equation on closed manifolds yields
		\[F(x)\le F(x_0)\le  C\bigg(1+\frac{1}{r^2}+(f+|\nabla f|+|\nabla^2 f|)\sup_{B_r^+}e^{-2u}\bigg).\]
		
		When $x_0 \in \pa M$, we know from  Corollary \ref{Delta boundary estimates} that $F(x)\le F(x_0)\le C$ and thus complete the proof of Theorem \ref{Thm:local estimates}.
	\end{proof}

\subsection{Counter-examples to boundary local $C^2$ estimates}\label{subsect:counter-examples}
	As mentioned in the introduction, the boundary local $C^2$ estimates for  \eqref{PDE:prescribing_sigma_2} fail in general if $c$ is allowed to be negative somewhere. To be consistent with our notations, we first state Li-Nguyen's counter-examples in \cite{Li-Luc4}.
\begin{proposition}\label{prop:counterexample_Li-Luc}
For $2\leq k\leq n$ and any  constant $c<0$, there exist constants $C_{0}>0,R_{0}>1$
and a sequence of functions $\left\{ u_{j}\right\} \subset C^{\infty}\left(\overline{B_{R_0}}\backslash B_{1}\right)$
satisfying 
\[
\begin{cases}
u_j^{-\frac{4}{n-2}}\sigma_{k}^{1/k}\big(\lambda\big(A_{u_j^{4/(n-2)}|\ud x|^2}\big)\big)=1&\mathrm{~~in~~}B_{R_{0}}\backslash \overline{B_1},\\
u_{j}>0\mathrm{~~and~~}\lambda\big(A_{u_j^{4/(n-2)}|\ud x|^2}\big)\in\Gamma_{k}^+ &\mathrm{~~ in~~}\overline{B_{R_0}}\backslash B_{1},\\
\frac{\partial u_{j}}{\partial r}+\frac{n-2}{2}u_{j}=-\frac{n-2}{2}cu_{j}^{\frac{n}{n-2}} &\mathrm{~~on~~}\partial B_{1},
\end{cases}
\]
 and 
\[
\left|u_{j}\right|+|u_{j}^{-1}|+|\nabla u_{j}|\leq C_{0}\qquad \mathrm{~~in~~} \overline{B_{R_0}}\backslash B_{1},
\]
 such that 
\[
\lim_{j\rightarrow\infty}\inf_{\partial B_{1}}\left|\nabla^{2}u_{j}\right|=\infty.
\]
\end{proposition}

Moreover, the functions $u_j$ are singular radial solutions to the $\sigma_k$-Yamabe equation on annuli, which had been classified by Chang-Hang-Yang \cite{CHY}.

We first need some basic facts about the set of possible umbilic points for ellipsoids, which are of independent interest. For an ellipsoid in $\Rn$ with pairwise distinct radii, it has exactly four umbilic points when $n=3$; however as we shall see, it is a totally non-umbilic hypersurface when $n \geq 4$.
\begin{proposition}\label{prop:ellipsoid}
For $n \geq 3$, let 
\begin{equation}\label{def:ellipsoid}
\Sigma=\left\{x \in \Rn;~~ \sum_{i=1}^{n}\frac{x_{i}^{2}}{a_{i}^{2}}=1\right\}, \qquad a_i \in \R_+ \mathrm{~~for~~} 1 \leq i \leq n
\end{equation}
be an ellipsoid. Then the mean curvature is 
$$H=\frac{1}{\sqrt{\sum\limits_{i=1}^n a_i^{-4} x_i^2}}\left(\sum\limits_{i=1}^n a_i^{-2}-\frac{\sum\limits_{i=1}^n a_i^{-6} x_i^2 }{\sum\limits_{i=1}^n a_i^{-4} x_i^2}\right).$$
In addition, assume $a_1<a_2<\cdots<a_n$, then
\begin{itemize}
\item [(i)] For $n=3$, $\Sigma$ has exactly four umbilic points, which are
$$\left\{x \in \R^3; x_1=\pm a_1 \sqrt{\frac{a_2^2-a_1^2}{a_3^2-a_1^2}},x_2=0,x_3=\pm a_3\sqrt{\frac{a_3^2-a_2^2}{a_3^2-a_1^2}}\right\}.$$

\item [(ii)] For $n \geq 4$, $\Sigma$ is totally non-umbilic.
\end{itemize}
\end{proposition}

For a general ellipsoid, we give a sufficient and necessary condition of its umbilic point (if exists).
\begin{proposition}\label{prop:umbilic_pts_ellipsoid}
For an ellipsoid $\Sigma$ as in \eqref{def:ellipsoid} with generic $a_i \in \R_+$, then $p=(x_1,\cdots,x_n) \in \Sigma$ is umbilic if and only if the following  holds:
\begin{align}\label{eqn:umbilic_generic}
&\sum\limits_{i=1}^n a_i^{-4}-2\frac{\sum\limits_{i=1}^n a_i^{-8}x_i^2}{\sum\limits_{i=1}^n a_i^{-4} x_i^2}+\frac{(\sum\limits_{i=1}^n a_i^{-6}x_i^2)^2}{(\sum\limits_{i=1}^n  a_i^{-4} x_i^2)^2}\no\\
=&\frac{1}{n-1}\left(\sum\limits_{i=1}^n a_i^{-2}-\frac{\sum\limits_{i=1}^n a_i^{-6} x_i^2 }{\sum\limits_{i=1}^n a_i^{-4} x_i^2}\right)^2.
\end{align}
\end{proposition}

The proof of Propositions \ref{prop:ellipsoid} and \ref{prop:umbilic_pts_ellipsoid} involves the standard techniques in differential geometry.  For clarify, we leave it to Appendix \ref{Appendix:A}.

\medskip

We consider the above sequence $\{u_j\}$ in Li-Nguyen's counterexamples on $B_{R_0}\backslash \overline{\Omega_\ve}$, where 
$$\Omega_\ve=\left\{x \in \Rn;\sum_{i=1}^{n}\frac{x_{i}^{2}}{a_{i}^{2}}< 1\right\}$$
with $a_{i}=1+(n-i)\varepsilon$ for $1 \leq i \leq n$ and some $0<\varepsilon \ll1$.
 Clearly, $\Sigma_\ve:=\pa \Omega_\ve$ touches $\pa B_1$
at two points $x_{0}^{\pm}=(0,\cdots,0,\pm1)$ such that $B_{R_{0}}\supset\Omega_\ve\supset B_{1}$. Moreover, by Proposition \ref{prop:ellipsoid} we know that  the mean curvature $h_\ve$ on $\Sigma_\ve \subset \pa(B_{R_0}\backslash \overline\Omega_\ve)$ is 
$$h_\ve=-\frac{1}{n-1}\frac{1}{\sqrt{\sum\limits_{i=1}^n a_i^{-4} x_i^2}}\left(\sum\limits_{i=1}^n a_i^{-2}-\frac{\sum\limits_{i=1}^n a_i^{-6} x_i^2 }{\sum\limits_{i=1}^n a_i^{-4} x_i^2}\right) \to -1,  \qquad \mathrm{as~~} \ve \to 0$$
and $x_0^{\pm}$ are  non-umbilic points.

We claim that on $\partial\Omega_\ve \subset \pa(B_{R_0}\backslash \overline{\Omega_\ve})$, there exists a small enough $\ve$ such that $c_j(x_0^\pm)<0$ and
\[
\frac{\partial u_{j}}{\partial \vec n}-\frac{n-2}{2}h_{\ve}u_{j}=-\frac{n-2}{2}c_ju_{j}^{\frac{n}{n-2}},
\]
  where $\vec n=(n_1,\cdots,n_n)$ is the inward unit normal vector on $\partial(B_{R_0}\backslash \overline\Omega_\ve)$ and $\vec n(x_{0}^{\pm})=x_0^{\pm}$.

To this end, notice that at $x_{0}^{\pm}$, $\frac{\partial u_{j}}{\partial \vec n}=\frac{\partial u_{j}}{\partial r}$
and $\nabla^2 u_{j}(\vec n,\vec n)=\partial_r^{2}u_{j}$, furthermore
\begin{align*}
\frac{\partial u_{j}}{\partial \vec n}& =\frac{\partial u_{j}}{\partial r}=-\frac{n-2}{2}u_{j}-\frac{n-2}{2}cu_{j}^{\frac{n}{n-2}}\\
 & =-\frac{n-2}{2}\left[u_{j}(h_\ve+1)+cu_{j}^{\frac{n}{n-2}}\right]+\frac{n-2}{2}h_\ve u_{j}.
\end{align*}

For the constant $c<0$ and 
\[
\left|u_{j}\right|+|u_{j}^{-1}|+\left|\nabla u_{j}\right|\leq C_{0} \qquad \mathrm{~~in~~} B_{R_{0}}\backslash \overline{B_1},
\]
 we can choose sufficiently small $\varepsilon$  such that 
\[
u_{j}(h_\ve+1)+cu_{j}^{\frac{n}{n-2}}<0.
\]
This implies that the prescribed mean curvature $c_j(x_0^{\pm})<0$ for some sufficiently small $\ve$. Meanwhile, it follows from Proposition \ref{prop:counterexample_Li-Luc} that $\lim_{j \to \infty}|\nabla^{2}u_{j}(\vec n,\vec n)|(x_{0}^{\pm})=\infty$.

Let $\xi_{k}=\frac{\pa}{\pa {x_k}}-n_k \vec n, 1 \leq k \leq n$ denote
the vector fields on $\Sigma_{\ve}$. Take $x_0^+$ for example. In particular, $T_{x_0^+}\Sigma_{\ve}=\mathrm{span}\{\xi_\alpha; 1 \leq \alpha \leq n-1\}$, then given $\varphi \in C^2(\Sigma_{\ve})$, a direct computation yields
\begin{align*}
\xi_{\alpha}\xi_{\beta} \varphi=\pa_\alpha\pa_\beta\varphi-\frac{\partial n_\beta}{\partial x_{\alpha}}\frac{\pa \varphi}{\pa \vec n} \qquad \mathrm{~~at~~} x_0^+.
\end{align*}
Similarly, let $\zeta_{k}=(\delta_{kl}-x_{k} x_{l})\frac{\pa}{\partial x_l}, 1 \leq k\leq n$ denote
the vector fields on $\pa B_1$, then given $\varphi \in C^2(\pa B_1)$, at $x_0^+$ we have
\begin{align*}
\zeta_{\alpha}\zeta_{\beta} \varphi=\pa_\alpha\pa_\beta\varphi-\delta_{\alpha \beta}\frac{\pa \varphi}{\pa r}.
\end{align*}

We now start to compute 
$$\tilde  \nabla^2 c_j(\xi_\alpha,\xi_\beta) \qquad \mathrm{~~at~~}\quad x_0^+.$$

To continue, we need some observations:
\begin{itemize}
\item At $x_0^+$,  $\xi_{\alpha}(\pa_r u_j)=\zeta_{\alpha}(\pa_r u_j)$ for all  $\alpha$, and
\begin{align*}
\tilde  \nabla^2 c_j(\xi_\alpha,\xi_\beta)=\xi_\alpha \xi_\beta c_j-\tilde \nabla_{\xi_\alpha}\xi_\beta c_j=\xi_\alpha \xi_\beta c_j,
\end{align*}
where the last identity follows from
\begin{align*}
\tilde \nabla_{\xi_\alpha}\xi_\beta=&(\nabla_{\xi_\alpha}\xi_\beta)^{\top}=\big(\nabla_{\xi_\alpha}(\pa_\beta-n_\beta \vec n)\big)^{\top}=-n_\beta \nabla_{\xi_\alpha}\vec n=0.
\end{align*}
Here $\top$ denotes the orthogonal projection on $T_{x_0^+} \Sigma_\ve$.

\item For a smooth radial function $v$, on $\Sigma_\ve$ we have 
\[
\frac{\partial v}{\pa \vec n}=\frac{\sum_{i}\frac{x_{i}}{a_{i}^{2}}\frac{\partial v}{\partial x_{i}}}{\sqrt{\sum_{i=1}^{n}\frac{x_{i}^{2}}{a_{i}^{4}}}}=\frac{\frac{1}{r}v_{r}}{\sqrt{\sum_{i=1}^{n}\frac{x_{i}^{2}}{a_{i}^{4}}}}.
\]

\end{itemize}

Then,  at $x_0^+$ we have

\begin{align*}
\frac{n-2}{2}\tilde \nabla^2 c_j(\xi_\alpha,\xi_\beta) =&\frac{n-2}{2}\xi_\alpha \xi_\beta c_j\\
 =&\xi_{\alpha}\xi_{\beta}\left[u_j^{-\frac{n}{n-2}}\big(-\frac{\partial u_{j}}{\partial \vec n}+\frac{n-2}{2}h_{\ve} u_j\big)\right]\\
=&
-u_j^{-\frac{n}{n-2}}\frac{\xi_{\alpha}\xi_{\beta}(\pa_r u_{j})}{\sqrt{\sum_{i=1}^{n}\frac{x_{i}^{2}}{a_{i}^{4}}}}+O(1)\nonumber \\
=&-u_j^{-\frac{n}{n-2}}\left(\pa_\alpha\pa_\beta(\pa_r u_{j})-\frac{\partial n_{\beta}}{\partial x_{\alpha}}\right)(u_j)_{rr}+O(1).
\end{align*}
Meanwhile, on $\partial B_{1}$, $\pa_r u_{j}=\frac{n-2}{2}u_j-\frac{n-2}{2}cu_j^{\frac{n}{n-2}}$. Then at $x_0^+$ we have 
\[
O(1)=\zeta_{\alpha}\zeta_{\beta}(\frac{n-2}{2}u_j-\frac{n-2}{2}cu_j^{\frac{n}{n-2}})=\zeta_{\alpha}\zeta_{\beta}(\pa_r u_{j})=\pa_\alpha\pa_\beta(\pa_r u_{j})-\delta_{\alpha \beta}(u_j)_{rr},
\]
whence,
\begin{equation*}
\pa_\alpha\pa_\beta(\pa_r u_{j})=\delta_{\alpha \beta}(u_j)_{rr}+O(1).\label{eq:third derivaties}
\end{equation*}
Hence, we conclude that at $x_0^+$,
\begin{align*}
\frac{n-2}{2}\tilde \nabla^2 c_j(\xi_\alpha,\xi_\beta) =&-u_j^{-\frac{n}{n-2}}(\delta_{\alpha \beta}-\frac{\partial n_\beta}{\partial x_{\alpha}})(u_j)_{rr}+O(1)\\
 =&-u_j^{-\frac{n}{n-2}}\delta_{\alpha \beta}\left(-\frac{1}{a_{\beta}^{2}}+1\right)(u_j)_{rr}+O(1)\\
=&-u_j^{-\frac{n}{n-2}}\frac{2(n-\beta)\varepsilon+(n-\beta)^{2}\varepsilon^{2}}{a_{\beta}^{2}}\delta_{\alpha \beta}(u_j)_{rr}+O(1).
\end{align*}
Based on the above estimate, we can choose a smaller $\ve_j=O(\frac{1}{| \nabla^2 u_j(\vec n,\vec n)(x_0^+)|})$  such that 
$$|\tilde \nabla^2 c_j(x_0^+)|\leq 1.$$

This produces counter-examples for boundary local $C^2$ estimates while $c_j$ is negative somewhere.
\begin{proposition}\label{prop:counterexample_new}
For any $2\leq k\leq n$ and $\ve>0$, let $\Omega_\ve=\{x \in \Rn; \sum_{i=1}^{n}x_{i}^{2}/a_{i}^{2}<1\}$
with $a_{i}=1+(n-i)\varepsilon$ for $1 \leq i \leq n$, and $x_0^{\pm}=( 0, \cdots,0,\pm 1) \in \pa \Omega_\ve$. Then there exist $C_{0}>0,R_{0}>1$,  sequences of positive constants $\{\ve_j\}$,  smooth functions $\{c_j\}$ on $ \pa \Omega_{\ve_j}$ and $\left\{ u_{j}\right\} \subset C^{\infty}\left(\overline{B_{R_0}}\backslash\Omega_{\ve_j}\right)$
satisfying 
\[
\begin{cases}
u_j^{-\frac{4}{n-2}}\sigma_{k}^{1/k}\big(\lambda\big(A_{u_j^{4/(n-2)}|\ud x|^2}\big)\big)=1 &\quad\mathrm{~~in~~}B_{R_{0}}\backslash \overline{\Omega_{\ve_j}},\\
u_{j}>0\mathrm{~~and~~}\lambda\big(A_{u_j^{4/(n-2)}|\ud x|^2}\big)\in\Gamma_{k}^+ &\quad \mathrm{~~in~~}\overline{B_{R_0}}\backslash\Omega_{\ve_j},\\
\frac{\partial u_{j}}{\partial \vec n}-\frac{n-2}{2}h_{\ve_j}u_{j}=-\frac{n-2}{2}c_j u_{j}^{\frac{n}{n-2}} &\quad \mathrm{~~on~~} \pa \Omega_{\ve_j},
\end{cases}
\]
 with the following properties:
 \begin{itemize}
\item [(1)]
$\left|u_{j}\right|+|u_{j}^{-1}|+\left|\nabla u_{j}\right|\leq C_{0}$ in $\overline{B_{R_0}}\backslash\Omega_{\ve_j}$;

\item [(2)] $\lim_{j\rightarrow\infty}\left|\nabla^{2}u_{j}(\vec n,\vec n)(x_0^{\pm})\right|=\infty$ and $x_0^{\pm}$ are non-umbilic points on each $\pa \Omega_{\ve_j}$;
\item [(3)] $\ve_j\leq \frac{C_0}{| \nabla^2 u_j(\vec n,\vec n)(x_0^\pm)|}$;
\item [(4)] $c_j(x_0^{\pm})<0, |\tilde \nabla c_j|\leq C_0, |\tilde \nabla^2 c_j(x_0^\pm)|\leq 1$.
\end{itemize}
Here $h_{\ve_j}$ is the mean curvature. 
\end{proposition}

	\section{First boundary $\sigma_{2}$-eigenvalue problem}\label{Sect:first_bdry_eigenvalue}

	One of two types of the first eigenvalue problem associated with boundary
	Yamabe problem is to find positive solutions of 
	\[
	\begin{cases}
		-\Delta_{g}u+\frac{n-2}{4(n-1)}R_{g}u=\lambda_{1}u & \qquad\mathrm{~~in~~}M,\\
		-\frac{\pa u}{\pa\vec{n}_{g}}+\frac{n-2}{2}h_{g}u=0 & \qquad\mathrm{~~on~~}\pa M,
	\end{cases}
	\]
	where $\lambda_1 \in \R$. This is equivalent to 
	\[
	\begin{cases}
		\sigma_{1}(A_{u^{4/(n-2)}g})=\frac{n-2}{2}\lambda_{1} & \qquad\mathrm{~~in~~}M,\\
		h_{u^{4/(n-2)}g}=0 & \qquad\mathrm{~~on~~}\pa M.
	\end{cases}
	\]
	As was shown in \cite{escobar4}, the sign of $\lambda_{1}$
	coincides with the sign of boundary Yamabe constant $Y(M,\pa M,[g])$ in \eqref{def:bdry_Yamabe_constant}.
	
	The above eigenvalue problem motivates the
	first boundary $\sigma_{2}$-eigenvalue problem under constraint
	that $\lambda(A_{g})\in\Gamma_{2}^{+}$ and $h_{g}\ge0$ on $\pa M$. We would like to mention that the first eigenvalue problem for fully nonlinear PDEs on closed manifolds with constraint of positive cone  has appeared in Ge-Lin-Wang \cite{Ge-Lin-Wang}, Guan-Wang \cite[Section $4$ in the arXiv version]{Guan-Wang}, Ge-Wang \cite{Ge-Wang1}, Duncan-Nguyen \cite{Duncan-Luc} etc.
	
	\begin{proposition}\label{prop:special_g} Let $(M,g)$ be a smooth compact
		manifold of dimension $n\geq3$ with boundary. Suppose $\lambda(A_{g})\in\Gamma_{2}^{+}$
		and $h_{g}\ge0$ on $\pa M$. Then there exists a solution $u\in C^{\infty}(\overline{M})$
		to 
		\[
		\begin{cases}
			\sigma_{2}(A_{u})=f(x)e^{u} & \qquad\mathrm{~~in~~}M,\\
			h_{g_{u}}=0 & \qquad\mathrm{~~on~~}\partial M,
		\end{cases}
		\]
		where $f(x)=\sigma_{2}(A_{g})>0$ and $g_{u}=e^{-2u}g$. 
	\end{proposition} 
	\begin{proof}
		Consider a family of solutions to 
		\[
		\begin{cases}
			\sigma_{2}(A_{u})=f(x)e^{u} & \qquad\mathrm{~~in~~}M,\\
			\frac{\partial u}{\partial\vec{n}}=-th_{g} & \qquad\mathrm{~~on~~}\partial M.
		\end{cases}
		\]
		
		At $t=0$, $u=0$ is the only solution and the linearized operator
		is elliptic and invertible. As the same reason above, we only need
		to show the closeness by a priori estimates and the proof follows
		by the continuity method.
		
		Choose a smooth function $l$ such that $\int_{M}l\ud\mu_{g}=\int_{\partial M}(h_{g}+1)\ud\sigma_{g}$
		and find a unique (up to a constant) smooth solution $v$  to the following PDE with Neumann boundary condition
		\[
		\begin{cases}
			\Delta v=tl & \qquad\mathrm{~~in~~}M,\\
			\frac{\partial v}{\partial\vec{n}}=-th_{g}-t & \qquad\mathrm{~~on~~}\partial M.
		\end{cases}
		\]
		
		We first observe that $\min_{\overline{M}}(v-u)$ can not happen on
		boundary $\partial M$ due to $\frac{\partial(v-u)}{\partial\vec{n}}=-t$
		on $\partial M$ for $t>0.$ Thus, we may assume $\min_{\overline{M}}(v-u)=(v-u)(x_{0})$
		for some $x_{0}\in M$ and $\nabla^{2}(v-u)(x_{0})\ge0$ and $\nabla v(x_{0})=\nabla u(x_{0}).$
		Then at $x_{0}$ we have 
		\[
		\big(\Delta v+\sigma_{1}(A_{g})\big)^{2}\frac{\binom{n}{2}}{n^{2}}=\sigma_{1}^{2}(A_{v})\frac{\binom{n}{2}}{n^{2}}\ge\sigma_{2}(A_{v})\ge\sigma_{2}(A_{u})=fe^{u}
		\]
		whence, 
		\[
		e^{u(x_{0})}\le C,
		\]
		where $C$ is a constant independent of $t$.
		
		Notice that 
		\[
		(v-u)(x)\ge(v-u)(x_{0})
		\]
		implies 
		\[
		u(x)\le u(x_{0})+v(x)-v(x_{0}).
		\]
		
		If $h_{g}\ge0$, then we let $\min_{\overline{M}}u=u(x_{1})$ for
		some $x_{1}\in M$ and for $t>0$, $\nabla u(x_{1})=0$ and $\nabla^{2}u(x_{1})\ge0.$
		Then at $x_{1}$, we have
		
		\[
		fe^{u}=\sigma_{2}(A_{u})\ge\sigma_{2}(A_{g}),
		\]
		which shows 
		\[
		u(x_{1})\ge\log[\sigma_{2}(A_{g})/\sigma_{2}(A_{g})]=0.
		\]
		
		Consequently, $u$ is uniformly bounded in $t\in[0,1]$. As before,
		we can obtain uniform a priori $C^{1}$- and $C^{2}$-estimates in
		$t\in[0,1]$ and thus $u\in C^{\infty}(\overline{M})$. 
	\end{proof}
	\begin{proposition}\label{prop:sigma_2_1st_bdry_eigenvalue} Let
		$(M,g)$ be a smooth compact manifold of dimension $n\geq3$ with
		boundary. Suppose $\lambda(A_{g})\in\Gamma_{2}^{+}$ and $h_{g}\ge0$ on $\pa M$.
		Then there exists a unique constant $\lambda>0$ and $u\in C^{\infty}(\overline{M})$
		such that $g_{u}=e^{-2u}g$ satisfying 
		\[
		\begin{cases}
			\sigma_{2}(A_{u})=\lambda & \quad\mathrm{~~in~~}M\\
			h_{g_{u}}=0 & \quad\mathrm{~~on~~}\partial M.
		\end{cases}
		\]
	\end{proposition} 
	\begin{proof}
		A direct consequence of Proposition \ref{prop:special_g} is that there exists
		a conformal metric, still denoted by $g$, such that $\lambda(A_{g})\in\Gamma_{2}^{+}$
		and $h_{g}=0$ on $\pa M$.
		
		We first prove that for every $\varepsilon>0$, there exists a solution
		$u$ such that for 
		\[
		\begin{cases}
			\sigma_{2}^{1/2}(A_{u})=e^{\varepsilon u} & \quad\mathrm{~~in~~}M,\\
			\frac{\partial u}{\partial\vec{n}}=0 & \quad\mathrm{~~on~~}\partial M.
		\end{cases}
		\]
		
		Consider the following equation 
		\begin{equation}
			\begin{cases}
				\sigma_{2}^{1/2}(A_{u})=\left(t+(1-t)f(x)\right)e^{\varepsilon u} & \quad\mathrm{~~in~~}M,\\
				\frac{\partial u}{\partial\vec{n}}=0 & \quad\mathrm{~~on~~}\partial M,
			\end{cases}\label{eq:continuity of equation}
		\end{equation}
		where $t\in[0,1]$ and $f(x)=\sigma_{2}^{1/2}(A_{g}).$
		
		Denote 
		\[
		I=\{t\in[0,1];\mathrm{~~there~~exists~~a~~solution~~}u\mathrm{~~satisfying~~(\ref{eq:continuity of equation})~~with~~}\lambda(A_{g_{u}})\in\Gamma_{2}^{+}\}.
		\]
		
		At $t=0$, $u=0$ is the unique solution. Actually, by the classical
		elliptic theory, the linearized operator of $\sigma_{2}^{1/2}(A_{u})-[t+(1-t)f(x)]e^{\varepsilon u}$
		at $u=0$ is uniformly elliptic when $t\in I$ and the openness follows.
		To ensure $I=[0,1]$, we only need to obtain the uniform $C^{2}$-estimate
		independent of $t\in[0,1].$
		
		Suppose $\max_{\overline{M}}u=u(x_{1})$, then 
		\[
		\left(t+(1-t)f(x)\right)e^{\varepsilon u}=\sigma_{2}^{1/2}(g^{-1}A_{u})(x_{1})\le\sigma_{2}^{1/2}(A_{g})(x_{1})=f(x_{1}),
		\]
		thereby, 
		\[
		\varepsilon u(x_{1})\le t\log f(x_{1})\le \frac{1}{2}|\log\max_{\overline{M}}\sigma_{2}(A_{g})|.
		\]
		Similarly, we have
		
		\[
		\varepsilon\min_{\overline{M}}u\ge\log\frac{\min_{\overline{M}}\sigma_{2}^{\frac 1 2}(A_{g})}{1+\min_{\overline{M}}\sigma_{2}^{\frac 1 2}(A_{g})}.
		\]
		Thus, we arrive at 
		\begin{equation}
			\log\frac{\min_{\overline{M}}\sigma_{2}^{\frac 1 2}(A_{g})}{1+\min_{\overline{M}}\sigma_{2}^{\frac 1 2}(A_{g})}\le\varepsilon u\le |\log\max_{\overline{M}}\sigma_{2}(A_{g})|.\label{eq:bounded co}
		\end{equation}
		
		Thus by the a priori $C^{1}$- and $C^{2}$-estimates we have 
		\begin{align}
			& \|\nabla u\|_{C^{0}(\overline{M})}\no\\
			\le & C\left(1+\sup_{M}(t+(1-t)f)e^{\varepsilon u}+\sup_{M}\left((t+(1-t)f)\varepsilon+(1-t)|\nabla f|\right)e^{\varepsilon u}\right)\label{eq:c1 estimate}
		\end{align}
		and 
		\begin{align}
			& \|\nabla^{2}u\|_{C^{0}(\overline{M})}\no\\
			\le & C\bigg(1+\frac{1}{\inf_{M}(t+(1-t)f)e^{\varepsilon u}}+\sup_{M}(t+(1-t)f)e^{\varepsilon u}\no\\
			& \qquad+\sup_{M}\left((t+(1-t)f)\varepsilon+(1-t)|\nabla f|\right)e^{\varepsilon u}\no\\
			& \qquad+\sup_{M}((t+(1-t)f)\varepsilon^{2}+(1-t)|\nabla^{2}f|)e^{\varepsilon u}\bigg).\label{eq:c2 estimate}
		\end{align}
		Therefore, we have 
		\[
		\|\nabla u\|_{C^{0}(\overline{M})}+\|\nabla^{2}u\|_{C^{0}(\overline{M})}\le C,
		\]
		where $C$ is independent of $t$ and $\ve$.
		
		By the continuity method, we know that for every $\varepsilon>0$,
		there exists a solution $u$ to 
		\[
		\begin{cases}
			\sigma_{2}(g^{-1}A_{u})=e^{\varepsilon u} & \qquad\mathrm{~~in~~}M,\\
			\frac{\partial u}{\partial\vec{n}}=0 & \qquad\mathrm{~~on~~}\partial M,
		\end{cases}
		\]
		satisfying 
		\[
		\|\nabla u\|_{C^{0}(\overline{M})}+\|\nabla^{2}u\|_{C^{0}(\overline{M})}\le C,
		\]
		where $C$ is independent of $\varepsilon.$
		
		Rewrite the above equation as 
		\[
		\begin{cases}
			\sigma_{2}(A_{u})=e^{\varepsilon(u-\bar{u})+\varepsilon\bar{u}} & \qquad\mathrm{~~in~~}M,\\
			\frac{\partial u}{\partial\vec{n}}=0 & \qquad\mathrm{~~on~~}\partial M,
		\end{cases}
		\]
		where $\bar{u}=\fint_{M}u\ud\mu_{g}$. Then, for any sequence $\{\ve_{i}\}$
		with $\ve_{i}\to0$, there holds $\|\varepsilon_{i}\nabla u\|_{C^{0}(\overline{M})}\rightarrow0$
		and thereby up to a subsequence, $\ve_{i}u=\varepsilon_{i}(u-\bar{u})+\varepsilon_{i}\bar{u}\rightarrow\lambda\in\R.$
		
		Let $v_{\varepsilon_{i}}=u-\bar{u}$ satisfy 
		\[
		\begin{cases}
			\sigma_{2}(A_{v_{\varepsilon}})=e^{\varepsilon_{i}v_{\varepsilon_{i}}+\varepsilon_{i}\bar{u}} & \qquad\mathrm{~~in~~}M,\\
			\frac{\partial v_{\varepsilon_{i}}}{\partial\vec{n}}=0 & \qquad\mathrm{~~on~~}\partial M,
		\end{cases}
		\]
		by virtue of the fact that $A_{u}=A_{v_{\ve_{i}}+\bar{u}}=A_{v_{\ve_{i}}}$,
		and $\|v_{\varepsilon_{i}}\|_{C^{2}(\overline{M})}\le C$. Thus for
		any $\alpha\in(0,1)$, there holds 
		\[
		\|v_{\varepsilon_{i}}\|_{C^{2,\alpha}(\overline{M})}\le C
		\]
		and for $\gamma\in(0,\alpha)$, as $\ve_{i}\to0$ we have 
		\[
		v_{\varepsilon_{i}}\rightarrow v\qquad\mathrm{~~in~~}C^{2,\gamma}(\overline{M}),
		\]
		where $v$ satisfies
		\[
		\begin{cases}
			\sigma_{2}(A_{v})=e^{\lambda} & \qquad\mathrm{~~in~~}M,\\
			\frac{\partial v}{\partial\vec{n}}=0 & \qquad\mathrm{~~on~~}\partial M.
		\end{cases}
		\]
		
		Moreover, the uniqueness of the constant $e^{\lambda}$ follows from
		the strong maximum principle and Hopf lemma for the elliptic equation. 
	\end{proof}
	With the help of Proposition \ref{prop:sigma_2_1st_bdry_eigenvalue}, it is not
	hard to show the following:
	
	Let $(M,g)$ be a compact manifold of dimension $n\geq3$ with boundary.
	Suppose $\lambda(A_{g})\in\Gamma_{2}^{+}$ and $h_{g}\ge0$ on $\pa M$. Then
	there exists no metric $g_{u}\in[g]$ such that $\lambda(A_{u})\in\Gamma_{2}^{+}$, $\sigma_{2}(A_{u})=0$
	in $M$ and $h_{g_{u}}=0$ on $\pa M$.

	\section{Blow up analysis}\label{Sect:blowup_analysis}

	In this section, we mostly denote a conformal metric by $g_{u}=u^{4/(n-2)}g$. Given a positive $C^2$ function $u$, by an abuse of notation we define 
	\begin{equation}\label{def:new_A_u}
		A_{u}:=-\frac{2}{n-2}\frac{\nabla^{2}u}{u}+\frac{2n}{(n-2)^{2}}\frac{\ud u\otimes\ud u}{u^{2}}-\frac{2}{(n-2)^{2}}\frac{|\nabla u|^{2}}{u^{2}}g+A_{g}
	\end{equation}
	for convenience. We sometimes denote $A_u$ as $A_{g_u}$ or $A_{u^{4/(n-2)}g}$  to indicate the metric $g_u$.

	To obtain the existence for boundary  $\sigma_{2}$-curvature equation \eqref{PDE:prescribing_sigma_2}, we adopt the method of Leray-Schauder degree. To this end, we first study  the solvability of positive solutions to  a slightly more general boundary $\sigma_2$-curvature equation:
	\begin{equation}
		\begin{cases}
			\sigma_{2}^{1/2}(A_{u})=fu^{\frac{4}{n-2}}+f_{0}, \quad \lambda(A_u) \in \Gamma_2^+ & \qquad\mathrm{~~in~~}M,\\
			\frac{2}{n-2}\frac{\pa u}{\pa\vec{n}}=-cu^{\frac{n}{n-2}}+h_{g}u, & \qquad\mathrm{~~on~~}\pa M,
		\end{cases}\label{PDE:sigma_2_new}
	\end{equation}
	where $f_{0}$ is a non-negative function in $\overline{M}$. For
	simplicity, we just take $f_{0}$ as a constant function, which is enough to our use. Clearly,  $w=-2\log u/(n-2)$ for another conformal factor $e^{2w}$.

	Our perspective is to develop the blow-up theory for  boundary $\sigma_2$-curvature equation \eqref{PDE:sigma_2_new}. The results of this section will be applied in the next Section \ref{Sect:degree_theory}.
	
	\begin{convention} 
		For a nonnegative constant $T$  we define
		\[
		\mathbb{R}_{-T}^{n}:=\{y=(y',y_{n})\in\Rn;y_{n}>-T\}.
		\]
		In particular, $\Rn_{+}:=\{y=(y',y_{n})\in\Rn;y_{n}>0\}$
		denotes the upper half space. We let
		$$B_\rho=B_\rho(0),\quad B_\rho^+=B_\rho \cap \Rn_+,  \quad \pa^+ B_\rho^+ = \pa B_\rho^+ \cap \Rn_+,\quad D_\rho=\pa B_\rho^+\backslash \pa^+ B_\rho^+.$$
		Denote by $\mathbf{e}_{n}$ the unit
		direction vector in the $n$-th coordinate and by $g_{\mathrm{E}}$ the flat metric in $\Rn$.
	\end{convention}

	\subsection{Prior results}

	We collect some known results together and present some preliminaries relevant to the forthcoming blow-up analysis.  
	
	Let us start with the definition of viscosity solutions to fully nonlinear PDEs in  a Euclidean domain endowed with a background metric $g_{\mathrm{E}}$. See Y. Y. Li \cite{LiYY,LiYY0} and Li-Nguyen \cite{Li-Luc1} etc. Via geodesic normal/Fermi coordinates around some fixed point, it is transparent that a viscosity solution to fully nonlinear PDEs  on manifolds can be similarly defined. 	
	
	\begin{definition}\label{vis-definition} 
		
		Let $\Omega$ be an open subset in $\Rn$ for $n\geq 3$. A positive function
		$w \in C^0(\Omega)$  is called a viscosity supersolution (resp.~ subsolution)
		of $\lambda(A_{w})\in\partial \Gamma_2^+$ when the following holds: if
		$x_{0}\in\Omega$, $\varphi\in C^{2}(\Omega),(w-\varphi)\left(x_{0}\right)=0$,
		and $w-\varphi\geq0$ near $x_{0}$, then 
		\[
		\lambda\left(A_{\varphi}\right)(x_{0})\in \overline{\Gamma_2^+}
		\]
		(resp.~ if $(w-\varphi)\left(x_{0}\right)=0$ and $w-\varphi\leq0$
		near $x_{0}$, then $\lambda\left(A_{\varphi}\right)(x_{0})\in \mathbb{R}^{n}\backslash\Gamma_2^+)$.
		We say that $w$ is a viscosity solution of $\lambda(A_{w})\in\partial \Gamma_2^+$
		if it is both a viscosity supersolution and a viscosity subsolution. \end{definition}

	\begin{theorem}[\protect{\cite[Theorems 1.3 and 1.4]{Li-Luc1}}]
		\label{thm:(Bocher-Type-theorems--Li-Nguyen} Let $u\in\operatorname{LSC}\left(B_1\backslash\{0\}\right)\cap L_{\mathrm{loc}}^{\infty}\left(B_1\backslash\{0\}\right)$
		be a positive viscosity solution of $\sigma_{2}(\lambda(A_{u}))=0$ and $\lambda\left(A_{u}\right)\in\overline{\Gamma_2^+}$
		in $B_1\backslash\{0\}$. Here we use `LSC' instead of `lower semicontinuous' for short. Then
		\begin{itemize}
			\item[(1)] for $n=3$, either $u(x)=C |x|^{2-n}$
			for some $C\in \R_+$, or $u$ can be extended to a positive function in
			$C_{\mathrm{loc}}^{1/2}(B_1)$. Moreover,
			there holds 
			\[
			\big\| u^{-\frac 1 2}\big\|_{C^{\frac 1 2}(B_{1/2})}\leq C(\Gamma_2^+)\sup_{\partial B_{3/4}}u^{-\frac 1 2};
			\]	
			\item[(2)] for $n=4$, \[\log u(x)=-\alpha\log |x|+\mathring{w}(x),\]
			where $\alpha\in [0,n-2]$ and $\mathring{w}\in L^\infty_{\mathrm{loc}}(B_1)$ satisfies
			\[\min_{\partial B_r}\mathring{w}\le \mathring{w}(x) \le \max_{\partial B_r}\mathring{w} \qquad \mathrm{~~in~~} B_r\backslash\{0\}, ~~ \forall~ r \in (0,1). \]
			In particular, if $\alpha=n-2$, then $\mathring{w}$ is constant. 
		\end{itemize}
	\end{theorem}

	\begin{theorem}[\protect{\cite[Theorem B]{Han-Li-Tei}}]
		\label{thm:Han-Li-Texi} Assume that $\sigma_{2}(\lambda(g_u^{-1}A_{u}))=1$
		in $B_1\backslash\{0\}$ and $\lambda\left(A_{u}\right)\in\overline{\Gamma_2^+}$. Then for the radial solution $\bar{u}$
		to $\sigma_{2}(\lambda(g_{\bar{u}}^{-1}A_{\bar{u}}))=1$ in $\mathbb{R}^{n}\backslash\{0\}$,
		\[
		|u(x)-\bar{u}(|x|)|\le C|x|^{\alpha}\bar{u}(|x|),
		\]
		where for $n=3$,  $\bar{u}$ is positive, finite limit at $x=0.$ 
		
	\end{theorem}
	
	\begin{theorem}[\protect{\cite[Theorem 1.2]{LiYY0}}] \label{thm:(Li--JFA-2006}Let
		$u\in C^{2}(\mathbb{R}^{n}\backslash\{0\})$ be the solution to $\sigma_{2}(\lambda(g_u^{-1}A_{u}))=1$
		on $\mathbb{R}^{n}\backslash\{0\}$ with $\lambda(A_{u})\in\Gamma_{2}^{+}$.
		Assume that $u$ can not be extended as a $C^{2}$ positive function
		satisfying $\lambda(A_{u})\in\Gamma_{2}^{+}$ near zero. Then $u$ is radially
		symmetric about the origin and $u'(r)<0$ for all $0<r<\infty.$ 
	\end{theorem}
	
	\begin{theorem}[\cite{Li-Li2}] \label{Thm:Acta}
		
		For $n\geq 3$, let
		$\lambda({A_u})\in \Gamma_2^+$ and 
		$F=\sigma_{2}^{1/2}$. Assume that $u\in C^{2}(\Rn)$
		is a superharmonic solution of $F(\lambda(g_u^{-1}A_{u}))=1$.
		Then for some
		$\bar{x}\in\Rn$ and some positive constants $a$ and $b$ satisfying
		$\lambda(2b^{2}a^{-2}I)\in \Gamma_2^+$ and $F(2b^{2}a^{-2}I)=1$, 
		\[
		u(x)=\left(\frac{a}{1+b^{2}|x-\bar{x}|^{2}}\right)^{\frac{n-2}{2}},\quad x\in\mathbb{R}^{n}.
		\]
	\end{theorem}
	
	\begin{theorem}[\cite{Li-Li3}]\label{half liouville theorem}
		For $n\geq 3$, let
		$\lambda({A_u})\in \Gamma_2^+$ and 
		$F=\sigma_{2}^{1/2}$. Assume that $u\in C^{2}(\Rn_+)$
		is a superharmonic solution of $F(\lambda(g_u^{-1}A_{u}))=1$.
		Then 
		\[ u(x) =\left(\frac{a}{1+b\left|x-\bar{x}\right|^{2}}\right)^{\frac{n-2} {2}} \quad \mathrm{~~for~~} x=(x^{\prime}, x_{n}) \in \mathbb{R}_{+}^{n},\]
		where $\bar{x}=\left(\bar{x}^{\prime}, \bar{x}_{n}\right) \in \mathbb{R}^{n}, a>0$ and $b+\min\left\{ \bar{x}_{n}, 0\right\}^{2}>0$ are two constants satisfying $\lambda(2b^{2}a^{-2}I)\in \Gamma_2^+, F\left(2 a^{-2} b I\right)=1$ and $(n-2) a^{-1} b \bar{x}_{n}=c$.
	\end{theorem}
	
	The following are concerned with degenerate $\sigma_2$-curvature equations. 
	
	\begin{theorem}[\protect{\cite[Theorem 1.18]{LiYY0}}] \label{thm:(Li--CPAM09} Let
		$u\in C^{0,1}(\mathbb{R}^{n}\backslash\{0\})$ be a viscosity solution
		to $\sigma_{2}(A_u)=0$ in $\mathbb{R}^{n}\backslash\{0\}$
		with $\lambda(A_u)\in\Gamma_{2}^{+}$. 
		Then $u$ is radially symmetric about the origin and $u'(r)<0$ for
		all almost $0<r<\infty.$ \end{theorem}
	
	As in \cite{{Li-Luc1}}, we denote by 
	$\mu_{\Gamma}^+ \in [0,n-1]$ the unique number such that $(-\mu_\Gamma^+,1,\cdots,1) \in \pa \Gamma_2^+$. It is clear that $\mu_{\Gamma}^{+}=1/2$ for $n=3$ and $\mu_{\Gamma}^{+}=1$ for $n=4$. 
	\begin{theorem}[\protect{\cite[Theorem 2.2]{Li-Luc1}}] \label{thm:(Li-Nguyen2015-CPAM-Prop2.2)radial solution}
		For $0\leq a<b\leq \infty$, every radially symmetric positive viscosity solution $u$ of $\sigma_{2}(\lambda(A_u))=0$ with $\lambda(A_u)\in \Gamma_{2}^{+}$ in an annulus $\{a<|x|<b\}$ is one of the following smooth solutions:
		
		(a) $u(x)=C_{1}|x|^{-C_{2}}$ with $C_{1}>0,0\leq C_{2}\leq n-2$
		when $\mu_{\Gamma}^{+}=1$,
		
		(b) $u(x)=\left(C_{3}|x|^{-\mu_{\Gamma}^{+}+1}+C_{4}\right)^{\frac{n-2}{\mu_{\Gamma}^{+}-1}}$
		with $C_{3}\geq0,C_{4}\geq0,C_{3}+C_{4}>0$ when $\mu_{\Gamma}^{+}\neq1$.
		
		%
	\end{theorem}


	For later use, we consider the following fully nonlinear degenerate elliptic equation with Neumann boundary condition. Let $\partial_1\Omega \subset \pa M$ be a partial boundary of a domain $\Omega \subset \overline M$ and $b(x)$ be a smooth function defined on $\pa \Omega$. As we shall see shortly, the following boundary degenerate elliptic equation naturally arises as a rescaled limit equation of a  blow-up sequence of solutions to \eqref{PDE:sigma_2_new}:
	\begin{align}\label{visocosity solution}
		\begin{cases}
			\lambda(A_{u})\in \partial{\Gamma_2^+} &\qquad\mathrm{~~in~~} \Omega,\\
			-\frac{\partial u}{\partial \vec n}+bu= 0 &\qquad\mathrm{~~on~~} \partial_1\Omega.
		\end{cases}
	\end{align}
	Readers are referred to Appendix \ref{Appendix:C} for a thorough discussion.

	We introduce the following definition of viscosity solutions to \eqref{visocosity solution},  which borrows some ideas from  that of Y. Y. Li and collaborators above and the classical definition used in the study of nonlinear Neumann boundary problem of fully nonlinear elliptic PDEs (see e.g. Crandall-Ishii-Lions \cite{Crandall-Ishii-Lions} and references therein for a survey).
	
	\begin{definition} \label{vis-definition boundary}
		A positive function $u\in C^{0}\left(\Omega\cup \partial_1\Omega \right)$ is
		called a viscosity supersolution (resp. subsolution) of (\ref{visocosity solution})
		if the following hold: Let $\psi\in C^{2}\left(\Omega\cup \partial_1\Omega\right)$
		and $u-\psi$ has a local minimum (resp. a local maximum) at
		$\bar{x}$. If $\bar{x}\in\partial_1\Omega$, then 
		\begin{align}\label{vis-boundary}
		\begin{split}
				&\lambda(A_{\psi})(\bar{x})\in \overline{\Gamma_2^+} \quad \mathrm{or}\quad [-\frac{\partial\psi}{\partial \vec n}+b\psi](\bar{x})\ge0,\\
				&\qquad\qquad\qquad\left(\mathrm{resp.}~~ \lambda(A_{\psi})(\bar{x})\in\mathbb{R}^n\backslash \Gamma_2^+ \quad \mathrm{or}\quad [-\frac{\partial\psi}{\partial \vec n}+b\psi](\bar{x})\le 0\right).				\end{split}
		\end{align}
		If $\bar x\in \Omega $, then 
		\[\lambda(A_{\psi})(\bar{x})\in\overline{\Gamma_2^+},\quad (\mathrm{resp.}~~ \lambda(A_{\psi})\in \mathbb{R}^n\backslash \Gamma_2^+ ).\]
		
		We denote the viscosity supersolution (resp. subsolution) by 
		\[\begin{cases}
			\lambda(A_{u})\in\overline{\Gamma_2^+} & \mathrm{in~~} \Omega\\
			-\frac{\partial u}{\partial \vec n}+bu\ge 0 & \mathrm{on~~} \partial_1\Omega
		\end{cases},
		\quad\bigg(\mathrm{resp.}	 \begin{cases}
			\lambda(A_{u})\in \mathbb{R}^n\backslash \Gamma_2^+ & \mathrm{in~~}\Omega\\
			-\frac{\partial u}{\partial \vec n}+bu\le 0 & \mathrm{on~~} \partial_1\Omega
		\end{cases}.\bigg)\]

		Similarly, we define a \emph{strict} viscosity supersolution (resp. subsolution) by  
		\[\begin{cases}
			\lambda(A_{u})\in\overline{\Gamma_2^+} & \mathrm{in~~} \Omega\\
			-\frac{\partial u}{\partial \vec n}+bu> 0 &\mathrm{on~~} \partial_1\Omega
		\end{cases},
		\quad\bigg(\mathrm{resp.}	 \begin{cases}
			\lambda(A_{u})\in \mathbb{R}^n\backslash \Gamma_2^+ &\mathrm{in~~} \Omega\\
			-\frac{\partial u}{\partial \vec n}+bu<0 &\mathrm{on~~} \partial_1\Omega
		\end{cases},\bigg)\] 
		if the inequalities in (\ref{vis-boundary}) are strict.
		
		We say $u \in C^{0}\left(\Omega\cup \partial_1\Omega \right)$ is a viscosity solution to 
		\[\begin{cases}
			\lambda(A_{u})\in \partial{\Gamma_2^+} &\qquad\mathrm{~~in~~} \Omega,\\
			-\frac{\partial u}{\partial \vec n}+bu= 0 &\qquad\mathrm{~~on~~} \partial_1\Omega.
		\end{cases}\]
		if $u$ is both a viscosity supersolution and a viscosity subsolution. \end{definition}

	Let $\Omega^{+}\subset\mathbb{R}_{+}^{n}$
	be an open set and denote 
	\[
	\partial^{\prime\prime}\Omega^{+}=\overline{\partial\Omega^{+}\cap\mathbb{R}_{+}^{n}},\quad\partial^{\prime}\Omega^{+}=\partial\Omega^{+}\backslash\partial^{\prime\prime}\Omega^{+}.
	\]
	
	As shown in Y. Y. Li \cite[Proposition 5.1]{LiYY}, an  even reflection of a positive Lipschitz viscosity solution $u^+$ to  $\lambda(A_{u^+}) \in \pa \Gamma_2^+$ in $B_1^+$ with vanishing Neumann boundary condition on $\pa' B_1^+$ keeps the equation invariant. For our purpose, we aim to show that the same result still holds for a viscosity solution as Definition \ref{vis-definition boundary}, which is slightly different from that of \cite{LiYY}. For clarity, we restate it in the following proposition, whose proof  we present in Appendix \ref{Appendix:B},  in the spirit is  very similar to that of   \cite[Proposition 5.1]{LiYY}.
	
	\begin{proposition}\label{even reflection in Appendix}
		Let $u^{+}\in C^{0,1}(\overline{B_{1}^{+}})$ and $u^{-}\in C^{0,1}(\overline{B_{1}^{-}})$
		be two positive functions satisfying $u^{+}=u^{-}$ on $\partial^{\prime}B_{1}^{+}$. Here $B_1^-=B_1 \cap\{x_n<0\}$.
		Suppose $u^{+}$ is a viscosity supersolution of 
		\[
		\begin{cases}
			\lambda(A_{u})\in\partial\Gamma_2^+ & \quad \mathrm{~~in~~} B_{1}^{+},\\
			-\frac{\partial u}{\partial x_n}=0 &\quad \mathrm{~~on~~}  \partial'B_{1}^{+},
		\end{cases}
		\]
		and $u^{-}$ is a viscosity supersolution of
		\[
		\begin{cases}
			\lambda(A_{u})\in\partial\Gamma_2^+ &\quad \mathrm{~~in~~}  B_{1}^{-},\\
			\frac{\partial u}{\partial x_n}=0 &\quad \mathrm{~~on~~}  \partial' B_{1}^{-}.
		\end{cases}
		\]
		Then 
		\[
		u\left(x^{\prime},x_{n}\right):=\begin{cases}
			u^{+}\left(x^{\prime},x_{n}\right) & \mathrm{~~if~~}x_{n}\geq0,\\
			u^{-}\left(x^{\prime},x_{n}\right) & \mathrm{~~if~~}x_{n}<0,
		\end{cases}
		\]
		is a $C^{0,1}$ viscosity supersolution of $\lambda(A_{u})\in\partial\Gamma_2^+$
		in $B_{1}$. 
	\end{proposition}

	We would like to mention that most theorems just described hold for more general equations in \cite{Li-Luc1,Li-Li2,Li-Li3,Han-Li-Tei,LiYY0,LiYY} and state them as $\sigma_2$-curvature equation for our purpose.
	
	We apply the  local $C^1$ estimates for $\sigma_2$-curvature equation in \cite{Jin-Li-Li} to \eqref{PDE:sigma_2_new} just in a different conformal factor.
	\begin{theorem}[\protect{\cite[Theorem 1.4]{Jin-Li-Li}}]\label{Thm:local $C^1$ estimates}
		Let $(M,g)$ be a smooth compact Riemannian manifold of dimension
		$n\geq3$ with boundary $\pa M$ and assume $0<f \in C^\infty(M)$. Let $\mathcal{O}$ is an open set in $\overline{M}$. Let $u \in C^4(\overline M)$ be a solution to \eqref{PDE:sigma_2_new}. 
		If there exists a positive constant $C_0$ such that
		$$u\le C_0\qquad \mathrm{~~in~~}\quad \mathcal{O},$$
		 then for any open subset $\mathcal{O}' \Subset \mathcal{O}$, there exists a positive constant $C$ depending on $\mathrm{dist}_g(\mathcal{O}', \mathcal{O})$, $n$, $\sup_{\mathcal{O}} (f+f_0+|f_x|)$, $\|c\|_{C^2(\overline{\mathcal O})}$,$\|g\|_{C^3(\overline {\mathcal O})}$ and $C_0$, such that
		$$|\nabla (\log u)| \leq C \qquad \mathrm{~~in~~} \mathcal{O}'.$$
	\end{theorem}

	\subsection{Basic results for all dimensions}
	
	\label{Subsect:Basic facts in all dimensions}
	
	The first goal is to establish Lemma \ref{lem fudament-1-1} for the boundary $\sigma_2$-curvature equation \eqref{PDE:sigma_2_new} when $n\geq 3$, which is the counterpart corresponding to those of Schoen and Zhang \cite[Lemma 3.1]{Schoen-Zhang} in scalar curvature equation and Han-Li \cite[Proposition 1.1]{han-li2} in the boundary Yamabe problem. 
	
	As a first step, we present some common ingredients for dimensions $n\geq 3$. 
	\begin{lemma}
		\label{lem:preparation lemma}
		 For any $R\geq1$ and $0<\ve<1$,
		there exists a positive constant $C_{0}$ depending on $M,g,f,c,R$ and $\ve$, such
		that for any compact set $K\subset\overline{M}$ and any positive smooth solution
		$u$ to $(\ref{PDE:sigma_2_new})$ with 
		\[
		\max_{q\in\overline{M\backslash K}}d_{g}^{\frac{n-2}{2}}(q,K)u(q)\geq C_{0},
		\]
		where $d_{g}(q,K):=1$ for $K=\emptyset$, and the following properties hold. 
		
		There exist $q_{0}\in\overline{M}\backslash K$, which is a local
		maximum point of $u$ in $\overline{M}$, and $r_\ast \in (0,1)$ such that one of the three cases happens:
		\begin{itemize}
			
			\item[($a_1$)] If $q_{0}\in M$ and $d_g(q_0,\partial M)\ge r_*$, then 
			\[
			\left\Vert u^{-1}(q_{0})u\circ\exp_{q_{0}}\left(u(q_{0})^{-\frac{2}{n-2}}y\right)-\left(1+\frac{f(q_{0})}{2\sqrt{\binom{n}{2}}}|y|^{2}\right)^{\frac{2-n}{2}}\right\Vert _{C^{2}\left(B_{2R}\right)}<\ve.
			\]
			
			\item[($a_2$)] If $q_{0}\in \overline{M}$ and $d_g(q_0,\partial M)<r_*$, then 
			\[
			\left\Vert u^{-1}(q_{0})u\circ G_{q_0'}\left(u^{-\frac{2}{n-2}}(q_{0})y+{d}_{0}\mathbf{e}_{n}\right)-\left(1+\frac{f(q_{0})}{2\sqrt{\binom{n}{2}}}|y|^{2}\right)^{\frac{2-n}{2}}\right\Vert _{C^{2}\left(B_{2R}^{+}\right)}<\ve,
			\]	
			where $d_0:=d_g(q_0, \partial M)=d_g(q_0, q_0')$ for some $q_0'\in \partial M$.	
			
			\item[($a_3$)] If $q_{0}\in\partial M$, then 
			\[
			\left\Vert u^{-1}(q_{0})u\circ G_{q_{0}}\left(u(q_{0})^{-\frac{2}{n-2}}y\right)-(\frac{2b\sqrt{\binom{n}{2}}}{f(q_0)})^{\frac{n-2}{4}}\left(1+b|y+\tilde T_c \mathbf{e}_{n}|^{2}\right)^{\frac{2-n}{2}}\right\Vert _{C^{2}\left(B_{2R}^{+}\right)}<\ve,
			\]
			where 
			\[
			b=\frac{f(q_0)}{2\sqrt{\binom{n}{2}}}\left(1+\frac{c(q_0)^2}{2f(q_0)}\sqrt{\binom{n}{2}}\right)^2\quad\mathrm{~~and~~}\quad \tilde T_c=\frac{c(q_{0})}{\sqrt{2f(q_{0})}}\frac{\binom{n}{2}^{1/4}}{\sqrt{b}}.
			\]
		\end{itemize}
	\end{lemma} 
	\begin{proof}
		Suppose not, there exist positive constants $R,\ve$ such that for any $i \in \N$,
		there exist a sequence of compact sets $K_{i}\subset\overline{M}$
		and a sequence of positive solutions $u_{i}$ to (\ref{PDE:sigma_2_new}) with
		\[
		\max_{q\in\overline{M\backslash K_{i}}}d_{g}^{\frac{n-2}{2}}(q,K_{i})u_{i}(q)\geq i,
		\]
		and there is no $q_{0}\in\overline{M}\backslash K_i$, which is a local maximum point of $u_i$ in $\overline M$, satisfying any of $(a_1)$-$(a_3)$.
		
		Let $q_{i}\in\overline{M}\backslash K_i$ be the point such that
		\[
		d_{g}^{\frac{n-2}{2}}(q_{i},K_{i})u_{i}(q_{i})=\max_{q\in\overline{M\backslash K_{i}}}d_{g}^{\frac{n-2}{2}}(q,K_{i})u_{i}(q)\to\infty\qquad\mathrm{as~~}i\to\infty.
		\]
		This indicates that $u_{i}(q_{i})\to\infty$.
		
		We abbreviate $\lambda_{i}=u_{i}^{-2/(n-2)}(q_{i})$, $d_{i}:=d_{g}(q_{i},q_{i}')=d_{g}(q_{i},\pa M)$ for some $q_{i}'\in\pa M$
		and $T_{i}=\lambda_{i}^{-1}d_{i}$.
		
		Fix a small $\delta>0$ and assume that up to a subsequence, $d_{i}<\delta$
		for all $i$. Under Fermi coordinates around $q_{i}'$, we introduce
		\[
		x=\Psi_{i}(y)=G_{q_{i}'}\left(\lambda_{i}y+d_{i}\mathbf{e}_{n}\right)\quad\Rightarrow\quad q_{i}=\Psi_{i}(0), ~~q_i'=\Psi_i(-T_i \mathbf{e}_n).
		\]
		Define 
		\[
		v_{i}(y):=\frac{1}{u_{i}(q_{i})}u_{i}\circ\Psi_{i}(y)=\lambda_i^{\frac{n-2}{2}}u_{i}\circ\Psi_{i}(y) \qquad\mathrm{~~and~~}\qquad g_{i}=\lambda_i^{-2}\Psi_{i}^{\ast}(g),
		\]
		then
		$$v_i^{\frac{4}{n-2}}g_i=\Psi_i^\ast(u_i^{\frac{4}{n-2}}g).$$
		
		For each $i$, we consider $u_i$ in a domain $\Psi_{i}(\Omega_{i})$ with
		\[
		\Omega_{i}:=\{y\in\Rn_{-T_{i}};d_{g}(\Psi_{i}(y),q_{i})<\frac{1}{4}d_{g}(q_{i},K_{i})\}.
		\]
		Then for any $x\in\Psi_{i}(\Omega_{i})$, there holds
		\begin{equation}\label{est:distance lower}
			d_{g}(x,K_{i})\geq\frac{1}{2}d_{g}(q_{i},K_{i}) \quad \Longrightarrow \quad \Psi_{i}(\Omega_{i})\cap K_{i}=\emptyset.
		\end{equation}
		It follows that for any $x\in\Psi_{i}(\Omega_{i})$, by definition of $q_{i}$
		we have 
		\begin{align*}
			d_{g}^{\frac{n-2}{2}}(q_{i},K_{i})u_{i}(q_{i}) & \ge d_{g}^{\frac{n-2}{2}}(x,K_{i})u_{i}(x)\\
			& \ge\bigg(\frac{1}{2}d_{g}(q_{i},K_{i})\bigg)^{\frac{n-2}{2}}u_{i}(x).
		\end{align*}
		This means
		\[
		u_{i}(x)\le2^{\frac{n-2}{2}}u_{i}(q_{i})
		\quad \Longrightarrow \quad
		v_{i}(y)\le2^{\frac{n-2}{2}},\qquad\forall~y\in~\Omega_{i}.
		\]
		
		Notice that $v_{i}$ satisfies 
		\begin{align*}
			& \sigma_{2}^{1/2}\left(g_{i}^{-1}(-\frac{2}{n-2}\frac{\nabla_{g_{i}}^{2}v_{i}}{v_{i}}+\frac{2n}{(n-2)^{2}}\frac{\ud v_{i}\otimes\ud v_{i}}{v_{i}^{2}}-\frac{2}{(n-2)^{2}}\frac{|\nabla v_{i}|_{g_{i}}^{2}}{v_{i}^{2}}g_{i}+A_{g_{i}})\right)\\
			= & f\circ\Psi_{i}v_{i}^{\frac{4}{n-2}}+\lambda_{i}^{2}f_{0}\qquad\mathrm{~~in~~}\quad|y+T_{i}\mathbf{e}_{n}|<\frac{d_{g}(q_{i},K_{i})}{4\lambda_{i}},~~y_{n}>-T_{i}
		\end{align*}
		and the boundary condition 
		\[
		\frac{2}{n-2}\frac{\partial v_{i}}{\partial\vec{n}_{g_{i}}}=-c\circ\Psi_{i}v_{i}^{\frac{n}{n-2}}+h_{g_{i}}v_{i}\qquad\mathrm{~~on~~}\quad y_{n}=-T_{i}.
		\]
		
		Up to a further subsequence, we assume $\lim_{i\rightarrow\infty}T_{i}:=T\in[0,\infty]$.
		It suffices to consider $T<\infty$, since $T=\infty$
		can be similarly handled as  \emph{Case 1} below.
		
		By  $C^1$ estimates in \cite{Jin-Li-Li} and Theorem \ref{Thm:local estimates} we have
		\[
		\|v_{i}\|_{C^{2}(\tilde{K})}\le C_{\tilde{K}}\qquad\mathrm{~~for~~any~~compact~~set~~}\tilde{K}\subset\Omega_{i},
		\]
		which can be improved by Evans-Krylov theory to 
		\[
		\|v_{i}\|_{C^{2,\alpha}(\tilde{K})}\le C_{\tilde{K}}\qquad\mathrm{for~~any~~}0<\alpha<1.
		\]
		
		Notice that $T<\infty$ and $u_{i}(q_{i})\to\infty$, this forces $q_{i}\rightarrow q_{\infty}\in\pa M$.
		Hence, a subsequence $v_{i}$ converges in $C_{\mathrm{loc}}^{2}(\overline{\mathbb{R}_{-T}^{n}})$
		to $v_{\infty}$ with $v_{\infty}(0)=1$ by virtue of $v_i(0)=1, ~\forall~i \in \N$, and $v_{\infty}$ satisfies
		\begin{equation}\label{eqns:limit}
		\begin{cases}
			\displaystyle \quad \sigma_{2}^{1/2}(-\frac{2}{n-2}\frac{\nabla_{g_{\mathrm{E}}}^{2}v_{\infty}}{v_{\infty}}+\frac{2n}{(n-2)^{2}}\frac{\ud v_{\infty}\otimes\ud v_{\infty}}{v_{\infty}^{2}}-\frac{2}{(n-2)^{2}}\frac{|\nabla v_{\infty}|_{g_{\mathrm{E}}}^{2}}{v_{\infty}^{2}}g_{\mathrm{E}})\\
		\displaystyle	=f(q_{\infty})v_{\infty}^{\frac{4}{n-2}} \qquad \qquad \quad \mathrm{~~in~~}\quad \mathbb{R}_{-T}^{n},\\
		\displaystyle \frac{\partial v_{\infty}}{\partial y_{n}}=-\frac{n-2}{2} c(q_{\infty})v_{\infty}^{\frac{n}{n-2}} \qquad \mathrm{~~on~~}\quad y_{n}=-T.
		\end{cases}
		\end{equation}
		Then we assert  that  the unique solution is
		\begin{equation}\label{standard bubble}
			v_{\infty}(y)=\frac{1}{f(q_{\infty})^{\frac{n-2}{4}}}\bigg(\frac{(2b\binom{n}{2}^{\frac{1}{2}})^{\frac{1}{2}}}{1+b|y-(0,\bar{y}_{n})|^{2}}\bigg)^{\frac{n-2}{2}}\quad\mathrm{~~for~~}y=(y',y_{n})\in\mathbb{R}_{-T}^{n},
		\end{equation}
		where $b\in\R_{+}, \bar{y}_{n}\in\R$ satisfy
		\begin{subequations}
		\begin{align}
			\bar{y}_{n}+T=-\frac{c(q_{\infty})}{\sqrt{2f(q_{\infty})}}\frac{\binom{n}{2}^{1/4}}{\sqrt{b}},\label{eq:bdry_T}\\
		(\bar{y}_{n})^{2}b-\sqrt{\frac{2}{f(q_{\infty})}}\binom{n}{2}^{1/4}\sqrt{b}+1=0.\label{eq:bdry_fixed_pt}
		\end{align}
		\end{subequations}
		
		Since the derivation of formulae \eqref{standard bubble}, \eqref{eq:bdry_T} and \eqref{eq:bdry_fixed_pt} is very delicate, it is reasonable to give a geometric interpretation. For brevity, we denote by $A_{v_\infty}$ the Schouten tensor of metric $g_{v_\infty}:=v_\infty^{4/(n-2)}|\ud y|^2$ and define $T_c:=\sqrt{b}(\bar y_n +T)$ in consistent of the notations in \cite{Chen-Ruan-Sun,Chen-Sun}. Denote by a mapping $\pi: \Sn((-T+T_c)\mathbf{e}_n) \setminus\{(-T+T_c) \mathbf{e}_n+\mathbf{e}_{n+1}\} \to \{(y,0)+(-T+T_c) \mathbf{e}_n \in \mathbb{R}^{n+1}\}\simeq \mathbb{R}^n$ the stereographic projection from the unit sphere $\Sn((-T+T_c)\mathbf{e}_n)$ in $\mathbb{R}^{n+1}$ centered at $(-T+T_c) \mathbf{e}_n$. By  the classification theorem of A. Li and Y.Y. Li \cite{Li-Li1}  we know
		$$(\pi^{-1})^\ast(\lambda g_{\Sn})=v_\infty^{\frac{4}{n-2}}(y) |\ud y|^2:=(4 \lambda) W^{\frac{4}{n-2}} (y)|\ud y|^2 \qquad \mathrm{for~~some~~}\lambda \in \R_+,$$
		where $W$ is the \emph{bubble function} up to horizontal translations, like
		$$W(y)=\left(\frac{\sqrt{b}}{1+b|y-\bar y_n \mathbf{e}_{n} |^2}\right)^{\frac{n-2}{2}},$$
		see Chen-Ruan-Sun \cite{Chen-Ruan-Sun} and Chen-Sun \cite[pp.8-9]{Chen-Sun} for its geometric meaning.  We now interpret the limit equations \eqref{eqns:limit} as
	$$\begin{cases}
	\displaystyle f(q_\infty)=\sigma_2^{1/2}(g_{v_\infty}^{-1}A_{v_\infty})=\lambda^{-1}\sigma_2^{1/2}(A_{\Sn})=\frac{\sqrt{\binom{n}{2}}}{2\lambda} &\quad \Longrightarrow \quad \lambda=\frac{\sqrt{\binom{n}{2}}}{2f(q_\infty)},\\
	\displaystyle c(q_\infty)=h_{g_{v_\infty}}= \frac{1}{\sqrt{4\lambda}}(-2 T_c)=-\frac{T_c}{\sqrt{\lambda}} &\quad \Longrightarrow \quad T_c=-c(q_\infty) \sqrt{\lambda}.
	\end{cases}$$
	Together with $v_\infty(0)=1$:
	$$\frac{2\sqrt{\lambda b}}{1+b |\bar y_n|^2}=1,$$
	we obtain the desired formulae.

		\vskip 4pt
		\emph{Case 1:} $c(q_{\infty})=0$.
		\vskip 4pt
		
		Denote $\hat{y}=-T\mathbf{e}_{n}$ and 
		\[
		\max_{\overline{\R_{-T}^{n}}}v_{\infty}=v_{\infty}(\hat{y})=\left(\frac{2b\sqrt{\binom{n}{2}}}{f(q_{\infty})}\right)^{\frac{n-2}{4}}.
		\]
		There exists a sequence $y_{i}\rightarrow\hat{y},$ and each $y_{i}$ is a
		local maximum point of $v_{i}$ such that $v_{i}(y_{i})\rightarrow v_{\infty}(\hat{y})$. This together with (\ref{est:distance lower}) implies
		that every point $p_{i}=\Psi_{i}(y_{i})\in\overline{M}\backslash K_{i}$
		is also a local maximum point of $u_{i}$. Up to a subsequence, $p_i \to  p_\infty \in \pa M$. This enables us to repeat
		the previous procedure with $p_{i}$ replacing $q_{i}$.  Eventually we 
		conclude that given any $0<\ve<1$ and $R\geq 1$, for all sufficiently large $i$ we have
		\[
		\left\Vert u_{i}^{-1}\left(p_{i}\right)u\circ G_{p_{i}'}\left(u_{i}^{-\frac{2}{n-2}}(p_{i})y+\tilde{d}_{i}\mathbf{e}_{n}\right)-\tilde v_\infty(y)\right\Vert _{C^{2}\left(B_{2R}^{+}\right)}<\ve.
		\]
		Here  $\tilde{d}_{i}=d_{g}(p_{i},\pa M)=d_{g}(p_{i},p_{i}')$ for some
		$p_{i}'\in\pa M$ and up to horizontal translations we have
		\begin{equation*}
			\tilde v_{\infty}(y)=\frac{1}{f(p_{\infty})^{\frac{n-2}{4}}}\bigg(\frac{(2\tilde b\sqrt{\binom{n}{2}})^{\frac{1}{2}}}{1+\tilde b|y-(0,\bar{y}_{n})|^{2}}\bigg)^{\frac{n-2}{2}}
		\end{equation*}
		for $y=(y',y_{n})\in\mathbb{R}_+^{n}$ and $\tilde b \in \R_+$, $\bar{y}_{n}\in\R$ satisfying
		\begin{equation}\label{bdry_T-1}
			\bar{y}_{n}=-\frac{c(p_{\infty})}{\sqrt{2f(p_{\infty})}}\frac{\binom{n}{2}^{1/4}}{\sqrt{\tilde b}}
		\end{equation}
		and
		\begin{equation}\label{value_b}
			(\bar{y}_{n})^{2}\tilde b-\sqrt{\frac{2}{f(p_{\infty})}}\binom{n}{2}^{1/4}\sqrt{\tilde b}+1=0.
		\end{equation}
		For brevity, we let $\tilde{v}_i(y):=u_{i}^{-1}\left(p_{i}\right)u\circ G_{p_{i}'}(u_{i}^{-2/(n-2)}(p_{i})y+\tilde{d}_{i}\mathbf{e}_{n})$.
		Notice that $ \nabla \tilde v_i(0)\rightarrow 0$, this forces $\bar y_n=0$. This together with \eqref{bdry_T-1} implies $c(p_\infty)=0$. By \eqref{value_b} we obtain
		$$\tilde b=\frac{f(p_{\infty})}{2\sqrt{\binom{n}{2}}},$$
		whence,
		$$\tilde v_\infty(y)=\left(1+\frac{f(p_{\infty})}{2\sqrt{\binom{n}{2}}}|y|^{2}\right)^{\frac{2-n}{2}}.$$

		However, this violates the contradiction hypothesis $(a_2)$.

		\vskip 4pt
		\emph{Case 2:} $c(q_{\infty})>0$.
		\vskip 4pt
		
		This means $\bar{y}_{n}<-T$ by virtue of \eqref{eq:bdry_T}. It follows
		that there exists $y_{i}\in\Psi_{i}^{-1}(\partial M)$ that is a local
		maximum point of $v_{i}$ such that $y_{i}\rightarrow-T\mathbf{e}_{n}$.
		The corresponding point $p_{i}=\Psi_{i}(y_{i})\in\pa M$ is also a
		local maximum point of $u_{i}$. Noticing that $d_g(p_{i},\partial M)u_{i}^{2/(n-2)}(p_{i})=0$,
		we prefer to use $G_{p_{i}}$ instead of $\Psi_{i}$ above and repeat
		the previous procedure. Thus up to a subsequence, for all sufficiently
		large $i$ we obtain 
		\[
		\left\Vert u_{i}^{-1}\left(p_{i}\right)u_i\circ G_{p_{i}}\left(u_{i}^{-\frac{2}{n-2}}(p_{i})y\right)-\left(\frac{2b\sqrt{\binom{n}{2}}}{f(p_{\infty})}\right)^{\frac{n-2}{4}}(1+b|y+\tilde T_c\mathbf{e}_{n}|^{2})^{\frac{2-n}{2}}\right\Vert _{C^{2}\left(B_{2R}^{+}\right)}<\ve,
		\]
		where $p_{i}\to p_{\infty}\in\pa M$ and 
		\[
		\tilde T_c=\frac{c(p_{\infty})}{\sqrt{2f(p_{\infty})}}\frac{\binom{n}{2}^{1/4}}{\sqrt{b}}
		\]
		with $b\in\R_{+}$ being determined by 
		\[
		\frac{1}{f(p_{\infty})^{\frac{n-2}{4}}}\left(\frac{(2b\sqrt{\binom{n}{2}})^{\frac{1}{2}}}{1+b|\tilde T_c|^{2}}\right)^{\frac{n-2}{2}}=1.
		\]
		
		As $f\in C^{2}(\overline{M})$ and $c\in C^{1}(\pa M)$, we know that
		for sufficiently large $i$
		
		\[
		\left\Vert u_{i}^{-1}\left(p_{i}\right)u_i\circ G_{p_{i}}\left(u_{i}^{-\frac{2}{n-2}}(p_{i})y\right)-(\frac{2b\sqrt{\binom{n}{2}}}{f(p_{i})})^{\frac{n-2}{4}}(1+b_{i}|y+\tilde T_c^{(i)}\mathbf{e}_{n}|^{2})^{\frac{2-n}{2}}\right\Vert _{C^{2}\left(B_{2R}^{+}\right)}<\ve.
		\]
		Here 
		\[
		\tilde T_c^{(i)}=\frac{c(p_{i})}{\sqrt{2f^(p_{i})}}\frac{\binom{n}{2}^{1/4}}{\sqrt{b_i}}
		\]
		with $b_{i}\in\R_{+}$ being determined by 
		\[
		\left(\frac{2b_i\sqrt{\binom{n}{2}}}{f(p_{i})}\right)^{\frac{n-2}{4}}(1+b_{i}|\tilde T_c^{(i)}|^{2})^{\frac{2-n}{2}}=1.
		\]
		It is not hard to see that 
		\begin{align*}
			b_{i}=\frac{f(p_{i})}{2\sqrt{\binom{n}{2}}}\left(1+\frac{c(p_{i})^2}{2f(p_i)}\sqrt{\binom{n}{2}}\right)^2.
		\end{align*}
		However, it again violates the contradiction hypothesis $(a_3)$. 
	\end{proof}

	Now we give a lemma only relying on the Liouville theorems.

		%

	\begin{lemma}\label{lem fudament-1-1}
		 For any $R\geq1$ and $0<\ve<1$, there exist two positive constants $C_{0}, C_{1}$ depending on $M,g,f,c,R$ and $\ve$, such that for all positive solutions $u$
		to (\ref{PDE:sigma_2_new}) with $\max_{\overline{M}}u\geq C_{0}$,
		and there exists a finite set $\left\{ q_{1},\cdots,q_{N}\right\} \subset\overline{M}$
		with $N\geq1$ and the following are true:
		
		(i) Each $q_{i}$ is a local maximum point of $u$ in $\overline{M}$
		and 
		\[
		\overline{B_{{r}_{i}}\left(q_{i}\right)}\cap\overline{B_{{r}_{j}}\left(q_{j}\right)}=\emptyset \qquad\mathrm{~~for~~}i\neq j,
		\]
		where ${r}_{i}=Ru^{-2/(n-2)}(q_{i})$ and $B_{r_i}(q_i):=\{q \in \overline M; d_g(q,q_i)<r_i\}$.
		
		(ii) One of the following three cases happens:
		\begin{itemize}
			\item[($a_1$)] $q_{i}\in M$, 
			\[
			\left\Vert u^{-1}\left(q_{i}\right)u\circ\exp_{q_{i}}\left(u^{-\frac{2}{n-2}}\left(q_{i}\right)y\right)-\left(1+\frac{f(q_{i})}{2\sqrt{\binom{n}{2}}}|y|^{2}\right)^{\frac{2-n}{2}}\right\Vert _{C^{2}\left(B_{2R}\right)}<\varepsilon.
			\]
			\item[($a_2$)] $q_{i}\in \overline{M}$ near $\pa M$,
			\[
			\left\Vert u^{-1}\left(q_{i}\right)u\circ G_{q_i'}\left(u^{-\frac{2}{n-2}}(q_{i})y+{d}_{i}\mathbf{e}_{n}\right)-\left(1+\frac{f(q_{i})}{2\sqrt{\binom{n}{2}}}|y|^{2}\right)^{\frac{2-n}{2}}\right\Vert _{C^{2}\left(B_{2R}^{+}\right)}<\ve,
			\]	
			where $d_i:=d_g(q_i, \partial M)=d_g(q_i, q_i')$ for some $q_i'\in \partial M$.
			
			\item[($a_3$)] $q_{i}\in\partial M$, 
			
			\begin{align*}
				\left\Vert u^{-1}(q_{i})u\circ G_{q_{i}}\left(u(q_{i})^{-\frac{2}{n-2}}y\right)-(\frac{2b\sqrt{\binom{n}{2}}}{f(q_{i})})^{\frac{n-2}{4}}\left(1+b_i|y+ \tilde T_c^{(i)}\mathbf{e}_{n}|^{2}\right)^{\frac{2-n}{2}}\right\Vert _{C^{2}\left(B_{2R}^{+}\right)}<\ve,
			\end{align*}
			where 
			\[
			b_i=\frac{f(q_{i})}{2\sqrt{\binom{n}{2}}}\left(1+\frac{c(q_{i})^2}{2f(q_i)}\sqrt{\binom{n}{2}}\right)^2\quad\mathrm{~~and~~}\quad \tilde T_c^{(i)}=\frac{c(q_{i})}{\sqrt{2f(q_{i})}}\frac{(\binom{n}{2})^{1/4}}{\sqrt{b}}.
			\]
		\end{itemize}
		(iii) For $q\in\overline{M}$ there holds
		\[
		u(q)\le C_{1}d_{g}(q,\{q_{1},\cdots,q_{N}\})^{-\frac{n-2}{2}}.
		\]
	\end{lemma} 
	\begin{proof}
		A procedure is set up to seek such a finite set of points in $\overline{M}$.
		We start with $d_{g}(q,K)=1$ for $K=\emptyset$, then 
		\[
		\max_{q\in\overline{M\backslash K}}d_{g}^{\frac{n-2}{2}}(q,K)u(q)=\max_{\overline{M}}u=u(q_{1})\geq C_{0}.
		\]
		Then the assertion  \textit{(ii)} follows from Lemma \ref{lem:preparation lemma}.
		
		Take $r_{1}=Ru^{-2/(n-2)}(q_{1})$ and $K_{1}=\overline{B_{r_{1}}(q_{1})}$.
		Now if 
		\[
		\max_{q\in\overline{M\backslash K_{1}}}d_{g}^{\frac{n-2}{2}}(q,K_{1})u(q)<C_{0},
		\]
		then we shall show that \textit{(iii)} holds, which will be seen in the proof
		below soon, and this procedure stops. If not, that is, 
		\[
		\max_{q\in\overline{M\backslash K_{1}}}d_{g}^{\frac{n-2}{2}}(q,K_{1})u(q)\geq C_{0},
		\]
		then we can apply Lemma \ref{lem:preparation lemma} to find
		a local maximum point  of $u$,  $q_{2} \in \overline{M}\backslash K_{1}$ satisfying \textit{(ii)}, and define $r_{2}=Ru^{-2/(n-2)}(q_{2})$.
		
		For any fixed $\ve\in(0,1)$, it follows from \textit{(ii)} that there exists
		a positive constant $c_{1}=c_{1}(\ve)\to0$ as $\ve\searrow0$
		such that one of three cases would occur corresponding to $(a_1),(a_2),(a_3)$:
		\begin{itemize}
			\item[($b_1$)] $u\circ\exp_{q_{2}}(u^{-2/(n-2)}(q_{2})y)$
			has no local maximum point in $|y|\ge c_{1}$; 
			\item[($b_2$)] $u\circ G_{q_{2}'}(u^{-2/(n-2)}(q_{2})y+d_2 e_n)$ 
			has no local maximum point in $|y|\ge c_{1}$; 
			\item[($b_3$)] $u\circ G_{q_{2}}(u^{-2/(n-2)}(q_{2})y)$ 
			has no  local maximum point in $|y+\tilde T_c^{(2)} \mathbf{e}_n|\ge \tilde T_c^{(2)}+ c_{1}$.
		\end{itemize}
		
		To see $\overline{B_{r_{2}}(q_{2})}\cap\overline{B_{r_{1}}(q_{1})}=\emptyset$, it suffices to show that $q_{1}\notin \overline{B_{2r_{2}}(q_{2})}$. 
		
		Suppose not, there holds $d_{g}(q_{1},q_{2})\le2Ru^{-2/(n-2)}(q_{2})$. For cases $(b_1)$ and $(b_2)$, noticing that $q_{1}$ is also a local maximum point of $u$, we have $d_{g}(q_{1},q_{2})\le c_{1}u^{-2/(n-2)}(q_{2})$. This together with \textit{(ii)} (either $(a_1)$ or $(a_2)$)
		implies that $u(q_{2})\le u(q_{1})\le2u(q_{2})$ for sufficiently small $\ve$. Thus,  we obtain $d_{g}(q_{1},q_{2})\le c_{1}u^{-2/(n-2)}(q_{2})\le 2^{2/(n-2)}c_1u^{-2/(n-2)}(q_{1})<r_{1}$,
		which contradicts $q_{2}\notin\overline{B_{r_{1}}(q_{1})}:=K_1$. For case $(b_3)$,  $q_2\in \partial M$. By choosing $\ve$ small enough we have $|y+\tilde T_c^{(2)} \mathbf{e}_n|>\tilde T_c^{(2)}+c_1$ for $y=u(q_2)^{2/(n-2)} G_{q_2}^{-1}(q_1)$. However, by $(b_3)$ we know that $q_1$ is not a local maximum of $u$. This gives us a contradiction.
		
		We may continue to apply the similar argument to find other points $q_{j}$. 
		
		Suppose that we have found a finite set of disjoint balls $\{\overline{B_{r_j}(q_j)}; 1 \leq j \leq i\}$ with each $q_j$ satisfying \textit{(ii)}. Let $K_i=\cup_{j=1}^{i}\overline{B_{r_j}(q_j)}$ and assume 
		$$
		\max_{q\in\overline{M\backslash K_{i}}}d_{g}^{\frac{n-2}{2}}(q,K_{i})u(q)\geq C_{0},
		$$ 
		then we use Lemma \ref{lem:preparation lemma} to find a local maximum point of $u$, $q_{i+1} \in \overline M \backslash K_i$ satisfying \textit{(ii)} and define  $r_{i+1}=Ru^{-2/(n-2)}(q_{i+1})$. 	To see $K_i \cap \overline{B_{r_{i+1}}(q_{i+1})}=\emptyset$. If not, there exists some $1 \leq j_0\leq i$ such that
		$$B_{r_{i+1}}(q_{i+1})\cap B_{r_{j_0}}(q_{j_0})\neq \emptyset.$$
		Now we choose 
		$$u(q_{i_0})=\min\{u(q_{j_0}),u(q_{i+1})\}.$$
		It is enough to replace $q_2,q_1$ in the above argument by $q_{i_0}$ and the other point, respectively. This give us a contradiction.

		We continue this procedure and it will stop after finitely many steps,
		since there exists a positive constant $a(n)$ such that $\int_{B_{r_{i}}(q_{i})}|\nabla u|^{2}\ud\mu_{g}\ge a(n)$ by virtue of \textit{(ii)}.
		Finally, we can find a positive integer $N$ and $\{q_{1},\cdots,q_{N}\}\subset\overline{M}$
		satisfying \textit{(ii)}. Meanwhile, 
		\[
		d_{g}^{\frac{n-2}{2}}(q,\cup_{i=1}^{N}\overline{B_{r_{i}}(q_{i})})u(q)\le C_{0}.
		\]
		
		If $q\in\overline{B_{2r_{i}}(q_{i})}$ for some $i$, then 
		\[
		d_{g}(q,\{q_{1},\cdots q_{N}\})\le d_{g}(q,q_{i})\le2r_{i}.
		\]
		By \textit{(ii)}, choosing $\ve$ sufficiently small we have
		\[
		u(q)\le2u(q_{i})=2R^{\frac{n-2}{2}}r_{i}^{-\frac{n-2}{2}}.
		\]
		Thus, we obtain 
		\[
		d_{g}^{\frac{n-2}{2}}(q,\{q_{1},\cdots,q_{N}\})u(q)\le2(2R)^{\frac{n-2}{2}}.
		\]
		
		If $q\notin\cup_{i=1}^{N}\overline{B_{2r_{i}}(q_{i})}$, then it is
		not hard to check that 
		\[
		d_{g}(q,\{q_{1},\cdots q_{N}\})\le2d_{g}(q,\cup_{i=1}^{N}\overline{B_{r_{i}}(q_{i})}),
		\]
		whence, 
		\[
		d_{g}^{\frac{n-2}{2}}(q,\{q_{1},\cdots, q_{N}\})u(q)\le2^{\frac{n-2}{2}}C_{0}.
		\]
		
		Hence, we choose $C_{1}=\max\{2(2R)^{\frac{n-2}{2}},2^{\frac{n-2}{2}}C_{0}\}$
		such that 
		\[
		u(q)\le C_{1}d_{g}^{-\frac{n-2}{2}}(q,\{q_{1},\cdots q_{N}\}),\qquad\forall~q\in\overline{M}.
		\]
		
		We finish the whole procedure. 
	\end{proof}
	
	Lemma \ref{lem fudament-1-1} gives a rough description of large solutions to \eqref{PDE:sigma_2_new}. Particularly emphasize that the positive integer $N$ may depend on $u$. We pause here for some historical notes. Schoen initiated a blow-up approach for the scalar curvature equation on closed manifolds, the statement similar to above Lemma \ref{lem fudament-1-1} for scalar curvature equation first appeared in Schoen and Zhang \cite[Lemma 3.1]{Schoen-Zhang}. 
	In the study of compactness of boundary Yamabe problem on locally conformally flat manifolds with umbilic boundary,  Han and Li \cite[Proposition 1.1]{han-li2} set up a more delicate procedure of seeking finitely many concentration points for each large solution meanwhile coming with more deep insights. 
		
	Suppose $\{u_{i}\}$ is a sequence of positive smooth solutions to \eqref{PDE:sigma_2_new}, more delicate characterization demands a precise definition of isolated blow-up points for $u_i$.

	\begin{definition}\label{def:isolated blow up} 
		
		A blow-up point $x_{i}\rightarrow x_{0}\in\overline{M}$
		for $u_{i}$ is isolated, if there exist positive constants $\delta,C$ such that 
		\[
		u_{i}(x)\leq Cd_{g}\left(x,x_{i}\right)^{-\frac{n-2}{2}}\qquad\mathrm{~~for~~all~~}x\in\overline{M},~d_{g}\left(x,x_{i}\right)<\delta,
		\]
		where $x_{i}$ is a local maximum point of $u_{i}$ and $u_{i}(x_{i})\rightarrow+\infty.$
	\end{definition}

	\begin{lemma} \label{lem fudament} 
	 Let $x_{i}\rightarrow x_0 \in\overline{M}$
		be an isolated blow-up point for ${u_{i}}$. Given two sequences of positive real
		numbers $R_{i}\rightarrow\infty$, $\varepsilon_{i}\rightarrow0$,
		there exists a subsequence of  $u_{i}$ such that
		\[
		r_{i}:=R_{i}u_{i}^{-\frac{2}{n-2}}\left(x_{i}\right)\rightarrow0,
		\]
		and either Case 1: for $x_{i}\in M$, 
		\[
		\left\Vert u_{i}^{-1}\left(x_{i}\right)u_{i}\circ \exp_{x_{i}}(u_{i}^{-\frac{2}{n-2}}\left(x_{i})y\right)-\left(1+\frac{f(x_0)}{2\sqrt{\binom{n}{2}}}|y|^{2}\right)^{\frac{2-n}{2}}\right\Vert _{C^{2}(B_{2R_{i}})}<\varepsilon_{i}
		\]
		or Case 2: for $x_{i}$ near  $\partial M$,
		\[
		\left\Vert u_{i}^{-1}\left(x_{i}\right)u_{i}\circ G_{x_{i}'}(\frac{y}{u_{i}^{\frac{2}{n-2}}\left(x_{i}\right)}+d_i\mathbf{e}_n)-(\frac{2b\sqrt{\binom{n}{2}}}{f(x_0)})^{\frac{n-2}{4}}\left(1+b|y+\tilde T_c\mathbf{e}_{n}|^{2}\right)^{\frac{2-n}{2}}\right\Vert _{C^{2}(B^+_{2R_{i}})}<\varepsilon_{i},
		\]
		where 
		
		\[
		b=\frac{f(x_0)}{2\sqrt{\binom{n}{2}}}\left(1+\frac{c(x_0)^2}{2f(x_0)}\sqrt{\binom{n}{2}}\right)^2\quad\mathrm{~~and~~}\quad \tilde T_c=\frac{c(x_0)}{\sqrt{2f(x_0)}}\frac{(\binom{n}{2})^{1/4}}{\sqrt{b}},
		\] and  $d_i=d_g(x_i, \partial M)=d_g(x_i,x_i')$ for some $x_i' \in \pa M$. 
		
		If $T_i=d_iu_i(x_i)^{2/(n-2)}\rightarrow \infty$, then Case 1 happens. If $T_i=d_iu_i(x_i)^{2/(n-2)}$ is bounded, then Case 2 happens.
	\end{lemma}
	
	The proof of Lemma \ref{lem fudament} is a standard argument of limitation as before, so we omit it.
	\medskip

	\begin{definition} \label{def:interior_isolated_simple}
		
		An isolated blow-up point $x_{i}\rightarrow x_{0}\in M$
		for $u_{i}$ is simple, if there exist constants $\tilde{r},C>0$ such that
		\begin{equation}
			\overline{w}_{i}'(r)<0\qquad\mathrm{~~for~~all~~}Cu_{i}(x_{i})^{-\frac{2}{n-2}}\le r\le\tilde{r},\label{eq:isolated blow up-1}
		\end{equation}
		where $\overline{w}_{i}(r):=r^{\frac{n-2}{2}}\bar{u}_{i}(r)$ and
		\[
		\bar{u}_{i}(r)=\frac{1}{r^{n-1}}\int_{\exp_{x_{i}}(\partial B_{r})}u_{i}\ud\sigma_{g}.
		\]

	\end{definition}
	
	A `\emph{fine}' property for $\bar u_i$ is as follows: If we let $v_{i}(z)=s^{(n-2)/2}u_{i}(\exp_{x_{i}}(sz)):=s^{(n-2)/2}u_i\circ\mathcal{T}_{i}(z)$
	for $s\in\R_{+}$, then 
	\[
	r^{\frac{n-2}{2}}\bar{v}_{i}(r)=(sr)^{\frac{n-2}{2}}\bar{u}_{i}(sr),
	\]
	where 
	$$\bar v_i(r)=\frac{1}{r^{n-1}}\int_{\pa B_r} v_i(z) \ud \sigma_{s^{-2}\mathcal{T}_{i}^\ast(g)}.$$

	\begin{definition}[New boundary  integral average]\label{def:bdry_blowup_pts}

		An isolated blow-up point $x_{i}\rightarrow x_{0}\in\partial M$ for
		$u_{i}$ is simple, if there exist constants $\tilde{r},C>0$ such that 
		\[
		\overline{w}_{i}'(r)<0,\qquad\mathrm{~~for~~all~~}Cu_{i}(x_{i})^{-\frac{2}{n-2}}\le r\le\tilde{r},
		\]
		where $\overline{w}_{i}(r):=r^{\frac{n-2}{2}}\bar{u}_{i}(r)$ and fix $\theta\in(0,\pi/2)$,
		\[
		\bar{u}_{i}(r)=\frac{1}{r^{n-1}}\int_{G_{x_i'}(\pa^{+}B_{r,\theta}^{+}+d_i\mathbf{e}_n)}u_{i}\ud\sigma_{g},
		\]
		where \footnote{Instead of the average over the hemi-sphere in the boundary Yamabe problem, we use the average over a small spherical cap, since we can not obtain $C^2$ estimates up to the boundary for  degenerate boundary $\sigma_2$-curvature equations as indicated by Lemma \ref{lem:double_normal_derivatives}.  Otherwise, it will bring us additional trouble verifying two properties: \eqref{average_limit_fcn1} and \eqref{average_limit_fcn2}, provided that the average over the hemi-sphere was used.}$\partial^{+}B_{r,\theta}^{+}=\{z\in\Rn;|z|=r,\mathrm{dist}_{\Sp^{n-1}}(N,z/r)<\theta\}$
		with the north pole $N$, and
		\,$d_i=d_g(x_i, \partial M)=d_g(x_i,x_i')$ for some $x_i' \in \pa M$.
	\end{definition}
	Moreover, if we let
	$$v_{i}(z)=s^{(n-2)/2}u_{i}(G_{x_i'}(sz+d_i\mathbf{e}_n)):=s^{(n-2)/2}u_i\circ \Psi_{i}(z)$$
	and
	$$\bar{v}_{i}(r)=\frac{1}{r^{n-1}}\int_{\pa^{+}B_{r,\theta}^{+}}v_i(z)\ud \sigma_{s^{-2}\Psi_{i}^\ast(g)}\qquad \mathrm{~~for~~}\quad  r\in\R_{+},$$
	then a `\emph{fine}' property for $\bar{v}_{i}$ is given by
	\begin{equation}\label{property_bdry_blowup_pt}
		r^{\frac{n-2}{2}}\bar{v}_{i}(r)=(sr)^{\frac{n-2}{2}}\bar{u}_{i}(sr)\qquad \mathrm{for~~} s \in \R_+.
	\end{equation}
	Indeed, this follows from
	\begin{align*}
		(sr)^{\frac{n-2}{2}}\bar{u}_{i}(sr)=&\frac{1}{(sr)^{n-1}}\int_{\pa^{+}B_{sr,\theta}^{+}}(sr)^{\frac{n-2}{2}}u_i(G_{x_i'}(y+d_i\mathbf{e}_n))\ud \sigma_{G_{x_i'}^\ast(g)}\\
		=&\frac{1}{r^{n-1}}\int_{\pa^{+}B_{r,\theta}^{+}}(sr)^{\frac{n-2}{2}}u_i(G_{x_i'}(sz+d_i\mathbf{e}_n))\ud \sigma_{s^{-2}\Psi_{i}^\ast(g)}\\
		=&r^{\frac{n-2}{2}}\frac{1}{r^{n-1}}\int_{\pa^{+}B_{r,\theta}^{+}}v_i(z)\ud \sigma_{s^{-2}\Psi_{i}^\ast(g)}=r^{\frac{n-2}{2}}\bar{v}_{i}(r).
	\end{align*}
	
	\begin{lemma}[Harnack inequality] \label{lem:harnack inequality} 
		Let $x_{i}\rightarrow x_{0}\in\pa M$ be an isolated blow-up point for $u_i$.
		For $0<r\le\rho/2$, there holds
		
		\[
		\max_{G_{x_i'}(\overline{B_{2r}^{+}}\backslash B_{r/2}^{+})}u_{i}\leq C\min_{G_{x_i'}(\overline{B_{2r}^{+}}\backslash B_{r/2}^{+})}u_{i},
		\]
		where $C$ is a positive constant independent of $i,r$. \end{lemma} 
	\begin{proof}
		We define 
		\[
		\Psi_{i}(z)=G_{x_i'}(rz+d_i\mathbf{e}_n)
		\]
		and 
		\[
		\tilde{u}_{i}(z)=r^{\frac{n-2}{2}}u_{i}\circ\Psi_{i}(z).
		\]
		
		Let $g_{i}=r^{-2}\Psi_{i}^{\ast}g$, then it is not hard to see that
		\[
		g_{\tilde{u}_{i}}=\tilde{u}_{i}^{\frac{4}{n-2}}g_{i}=\Psi_{i}^{\ast}(u_{i}^{\frac{4}{n-2}}g)=\Psi_{i}^{\ast}(g_{u_{i}}).
		\]
		
		Denote $A^{+}_i(1/2,2):=\overline{B_{2}^{}}\backslash B_{1/2}^{}(-\frac{d_i\mathbf{e}_n}{r})\cap \{z_n\ge -\frac{d_i}{r}\}$
		for simplicity. Thus, we obtain 
		\[
		\sigma_{2}^{1/2}(g_{i}^{-1}A_{\tilde{u}_{i}})=f\circ\Psi_{i}\tilde{u}_{i}^{\frac{4}{n-2}}\quad\mathrm{~~in~~}\quad A^{+}_i(1/2,2)
		\]
		and on $\{z_{n}=-\frac{d_i}{r}\}\cap\pa A_i^{+}(1/2,2)$, 
		\[
		\frac{2}{n-2}\frac{\partial\tilde{u}_{i}}{\partial\vec{n}_{g_{i}}}=-c\circ\Psi_{i}\tilde{u}_{i}^{\frac{n}{n-2}}+h_{g_{i}}\tilde{u}_{i}.
		\]
		Here $\vec{n}_{g_{i}}=r\vec{n}\circ\Psi_{i}$.
		
		By definition (\ref{def:isolated blow up}) of isolated blow-up point,
		for $0<r<\rho/2$ and $z\in B_{2}(-\frac{d_i\mathbf{e}_n}{r})\cap \{z_n\ge -\frac{d_i}{r}\}$, we have 
		\begin{align*}
			\tilde{u}_{i}(z) & \le r^{\frac{n-2}{2}}u_{i}(G_{i}(rz))\\
			& \le Cr^{\frac{n-2}{2}}d_{g}\left(G_{i}(rz),x_{i}\right)^{-\frac{n-2}{2}}\\
			& \le C|z|^{-\frac{n-2}{2}}.
		\end{align*}
		
		By Theorem \ref{Thm:local $C^1$ estimates} we have
		\[
		\|\nabla\log\tilde{u}_{i}\|_{C^{0}(A_i^{+}(1/2,2))}\le C(\sup\tilde{u}_{i}^{\frac{2}{n-2}}+1)\le C.
		\]
		
		Let		\[
		\max_{A_i^{+}(1/2,2)}\tilde{u}_{i}=\tilde{u}_{i}(z_{0})\qquad\mathrm{~~and~~}\qquad\min_{A_i^{+}(1/2,2)}\tilde{u}_{i}=\tilde{u}_{i}(z_{1}).
		\]
		Then for a geodesic $\gamma(t):[0,1]\to\overline{M}$ joining $z_{0}$
		and $z_{1}$, with a positive constant $C_{1}$ there holds $\int_{0}^{1}|\gamma'(t)|\ud t\le C_{1}$.
		Observe that 
		\begin{align*}
			\log\big(\max_{A_i^{+}(1/2,2)}\tilde{u}_{i}\big)-\log\big(\min_{A_i^{+}(1/2,2)}\tilde{u}_{i}\big)=  \int_{0}^{1}\langle\nabla\log\tilde{u}_{i},\gamma'(t)\rangle\ud t\le CC_{1}.
		\end{align*}
		This means 
		\[
		\frac{\max_{A_i^{+}(1/2,2)}\tilde{u}_{i}}{\min_{A_i^{+}(1/2,2)}\tilde{u}_{i}}\le C
		\]
		and the desired assertion follows. 
	\end{proof}

	Since our definition of isolated blow-up point is slightly different from others, Lemma \ref{lem fudament} enables us to give a proof of the following elementary property for each isolated blow-up point, which might be well-known in other classical problems and thus is left to Appendix \ref{Appendix:B}.
	
	\begin{lemma} \label{one critical point}
		Let
		$x_{i}\rightarrow x_{0}\in\overline{M}$ be an isolated blow-up point
		for ${u_{i}}$. Then for sufficiently large $i$, $r^{(n-2)/2}\overline{u}_{i}(r)$ has only
		one critical point in $(0,r_{i})$, where $r_i$ is given in Lemma \ref{lem fudament}.
	\end{lemma}

	\subsection{Isolated singularity behavior}\label{Subsect:isolated_singularity}

	From now on, we restrict consideration to the lower dimensions $n=3,4$. 
	We first set up some notation.  If $x_0$ is an interior point,  then we take $r(x)=d_{g}(x,x_0)$ and computations are employed under geodesic normal coordinates
	around $x_0$; if $x_0\in\partial M$, then we take $r(x)^{2}=|x'|^{2}+x_{n}^{2}$ under
	Fermi coordinates around $x_0$. To emphasize the fixed point $x_0$, we sometimes denote the above $r(x)$ by $r_{x_0}(x)$. 
	
	Near an interior isolated blow-up point, the following interior singularity behavior has been described  in  Li-Nguyen \cite{Li-Luc2}.

	\begin{lemma}[\protect{\cite[Lemma 3.4]{Li-Luc2}}]
		Suppose $n=3,4$ and  $x_0\in M$. There exists a positive constant  $r_{1}$ depending
		on $g$ such that for all $0<\delta<1/2$, the function
		$\underline{v}_{\delta}(r):=r^{-(n-2-\delta)}e^{r}$ satisfies 
		\[
		\lambda\left(A_{\underline{v}_{\delta}^{\frac{4}{n-2}}g}\right)\in\mathbb{R}^{n}\backslash\overline{\Gamma_2^{+}}\quad\mathrm{~~in~~}\left\{ 0<r<r_{1}\right\}.
		\]
	\end{lemma}
	We  next give boundary singularity behavior near a boundary blow-up point under the boundary condition that $\frac{2}{n-2}\frac{\partial u}{\partial\vec{n}}=-cu^{\frac{n}{n-2}}+h_g u$, which becomes more subtle due to  the effect of the second fundamental form.
	
	\begin{lemma} \label{lem:four-dim lower bound-boundary}
		Suppose $n=3, 4$ and   $x_0\in \partial M$. For any $0<\delta<1/2$, there exists a positive constant 
		$r_{1}:=r_1$ depending on $\delta,g$ and $L_{\alpha \beta}(x_0)$  such that
		the function 
		
		\[
		\underline{v}_{\delta}(r)=r^{-(n-2-\delta)}e^{br}, \quad  b\in \R_+
		\]
		satisfies
		\[
		\lambda\left(g^{-1}A_{\underline{v}_{\delta}^{\frac{4}{n-2}}g}\right)\in\mathbb{R}^{n}\backslash\overline{\Gamma_2^{+}}\qquad\mathrm{~~in~~}\left\{ 0<r(x)<r_{1}\right\} 
		\]
		and $\frac{\partial v_{\delta}}{\partial\vec{n}}=0$
		on $\partial M\cap\left\{ 0<r(x)<r_{1}\right\}$. Here  $b=1$ for $n=3$ and  $b$ is a positive constant depending on $g$ and $L_{\alpha \beta}(x_0)$ for $n=4$.
	\end{lemma}
	
	For clarity, we leave the proof of the above lemma to Appendix \ref{Appendix:B}.

	\begin{remark}
		For $x_0$ close to the boundary, we take $r=r_{x_i'}$ with $d_g(x_i, x_i')=d_g(x_i, \partial M)$ for some $x_i' \in \pa M$, then Lemma \ref{lem:four-dim lower bound-boundary} still holds.
	\end{remark}


	\begin{lemma}\label{lem:superharmonic}
		Let $\Omega$ be an open neighborhood of $p\in\overline{M}$ and  $L_{g}=\Delta_{g}-(n-2)R_{g}/[4(n-1)]$ denote the conformal Laplacian
		of $g$. Suppose $R_g\geq 0$ in $M$ and $h_{g}\le 0$ on $\partial M$. If $w\in LSC(\Omega\backslash\{p\})$ is a non-negative function in $\Omega\backslash\{p\}$ and satisfies $L_{g}w\leq0$ in the viscosity
		sense in $\Omega\backslash\{p\}$  and $\frac{2}{n-2}\frac{\partial w}{\partial\vec{n}}=h_{g} w-cw^{\frac{n}{n-2}}$
		on $\Omega\cap\partial M$ (if nonempty) away from $p$, then with $r(x)=d_{g}(x,p)$ for $p\in M$, there holds
		\[
		\lim_{r\rightarrow0}\left(r^{n-2}\min_{\partial B_{r}(p)\cap M}w\right)\quad\mathrm{~~exists~~and~~is~~finite;}
		\]
		with $r(x)^{2}=|x'|^{2}+x_{n}^{2}$ under Fermi coordinates around
		$p\in\pa M$, there holds
		\[
		\lim_{r\rightarrow0}\left(r^{n-2}\min_{G_{p}(\pa^{+}B_{r}^{+})}w\right) \quad \mathrm{~~exists~~and~~is~~finite}.
		\]
		
	\end{lemma} 
	\begin{proof}
		The first assertion has been available in Li-Nguyen \cite[Lemma 3.3]{Li-Luc2}.
		For $p\in\pa M$, we introduce an auxiliary function $\underline{w}(\rho)G(r)/G(\rho)$,
		where 
		\[
		G(r)=r^{2-n}-Kr^{\frac{5}{2}-n}-(\delta^{2-n}-K\delta^{\frac{5}{2}-n}) \qquad \mathrm{~~in~~}B_\delta^+
		\]
		for some $K \in \R_+$ and $\underline{w}(r)=\min_{G_{p}(\pa^{+}B_{r}^{+})}w$.
		
		Under Fermi coordinates around $p$, for sufficiently small $\delta$ and sufficiently large $K$,
		a direct computation yields 
		\[
		L_{g}G=\frac{1}{2}K(n-\frac{5}{2})r^{\frac{1}{2}-n}+O(r^{1-n})>0, ~~G>0  \qquad \mathrm{~~in~~} B_{\delta}^{+}.
		\]
		Moreover, we notice that
		\begin{align*}
			G=0 & \quad\mathrm{~~on~~}\quad\pa^{+}B_{\delta}^{+};\\
			\frac{\underline{w}(\rho)}{G(\rho)}\frac{\partial G(r)}{\partial x_{n}}=0\geq\frac{\partial w}{\partial x_{n}}=-\frac{n-2}{2}cw^{\frac{n}{n-2}}+\frac{n-2}{2}h_g w & \quad\mathrm{~~on~~}\quad D_\delta \backslash\{0\}.
		\end{align*}
		
		By the strong maximum principle together with Hopf lemma, we know
		that
		\[
		w(x)\ge\frac{\underline{w}(\rho)}{G(\rho)}G(r)\qquad \mathrm{~~for~~all~~} \rho\leq r(x)\leq \delta.
		\]
		This indicates that $\underline{w}/G$ is non-decreasing and also
		nonnegative in $(0,\delta)$. Thus, $\lim_{r\rightarrow0}\underline{w}(r)/G(r)$
		exists and is finite. The second assertion follows. 
	\end{proof}
	
	In the following, we shall describe the behavior of viscosity
	solutions to a degenerate $\sigma_{2}$-curvature equation \eqref{PDE:degen_bdry_sigma_2} below
	with certain Neumann boundary condition. Although the condition $h_g=0$ is enough to our use,  we relax it here to $h_g \leq 0$.
	
	\begin{theorem}\label{Thm:isolated behavior}
		For $n=3,4$, let $u_{i}$ be a sequence of positive smooth solutions to \eqref{PDE:sigma_2_new}
		and $x_{i}\rightarrow x_{0}$ be an isolated blow-up
		point for $u_i$. Assume that $c\ge 0$ and $h_g\le 0$ on $\partial M$. Then there exists $v_{\infty}\in C_{\mathrm{loc}}^{1,1}(B_{r_{0}}(x_{0})\backslash\{x_{0}\})$
		such that for any fixed $y_{0}\in B_{r_{0}}(x_{0})\backslash\{x_{0}\}$,
		up to a subsequence we have for any $\ensuremath{0<\gamma<1}$,
		\begin{align*}
			u_{i}(x)\rightarrow0 & \quad\mathrm{~~in~~}\quad C_{\mathrm{loc}}^{1,\gamma}(M\cap B_{r_{0}}(x_{0})\backslash\{x_{0}\})\cap C_{\mathrm{loc}}^{\gamma}(\overline{M}\cap B_{r_{0}}(x_{0})\backslash\{x_{0}\}),\\
			\frac{u_{i}(x)}{u_{i}(y_{0})}\rightarrow v_{\infty} & \quad\mathrm{~~in~~}\quad C_{\mathrm{loc}}^{1,\gamma}(M\cap B_{r_{0}}(x_{0})\backslash\{x_{0}\})\cap C_{\mathrm{loc}}^{\gamma}(\overline{M}\cap B_{r_{0}}(x_{0})\backslash\{x_{0}\}).
		\end{align*}
		Moreover, $v_{\infty}$ is a positive Lipschitz viscosity solution
		to 
		\begin{equation}
			\begin{cases}
				\sigma_{2}(A_{v_{\infty}^{\frac{4}{n-2}}g})=0 & \quad\mathrm{~~in~~}\quad B_{r_{0}}(x_{0})\backslash\{x_{0}\},\\
				\frac{\partial v_{\infty}}{\pa\vec{n}}=\frac{n-2}{2}h_{g}v_{\infty} & \quad\mathrm{~~on~~}\quad (B_{r_{0}}(x_{0})\cap\partial M)\backslash\{x_{0}\}~~ (\mathrm{if~~} x_0 \in \pa M).
			\end{cases}\label{PDE:degen_bdry_sigma_2}
		\end{equation}
		
		If $x_0\in \partial M$, then there exists a positive constant $\mathcal{N}=\mathcal{N}(\delta)$ such that
		\[v_{\infty}(x)\geq  \frac{\mathcal{N}}{r_{x_0}^{n-2-\delta}(x)}
		\qquad\mathrm{~~in~~} \left\{ 0<r_{x_0}(x)\leq C_3\delta\right\}, \forall~~
		0<\delta<\frac{1}{2}.\]
		
		If $x_0\in M$, then $$\lim_{x\rightarrow x_{0}}v_{\infty}(x)r_{x_0}^{n-2}(x)=a \qquad \mathrm{for~~some~~}a \in \R_+.$$
	\end{theorem} 
	\begin{proof}
		Since $x_{i}\rightarrow x_{0}$ is an isolated blow-up point, by definition we have
		\begin{equation}
			u_{i}(x)\le C_{1}d_{g}(x,x_{i})^{-\frac{n-2}{2}},\quad\forall~x\in\overline{M},~~d_{g}(x,x_{i})<\rho.\label{eq:growth estimate for use-1}
		\end{equation}
		
		\emph{Step 1:} There holds $u_{i}\rightarrow0$ in $B_{r_{0}}(x_{0})\backslash\{x_{0}\}.$
		\vskip 8pt
		
		By contradiction,  using local $C^1$ estimates we know that up to a subsequence, $u_{i}\rightarrow u_{\infty}>0$ in $C_{\mathrm{loc}}^{\alpha}(B_{r_{0}}(x_{0})\backslash\{x_{0}\}\cap\overline{M})$ for any $0<\alpha<1$
		and thus
		\begin{equation}\label{limit decay}
			u_{\infty}(x)\le C_{1}d(x,x_{0})^{-\frac{n-2}{2}},\quad~\forall~ x\in\overline{M},~~d_{g}(x,x_0)<\rho.
		\end{equation}
		\begin{claim}
			(1) If $x_{0}\in M$, then for any $0<\delta<1/2$, there
			exist positive constants $\mathcal{N}$ and $r_{1}$ independent of $i$, and $K=K(\delta)$
			and  $i_0=i_0(\delta)\in \N$ such that for all $i\geq i_0$, we have 
			\[
			u_{i}(x)\geq\frac{\mathcal{N}}{d_{g}^{n-2-\delta}(x,x_{i})}\qquad\mathrm{~~in~~}\left\{ Ku_{i}\left(x_{i}\right)^{-\frac{2}{n-2}}\leq d_{g}(x,x_{i})\leq r_{1}\right\} ,
			\]
			where $r_{1}$ is given in Lemma \ref{lem:four-dim lower bound-boundary}.
			
			(2) If $x_{0}\in\partial M$, then for any $0<\delta<1/2$, there
			exist positive  constants $\mathcal{N}=\mathcal{N}(\delta)$ and $r_{1}$ independent of $i$, $K=K(\delta)>0$
			and $i_0=i_0(\delta)\in \N$ such that for all $i\geq i_0$, we have 
			\[
			u_{i}(x)\geq\frac{\mathcal{N}}{r_{x_i'}^{n-2-\delta}(x)}\quad\mathrm{~~in~~}\left\{ Ku_{i}\left(x_{i}\right)^{-\frac{2}{n-2}}\leq r_{x'_i}(x)\leq r_{1},z_{n}+d_i\ge 0, x=G_{x_i'}(z+d_i\mathbf{e}_n)\right\} ,
			\]
			where $d_i:=d_g(x_i, \partial M)=d_g(x_i',x_i)$ for some $x_i' \in \pa M$,  $r_{x'_i}(x)=|z+d_i\mathbf{e}_n|$ and $r_{1}$ is given in Lemma \ref{lem:four-dim lower bound-boundary}.  
		\end{claim}
		
		It suffices to prove Claim \emph{(2)}, since  the proof of Claim \emph{(1)} is identical to Li-Nguyen's. 
		
		We can find a positive constant $\mathcal{N}$ independent of $i$ such that
		\[
		u_{i}\geq \mathcal{N}\underline{v}_{\delta}\quad\mathrm{~~on~~}\left\{ r_{x_i'}=r_{1}\right\},
		\]
		where 	
		$$
		\underline{v}_{\delta}(x)=r_{x_i'}^{-(n-2-\delta)}(x)e^{br_{x_i'}(x)} 
		$$
		is given in Lemma \ref{lem:four-dim lower bound-boundary}.
		
		Fix a $K\in \R_+$, define
		\[
		\bar{\beta}=\sup\left\{ 0<\beta<\mathcal{N};u_{i}(x)\geq\beta\underline{v}_{\delta}(x)\quad\mathrm{in~~}\left\{ r_{i}:=Ku_{i}\left(x_{i}\right)^{-\frac{2}{n-2}}<r_{x_i'}(x)<r_{1},x_{n}\ge0\right\} \right\} .
		\]
		
		Our goal is to show that for all sufficiently large $i$, there exists a large $K=K(\delta) \in \R_+$ such that  $\bar{\beta}=\mathcal{N}$.
		
		By contradiction, suppose $\bar{\beta}<\mathcal{N}$. We can find
		a point $y_{i} \in \{r_{i}\le r_{x_i'}(x)\le r_{1},x_n\geq 0\}$ such that $u_{i}$
		touches $\bar \beta\underline{v}_{\delta}$ at $\hat x_{i}$, that is, 
		\begin{equation}\label{first touch}
			u_{i}\ge\bar \beta\underline{v}_{\delta} \quad \mathrm{~~and~~}\quad
			u_{i}(\hat x_{i})=\bar \beta\underline{v}_{\delta}(\hat x_{i}).
		\end{equation}
		
		First notice that $\hat x_{i} \notin \{r_{x_i'}=r_{1},x_n\geq 0\}$ due to the definition of $\bar{\beta}$. The following is to show that $r_{x_i'}(\hat x_i)=r_{i}$. 
		
		To this end, we assume $r_{i}< r_{x_i'}(\hat x_i)< r_{1}$ and need to exclude two possibilities:
		\begin{itemize}
			\item If $\hat x_{i}\in M$, then we have 
			$\nabla(u_{i}-\bar \beta\underline{v}_{\delta})=0$
			and $\nabla^{2}(u_{i}-\bar\beta\underline{v}_{\delta})\ge0$ at $\hat x_{i}$. 
			\item If $\hat x_{i}\in\partial M$, then $\frac{\partial}{\partial\vec{n}}(u_{i}-\bar\beta\underline{v}_{\delta})(\hat x_{i})\ge0$.
			Notice that on $\pa M$, 
			$$\frac{2}{n-2}\frac{\partial}{\partial\vec{n}}(u_{i}-\bar \beta\underline{v}_{\delta})=-cu_{i}^{\frac{n}{n-2}}+h_g u_i\le0$$
			by virtue of $\frac{\partial\underline{v}_{\delta}}{\partial\vec{n}}=0$ and $h_g\leq 0$,
			in particular there holds $\frac{\partial}{\partial\vec{n}}(u_{i}-\bar \beta\underline{v}_{\delta})(\hat x_{i})=0$,
			whence 
			\begin{equation}\label{first derivative vanish}
				\nabla(u_{i}-\bar \beta\underline{v}_{\delta})(\hat x_{i})=0.
			\end{equation}
			This
			forces $c(\hat x_{i})=0$, $h_g(\hat x_i)=0$ and then $\frac{\partial u_{i}}{\partial\vec{n}}(\hat x_{i})=0$.
			Otherwise, if $c(\hat x_{i})>0$, then $u_{i}(\hat x_{i})=0$, which contradicts
			the positivity of $u_{i}$. Furthermore, since $\hat x_{i}$ is also a
			maximum point of $\frac{\partial}{\partial\vec{n}}(u_{i}-\bar \beta\underline{v}_{\delta})$ on the boundary,
			we know that $\nabla_{\alpha}(u_{i}-\bar \beta\underline{v}_{\delta})_{n}(\hat x_{i})=0$  and
			it implies that $\nabla_{\alpha n}(u_{i}-\bar \beta\underline{v}_{\delta})(\hat x_{i})=0$ by \eqref{first derivative vanish}.
			By \eqref{first touch} and $\frac{\partial}{\partial\vec{n}}(u_{i}-\bar \beta\underline{v}_{\delta})(\hat x_{i})=0$, we have  $(u_{i}-\bar \beta\underline{v}_{\delta})_{nn}(\hat x_{i})\ge0$. Thus $\nabla^{2}(u_{i}-\bar \beta\underline{v}_{\delta})(\hat x_{i})\ge0$. 
		\end{itemize}
		In both cases, at $\hat x_{i}$ we obtain 
		\begin{align*}
			& -\frac{2}{n-2}\frac{\nabla_{g}^{2}u_{i}}{u_{i}}+\frac{2n}{(n-2)^{2}}\frac{\ud u_{i}\otimes\ud u_{i}}{u_{i}^{2}}-\frac{2}{(n-2)^{2}}\frac{|\nabla u_{i}|^{2}}{u_{i}^{2}}g+A_{g}\\
			\le & A_{\bar \beta\underline{v}_{\delta}}=A_{\underline{v}_{\delta}}=-\frac{2}{n-2}\frac{\nabla_{g}^{2}\underline{v}_{\delta}}{\underline{v}_{\delta}}+\frac{2n}{(n-2)^{2}}\frac{\ud\underline{v}_{\delta}\otimes\ud\underline{v}_{\delta}}{\underline{v}_{\delta}^{2}}-\frac{2}{(n-2)^{2}}\frac{|\nabla\underline{v}_{\delta}|^{2}}{\underline{v}_{\delta}^{2}}g+A_{g}.
		\end{align*}
		This reaches a contradiction between $\lambda(A_{u_{i}})\in\Gamma_{2}^{+}$
		and $\lambda(A_{\underline{v}_{\delta}})\in\Rn\backslash\overline{\Gamma_{2}^{+}}$.

		Consequently, we can find $\hat{x}_{i}$ with $r_{x_i'}(\hat{x}_{i})=r_{i}$
		such that $u_{i}$ touches $\beta\underline{v}_{\delta}$ at $\hat{x}_{i}$.
		
		By Lemma \ref{lem fudament}, we know either if $\lim_{i\rightarrow \infty}d_iu_i(x_i)^{2/(n-2)}=+\infty$, then
		\[
		\frac{1}{u_{i}\left(x_{i}\right)}u_{i}\circ \exp_{x_i}(\frac{y}{u_{i}\left(x_{i}\right)^{\frac{2}{n-2}}})\rightarrow\left(1+\frac{f(x_{0})}{2\sqrt{\binom{n}{2}}}|y|^{2}\right)^{\frac{2-n}{2}}\quad\mathrm{~~in~~}\quad C_{\mathrm{loc}}^{2,\alpha}(\mathbb{R}^{n});
		\]
		\\
		or if $\lim_{i\rightarrow \infty}d_iu_i(x_i)^{2/(n-2)}$ is bounded, then
		\[
		\frac{1}{u_{i}(x_{i})}u_{i}\circ G_{x_i'}(\frac{y}{u_{i}(x_{i})^{\frac{2}{n-2}}}+d_i\mathbf{e}_n)\rightarrow\left(\frac{2b\sqrt{\binom{n}{2}}}{f(x_{0})}\right)^{\frac{n-2}{4}}(1+b|y+\tilde T_c\mathbf{e}_{n}|^{2})^{\frac{2-n}{2}}\quad\mathrm{~~in~~}\quad C_{\mathrm{loc}}^{2,\alpha}(\overline{\mathbb{R}_{+}^{n}}).
		\]
		This means that given $R\geq  1,\varepsilon>0$, for all sufficiently large $i$ there holds 
		\[	
		\|\frac{1}{u_{i}\left(x_{i}\right)}u_{i}\circ G_{x_i'}(\frac{y}{u_{i}\left(x_{i}\right)^{\frac{2}{n-2}}}+d_i\mathbf{e}_n)-\left(\frac{2b\sqrt{\binom{n}{2}}}{f(x_{0})}\right)^{\frac{n-2}{4}}(1+b|y+\tilde T_c\mathbf{e}_{n}|^{2})^{\frac{2-n}{2}}\|_{C^{2}(B_{R}^{+})}<\varepsilon.
		\]
		
		It suffices to consider the latter case $\lim_{i\rightarrow \infty }d_iu_i(x_i)^{2/(n-2)}$ is bounded, since the front case is relatively easier.

		For brevity, we denote 
		\[
		\tilde{u}_{i}(y)=	\frac{1}{u_{i}\left(x_{i}\right)}u_{i}\circ G_{x_i'}(\frac{y}{u_{i}\left(x_{i}\right)^{\frac{2}{n-2}}}+d_i\mathbf{e}_n)
		\]
		and 
		\[
		\varphi(y)=\left(\frac{2b\sqrt{\binom{n}{2}}}{f(x_{0})}\right)^{\frac{n-2}{4}}(1+b|y+\tilde T_c\mathbf{e}_{n}|^{2})^{\frac{2-n}{2}}.
		\]
		
		For any $\ve_{0}\in(0,1)$,
		there exists a positive constant $K\gg 1$ such that for all $|y|>K/2$, 
		\[
		\frac{b(|y|^{2}+\tilde T_c y_{n})}{1+b(|y'|^{2}+(y_{n}+\tilde T_c)^{2})}\ge1-\ve_{0}
		\]
		and \[ 	\frac{2b(|y'|+|y_{n}+\tilde T_c|)}{1+b(|y'|^{2}+(y_{n}+\tilde T_c)^{2})}\le \frac{4}{K}.
		\]
		
		Then for all $K/2<|y|<K$ with a larger $K$,  with a positive constant
		$C_{1}$ we have 
		\begin{equation}\label{est1:D_log_bubble}
			|\nabla_y \log\varphi(y)|\le \frac{2(n-2)}{K},
		\end{equation}
		and 
		\begin{align*}
			& |\partial_{y_k}\log\tilde{u}_{i}(y)-\partial_{y_k}\log\varphi(y)|\\
			= & |\frac{\partial_{y_k}\tilde{u}_{i}}{\widetilde{u}_{i}}-\frac{\partial_{y_k}\varphi}{\varphi}|=|\frac{\varphi\partial_{y_k}\tilde{u}_{i}-\tilde{u}_{i}\partial_{y_k}\varphi}{\tilde{u}_{i}\varphi}|\le\frac{|\varphi-\tilde{u}_{i}||\partial_{y_k}\varphi|+\varphi|\partial_{y_k}\tilde{u}_{i}-\partial_{y_k}\varphi|}{\tilde{u}_{i}\varphi}\\
			\le & C_{1}\varepsilon(\frac{1}{K}K^{n-2}+K^{n-2})\le  \frac{\delta}{4nK},
		\end{align*}
		where we take $\varepsilon, K$ such that $C_{1}\varepsilon(\frac{1}{K}K^{n-2}+K^{n-2})\leq\frac{\delta}{4nK}$. Thus for all $K/2<|y|<K$  we obtain
		\begin{equation}\label{upper bound-approximation}
			|\nabla_y \log\tilde{u}_{i}(y)|\le \frac{2(n-2)+\frac{\delta}{4n}}{K}.
		\end{equation}
		Notice that 
		\begin{align}\label{est2:D_log_bubble}
			& (y_k\partial_{y_k})\log \varphi(y)\no\\
			=&- \frac{n-2}{2}(y_{i}\partial_{y_{i}})\log[1+b(|y'|^{2}+(y_{n}+\tilde T_c)^{2})]\no\\
			= & -\frac{n-2}{2}\frac{2b(|y|^{2}+\tilde T_c y_{n})}{1+b(|y'|^{2}+(y_{n}+\tilde T_c)^{2})}\no\\
			\le & -(n-2)(1-\ve_{0}),
		\end{align}
		and thereby, for $K/2<|y|<K$,
		\begin{equation}\label{est:D_logu}
			y_k\partial_{y_k}\left[\log \tilde{u}_{i}(y)\right]<\frac{\delta }{4}-(n-2)(1-\varepsilon_{0}).
		\end{equation}

		Recall that $x=G_{x_i'}(z+d_i\mathbf{e}_n)$ and $r_{x'_i}^2(x)=|z+d_i\mathbf{e}_n|^2$. Let us denote  
		$$\hat{x}_i:=G_{x_i'}(\hat{z}_i+d_i\mathbf{e}_n)  \quad \mathrm{~~and~~then~~} \quad |\hat{z}_i+d_i\mathbf{e}_n|^2=r_i^2.$$ 
		If we let $z:=u_{i}\left(x_{i}\right)^{-2/(n-2)}y$ and $\hat{z}_i:=u_{i}\left(x_{i}\right)^{-2/(n-2)}\tilde{y}_{i}$ with $\tilde y_i=(\tilde y_i', (\tilde y_i)_n)$, then

		$$|\frac{\tilde{y}'_{i}}{u_{i}\left(x_{i}\right)^{\frac{2}{n-2}}}|^2+|\frac{(\tilde{y}_i)_n}{u_{i}\left(x_{i}\right)^{\frac{2}{n-2}}}+d_i|^2=\frac{K^2}{u_{i}\left(x_{i}\right)^{\frac{4}{n-2}}}.$$ 
		This in turn implies that $K>|\tilde{y}_i|>K/2$.

		The inner unit  normal with respect to $g$ at $\hat{x}_{i}$ is given by 
		\[\vec n(\hat{x}_{i})=\frac{g^{jk}(z_{k}+d_i\delta_{nk})\partial_{z_{j}}}{\|x\|}|_{x=\hat{x}_{i}}.\] 
		Here near $x_i'$,
		\[\|x\|^{2}:=g^{kl}(z_{k}+d_i\delta_{nk})(z_{l}+d_i\delta_{nl})=|z+d_i\mathbf{e}_n|^2(1+O(|z+d_i\mathbf{e}_n|)\] and
		\[g^{jk}(x)=\delta_{jk}+2\tilde{L}_{jk}(x'_i)(z_n+d_i)+O(|z+d_i\mathbf{e}_n|^2)|.\]

		Since $u_i \geq  \bar\beta\underline{v}_{\delta}$ in $\{r_i< r_{x_i'}(x)<r_1\}$ and $r_{x_i'}(\hat{x}_i)=r_i$, we have
		\begin{equation}\label{comparison}
			\nabla_{\vec n}\log u_{i}|_{\hat{x}_{i}}\ge \nabla_{\vec n}\log(\bar \beta\underline{v}_{\delta})|_{\hat{x}_{i}}.
		\end{equation}
		At $\hat x_i$, by \eqref{est:D_logu} we have
		\begin{align*}
			\nabla_{\vec n}(\log u_{i})=&\frac{(z_k+d_i\delta_{nk})\partial_{z_k}+O(|z+d_i\mathbf{e}_n|^2)\partial_{z_k}}{|z+d_i\mathbf{e}_n|(1+O(r_i))}\log u_{i}\nonumber\\
			=&\frac{u_{i}\left(x_{i}\right)^{-\frac{2}{n-2}}y_k+d_i\delta_{nk}}{r_i}\partial_{z_k} \log u_i+O(r_i)|\nabla_z \log u_i|\nonumber\\
			=&\frac{y_k\partial_{y_k}\log \tilde{u}_i(y)}{r_i}+d_iu_i(x_i)^{\frac{2}{n-2}}\frac{\partial_{y_n} \log \tilde{u}_i(y)}{r_i}+O(r_i) u_i(x_i)^{\frac{2}{n-2}}|\nabla_y \log \tilde u_i|\nonumber\\
			=&\frac{y_k\partial_{y_k}\log \tilde{u}_i(y)}{r_i}+O(\frac{1}{Kr_i})+O(1)\nonumber\\
			\le& \frac{\frac{\delta }{4}-(n-2)(1-\varepsilon_{0})}{r_i}+O(\frac{1}{Kr_i})+O(1),
		\end{align*}
		where the last equality follows from the boundness of $d_iu_i(x_i)^{2/(n-2)}$ and \eqref{upper bound-approximation}.
		Similarly, at $\hat x_i$, by \eqref{est1:D_log_bubble} and \eqref{est2:D_log_bubble} we have
		\begin{align*}
			\nabla_{\vec n}\log(\bar \beta\underline{v}_{\delta})
			&=[\frac{(z_k+d_i\delta_{nk})\partial_{z_k}+O(|z+d_i\mathbf{e}_n|^2)\partial_{z_k}}{|z+d_i\mathbf{e}_n|(1+O(r_i))}]\log(\bar\beta\underline{v}_{\delta})\nonumber\\
			&=\frac{z_k+d_i\delta_{nk}}{r_i}\partial_{z_k} \log \underline{v}_{\delta}+O(r_i)|\nabla_z \log \underline{v}_{\delta}|\nonumber\\
			&=-\frac{n-2-\delta}{r_i}+b+O(1).
		\end{align*}

		We now fix $\ve_{0}<\delta/[4(n-2)]$. It follows from (\ref{comparison}) that 
		$$\frac{\frac{\delta }{4}-(n-2)(1-\varepsilon_{0})}{r_i}+O(\frac{1}{Kr_i})+O(1)\ge -\frac{n-2-\delta}{r_i}+b+O(1).$$
		We may choose a larger $K$ depending on $\delta$ such that for all sufficiently large $i$
		$$O(\frac{1}{K})+O(r_i)+b r_i \leq \frac{\delta}{2}.$$
		Thus, we conclude that
		$$\frac{\delta }{4}-(n-2)(1-\varepsilon_{0})\ge -(n-2-\delta/2),$$ which contradicts the choice of $\varepsilon_{0}$.

		Up to now, we have established Claim \emph{(2)}, which clearly contradicts (\ref{limit decay}) as $i\rightarrow \infty$.
		
		Therefore, we conclude that $\bar \beta=\mathcal{N}$ and $u_{i}\rightarrow0$ in $C^0_{\mathrm{loc}}(B_{r_{0}}(x_{0})\backslash\{x_{0}\}\cap \overline{M}).$
		\vskip 8pt
		\emph{Step 2:} By local $C^{1}$- estimates,  we know 
		\[
		v_{i}(x)=\frac{u_{i}(x)}{u_{i}(y_{0})}\rightarrow v_{\infty}\quad\mathrm{~~in~~}\quad C_{\mathrm{loc}}^{\gamma}(B_{r_{0}}(x_{0})\backslash\{x_{0}\}\cap \overline{M}), \forall~ \gamma \in (0,1)
		\]
		and $v_{\infty}\in C_{\mathrm{loc}}^{0,1}(B_{r_{0}}(x_{0})\backslash\{x_{0}\})$ is the viscosity solution
		to the following PDE problem: 
		\begin{align*}
			\begin{cases}
				\sigma_{2}(A_{v_{\infty}^{\frac{4}{n-2}}g})=0 & \quad\mathrm{~~in~~}B_{r_{0}}(x_{0})\backslash\{x_{0}\}\cap \overline{M},\\
				\frac{\partial v_{\infty}}{\partial\vec{n}}=\frac{n-2}{2}h_{g}v_{\infty} & \quad\mathrm{~~on~~}\partial M\backslash \{x_0\},\quad\mathrm{in~~the~~viscosity~~sense}.
			\end{cases}
		\end{align*}
		Meanwhile, if $x_{0}\in M$, then $\|\log v_{i}\|_{C^{2}(K)}\le C$
		for any compact set $K\subset M\cap B_{r_{0}}(x_{0})\backslash\{x_{0}\}$.
		If $x_{0}\in\partial M$, then $\|\log v_{i}\|_{C^{1}(K)}^{2}\le C$
		for any compact set $K\subset\overline{M\cap B_{r_{0}}(x_{0})}\backslash\{x_{0}\}$.
		Thus $v_{\infty}\in C_{\mathrm{loc}}^{0,1}(\overline{M\cap B_{r_{0}}(x_{0})}\backslash\{x_{0}\})\cap C^{1,1}(M\cap B_{r_{0}}(x_{0})\backslash\{x_{0}\}).$
		\bigskip

		We just repeat the proof of \emph{Step 1} and find some postive
		constants $\mathcal{N}$ and $r_{1}$ independent of $i$, such that for any
		$\delta \in (0,1/2)$,  with some $K=K(\delta)\in \R_+$ and $i_0=i_0(\delta)\in \N$
	we have for all $i\geq i_0$, 
		\[
		v_{i}(x)\geq\frac{\mathcal{N}}{r_{x_i'}^{n-2-\delta}(x)}\quad\mathrm{~~in~~}\left\{ Ku_{i}\left(x_{i}\right)^{-\frac{2}{n-2}}\leq r_{x'_i}(x)\leq r_{1},z_{n}+d_i\ge 0, x=G_{x_i'}(z+d_i\mathbf{e}_n)\right\} .
		\]
		Clearly this implies that $v_{\infty}(x)\geq\frac{\mathcal{N}}{r_{x_0}^{n-2-\delta}(x)}$
		in $\left\{ 0<r_{x_0}(x)\leq r_{1}\right\} $.

		For $n=3,4$,  we have
		$v_{\infty}(x)\geq\frac{\mathcal{N}}{r_{x_0}^{n-2-\delta}(x)}$
		in $\left\{ 0<r_{x_0}(x)\leq C\sqrt{\delta}\right\} $ for some positive constant $C$ given in Lemma \ref{lem:four-dim lower bound-boundary}. 
	\end{proof}
	
	\begin{remark}
		Under the assumption that  $A_{g}\in\Gamma_{2}^{+}$ and $h_{g}\ge0$ on $\pa M$, Proposition \ref{prop:special_g} indicates that there exists a conformal metric of the background metric $g$, still denoted by itself, such that $A_{g}\in\Gamma_{2}^{+}$ and $h_{g}=0$ on $\pa M$.
	\end{remark}
	
	For future citation, we apply the same strategy above to separately state another theorem, where boundary blow-up points happening at a Euclidean domain,  and eventually can make an improvement compared with manifolds.
	
	\begin{theorem}\label{isolated behavior-1} 
		For $n=3,4$, assume that two sequences of constants $\ve_i \to 0$ and $D_i \rightarrow D\geq 0$. Let $B_{-D,r_0}^+:=\{x\in \Rn;|x-x_0|<r_0,x_{n}>-D\}$ and $x_0 \in \overline{B_{-D,r_0}^+}$~. Suppose a sequence of smooth metrics $g_{i}$ in $\overline {B_{-D_i,r_0}^+}$ converges to $g_{\mathrm{E}}$ in $C_{\mathrm{loc}}^{3}(\overline{B_{-D,r_{0}}^{+}})$ and $h_{g_i}\leq 0$ on $\{x_n=-D_i\}$. Let $u_{i}$ be a sequence
		of positive solutions to
		\begin{equation*}
			\begin{cases}
				\sigma_{2}^{1/2}(g_{i}^{-1}{A}_{u_{i}^{\frac{4}{n-2}}g_i})=f_i(x)u_{i}^{\frac{4}{n-2}}+\varepsilon_i f_{0},& \qquad\mathrm{~~in~~}B_{-D_i,r_0}^+,\\
				\frac{2}{n-2}\frac{\pa u_{i}}{\pa\vec{n}}=-c_i(x)u_{i}^{\frac{n}{n-2}}+h_{g_{i}}u_{i}, & \qquad\mathrm{~~on~~}x_{n}=-D_{i},
			\end{cases}
		\end{equation*}
		and $x_{i}\rightarrow x_0$ be an
		isolated blow-up point, where $f_i,c_i$ are smooth functions and $f_i \to f(x_0)$ in $C^3_{\mathrm{loc}}(\overline{B_{-D,r_{0}}^{+}})$ and $c_i \to c(x_0)$ in $C^3_{\mathrm{loc}}(\overline{B_{-D,r_{0}}^{+}})$ whenever $x_0 \in \{x_n=-D\}$.
		Then there exist $a\in\R_{+}$ and $v_{\infty}\in C_{\mathrm{loc}}^{0,1}(\overline{B_{-D,r_{0}}^{+}}\backslash\{x_0\})\cap C_{\mathrm{loc}}^{1,1}(B_{-D,r_{0}}^{+}\backslash\{x_0\})$
		such that for any fixed $y_{0}\in B_{-D,r_{0}}^{+}\backslash\{x_0\}$,
		we have $\text{for any \ensuremath{0<\gamma<1}}$,
		\begin{align*}
			u_{i}(x)\rightarrow0 & \quad\mathrm{~~in~~}\quad C_{\mathrm{loc}}^{\gamma}(\overline{B_{-D,r_{0}}^{+}}\backslash\{x_0\})\cap C_{\mathrm{loc}}^{1,\gamma}(B_{-D,r_{0}}^{+}\backslash\{x_0\}),\\
			\frac{u_{i}(x)}{u_{i}(y_{0})}\rightarrow v_{\infty} & \quad\mathrm{~~in~~}\quad C_{\mathrm{loc}}^{\gamma}(\overline{B_{-D,r_{0}}^{+}}\backslash\{x_0\})\cap C_{\mathrm{loc}}^{1,\gamma}(B_{-D,r_{0}}^{+}\backslash\{x_0\}),\\
			& \lim_{x\rightarrow 0}v_{\infty}(x)|x-x_0|^{n-2}=a \qquad \mathrm{for~~some~~} a \in \R_+.
		\end{align*}
		Moreover, $v_{\infty}$ is a positive Lipschitz viscosity solution
		to 
		\begin{equation}
			\begin{cases}
				\sigma_{2}(A_{v_{\infty}^{\frac{4}{n-2}}g_{\mathrm{E}}})=0 & \quad\mathrm{~~in~~}\quad B_{-D,r_{0}}^{+}\backslash\{x_0\},\\
				\frac{\partial v_{\infty}}{\pa\vec{n}}=0 & \quad\mathrm{~~on~~}\quad\{x; x_{n}=-D\}\backslash\{x_0\}.
			\end{cases}\label{PDE:degen_bdry_sigma_2-1}
		\end{equation}

	\end{theorem} 
	
	\begin{proof}
		
		Since the proof is in spirit similar to that of Theorem \ref{Thm:isolated behavior}, we only outline the proof. 
		
		Since $x_{i}\rightarrow x_{0}$ is an isolated blow-up point, we have
		for some $\rho>0$, $C_{1}>0$
		\begin{equation}
			u_{i}(x)\le C_{1}|x-x_{i}|^{-\frac{n-2}{2}},\quad\forall~x\in B_{i}:=B_{-D_i,r_0}^+,~~|x-x_{i}|<\rho.\label{eq:growth estimate for use-1-limit case}
		\end{equation}
		
		\emph{Step 1:} There holds $u_{i}\rightarrow0$ in $C^{0}(\overline{B_{-D,r_{0}}^{+}}\backslash\{x_0\}).$
		\medskip
		
		For simplicity, we assume that $D_{i}=0,x_0=0, D=0$.
		
		Suppose not, by gradient estimate, Harnack inequality and local $C^{2}$-estimate,
		we know that $u_{i}\rightarrow u_{\infty}>0$ in $C_{\mathrm{loc}}^{2,\alpha}(\overline{B_{r_{0}}^{+}}\backslash\{0\})$
		and thus $u_{\infty}(x)\le C_{1}|x|^{-\frac{n-2}{2}}$ in $\overline{B_{r_{0}}^{+}}\backslash\{0\}$.
		
		\begin{claim}
			For any $0<\delta<1/2$, there exist positive constants $\mathcal{N}$ and $r_{1}$
			independent of $i$, $K=K(\delta)$ and $i_0=i_0(\delta)\in \N$ such that
			for all $i\geq i_0$, we have 
			\[
			u_{i}(x)\geq\frac{\mathcal{N}}{|x|^{n-2-\delta}}\qquad\mathrm{~~in~~}\left\{x \in \overline{\Rn_+}; Ku_{i}\left(x_{i}\right)^{-\frac{2}{n-2}}\leq|x|\leq r_{1}\right\} ,
			\]
			where $r_{1}$ is given in Lemma \ref{lem:four-dim lower bound-boundary}.
		\end{claim}
		
		Assuming the above claim temporarily, we obtain a contradiction with (\ref{eq:growth estimate for use-1-limit case}) and thus the desired assertion follows.

		To prove the above claim, we nearly follow the same strategy as in Theorem \ref{Thm:isolated behavior}. 
		
		There exists a positive constant $\mathcal{N}$ independent of $i$ such that
		\[
		u_{i}\geq \mathcal{N} \underline{v}_{\delta}\qquad\mathrm{~~on~~}\left\{ r(x):= |x|=r_{1},x_n>0\right\},
		\]
		where $\underline v_\delta$ is defined in  Lemma \ref{lem:four-dim lower bound-boundary}.
		
		Fix a $K \in \R_+$, define 
		\[
		\bar{\beta}=\sup\left\{ 0<\beta<\mathcal{N}; u_{i}\geq\beta\underline{v}_{\delta}\quad\mathrm{~~in~~}\left\{ r_{i}:=Ku_{i}\left(x_{i}\right)^{-\frac{2}{n-2}}<|x|<r_{1},x_{n}\ge0\right\} \right\} .
		\]
		
		Our purpose is to show that for all sufficiently large $i$, there exists a large $K=K(\delta) \in \R_+$ such that  $\bar{\beta}=\mathcal{N}$.

		We argue by contradiction. Suppose $\bar{\beta}<\mathcal{N}$. With the help of Lemma \ref{lem:four-dim lower bound-boundary}, we also can establish that the first touch point of $u_{i}$ and $\bar \beta\underline{v}_{\delta}$ only occurs at some $\hat{x}_{i}$ with $\left|\hat{x}_{i}\right|=r_{i}$.
		
		The inner unit  normal with respect to $g_i$ at $\hat{x}_{i}$ is given by
		$$\vec n_{g_i}(\hat x_i)=\frac{g_{i}^{jk}x_{k}\partial_{x_{j}}}{\|x\|^{2}}|_{x=\hat{x}_{i}},$$
		where
		$\|x\|^{2}:=g_{i}^{kl}x_{k}x_{l}$.
		
		Since $u_i \geq  \bar\beta\underline{v}_{\delta}$ in $\{r_i<|x|<r_1,x_n>0\}$, we have
		\[
		\nabla_{\vec n_{g_i}}\log u_{i}|_{\hat{x}_{i}}\ge\nabla_{\vec n_{g_i}}\log(\bar \beta\underline{v}_{\delta})|_{\hat{x}_{i}}.
		\]
		Let  $\partial_{r}:=\frac{1}{r}x_{k}\partial_{x_{k}}$,  there
		exists a positive constant $C$ such that 
		\begin{equation}
			\partial_{r}\log u_{i}|_{\hat{x}_i}\ge\partial_{r}\log(\bar \beta\underline{v}_{\delta})|_{\hat{x}_i}+C\geq -\frac{n-2-\delta}{r_{i}}.\label{eq:lower bound-1-limit case}
		\end{equation}
		
		Write $\hat{x}_{i}:=u_{i}\left(x_{i}\right)^{-2/(n-2)}\tilde{y}_{i}$
		with $|\hat x_i|=r_{i}$, which means $|\tilde{y}_{i}|=K$.
		
		Without loss of generality, we just consider the case: $|x_i|u_i(x_i)^{\frac{2}{n-2}}<\infty$. 			Given $R,\varepsilon>0$, for all large $i$ there holds
		\[
		\|\frac{1}{u_{i}\left(x_{i}\right)}u_{i}(u_{i}\left(x_{i}\right)^{-\frac{2}{n-2}}y)-\left(\frac{2b\sqrt{\binom{n}{2}}}{f(x_0)}\right)^{\frac{n-2}{4}}(1+b|y+\tilde T_c\mathbf{e}_{n}|^{2})^{\frac{2-n}{2}}\|_{C^{2}(B_{R}^{+})}<\varepsilon,
		\]
		where $b,\tilde T_c$ are given in Lemma \ref{lem fudament}.
		
		For any $\ve_0 \in (0,1)$, there exists  a large $K$ such that for all $|y|=K$,
		\[
		\frac{b(|y|^{2}+\tilde T_c y_{n})}{1+b(|y'|^{2}+(y_{n}+\tilde T_c)^{2})}\ge1-\ve_{0}.
		\]
		
		For brevity, we denote 
		\[
		\tilde{u}_{i}(y)=\frac{1}{u_{i}\left(x_{i}\right)}u_{i}(u_{i}\left(x_{i}\right)^{-\frac{2}{n-2}}y)
		\]
		and 
		\[
		\varphi(y)=\left(\frac{2b\sqrt{\binom{n}{2}}}{f(x_0)}\right)^{\frac{n-2}{4}}(1+b|y+\tilde T_c\mathbf{e}_{n}|^{2})^{\frac{2-n}{2}}.
		\]
		Then for all $|y|=K$ with a larger $K$, with a positive constant
		$C_{1}$ we have 
		\begin{align*}
			& |\partial_{r}\log\tilde{u}_{i}(y)-\partial_{r}\log\varphi(y)|\\
			= & |\frac{\partial_{r}\tilde{u}_{i}}{\tilde{u}_{i}}-\frac{\partial_{r}\varphi}{\varphi}|=|\frac{\varphi\partial_{r}\tilde{u}_{i}-\tilde{u}_{i}\partial_{r}\varphi}{\tilde{u}_{i}\varphi}|\le\frac{|\varphi-\tilde{u}_{i}||\partial_{r}\varphi|+\varphi|\partial_{r}\tilde{u}_{i}-\partial_{r}\varphi|}{\tilde{u}_{i}\varphi}\\
			\le & C_{1}\varepsilon(\frac{1}{K}K^{n-2}+K^{n-2})\le \frac{\delta}{2K},
		\end{align*}
		where we take $C_{1}\varepsilon(\frac{1}{K}K^{n-2}+K^{n-2})\leq\frac{\delta}{2K}$
		. Notice that 
		\begin{align*}
			& \partial_{r}\big|_{y=\tilde{y}_{i}}\log\varphi(y)\\
			= & \frac{n-2}{2}(|y|^{-1}y_{i}\partial_{y_{i}})\big|_{y=\tilde{y}_{i}}\log[1+b_{i}(|y'|^{2}+(y_{n}+\tilde T_c)^{2})]\\
			= & -\frac{n-2}{2}\frac{1}{K}\frac{2b_{i}(|y|^{2}+\tilde T_c y_{n})}{1+b_{i}(|y'|^{2}+(y_{n}+\tilde T_c)^{2})}\\
			\le & -\frac{n-2}{K}(1-\ve_{0}),
		\end{align*}
		and thereby for all $|y|=K$,
		\[
		\partial_{r}[\log u_{i}^{-1}(x_{i})u_{i}(u_{i}(x_{i})^{-\frac{2}{n-2}}y)]<-\frac{n-2}{K}(1-\ve_{0})+\frac{\delta}{2K},
		\]
		that is,
		\[
		u_{i}\left(x_{i}\right)^{-\frac{2}{n-2}}\partial_{r}\log u_{i}\left(\hat{x}_{i}\right)\le-\frac{(n-2)(1-\ve_{0})-\frac{\delta}{2}}{K}.
		\]
		This in turn implies 
		\[
		\partial_{r}\log u_{i}\left(\hat{x}_{i}\right)\leq-\frac{(n-2)(1-\ve_{0})-\frac{\delta}{2}}{r_{i}}.
		\]
		If we choose $2(n-2)\ve_{0}<\delta$, then for all sufficiently large $i$, we again obtain a contradiction with (\ref{eq:lower bound-1-limit case}) by choosing $K$ large enough.
		
		\medskip

		\emph{Step 2:} By local $C^{1}$ and $C^{2}$ estimates, for any
		$0<\gamma<1$, we know 
		\[
		v_{i}(x)=\frac{u_{i}(x)}{u_{i}(y_{0})}\rightarrow v_{\infty}\quad\mathrm{~~in~~}C_{\mathrm{loc}}^{\gamma}(\overline{B_{r_{0}}^{+}}\backslash\{0\})
		\]
		and $v_{\infty}\in C_{\mathrm{loc}}^{0,1}(\overline{B_{r_{0}}^{+}}\backslash\{0\})\cap C_{\mathrm{loc}}^{1,1}(B_{r_0}^+)$ is the viscosity solution
		to  
		\begin{align*}
			\begin{cases}
				\sigma_{2}(A_{v_{\infty}^{\frac{4}{n-2}}g_{\mathrm{E}}})=0 & \quad\mathrm{~~in~~}B_{r_{0}}^{+},\\
				\frac{\partial v_{\infty}}{\partial x_n}=0 & \quad\mathrm{~~on~~}D_{r_0}\quad\mathrm{in~~the~~viscosity~~sense}.
			\end{cases}
		\end{align*}
		
		\medskip
		
		\emph{Step 3:} There holds $\lim_{x\rightarrow0}v_{\infty}(x)|x|^{n-2}=a$
		for some $a\in\R_{+}$.
		
		\medskip
		
		By Lemma \ref{lem:superharmonic}, we have 
		\begin{align*}
			a:= & \lim_{r\rightarrow0}r^{n-2}\min_{\partial B_{r}^{+}(0)}v_{\infty}\quad\mathrm{~~is~~finite}.
		\end{align*}
		Then the Harnack inequality in Lemma \ref{lem:harnack inequality} yields $A:=\lim\sup_{r\rightarrow0}[r^{n-2}v_{\infty}](x_{0})$
		is finite.
		
		Now it remains to show that $A=a$. 
		
		We argue by contradiction. If $A>a$, then we can find a sequence
		$z_{i}\rightarrow0$ and some $\varepsilon>0$ such that 
		\[
		A+\varepsilon\ge v_{\infty}(z_{i})r(z_{i})^{n-2}\ge a+2\varepsilon,
		\]
		while, 
		\[
		r(z_{i})^{n-2}\min_{r(x)=r(z_{i})}v_{\infty}(z_{i})\le a+\varepsilon.
		\]
		Let $r_{i}=r(z_{i}) \to 0$ and define 
		\[
		\hat{v}_{i}(y)=r_{i}^{n-2}v_{\infty}(r_{i}y).
		\]
		Then $\hat{v}_{i}\in C_{\mathrm{loc}}^{1,1}(B^{+}_{r_{0}r_{i}^{-1}})$
		satisfying 
		\begin{align*}
			\sigma_{2}(A_{\hat{v}_{i}^{\frac{4}{n-2}}g_{\mathrm{E}}})=0 & \quad\mathrm{~~in~~}\quad B_{r_{0}r_{i}^{-1}}^{+},\\
			\frac{\partial\hat{v}_{i}}{\partial y_n}=0 & \quad\mathrm{~~on~~}\quad D_{r_0r_i^{-1}}\backslash \{0\},\\
			\min_{r(y)=1}\hat{v}_{i}\le a+\varepsilon,\qquad & \max_{r(y)=1}\hat{v}_{i}\ge a+2\varepsilon.
		\end{align*}
		Thus, up to a subsequence, $\hat{v}_{i}$ converges to $v_{*}$ in
		$C_{\mathrm{loc}}^{1,1}({\mathbb{R}_{+}^{n}})\cap C^{\gamma}_{\text{loc}}(\overline{\mathbb{R}_{+}^{n}}\backslash\{0\})$ for any $0<\gamma<1$
		satisfying the following solution in viscosity sense:
		\begin{align}
			\sigma_{2}(A_{\hat{v}_{*}^{\frac{4}{n-2}}g_{\mathrm{E}}})=0\qquad & \mathrm{~~in~~}\quad\mathbb{R}_{+}^{n},\nonumber \\
			\frac{\partial\hat{v}_{*}}{\partial y_{n}}=0\qquad & \mathrm{~~on~~}\quad\pa\Rn_{+}\backslash\{0\},\nonumber \\
			\min_{r(y)=1}\hat{v}_{*}\le a+\varepsilon,\qquad & \max_{r(y)=1}\hat{v}_{*}\ge a+2\varepsilon.\label{eq:non-radial-1}
		\end{align}
		By an even reflection, it follows from Proposition \ref{even reflection in Appendix} and Theorem \ref{thm:(Li--CPAM09}
		that $\hat{v}_{*}$ is a radial solution with respect to the origin, which
		contradicts $(\ref{eq:non-radial-1})$. Hence $A=a$.
		
		To show $a>0$, we can just repeat the proof of \emph{Step 1} and find some
		positive constant $\mathcal{N}$ and $r_{1}$ independent of $i$ such that, for any 
		$0<\delta<1/2$, there holds for some large $K=K(\delta)\in \R_+$ and all sufficiently large $i$, 
		\[
		v_{i}(x)\geq\frac{\mathcal{N}}{|x|^{n-2-\delta}}\quad\mathrm{~~in~~}\left\{ Ku_{i}\left(x_{i}\right)^{-\frac{2}{n-2}}\leq |x|\leq r_{1},x_{n}\ge0\right\} .
		\]
		Clearly this implies that $v_{\infty}(x)\geq\frac{\mathcal{N}}{|x|^{n-2-\delta}}$
		in $\left\{ 0<|x|\leq r_{1}\right\} $.
		
		For $n=3$, by Li-Nguyen \cite[Theorem 1.2]{Li-Luc1} or Theorem \ref{thm:(Bocher-Type-theorems--Li-Nguyen} and even reflection  we have
		\[
		\lim_{x\rightarrow x_{0}}v_{\infty}(x)|x|^{n-2}=a>0.
		\]

		For $n=4$,  there exists a positive constant $\mathcal{N}$ depending on $\delta$ such that 
		$$v_{\infty}(x)\geq\frac{\mathcal{N}}{|x|^{n-2-\delta}}$$
		in $\left\{ 0<|x|\leq C_3\delta\right\} $ for all $0<\delta<1/2$. On the other hand, via an even reflection, we can apply Theorem \ref{thm:(Bocher-Type-theorems--Li-Nguyen} to conclude that $v_{\infty}(x)=e^{\mathring{w}}|x|^{-\alpha}$ for some bounded function $\mathring{w}$ and $\alpha \in (0,n-2]$.  Combing these two facts, we eventually have $v_{\infty}(x)=a |x|^{2-n}$ in $B_{r_{0}}^{+}$. This phenomenon is distinct from that in Theorem \ref{Thm:isolated behavior}.
	\end{proof}

	\subsection{Isolated blow-up points are simple}
	
	Thanks to Proposition \ref{prop:special_g}, here and below we can assume that the background metric satisfies $h_g=0$ on $\pa M$.

	Suppose that $\{u_{i}\}$ is a sequence of positive smooth solutions to \eqref{PDE:sigma_2_new}. Remind that
	$$r_{i}:=R_{i}u_{i}^{-\frac{2}{n-2}}(x_{i}), \quad R_i \to \infty$$
	are two sequences of positive real numbers given in Lemma \ref{lem fudament}.

	First we shall prove that every interior isolated blow-up point
	is simple. 
	
	\begin{lemma} \label{lem:interior blow up is simple}
		Let $x_{i}\rightarrow x_{0}\in M$ be an isolated blow-up point for
		$u_{i}$, then $x_{0}$ is simple. 
		\end{lemma} 
	\begin{proof}
		
		Let 
		$$\eta_i(z)=u_{i}^{-1}\left(x_{i}\right)u_{i}\circ \exp_{x_{i}}(u_{i}^{-\frac{2}{n-2}}\left(x_{i}\right)z)$$
		and 
		$$\eta(z)=\left(1+\frac{f(x_0)}{2\sqrt{\binom{n}{2}}}|z|^{2}\right)^{\frac{2-n}{2}}.$$
		For all sufficiently large $i$, similar to the proof of Lemma \ref{one critical point} and  \eqref{eq:first derivative negative upper bound},  we find some sufficiently small $\alpha>0$ such that 
		\[
		(s^{\frac{n-2}{2}}\overline{\eta}_{i}(s))'\le-Cs^{-\frac{n}{2}}\quad\mathrm{~~for~~}s_0(1+\alpha)\le s\le3R_{i},
		\]
		where $s_0=\sqrt{2\sqrt{\binom{n}{2}}/f(x_0)}$.
		This means 
		\[
		\overline{w}_{i}'(r)<0,\quad\mathrm{~~for~~all~~} Cu_{i}(x_{i})^{-\frac{2}{n-2}}\le r\le r_{i}.
		\]
		
		Suppose $x_0$ is not an isolated simple blow-up point and let $\mu_{i}$
		be the second critical point of $r^{(n-2)/2}\overline{u}_{i}(r)$,
		we have $\mu_{i}>r_{i}$ and $\lim_{i\rightarrow\infty}\mu_{i}=0.$
		Define 
		\[
		\tilde{u}_{i}(z)=\mu_{i}^{\frac{n-2}{2}}u_{i}\circ\mathcal{T}_{i}(z),\quad|z|<\delta/\mu_{i}
		\]
		and
		\[
		\overline{\tilde{u}_{i}}(r)=\frac{1}{r^{n-1}}\int_{|z|=r}\tilde{u}_{i}(z)\ud\sigma_{g_{i}},			\]
		where $\mathcal{T}_{i}(z)=\exp_{x_{i}}(\mu_{i}z)$. 
		
		Let $g_{i}:=\mu_{i}^{-2}\mathcal{T}_{i}^{\ast}(g)\to g_{\mathrm{E}}$ in $C^3_{\mathrm{loc}}(\Rn)$,
		then it is not hard to see that
		$$g_{\tilde{u}_{i}}:=\tilde{u}_{i}^{\frac{4}{n-2}}g_{i}=\mathcal{T}_i^\ast(u_i^{\frac{4}{n-2}}g).$$
		
		Furthermore, the following are true:
		\begin{itemize}
			\item $\sigma_{2}^{1/2}(g_i^{-1}A_{\tilde{u}_{i}})=f\circ\mathcal{T}_{i}\tilde{u}_{i}^{\frac{4}{n-2}}+\mu_{i}^{2}f_{0}\quad$
			in $~~|z|<\delta/\mu_{i}$;
			\item The isolated blow-up point implies 
			\begin{equation}
				|z|^{\frac{n-2}{2}}\tilde{u}_{i}(z)\le C\quad\mathrm{~~in~~}|z|<\delta/\mu_{i};\label{Isolated_2nd_critical pt}
			\end{equation}
			\item $(r^{\frac{n-2}{2}}\overline{\tilde{u}_{i}}(r))'|_{r=1}=0$;
			\item $r^{\frac{n-2}{2}}\overline{\tilde{u}_{i}}(r)$ has precisely one
			critical point in $0<r<1$;
			\item $(r^{\frac{n-2}{2}}\overline{\tilde{u}_{i}})'(r)<0$ for $4\tilde{u}_{i}(0)^{-\frac{2}{n-2}}<r<1$;
			\item $\lim_{i\rightarrow\infty}\tilde{u}_{i}(0)=\infty$ and $0$ is a
			local maximum point of $\tilde{u}_{i} $.
			
			This is due to 
			\[
			\tilde{u}_{i}(0)=\mu_{i}^{\frac{n-2}{2}}u_{i}(x_{i})\geq r_{i}^{\frac{n-2}{2}}u_{i}(x_{i})=R_{i}^{\frac{n-2}{2}}\to\infty\quad\mathrm{as~~}i\to\infty.
			\]
			\item $\tilde{u}_{i}(z)\rightarrow0$ in $C_{\mathrm{loc}}^{0}(\mathbb{R}^{n}\backslash\{0\})$.
			
			This is a direct consequence of Theorem \ref{isolated behavior-1}.
			
		\end{itemize}

		For any compact set $K\subset\mathbb{R}^{n}\backslash\{0\}$, the
		upper bound of $\tilde{u}_{i}$ in $K$ can also follow from \eqref{Isolated_2nd_critical pt}.
		Then the interior local estimate yields 
		\begin{equation}\label{priori estimates1}
			\|\nabla\log\tilde{u}_{i}\|_{C^{0}(K)}+\|\nabla^{2}\log\tilde{u}_{i}\|_{C^{0}(K)}\le C(K).
		\end{equation}
		Actually we have $\tilde{u}_{i}(z)\rightarrow0$ in $ C_{\mathrm{loc}}^{1,\gamma}({\mathbb{R}^{n}}\backslash\{0\})$
		for any $0<\gamma<1$.
		
		If we let $\omega_{i}(z)=\tilde{u}_{i}(z)/\tilde{u}_{i}(y_{0})$ for any fixed $y_0 \in \Rn \backslash \{0\}$,
		then  $\omega_{i}$ satisfies 
		\begin{align*}
			& \sigma_{2}^{1/2}(-\frac{2}{n-2}\frac{\nabla_{g_{i}}^{2}\omega_{i}}{\omega_{i}}+\frac{2n}{(n-2)^{2}}\frac{\ud\omega_{i}\otimes_{g_{i}}\ud\omega_{i}}{\omega_{i}^{2}}-\frac{2}{(n-2)^{2}}\frac{|\nabla\omega_{i}|_{g_{i}}^{2}}{\omega_{i}^{2}}g_{i}+A_{g_{i}})\\
			= & \omega_{i}^{\frac{4}{n-2}}(z)\tilde{u}_{i}^{\frac{4}{n-2}}(y_{0})f\circ\mathcal{T}_{i}(z)+\mu_{i}^{2}f_{0}\qquad\mathrm{~~in~~}B_{\delta/\mu_i}.
		\end{align*}
		
		By \eqref{priori estimates1} and $\omega_{i}(y_{0})=1$, with some $C_{K}\in\R_{+}$ we obtain  $\omega_{i}\le C_{K}$
		in $K$. Also by the gradient estimate,
		$\omega_{i}\geq c_{K}$ in $K$ for some $c_{K}\in\R_{+}$. Hence,
		the interior estimates show that $\|\log\omega_{i}\|_{C^{2}(K)}\le C.$
		Now there exists a subsequence $\omega_{i}$ converging to $\omega_{\infty}\in C^{1,1}(\mathbb{R}^{n}\backslash\{0\})$ satisfying
		\begin{align*}
			\sigma_{2}^{1/2}(-\frac{2}{n-2}\frac{\nabla_{g_{\mathrm{E}}}^{2}\omega_{\infty}}{\omega_{\infty}}+\frac{2n}{(n-2)^{2}}\frac{\ud\omega_{\infty}\otimes\ud\omega_{\infty}}{\omega_{\infty}^{2}}-\frac{2}{(n-2)^{2}}\frac{|\nabla\omega_{\infty}|_{g_{\mathrm{E}}}^{2}}{\omega_{\infty}^{2}}g_{\mathrm{E}}) & =0\quad\mathrm{~~in~~}\mathbb{R}^{n}\backslash\{0\},\\
			(r^{\frac{n-2}{2}}\overline{\omega}_{\infty})'(r)\leq0\quad\mathrm{~~for~~}0<r<1,\\
			(r^{\frac{n-2}{2}}\overline{\omega}_{\infty}(r))'|_{r=1}=0,
		\end{align*}
		where 
		\[	\overline{\omega}_{\infty}(r)=\frac{1}{r^{n-1}}\int_{\pa B_r}{\omega}_{\infty}(z)\ud\sigma.
		\]
		By the non-increasing and positivity of $r^{(n-2)/2}\overline{\omega}_{\infty}(r)$,
		we know that $\overline{\omega}_{\infty}$ must be singular at $0$.
		Hence, by Theorem \ref{thm:(Li--CPAM09} and Theorem
		\ref{thm:(Li-Nguyen2015-CPAM-Prop2.2)radial solution} (or \cite[Theorem 1.4]{Li-Luc1})
		we conclude that for $n=3$, $\omega_{\infty}(z)=b|z|^{2-n}$ for some $b\in\R_{+}$;
		for $n=4$, $\omega_{\infty}(z)=C_{1}|z|^{-\alpha}$ for some $C_{1}\in \R_+$
		and $0<\alpha \le n-2$. However, it clearly contradicts the property
		that $(r^{\frac{n-2}{2}}\overline{\omega}_{\infty}(r))'|_{r=1}=0$. 
	\end{proof}
	
	\begin{lemma} \label{lem:boundary bounded limit} Let  $x_{i}\rightarrow x_{0} \in \pa M$
		be an isolated blow-up point for $u_{i} $, then 
		\[
		d_{g}(x_{i},\partial M)u_{i}(x_{i})^{\frac{2}{n-2}}\mathrm{~~is~~bounded.}
		\]
	\end{lemma} 
	\begin{proof}
		For brevity, we let $d_i:=d_{g}(x_{i},\pa M)=d_{g}(x_{i},x_{i}')$
		for some $x_{i}'\in\pa M$.
		
		By contradiction, let $T_{i}:=d_i^{(n-2)/2}u_{i}(x_{i})\rightarrow\infty$. We define 
		\[
		v_{i}(z)=d_i^{\frac{n-2}{2}}u_{i}\circ\Psi_{i}(z) \qquad \mathrm{with~~}\quad\Psi_{i}(z)=G_{x_i'}(d_i(z+\mathbf{e}_n)).
		\]
		Let
		$$g_{i}=d_i^{-2}\Psi_{i}^{\ast}(g) \to g_{\mathrm{E}} \qquad \mathrm{in~~}C_{\mathrm{loc}}^3(
		\Rn_{-1})$$
		and
		\[
		g_{v_{i}}=v_{i}^{\frac{4}{n-2}}g_{i}=\Psi_{i}^{\ast}(u_{i}^{\frac{4}{n-2}}g)=\Psi_{i}^{\ast}(g_{u_{i}}).
		\]
		Denote $$B_{\frac{\rho}{d_i}}^{+}(-\mathbf{e}_n):=B_{\frac{\rho}{d_i}}^{}(-\mathbf{e}_n)\cap \{z_n> -1\}.$$
		Then $v_i$ satisfies
		\[
		\sigma_{2}^{1/2}(A_{v_{i}})=f\circ\Psi_{i}v_{i}^{\frac{4}{n-2}}+d_i^{2}f_{0}\qquad\mathrm{~~in~~}\quad B_{\frac{\rho}{d_i}}^{+}(-\mathbf{e}_n)\backslash\{0\}
		\]
		and 
		\[
		\frac{2}{n-2}\frac{\partial v_{i}}{\partial \vec n_{g_{i}}}=-c\circ\Psi_{i}v_{i}^{\frac{n}{n-2}}\quad\mathrm{~~on~~}\quad\{z_{n}=-1\}.
		\]
		Since $x_{0}$ is an isolated blow-up point, for any large $i$, we have 
		\[
		|z|^{\frac{n-2}{2}}v_{i}(z)\le C
		\]
		for any compact $K\subset\overline{\mathbb{R}_{-1}^{n}}$.
		
		Notice that $0$ is a local maximum point of $v_{i}$
		and $v_{i}(0)=T_{i}\rightarrow\infty.$ For any compact set $K\subset\overline{\mathbb{R}_{-1}^{n}}\backslash\{0\}$
		we have 
		\[
		\|\nabla\log v_{i}\|_{C^{0}(K)}\le C(K)
		\]
		and then 
		\[
		0\le\sup_{K_{1}}v_{i}\le c\inf_{K_{1}}v_{i},
		\]
		where $c$ is a positive constant independent of $i$ and $K_{1}$ is a compact set in
		$\overline{\mathbb{R}_{-1}^{n}}\backslash\{0\}.$
		
		By Lemma \ref{lem:interior blow up is simple}, $0$
		is a simple blow-up point of $v_{i}$ and thus there exists some $\tilde{r}\in\R$
		such that $\big(s^{(n-2)/2}\bar{v}_{i}(s)\big)'<0$ in $(\tilde{r_{i}},\tilde{r}),$
		where $s=|z|,\tilde{r}_{i}:=R_{i}v_{i}(0)^{-2/(n-2)}\rightarrow0$
		and $R_{i} \to \infty$.

		For any compact set $K\subset{\mathbb{R}_{-1}^{n}}\backslash\{0\}$, for large $i$,
		by the classical interior estimates (see \cite{Guan-Wang0,Wang,Chen2,LiYY}),
		we have 
		\[
		\|\nabla\log v_{i}\|_{C^{0}(K)}+\|\nabla^{2}\log v_{i}\|_{C^{0}(K)}\le C(K).
		\]
		By Theorem \ref{isolated behavior-1}, $v_{i}\rightarrow0$ in $C_{\mathrm{loc}}^{0}(\overline{\mathbb{R}_{-1}^{n}}\backslash\{0\})$. Actually we have $v_{i}(z)\rightarrow0$ in $C_{\mathrm{loc}}^{1,1}({\mathbb{R}_{-1}^{n}}\backslash\{0\})\cap C_{\mathrm{loc}}^{0,1}(\overline{\mathbb{R}_{-1}^{n}}\backslash\{0\})$.
		
		 We define $\omega_{i}(z)=v_{i}(z)/v_{i}(y_{0})$  for any fixed $y_0\neq 0$ with $\omega_{i}(y_{0})=1$.
		Then $\omega_{i}$ satisfies 
		$$
		\begin{cases}
			\sigma_{2}^{1/2}(-\frac{2}{n-2}\frac{\nabla_{g_{i}}^{2}\omega_{i}}{\omega_{i}}+\frac{2n}{(n-2)^{2}}\frac{\ud\omega_{i}\otimes_{g_{i}}\ud\omega_{i}}{\omega_{i}^{2}}-\frac{2}{(n-2)^{2}}\frac{|\nabla\omega_{i}|_{g_{i}}^{2}}{\omega_{i}^{2}}g_{i}+A_{g_{i}})&\\
			= \omega_{i}^{\frac{4}{n-2}}(z)v_{i}^{\frac{4}{n-2}}(y_{0})f\circ\Psi_{i}+d_i^{2}f_{0} &\qquad\mathrm{~~in~~}\quad  B_{\frac{\rho}{d_i}}^{+}(-\mathbf{e}_n)\backslash\{0\}\\
		\frac{2}{n-2}\frac{\partial\omega_{i}}{\partial \vec n_{g_i}}=-c(x_{0})\omega_{i}^{\frac{n}{n-2}}v_{i}(y_{0})^{\frac{2}{n-2}}
		&\qquad \mathrm{~~on~~} \quad z_{n}=-1.
		\end{cases}
		$$
		
		By $C^1$ estimates for $v_i$, we actually know $\omega_{i}(z)\le C_{K}$
		in every compact set $K\subset\overline{\mathbb{R}_{-1}^{n}}\backslash\{0\}$.
		Also by Harnack inequality, we have $\omega_{i}\geq c_{K}$ in $K$
		and thus $\|\log \omega_{i}\|_{C^{1}(K)}\le C$. Then there exists a subsequence
		$w_{i}$ converging to $\omega_{\infty}\in C_{\mathrm{loc}}^{1,1}(\mathbb{R}_{-1}^{n}\backslash\{0\})\cap C_{\mathrm{loc}}^{0,1}(\overline{\mathbb{R}_{-1}^{n}}\backslash\{0\})$
		satisfying the following equation in the viscosity
		sense 
		\[
		\begin{cases}
		\sigma_{2}(-\frac{2}{n-2}\frac{\nabla_{g_{\mathrm{E}}}^{2}\omega_{\infty}}{\omega_{\infty}}+\frac{2n}{(n-2)^{2}}\frac{\ud\omega_{\infty}\otimes_{g_{\mathrm{E}}}\ud\omega_{\infty}}{\omega_{\infty}^{2}}-\frac{2}{(n-2)^{2}}\frac{|\nabla\omega_{\infty}|_{g_{\mathrm{E}}}^{2}}{\omega_{\infty}^{2}}g_{\mathrm{E}})=0&\qquad\mathrm{in~~}\mathbb{R}_{-1}^{n}\backslash\{0\}\\
		\frac{\partial w_{\infty}}{\partial z_{n}}=0 &\qquad \mathrm{on~~} z_{n}=-1.
		\end{cases}
		\] 
		
		For small $r>0$, we define
		\[	\overline{\omega}_{\infty}(r)=\frac{1}{r^{n-1}}\int_{|z|=r}{\omega}_{\infty}(z)\ud\sigma.
		\]
		Moreover, $\big(s^{(n-2)/2}\overline \omega_{\infty}(s)\big)'\leq0$
		near zero.
		By the same reason as before, we know that $\overline \omega_{\infty}$ is singular
		at $0$. By Theorem \ref{isolated behavior-1} we have $\lim_{z\rightarrow 0}\omega_{\infty}(z)|z|^{n-2}=a$ for some $a \in \R_+$.
		
		On the other hand, it follows from Theorem \ref{thm:(Bocher-Type-theorems--Li-Nguyen}
		that when $n=3,4$, $\omega_{\infty}(z)=a|z|^{2-n}$. This yields 
		\[
		\frac{\partial\omega_{\infty}}{\partial z_{n}}\bigg|_{z_{n}=-1}=[(2-n)a|z|^{-n}z_{n}]\big|_{z_{n}=-1}=(n-2)a|z|^{-n}>0,
		\]
		which contradicts the boundary condition $\frac{\partial\omega_{\infty}}{\partial z_{n}}=0$
		on $\{z_{n}=-1\}$. 
	\end{proof}
	
	\begin{remark} As a direct consequence of Lemma \ref{lem:boundary bounded limit}, if
		$x_i \to x_{0} \in \overline M$ is an isolated blow-up point and $d_{g}(x_{i},\partial M)u_{i}(x_{i})^{2/(n-2)}\to \infty$, then $x_0$ has to be an interior blow-up point.
	\end{remark} 
	
	Now let us see that every boundary isolated blow-up point for $u_{i}$
	is simple.
	\begin{lemma} \label{lem:boundary blow up is simple-1}
		Let $x_{i}\rightarrow x_{0}\in\partial M$ be an isolated blow-up
		point for $u_{i}$, then $x_{0}$ is simple.
		 \end{lemma} 
		 
	\begin{proof}
	Using (\ref{eq:first derivative negative upper bound}) in the proof of Lemma \ref{one critical point},   we have 
		\[
		\overline{w}_{i}'(r)<0 \qquad\mathrm{~~for~~all~~} Cu_{i}(x_{i})^{-\frac{2}{n-2}}\le r\le r_{i}.
		\]
		
		Suppose it is not an isolated simple blow-up point and let $\mu_{i}$
		be the second critical point of $r^{(n-2)/2}\overline{u}_{i}(r)$,
		we know that $\mu_{i}>r_{i}$ and $\lim_{i\rightarrow\infty}\mu_{i}=0$. It is not hard to see that $d_i/\mu_i \to 0$ by Lemma \ref{lem:boundary bounded limit} and definition of $r_i$.
		
				Define 
		\[
		\Psi_{i}(z)=G_{x_i'}(\mu_{i}z+d_i \mathbf{e}_n)\quad\mathrm{~~and~~}\quad g_{i}=\mu_{i}^{-2}\Psi_{i}^{\ast}(g)\rightarrow g_E \quad \mathrm{~~in~~}  C^3_{\mathrm{loc}}(\Rn).
		\]
		We define 
		\[
		\tilde{u}_{i}(z)=\mu_{i}^{\frac{n-2}{2}}u_{i}\circ\Psi_{i}(z),\qquad z\in B_{\rho/\mu_{i}}^{+},
		\]
		then 
		\[
		g_{\tilde{u}_{i}}=\tilde{u}_{i}^{\frac{4}{n-2}}g_{i}=\Psi_{i}^{\ast}(u_{i}^{\frac{4}{n-2}}g)=\Psi_{i}^{\ast}(g_{u_{i}}).
		\]
		
		If we let 
		\[
		\overline{\tilde{u}_{i}}(r)=\frac{1}{r^{n-1}}\int_{\pa^{+}B^+_{r,\theta}}\tilde{u}_{i}(z)\ud\sigma_{g_{i}},
		\]
		then the following hold 
		\[
		\begin{cases}
			\sigma_{2}^{1/2}(g_i^{-1}A_{\tilde{u}_{i}})=f\circ\Psi_{i}\tilde{u}_{i}^{\frac{4}{n-2}}+\mu_{i}^{2}f_{0} &\mathrm{~~in~~}B_{\rho/\mu_{i}}\cap \{z_n\ge -\frac{d_i}{\mu_i}\},\\
			|z|^{\frac{n-2}{2}}\tilde{u}_{i}(z)\le C, &\mathrm{~~in~~}B_{\rho/\mu_{i}}\cap \{z_n\ge -\frac{d_i}{\mu_i}\},\\
			(r^{\frac{n-2}{2}}\overline{\tilde{u}_{i}}(r))'|_{r=1}=0,\\
			\frac{2}{n-2}\frac{\partial\tilde{u}_{i}}{\partial \vec n_{g_{i}}}=-c\circ\Psi_{i}\tilde{u}_{i}^{\frac{n}{n-2}} &\mathrm{~~on~~}B_{\rho/\mu_{i}}\cap \{z_n= -\frac{d_i}{\mu_i}\},\\
			r^{\frac{n-2}{2}}\overline{\tilde{u}_{i}}(r)\quad\mathrm{~~has~~precisely~~one~~critical~~point~~in~~}0<r<1,\\
			(r^{\frac{n-2}{2}}\overline{\tilde{u}_{i}}(r))'<0\quad\mathrm{~~in~~}C\tilde{u}_{i}(0)^{-\frac{2}{n-2}}<r<1,\\
			\lim_{i\rightarrow\infty}\tilde{u}_{i}(0)=\infty,~~0\mathrm{~~is~~a~~local~~maximum~~point~~of~~}\tilde{u}_i .
		\end{cases}
		\]
		Consequently, $0$ is an isolated blow up point for $\tilde{u}_i$. Meanwhile, it follows from
		Theorem \ref{isolated behavior-1} that $\tilde{u}_{i}(z)\rightarrow0$
		in $C_{\mathrm{loc}}^{0}(\overline{\mathbb{R}_{+}^{n}}\backslash\{0\})$.
		
		For any compact set $K\subset\overline{\mathbb{R}_{+}^{n}}\backslash\{0\}$,
		the  local $C^1$ estimates yield 
		\begin{align}\label{gradient estimate 5.22}
			\|\nabla\log\tilde{u}_{i}\|_{C^{0}(K)}\le C(K).
		\end{align}
		Actually we have $\tilde{u}_{i}(z)\rightarrow0$ in $C_{\mathrm{loc}}^{\gamma}(\overline{\mathbb{R}_{+}^{n}}\backslash\{0\})\cap C_{\mathrm{loc}}^{1,\gamma}({\mathbb{R}_{+}^{n}})$
		for any $0<\gamma<1$.
		
		Notice that for any fixed $y_0 \neq 0$, $\omega_{i}(z):=\tilde{u}_{i}(z)/\tilde{u}_{i}(y_{0})$ satisfies 
		\begin{align*}
			& \sigma_{2}^{1/2}(-\frac{2}{n-2}\frac{\nabla_{g_{i}}^{2}\omega_{i}}{\omega_{i}}+\frac{2n}{(n-2)^{2}}\frac{\ud\omega_{i}\otimes\ud\omega_{i}}{\omega_{i}^{2}}-\frac{2}{(n-2)^{2}}\frac{|\nabla\omega_{i}|_{g_{i}}^{2}}{\omega_{i}^{2}}g_{i}+A_{g_{i}})\\
			= & \omega_{i}^{\frac{4}{n-2}}\tilde{u}_{i}^{\frac{4}{n-2}}(y_{0})f\circ\Psi_{i}+\mu_{i}^{2}f_{0}\qquad\mathrm{~~in~~}\quad B_{\rho/\mu_{i}}^{+}
		\end{align*}
		and 
		$$\frac{2}{n-2}\frac{\partial\omega_{i}}{\partial z_{n}}=-c\circ \Psi_i\omega_{i}^{\frac{n}{n-2}}\tilde{u}_{i}(y_{0})^{\frac{2}{n-2}} \qquad \mathrm{~~on~~} D_{\rho/\mu_{i}}.$$
		
		By \eqref{gradient estimate 5.22} and $\omega_{i}(y_{0})=1$, we
		actually know $\omega_{i}(z)\le C_{K}$ in every compact set $K\subset\overline{\Rn_{+}}\backslash\{0\}$ and   $\omega_{i}\geq c_{K}$
		in $K$ and thus $\|\log\omega_{i}\|_{C^{1}(K)}\le C$. Then there
		exists a subsequence $\omega_{i}$ converging to $\omega_{\infty}\in C_{\mathrm{loc}}^{0,1}(\overline{\mathbb{R}_{+}^{n}}\backslash\{0\})\cap C_{\mathrm{loc}}^{1,1}({\mathbb{R}_{+}^{n}})$ satisfying the following equation
		\begin{equation*}
		\sigma_{2}\left(-\frac{2}{n-2}\frac{\nabla_{g_{\mathrm{E}}}^{2}\omega_{\infty}}{\omega_{\infty}}+\frac{2n}{(n-2)^{2}}\frac{\ud\omega_{\infty}\otimes\ud\omega_{\infty}}{\omega_{\infty}^{2}}-\frac{2}{(n-2)^{2}}\frac{|\nabla\omega_{\infty}|_{g_{\mathrm{E}}}^{2}}{\omega_{\infty}^{2}}g_{\mathrm{E}}\right)=0\quad\mathrm{in~~}\overline{\mathbb{R}_{+}^{n}}\backslash\{0\}
		\end{equation*}
		in viscosity
		sense, together with the following properties: 
		\begin{subequations}
\begin{align}
			(r^{\frac{n-2}{2}}\overline{\omega}_{\infty}(r))'\leq0 & \quad\mathrm{~~for~~}0<r<1,\label{average_limit_fcn1}\\
			(r^{\frac{n-2}{2}}\overline{\omega}_{\infty}(r))'|_{r=1}=0,\label{average_limit_fcn2}\\
			\frac{\partial \omega_{\infty}}{\partial z_{n}}=0, & \quad\mathrm{~~on~~}\quad\pa\Rn_{+}\backslash\{0\}.\no
		\end{align}
		\end{subequations}
		Here
		$$\overline{\omega}_\infty(r)=\frac{1}{r^{n-1}}\int_{\pa^{+}B^+_{r,\theta}}\overline{\omega}_\infty(z)\ud\sigma.$$
		
		By Proposition \ref{even reflection in Appendix} via an even reflection, we know that $\omega_{\infty}$  is a radial solution with singularity at the origin.   Again thanks to Theorem \ref{isolated behavior-1}, we
		have $\lim_{z\rightarrow0}\omega_{\infty}(z)|z|^{n-2}=a$ for some $a\in \R_+$. Hence, via an even reflection, we can apply Theorem
		\ref{thm:(Bocher-Type-theorems--Li-Nguyen} to conclude that for $n=3,4$, $\omega_{\infty}(z)=a|z|^{2-n}$, which contradicts the fact $(r^{(n-2)/2}\overline{\omega}_{\infty}(r))'|_{r=1}=0.$ 
	\end{proof}

	Though Lemma \ref{lem fudament-1-1} demonstrates each large solution $u$ can be well approximated under strong norms by standard bubbles in the distinct balls $B_{r_i}(p_i), 1\leq i \leq N(u)$, the next Lemma \ref{lem:pairwize_disjoint_sigular_pts} excludes the accumulation of these standard bubbles, which implies that only isolated blow-up points occur for any blow-up sequence of positive solutions to \eqref{PDE:sigma_2_new}.

	\begin{lemma}\label{lem:pairwize_disjoint_sigular_pts}
		Let $(M,g)$ be a smooth compact Riemannian manifold of dimension
		$n=3,4$ with boundary $\pa M$. Assume $h_g=0$ and $0\le c \in C^\infty(\pa M)$. Given any constants $R\geq1$ and $0<\ve<1$,
		there exist positive constants $\delta_{0}=\delta_{0}(M,g,c,R,\ve)$
		and $C_{0}=C_{0}(M,g,c,R,\ve)$ such that for all positive smooth solutions $u$ to
		\eqref{PDE:sigma_2_new} with $\max_{\overline{M}}u\geq C_{0}$, we
		have 
		\[
		\min\{d_{g}(q_{i},q_{j});1\le i\neq j\le N\}\ge\delta_{0},
		\]
		where $q_{1},\cdots,q_{N}$ are the points in $\overline{M}$ given
		in Lemma \ref{lem fudament-1-1}. 
	\end{lemma} 
	\begin{proof}
		If not, there exists a sequence of solutions $u_{i}$ to \eqref{PDE:sigma_2_new}
		such that $\min\{d_{g}(q_{i,a},q_{i,b});1\le a,b\le N_{i},a\neq b\}\rightarrow0$
		as $i\rightarrow\infty$, where each $q_{i,a}$ is a local maximum
		point of $u_{i}$. We may assume $\rho_{i}=\min_{a\neq b}d_{g}(q_{i,a},q_{i,b})=d_{g}(q_{i,1},q_{i,2})$
		and $q_{0}=\lim_{i\rightarrow\infty}q_{i,1}=\lim_{i\rightarrow\infty}q_{i,2}$.
		This implies $\lim_{i\to\infty}\rho_{i}=0$. Clearly,  we have either $q_{0}\in\partial M$ or $q_{0}\in M$.
		
		We focus on $q_{0}\in\partial M$, since the other case
		is relatively easier using the same idea.

		By Lemma \ref{lem fudament-1-1} (i) we obtain that for all $a\neq b$,
		\[
		Ru_{i}^{-\frac{2}{n-2}}(q_{i,a})\le d_{g}(q_{i,a},q_{i,b}),
		\]
		which implies that 
		\begin{equation}\label{blowtoinfty}
			\lim_{i\rightarrow\infty}u_{i}(q_{i,a})=+\infty.
		\end{equation}
		for $a=1,2$ and 
		\begin{equation}
			d_{g}(q_{i,1},q_{i,2})^{\frac{n-2}{2}}u_{i}(q_{i,a})\ge R\quad\mathrm{~~for~~}a=1,2.\label{eq:lower bound}
		\end{equation}

		We choose local coordinates near $q_{i,1}$ as follows: Let $d_i:=d_{g}(q_{i,1},\pa M)=d_{g}(q_{i,1},q_{i,1}')$
		for some $q_{i,1}'\in\pa M$ and denote by 
		\[
		\Psi_{i}(y):=G_{q_{i,1}'}(\rho_{i}y+d_{i}\mathbf{e}_{n}):B_{\rho/\rho_i}^{+}\to\overline{M}.
		\]
		We define 
		\[
		w_{i}(y)=\rho_{i}^{\frac{n-2}{2}}u_{i}\circ\Psi_{i}(y)
		\]
		and 
		\[
		q_{i,a}=\Psi_{i}(A_{i,a}).
		\]
		Under the above local coordinates we have 
		\[
		A_{i,1}=0,\quad|A_{i,2}|=1+O(\rho_i),\quad\lim_{i\rightarrow\infty}A_{i,2}=A \quad\mathrm{with}\quad|A|=1
		\]
		and 
		\[
		\min_{a\neq b}|A_{i,a}-A_{i,b}|\ge1+O(\rho_i).
		\]
		Notice that  for large 
		fixed $i$, all $A_{i,a} (1\leq a\leq N_{i})$ are distinct.
		
		We can deduce from (\ref{eq:lower bound}) that 
		\begin{equation}
			w_{i}(0),w_{i}(A_{i,2})\ge R \label{est:lower_bdd}
		\end{equation}
		and each $A_{i,a}$ is a local maximum point of $w_{i}$, for some $\rho>0$, 
		\begin{equation}
			\min_{a=1,2,\cdots,N_i}|y-A_{i,a}|^{\frac{n-2}{2}}w_{i}(y)\le C_{1},\quad |y|\le\frac{\rho}{2\rho_{i}},~~y_{n}\ge-\frac{d_{i}}{\rho_{i}}:=-T_{i}.\label{eq:decay est}
		\end{equation}
		If we let 
		\[
		g_{w_{i}}=w_{i}^{\frac{4}{n-2}}g_{i}\quad\mathrm{~~and~~}\quad g_{i}:=\rho_{i}^{-2}\Psi_{i}^{*}(g),
		\]
		then $w_{i}$ satisfies 
		\begin{align*}
			\begin{cases}
				\sigma_{2}^{1/2}(A_{g_{w_{i}}})=f\circ\Psi_{i}w_{i}^{\frac{4}{n-2}}+\rho_{i}^{2}f_{0} & \quad\mathrm{~~in~~}|y|<\frac{\rho}{2\rho_{i}},~~y_{n}>-T_{i},\\
				\frac{2}{n-2}\frac{\partial w_{i}}{\partial\vec{n}_{g_{i}}}(y)=-c\circ\Psi_{i}w_{i}^{\frac{n}{n-2}} & \quad\mathrm{~~on~~}y_{n}=-T_{i}.
			\end{cases}
		\end{align*}
		
		Up to a subsequence, we let $T=\lim_{i\rightarrow\infty}T_{i}\in[0,\infty]$. Though $T$ might be $T=\infty$ or $T\in\R_{+}$ or $T=0$, we manage to deal with them in a unified way.
		
		\begin{claim}
			As $i\to\infty$, there holds $w_{i}(0),w_{i}(A_{i,2})\rightarrow\infty$.
		\end{claim}
		
		If only one of them goes to infinity along a subsequence, say, $w_{i}(0)\rightarrow\infty$
		while $w_{i}(A_{i,2})$ is bounded, then $0$ is an isolated simple
		blow-up point by virtue of 
		Lemma \ref{lem:boundary blow up is simple-1}. Thus, there exists
		some positive constant $\tilde{r}$ such that $\big(s^{(n-2)/2}\overline{w}_{i}(s)\big)'<0$
		in $(\tilde{r}_{i},\tilde{r})$, where $$\overline{w}_{i}(r)=\frac{1}{r^{n-1}}\int_{\partial^{+}B_{r,\theta}^{+}}w_{i}\ud\sigma_{g_{i}}$$
		and $\tilde{r}_{i}:=Cw_{i}(0)^{-2/(n-2)}\rightarrow0$. By Theorem \ref{isolated behavior-1}, we know that $w_i(A_{i,2}) \rightarrow 0$, which is a contradiction with \eqref{est:lower_bdd}.

		If both $w_{i}(0)$
		and $w_{i}(A_{i,2})$ are bounded, then we first show that $\{w_{i}\}$
		is locally bounded.
		
		By negation, there exists a bounded sequence $\{A_{i,a_{i}}\}$ for
		some $1\leq a_{i}\leq N_{i}$ such that $w_{i}(A_{i,a_{i}})\rightarrow\infty.$
		By (\ref{eq:decay est}) we know that up to a subsequence, $\{A_{i,a_{i}}\}$
		is an isolated blow-up point for $w_{i}$ and thus also isolated simple.
		By Theorem \ref{isolated behavior-1}, we know that $w_i(0), w_{i}(A_{i,2}) \rightarrow 0$, which contradicts \eqref{est:lower_bdd}.
		Now $\{w_i\}$ is locally bounded
		and then, there exists a subsequence $w_{i}$ converging
		to $\tilde{w}_{\infty}^{1}$ satisfying 
		\[
		\begin{cases}
			\sigma_{2}^{1/2}(A_{g_{\tilde{w}_{\infty}^{1}}})=f(q_{0}) & \quad\mathrm{~~in~~}\quad\mathbb{R}_{-T}^{n},\\
			\frac{\partial\tilde{w}_{\infty}^{1}}{\partial y_{n}}=-\frac{n-2}{2}c(q_{0})(\tilde{w}_{\infty}^{1})^{\frac{n}{n-2}} & \quad\mathrm{~~on~~}\quad\{y_{n}=-T\}.
		\end{cases}
		\]
		However, it follows from Theorems \ref{half liouville theorem} and \ref{Thm:Acta} that the solution can not have two local maximum
		points, contradicting the fact that $0$ and $\lim_{i\rightarrow \infty} A_{i,2}$ are local
		maximum points of $w_{i}$. This finishes the proof of the claim.
		
		\vskip 8pt
		
		We now know that $\{0\}$ and $\{A_{i,2}\}$ are isolated blow-up
		points for $w_{i}$ and thus isolated simple. Denote by $\mathcal{S}$ the set of all blow-up points for $w_i$. Clearly, $\{0,A_{\infty,2}:=\lim_{i\rightarrow\infty}A_{i,2}\}\subset\mathcal{S}$.

		\vskip 8pt

		By Theorem \ref{isolated behavior-1}, we have 
		\[
		w_{i}(z)\rightarrow0\mathrm{~~in~~}C_{\mathrm{loc}}^{0}(\overline{\mathbb{R}_{-T}^{n}}\backslash\mathcal{S}).
		\]
		Moreover, for any compact set $K\subset\mathbb{R}_{-T}^{n}\backslash\mathcal{S}$,
		by interior estimates, we have 
		\[
		\|\nabla\log w_{i}\|_{C^{0}(K)}+\|\nabla^{2}\log w_{i}\|_{C^{0}(K)}\le C(K).
		\]
		By boundary local $C^1$ estimates, we have $w_{i}(z)\rightarrow0$ in $C_{\mathrm{loc}}^{\gamma}(\overline{\mathbb{R}_{-T}^{n}}\backslash\mathcal{S})\cap C_{\mathrm{loc}}^{1,\gamma}({\mathbb{R}_{-T}^{n}}\backslash\mathcal{S}), \forall~\gamma\in(0,1)$.
		
		It is not hard to see that for $y_0 \notin \mathcal{S}$, $\omega_{i}(z):=w_{i}(z)/w_{i}(y_{0})$ satisfies 
		\begin{align*}
		\begin{cases}
			\sigma_{2}^{1/2}(-\frac{2}{n-2}\frac{\nabla_{g_{i}}^{2}\omega_{i}}{\omega_{i}}+\frac{2n}{(n-2)^{2}}\frac{\ud\omega_{i}\otimes\ud\omega_{i}}{\omega_{i}^{2}}-\frac{2}{(n-2)^{2}}\frac{|\nabla\omega_{i}|_{g_{i}}^{2}}{\omega_{i}^{2}}g_{i}+A_{g_{i}})&\\
			=f\circ\Psi_{i}\omega_{i}(z)^{\frac{4}{n-2}}w_{i}^{\frac{4}{n-2}}(y_{0})+\rho_i^2 f_0 &\qquad\mathrm{in~~}\quad |z|<\frac{\rho}{2\rho_{i}},z_{n}>-T_{i}\\
		\frac{2}{n-2}\frac{\partial\omega_{i}}{\partial z_{n}}=-c\circ\Psi_{i}\omega_{i}^{\frac{n}{n-2}}w_{i}(y_{0})^{\frac{2}{n-2}} &\qquad \mathrm{on~~} \quad \{z_n=-T_i\}.
		\end{cases}
		\end{align*}
		
		By $C^1$ estimates and $\omega_{i}(y_{0})=1$, we actually know $\omega_{i}(z)\le C_{K}$
		on every compact set $K\subset\overline{\mathbb{R}_{-T}^{n}}\backslash\mathcal{S}$.
		This together with $C^1$ estimates implies $\omega_{i}\geq c_{K}>0$ in $K$ and $\|\log\omega_{i}\|_{C^{1}(K)}\le C$. Now there exists a subsequence
		$\omega_{i}$ converging to $\omega_{\infty}\in C_{\mathrm{loc}}^{1,1}(\mathbb{R}_{-T}^{n}\backslash\mathcal{S})\cap C_{\mathrm{loc}}^{0,1}(\overline{\mathbb{R}_{-T}^{n}}\backslash\mathcal{S})$
		which satisfies the following equation in the viscosity sense 
		\begin{align*}
		\begin{cases}
			\sigma_{2}\left(-\frac{2}{n-2}\frac{\nabla_{g_{\mathrm{E}}}^{2}\omega_{\infty}}{\omega_{\infty}}+\frac{2n}{(n-2)^{2}}\frac{\ud\omega_{\infty}\otimes\ud\omega_{\infty}}{\omega_{\infty}^{2}}-\frac{2}{(n-2)^{2}}\frac{|\nabla\omega_{\infty}|^{2}}{\omega_{\infty}^{2}}g_{\mathrm{E}}\right)=0 &\quad\mathrm{in~~}\mathbb{R}_{-T}^{n}\backslash\mathcal{S},\\
		\frac{\partial\omega_{\infty}}{\partial z_{n}}=0 &\quad \mathrm{on~~} \{z;z_{n}=-T\}\backslash\mathcal{S}.
		\end{cases}
		\end{align*}
		
		Denote
		\[
		\overline{\omega}_{\infty}(r)=\frac{1}{r^{n-1}}\int_{\pa^{+}B^+_{r,\theta}}\omega_{\infty}\ud\sigma
		\]
		and 
		\[
		\overline{\omega}_{\infty,1}(r)=\frac{1}{r^{n-1}}\int_{\pa^{+}B^+_{r,\theta}+A_{\infty,2}}\omega_{\infty}\ud\sigma.
		\]
		Since  $0$ and $A_{\infty,2}$ are isolated simple blow-up points for $w_i$, we have  $\big(r^{\frac{n-2}{2}}\overline{\omega}_{\infty}(r)\big)'\le0$ and
		$\big(r^{\frac{n-2}{2}}\bar{\omega}_{\infty,1}(r)\big)'\le0$ near
		zero. This implies that $\omega_{\infty}$ is singular at $0$ and $A_{\infty,2}$.
		Furthermore, by Theorem \ref{isolated behavior-1},  for $n=3,4$ we actually have \[\lim_{z\rightarrow0}\omega_{\infty}(z)|z|^{n-2}=a_{1},\]
		and \[\lim_{z\rightarrow A_{\infty,2}} \omega_{\infty}(z)|z-A_{\infty,2}|^{n-2}=a_{2}\]
		for some $a_{1},a_{2} \in \R_+$. However, by Theorem \ref{thm:(Bocher-Type-theorems--Li-Nguyen}
		and Proposition \ref{even reflection in Appendix} we know that $\omega_{\infty}$
		can not have two singular points. 							
	\end{proof}
	
	We would like to point out that the statement similar to above Lemma \ref{lem:pairwize_disjoint_sigular_pts}  in boundary Yamabe problem for dimension $n\geq 3$ is due to Han and Li \cite[Proposition 1.2]{han-li2}.

		\subsection{Compactness}\label{Subsect:compactness}
		
		For $\bar{H}\in\R$, a function $A_{n,\bar H}:[0,\infty) \to \R$ is introduced in \cite{Perales} by
		\[
		A_{n,\bar{H}}(R)=\begin{cases}
			(-\bar{H} R+n-1)^{n-1}/(n-1)^{n-1} & \mathrm{~~if~~}-\bar{H} R+n-1\geq0,\\
			0 & \mathrm{~~otherwise.}
		\end{cases}
		\]
		
		Fix a nonnegative constant $R$, let $M^{R}:=\{x\in M;d_{g}(x,\pa M)>R\}$
		and $\partial M^{R}$ be the boundary (as a metric subspace of
		$M$) of $M^{R}$.
		
		\begin{theorem}[\protect{Perales,  \cite[Theorem 1.3]{Perales}}]\label{Perales}
			Let $(M,g)$ be a connected and metrically complete Riemannian manifold
			of dimension $n\geq2$ with smooth boundary $\pa M$ and $\bar{H}\in\R$. Assume $\mathrm{Ric}_{g}(\overline{M}\backslash\partial M)\geq0,  \footnote{Remind that the mean curvature $H_{\pa M}$ in \cite{Perales} is with respect to the inward unit normal, it differs from ours by a minus sign. In other words, $H_{\pa M}=-H_g=-(n-1)h_g$.} H_{g}\geq-\bar{H}$
			on $\pa M$ and $\mathrm{Vol}(\partial M,g)<\infty$.
			Then, for $\mathcal{L}^{1}$-almost everywhere $R>0$, 
			\[
			\mathrm{Vol}\left(\partial M^{R},g\right)\leq\mathrm{Vol}(\partial M,g)A_{n,\bar{H}}(R),
			\]
			where $\mathcal{L}^{1}$ denotes the 1-dimensional Lebesgue measure.
		\end{theorem}			
		
		Suppose $\{u_i\}$ is a sequence of positive smooth solutions to \eqref{PDE:sigma_2_new} with $\max_{\overline M} u_i \to +\infty$ for dimensions $n=3,4$. It follows from Lemma \ref{lem:pairwize_disjoint_sigular_pts} that the singular set $\mathcal{S}=\mathcal{S}[\{u_i\}]$ is a finite set of isolated (simple) blow-up points, denoted by $\mathcal{S}=\{x_1,\cdots,x_N\} \subset \overline M$.

		As we shall see in Section \ref{Sect:degree_theory}, if the boundary is totally non-umbilic,   then $\mathcal{S}$ consists of only interior isolated blow-up points. The following assertion is crucial to our existence result.

		\begin{proposition}\label{prop:no blow-up pts}
			Assume $\mathcal{S}$ contains no boundary isolated blow-up points, then $\mathcal{S}=\emptyset$.
		\end{proposition}

		\begin{proof}
			As before, we assume $h_g=0$ on $\pa M$. The same proof of Theorem \ref{Thm:isolated behavior} can be used to show that for any `\emph{regular}' point $y_{0}\in \overline M \backslash \mathcal{S}$ and $\gamma \in (0,1)$,
			\begin{align*}
				u_{i}(x)\rightarrow0 & \quad\mathrm{~~in~~}\quad C_{\mathrm{loc}}^{\gamma}(\overline M \backslash \mathcal{S}),\\
				v_{i}(x):=\frac{u_{i}(x)}{u_{i}(y_{0})}\rightarrow v_{\infty} & \quad\mathrm{~~in~~}\quad C_{\mathrm{loc}}^{\gamma}(\overline M \backslash \mathcal{S})\cap C_{\mathrm{loc}}^{1,\gamma}(M \backslash \mathcal{S}),\\
				\lim_{x\rightarrow x_{k}}v_{\infty}(x)d_{g}^{n-2}(x,x_{k})=a_k>0 &\quad \mathrm{~~for~~}\quad 1 \leq k \leq N.
			\end{align*}
			Moreover,  $v_{\infty}$ is a positive viscosity solution to
			$\sigma_{2}(A_{v_{\infty}^{4/(n-2)}g})=0$ in $\overline M \backslash \mathcal{S}$
			and $\frac{\partial v_{\infty}}{\pa\vec{n}}=0$
			on $\partial M$.

			Let  $\tilde r$ be a smooth function in $\overline M$ and positive in $\overline M\backslash \mathcal{S}$ such that $\tilde r=r_{x_k}=d_g(x_k,\cdot)$ in a neighborhood of each $x_k \in \mathcal{S}$.
			Let $\{y\}$ be the geodesic normal coordinates around $x_k$, then $r_{x_k}(y)=|y|$ near $x_k$. Denote by $\{z=|y|^{-2}y\}$ the inverted coordinates in a component $N_k:=U_{x_k} \setminus\{x_k\}$. We will call each $N_k$ an \emph{end}.  Notice that
			$$g_{v_\infty}:=v_\infty^{\frac{4}{n-2}}g=(v_\infty \tilde r^{n-2} )^{\frac{4}{n-2}}\tilde r^{-4}g.$$
			Under inverted coordinates $\{z\}$, it is standard to verify that the metric 
			$$g_\ast=\Phi^{\ast}(\tilde r^{-4}g) \qquad \mathrm{with~~} y=\Phi(z)$$	
			is asymptotically flat of order two, see Lee-Parker  \cite[Definition 6.3]{Lee-Parker} or \cite{Gursky-Viaclovsky}; and thus
			$$\Phi^\ast(g_{v_\infty})=[(v_\infty \tilde r^{n-2})\circ \Phi]^{\frac{4}{n-2}} g_\ast.$$
			\begin{claim}
				The singular set $\mathcal{S}$ consists of one interior isolated blow-up point.
			\end{claim}

			We may choose $\delta$ sufficiently small such that $\mathcal{S} \subset M^\delta$ and define $M_{\mathrm{reg}}^\delta:=\overline{M^\delta}\backslash \mathcal{S}$. We would like to point out that part of Gursky-Viaclosky's argument in \cite[Section 7]{Gursky-Viaclovsky} is applicable to $(M_{\mathrm{reg}}^\delta,g_{v_\infty})$ and enables us to derive our conclusion. Since the integral estimate  in \cite[Theorem 3.5]{Gursky-Viaclovsky} is of independent interest and true for $\sigma_k$-curvature equation with $k\geq n/2$, with the help of our available estimates and Guan-Viaclovsky-Wang \cite[Theorem 1]{GVW}, we can see that their argument therein is still valid for $k \geq n/2$ corresponding to our case that $k=2, n=3,4$. We only give an outline of proof. Readers are referred to \cite[Section 7 on pp.511-524]{Gursky-Viaclovsky}  for details. 
			
			Let us  write
			$$M_{\mathrm{reg}}^\delta=M_0 \cup M_1$$
			and fix  a base point $y_0 \in M_0$, where $M_0$ is a smooth compact manifold with boundary containing $y_0$ and we use inverted coordinates $\{z\}$ on each \emph{end} $N_k$ to identify $N_k$ with $\Rn \backslash B_{R_0}(0)$  for some $R_0\in \R_+$ such that
			$$M_1 \approx \bigcup_{k=1}^N (\Rn \backslash B_{R_0}(0)).$$
			
			As in \cite{Gursky-Viaclovsky}, the first step is to show that $(M_{\mathrm{reg}}^\delta,g_{v_\infty})$ is geodesically complete. The next step is to establish  a `\emph{weak}' version of Bishop-Gromov type relative volume comparison\footnote{ The argument to derive the relative volume comparison is indeed a localized analysis around the interior point $y_0$, thus it still works well here.} for the $C^{1,1}$-metric $g_{v_\infty}$  in  $M_{\mathrm{reg}}^\delta$: For the above base point $y_0 \in M_0$, the ratio
			\begin{equation}\label{ineq:Gromov-Bishop volume comparison}
				\frac{\mathrm{Vol}_{g_{v_\infty}}(B_r^{g_{v_\infty}}(y_0))}{r^n}\leq |B_1(0)|
			\end{equation}
			is nonincreasing with respect to $r$, where $B_r^{g_{v_\infty}}(y_0)$ denotes a geodesic ball in the metric $g_{v_\infty}$ and $|B_1(0)|$ stands for the volume of the unit ball. In the third step, using the technique of tangent cone analysis developed  in \cite{Gursky-Viaclovsky}, we can prove the volume growth on each end: For each end $N_k$ and the above base point $y_0 \in M_0$, there holds
			\begin{equation}\label{volume_growth_in_ends}
				\lim_{r \to \infty}\frac{\mathrm{Vol}_{g_{v_\infty}}\left(B_r^{g_{v_\infty}}(y_0)\cap N_k\right)}{r^n}\geq |B_1(0)|.
			\end{equation}
			Finally, we combine \eqref{ineq:Gromov-Bishop volume comparison} and \eqref{volume_growth_in_ends} to conclude that $N=1$. This finishes the proof of the above claim.
			
			\vskip 8pt
			
			We next establish that the unique interior blow-up point is impossible either.

			Let $g_{i}=v_{i}^{4/(n-2)}g$ and define 
			\[
			M_{g_{i}}^{R}:=\{x\in M;d_{g_{i}}(x,\pa M)> R\} \qquad \mathrm{~~for~~} R \in \R.
			\]
		        Since $A_{u_{i}}\in\Gamma_{2}^{+}$, by Guan-Viaclovsky-Wang
			\cite[Theorem 1]{GVW} we have $\mathrm{Ric}_{g_{i}}\ge0$ in $\overline{M}$.
			Meanwhile,  on $\pa M$ we have
			\[
			H_{g_{i}}=u_{i}(y_{0})^{\frac{2}{n-2}}H_{u_{i}^{4/(n-2)}g}=(n-1)u_{i}(y_{0})^{\frac{2}{n-2}}c\geq 0.
			\]
			
			By definition we have $A_{n,0}(R)=1$ for $\bar H=0$. Then, on each fixed $(M,g_{i})$, we may apply Theorem \ref{Perales} to obtain
			\[
			\mathrm{Vol}\left(\partial M_{g_{i}}^{R},g_{i}\right)\leq\mathrm{Vol}(\partial M,g_{i}).
			\]
			
			On one hand, notice that 
			\[
			\mathrm{Vol}(\partial M,g_{i})=\int_{\pa M}v_{i}^{\frac{2(n-1)}{n-2}}\ud\sigma_{g}.
			\]
			Then there exists a positive constant $C_{\ast}$ independent of $i$
			such that for $\mathcal{L}^{1}$-almost everywhere $R>0$, 
			\begin{equation}
				\mathrm{Vol}(\partial M,g_{i})\leq C_{\ast}.\label{est:area_bdd}
			\end{equation}
			
			On the other hand, we always assume $\partial M_{g_{i}}^{R}\neq\emptyset$.
			For any $x\in\pa M_{g_{i}}^{R}$, let $\gamma(t)$ be a unit-speed
			minimizing geodesic to realize the distance $d_{g}(x,\pa M)$ starting
			at $x$ and ending at some $x'\in\pa M$. There exists a positive
			constant $C_{r_{0}}$ independent of $i$ such that for all $i\in\N$,
			\[
			\|v_{i}\|_{C^{0}(\overline{M}\backslash B_{r_{0}}(x_{0}))}\leq C_{r_{0}}.
			\]
			Let
			\[
			\mathrm{Diam}(M)=\max_{p,q\in\overline{M}}d_{g}(p,q)
			\]
			with 
			\[
			d_{g}(p,q)=\inf\{L_{g}(\gamma);\mathrm{~~all~~piecewise~~}C^{1}\mathrm{~~paths~~}\gamma:[0,\alpha]\to\overline{M}\mathrm{~~joining~~}p\mathrm{~~and~~}q\}
			\]
			and 
			\[
			L_{g}(\gamma)=\int_{0}^{\alpha}|\gamma'(t)|_{g}\ud t.
			\]
			With the choice of $R>2C_{r_{0}}^{2/(n-2)}\mathrm{Diam}(M)$, we
			claim that $x\in B_{r_{0}}(x_{0})\backslash\{x_{0}\}$. 
			
			Otherwise, if $x\notin B_{r_{0}}(x_{0})$, then we have 
			\begin{align*}
				R=d_{g_{i}}(x,\pa M)\leq d_{g_{i}}(x,x')\leq & \int_{0}^{d_{g}(x,\pa M)}|\gamma'(t)|_{g_{i}}\ud t\\
				\leq & C_{r_{0}}^{\frac{2}{n-2}}\mathrm{Diam}(M),
			\end{align*}
			which contradicts the choice of $R$.
			
			We define 
			\[
			t_{0}:=\inf\{t\in(0,d_{g}(x,\pa M)];\gamma(t)\notin B_{r_{0}}(x_{0})\}.
			\]
			By definition we know that $\gamma([t_{0},d_{g}(x,\pa M)])$ is outside
			of $B_{r_{0}}(x_{0})$.
			
			By letting $i\to\infty$ we have 
			\begin{align*}
				R \leq & \int_{0}^{d_{g}(x,\pa M)}|\gamma'(t)|_{g_{i}}\ud t\\
				= & \int_{0}^{t_{0}}|\gamma'(t)|_{g_{i}}\ud t+\int_{t_{0}}^{d_{g}(x,\pa M)}|\gamma'(t)|_{g_{i}}\ud t\\
				\leq & a^{\frac{2}{n-2}}\int_{0}^{t_{0}}d_{g}(\gamma(t),x_{0})^{-2}\ud t+C_{r_{0}}^{\frac{2}{n-2}}\mathrm{Diam}(M)\\
				\leq & a^{\frac{2}{n-2}}d_{g}(x,x_{0})^{-2}d_{g}(x_{0},\pa M)+\frac{R}{2}.
			\end{align*}
			This indicates that 
			\begin{equation}
				d_{g}(x,x_{0})\leq a^{\frac{1}{n-2}}\sqrt{\frac{2d_{g}(x_{0},\pa M)}{R}}.\label{est:r(x)}
			\end{equation}
			
			It is not hard to see that $M_{g_{i}}^{R}$ is star-shaped with
			respect to $x_{0}$. We may choose geodesic polar coordinates
			around $x_{0}$ such that 
			\[
			g=\ud r^{2}+h_{\alpha\beta}(r,\theta)\ud\theta_{\alpha}\otimes\ud \theta_{\beta},\qquad(r,\theta)\in(0,\mathrm{Diam}(M)]\times\Sp^{n-1},
			\]
			where $r(x)=d_{g}(x,x_{0})$. Moreover, there holds 
			\[
			\lim_{r\to0}\frac{h(r,\theta)}{r^{2}}=g_{\Sp^{n-1}}\qquad\mathrm{uniformly~~for~~all~~}\theta\in\Sp^{n-1}.
			\]
			
			If choosing anther larger $R$, by \eqref{est:r(x)}
			we have 
			\begin{align*}
				\mathrm{Vol}(\partial M_{g_{i}}^{R},g_{i})= & \int_{\partial M_{g_{i}}^{\bar \delta}}v_{i}^{\frac{2(n-1)}{n-2}}\ud\sigma_{g}\\
				\overset{i\to\infty}{\longrightarrow} &\int_{\Sp^{n-1}}v_\infty(x)\sqrt{\det h(r,\theta)}\ud\theta_{1}\wedge\cdots\wedge\ud\theta_{n-1}\\
				\geq& \left(\frac{a}{2}\right)^{\frac{2(n-1)}{n-2}}\int_{\Sp^{n-1}}d_{g}(x,x_{0})^{-2(n-1)}\sqrt{\det h(r,\theta)}\ud\theta_{1}\wedge\cdots\wedge\ud\theta_{n-1}\\
				\geq & C_{0}R^{\frac{n-1}{2}}
			\end{align*}
			for some positive constant $C_{0}$ independent of $i$ and $R$.
			This yields a contradiction with \eqref{est:area_bdd} provided that
			$R$ is large enough.
		\end{proof}

		\section{Existence of conformal metrics: Proof of Theorem \ref{Thm:existence}}\label{Sect:degree_theory}
		
		The method of Leray-Schauder degree theory will be applied to solve the boundary $\sigma_2$-curvature equation \eqref{PDE:prescribing_sigma_2} and complete the proof of Theorem \ref{Thm:existence} as well.

		We digress for a moment to give a brief summary of variational aspects of boundary $\sigma_2$-curvature equation. For the $\sigma_2$-curvature equation \eqref{PDE:prescribing_sigma_2} with $c=0$, a variational characterization was described in \cite{Chen}, where the following functional
	was introduced:
	\begin{equation}\label{def:F_2}
	\mathcal{F}_2[\tilde g]=\int_M \sigma_2(A_{\tilde g}) \ud \mu_{\tilde g}+\int_{\pa M} \mathcal{B}_{\tilde g}^2 \ud \sigma_{\tilde g}, \quad \tilde g \in [g],
	\end{equation}
	see \cite[formula $(2)$ on p.1030]{Chen}  for the definition of $\mathcal{B}_{\tilde g}^2$.\footnote{For $3 \leq k\leq n/2$, nowadays $\mathcal{B}_g^k$ is only known in locally conformally flat manifolds with boundary; see \cite{Chen}.}  The above $\mathcal{F}_2$ appears in the  Gauss-Bonnet formula \cite{Branson-Gilkey} on a four-dimensional compact manifold with boundary:
	$$32\pi^2\chi(M^4,\pa M^4)=\int_{M^4} |\mathcal{W}_g|^2 \ud\mu_g+16 \mathcal{F}_2[g],$$
	where $\chi$ is the Euler characteristic and $\mathcal{W}_g$ is the Weyl tensor.
	Then, for $\phi \in C^\infty(\overline M)$ and $n\neq 4$ there holds
	$$\frac{\ud}{\ud t}\bigg|_{t=0} \mathcal{F}_2[e^{-2t\phi}\tilde g]=(4-n)\left[\int_M\sigma_2(A_{\tilde g}) \phi \ud \sigma_{\tilde g}+\int_{\pa M}\mathcal{B}_{\tilde g}^2 \phi \ud \sigma_{\tilde g}\right].$$
	Later in critical dimension four,   a new conformal primitive (still denoted by $\mathcal{F}_2$) of $(\sigma_2,\mathcal{B}_g^2)$ is complemented by J. Case and Y. Wang \cite[Proposition 2.3]{Case-Wang}:
	$$\mathcal{F}_2[\tilde g]=\int_0^1\left[\int_{M^4}u\sigma_2(A_{\tilde g_s}) \ud \mu_{\tilde g_s}+\int_{\pa M^4} u\mathcal{B}_{\tilde g_s}^2 \ud \sigma_{\tilde g_s}\right] \ud s,$$
	where $\tilde g=e^{-2u}g$ and $\tilde g_s=e^{-2s u}g$. Then, for $\phi \in C^\infty(\overline M)$ there holds
	$$\frac{\ud}{\ud t}\bigg|_{t=0} \mathcal{F}_2[e^{-2t\phi}\tilde g]=\int_{M^4}\sigma_2(A_{\tilde g}) \phi \ud \sigma_{\tilde g}+\int_{\pa M^4}\mathcal{B}_{\tilde g}^2 \phi \ud \sigma_{\tilde g}.$$
	Thanks to a clever observation in \cite[Lemma 5]{Chen}, given $0<f \in C^\infty(\overline M)$ and $c=0$, on compact manifolds with umbilic boundary, the boundary $\sigma_2$-curvature equation \eqref{PDE:prescribing_sigma_2}   is equivalent to the Euler-Lagrange equation of a critical point $g_u$ with $\lambda(A_u) \in \Gamma_2^+$ of 
	$\mathcal{F}_2[\tilde g]$ under constraint that $\int_M  f \ud \mu_{\tilde g}=1$ for all $\tilde g \in [g]$. 
		
		Given a smooth compact manifold $(M, g)$ with boundary satisfying $\lambda(A_g) \in \Gamma_2^+$ and  $h_g\ge 0$ on $\partial M$, we may assume 
		\begin{equation}\label{positive matrix}
			\lambda_{n}V_{g}^{\frac{2}{n+1}}g-A_g>0,
		\end{equation}
		where $V_{g}=\mathrm{Vol}(M,g)$ and $\lambda_{n}=(\binom{n}{2})^{-1/2}$. Otherwise, we choose  $\tilde{g}=Bg$ for some sufficiently large $B \in \R_+$ as a background metric, such that
		$$\lambda(\tilde g^{-1} A_{\tilde{g}})=B^{-1}\lambda(A_g)\in \Gamma_2^+, \qquad h_{\tilde{g}}=\frac{1}{\sqrt{B}}h_g\ge 0$$
		and
		\begin{equation*}
			\lambda_{n}V_{\tilde{g}}^{\frac{2}{n+1}}\tilde{g}-A_{\tilde{g}}=\lambda_{n}B^{\frac{n}{n+1}+1}V_{g}^{\frac{2}{n+1}}g-A_g>0.
		\end{equation*}

		\begin{proof}[Proof of Theorem \ref{Thm:existence}]					
			We first  introduce a family of  boundary $\sigma_2$-curvature equations\footnote{To apply both local estimates and blow up analysis to \eqref{eq:path},  we may replace $\sigma_2^{1/2}(A_u)$ by $\sigma_2^{1/2}(A_u+S_g)$ with some smooth symmetric two-tensor $S_g\geq 0$ for the $\sigma_2$-curvature equations in \eqref{PDE:sigma_2_with_bdry} and \eqref{PDE:sigma_2_new}, due to the fact that  $\lambda(A_u) \in \Gamma_2^+$ implies $\lambda(A_u+S_g) \in \Gamma_2^+$.}
			\begin{equation}
				\begin{cases}
					\quad\sigma_{2}^{1/2}(\nabla^{2}u+\ud u\otimes\ud u-\frac{1}{2}|\nabla u|^{2}g+A_{g}+S_{g})\\
					=(1-t)(\int_{M}e^{-(n+1)u}\ud\mu_{g})^{\frac{2}{n+1}}+\zeta(t)fe^{-2u} & \quad\mathrm{~~in~~}M,\\
					\frac{\partial u}{\partial\vec{n}}=\zeta(t)(ce^{-u}-h_g) & \quad\mathrm{~~on~~}\pa M,
				\end{cases}\label{eq:path}
			\end{equation}
			where $S_{g}(t,\cdot)=(1-\zeta(t))(\lambda_{n}V_{g}^{2/(n+1)}g-A_{g})\ge 0$
			and $\zeta$ is a smooth nonnegative function in $[0,1]$ such that $\zeta(0)=0,$
			and $\zeta(t)=1$ for $1/2\le t\le1$. It is noticed that the resulting equation  \eqref{eq:path} at $t=1$ is exactly \eqref{PDE:prescribing_sigma_2}.
			
			Under assumption that $\lambda(A_g) \in \Gamma_2^+$ and $h_g\ge 0$,  thanks to Proposition \ref{prop:sigma_2_1st_bdry_eigenvalue} the background metric $g$ satisfies $\lambda(A_g)\in \Gamma_2^+$ and $h_g=0$ on $\pa M$. 
			
			Nowadays the strategy of the degree theory is more or less standard. The proof is accomplished by two steps.
			
			\emph{Step 1:} At $t=0$, the initial equation \eqref{eq:path} admits a unique solution $u=0$ and the linearized operator at $0$ is invertible. 
			Since the initial equation is identical to that of \cite[Theorem 2]{Chen},  we omit the proof.

			\emph{Step 2:} The Leray-Schauder degree is defined in a relatively compact set in $C^{4,\alpha}(\overline M)$. This reduces to establishing the uniform a priori $C^0$ estimates for solutions to \eqref{eq:path} with all $t \in [0,1]$. The $C^0$-bounds are enough to derive  a priori $C^1,C^2$ estimates and improve the regularity of solutions to $C^{4,\alpha}$ via the Evans-Krylov theory. By homotopy invariance the existence result follows.  See, for example,  \cite[pp.104-106]{Li-Luc3}  with slight modifications. 
			
			The key ingredient is to obtain the uniform $C^0$ estimates.

			\begin{lemma} \label{lem:finite supervolume with finite interval}Let
				$u$ be a solution to (\ref{eq:path}) for $t \in [0,1]$ with assumption that $h_g=0$ on $\pa M$. Then there exists a uniform positive constant $C$ such that
				\[
				(1-t)(\int_{M}e^{-(n+1)u}\ud\mu_{g})^{\frac{2}{n+1}}\le C.
				\]
			\end{lemma} 
			\begin{proof}
				Let $u(x_{0})=\max_{\overline{M}}u(x)$. If $x_0\in M$, then  we know $\nabla u(x_{0})=0$ and $\nabla^{2}u(x_{0})\le0$.
				Thus, we have 
				\[
				(1-t)(\int_{M}e^{-(n+1)u}\ud\mu_{g})^{\frac{2}{n+1}}\le\sigma_{2}^{1/2}(A_{g}(x_{0})+S_{g}(x_{0})).
				\]
				The desired estimate follows. 
				If $x_0\in \partial M$, then 
				$0\ge \frac{\partial u}{\partial\vec{n}}(x_0)=\zeta(t)ce^{-u(x_0)}\ge 0$. So  $\nabla u(x_{0})=0$ and $\nabla^{2}u(x_{0})\le0$ on $\pa M$. The left argument is similar as above. 
			\end{proof}

			\begin{lemma}\label{lem:totally_nonumbilic} Suppose $\partial M$
				is non-umbilic everywhere. Then for all $t \in [0,1]$
				there exists a uniform positive constant $C$ such that 
				\[
				\|u\|_{C^{0}(\overline{M})}\le C.
				\]
			\end{lemma} 
			\begin{proof}
				(1) $\inf_{M}u\ge-C$ for $t\in[0,1-\varepsilon]$.
				\vskip 4pt
				
				We argue by contradiction. Suppose that there exists a sequence of solutions $\{u_{i}\}$ to
				(\ref{eq:path}) corresponding to some $t_{i}\in[0,1-\varepsilon]$ such
				that $\inf_{M}u_{i}=u_{i}(p_{i})\rightarrow-\infty$ and $p_{i}\rightarrow p_{0}\in\overline{M}$.
				Let $\varepsilon_{i}=e^{\inf_{M}u_{i}}$ and $d_{i}=d_g(p_i,\pa M)$.
				
				Our discussion is divided into two cases.
				
				\emph{Case 1:} $d_{i}/\varepsilon_{i}\rightarrow\infty$.
				
				Under geodesic normal coordinates around $p_{i}$, let us denote
				\[
				\mathcal{T}_{i}:B_{d_i/\varepsilon_i}(0)\subset\mathbb{R}^{n}\rightarrow M \qquad \mathrm{with~~}\quad
				\mathcal{T}_{i}(x)=\exp_{p_{i}}(\varepsilon_{i}x).
				\]
				
				We define 
				\[
				\tilde{u}_{i}(x)=u_{i}\left(\mathcal{T}_{i}(x)\right)-\log\varepsilon_{i}\quad\mathrm{~~in~~}\quad B_{d_i/\varepsilon_i}(0)
				\]
				with  $\tilde{u}_{i}(x)\ge0$, and the metrics 
				\[
				g_{i}=\varepsilon_{i}^{-2}\mathcal{T}_{i}^{*}(g) \to g_{\mathrm{E}} \qquad \mathrm{in~~} C_{\mathrm{loc}}^3(\Rn).				\]
				Moreover, $\tilde u_i$ satisfies
				\begin{align*}
					& \sigma_{2}^{1/2}(\nabla_{g_{i}}^{2}\tilde{u}_{i}+\nabla_{g_{i}}\tilde{u}_{i}\otimes\nabla_{g_{i}}\tilde{u}_{i}-\frac{|\nabla\tilde{u}_{i}|_{g_{i}}^{2}}{2}g_{i}+A_{g_{i}}+S_{g_{i}})\\
					= & \ve_{i}^{2}(1-t_{i})(\int_{M}e^{-(n+1)u_{i}}\ud\mu_{g})^{\frac{2}{n+1}}+\zeta(t_{i})f(\mathcal{T}_{i}(x))e^{-2\tilde{u}_{i}} \qquad \mathrm{in~~} B_{d_i/\varepsilon_i}(0).
				\end{align*}
				
				Since $\tilde{u}_{i}\ge0$ in $B_{d_{i}/\ve_{i}}(0)$, we apply local $C^{1}$ estimates  to show
				\[
				\sup_{B_r(0)}|\nabla\tilde{u}_{i}|\le C_{r}.
				\]
				This together with $\tilde{u}_{i}(0)=0$ implies $\sup_{B_r(0)}|\tilde{u}_{i}|\le C$.
				On the other hand, it follows from  Lemma \ref{lem:finite supervolume with finite interval} that for $B_{1}(0)\subset B_{d_{i}/\ve_{i}}(0)$,
				\[
				\int_{B_{1}(0)}e^{-(n+1)\tilde{u}_{i}}\ud\mu_{g_{i}}=\varepsilon_{i}\int_{B_{\ve_{i}}(p_{i})}e^{-(n+1)u_{i}}\ud\mu_{g}\le C\varepsilon_{i}.
				\]
				Hence we arrive at a contradiction.
				
				\emph{Case 2:} $d_{i}/\varepsilon_{i}\le C$ for all $i$.
				
				Let $d_{i}:=d_{g}(p_{i},p_{i}')=d_g(p_i,\pa M)$ for some $p_{i}'\in\partial M$.
				Clearly, $d_{i}\rightarrow0$ as $i\rightarrow\infty.$ Under Fermi coordinates around $p_i'$, we let 
				\[
				\Psi_{i}(x):=G_{p_{i}'}(\varepsilon_{i}x)\qquad \mathrm{in~~} B_{\rho/\ve_i}^+
				\]
				
				Define
				\[
				\tilde{u}_{i}(x)=u_{i}(\Psi_{i}(x))-\log\varepsilon_{i}
				\]
				with  $\tilde{u}_{i}(x)\ge0$, and  the metrics 
				\[
				g_{i}=\varepsilon_{i}^{-2}\Psi_{i}^{*}(g)\rightarrow g_{\mathrm{E}}\quad\mathrm{~~in~~} C_{\mathrm{loc}}^{3}(\overline{\Rn_+}).
				\]
				Then $\tilde u_i$ satisfies
				\begin{align*}
					& \sigma_{2}^{1/2}(\nabla_{g_{i}}^{2}\tilde{u}_{i}+\nabla_{g_{i}}\tilde{u}_{i}\otimes\nabla_{g_{i}}\tilde{u}_{i}-\frac{|\nabla\tilde{u}_{i}|_{g_{i}}^{2}}{2}g_{i}+A_{g_{i}}+S_{g_{i}})\\
					= & \ve_{i}^{2}(1-t_{i})(\int_{M}e^{-(n+1){u}_{i}}\ud\mu_{g})^{\frac{2}{n+1}}+\zeta(t_{i})f\circ \Psi_{i}e^{-2\tilde{u}_{i}} \qquad \mathrm{~~in~~}   B_{\rho/\ve_i}^+
				\end{align*}
				together with the boundary condition that
				$$\frac{\partial\tilde{u}_{i}}{\partial \vec n_{g_i}}=\zeta(t_i) c\circ \Psi_ie^{-\tilde{u}_{i}} \qquad \mathrm{on~~} \quad D_{\rho/\ve_i}.$$

				Since $\tilde{u}_{i}\ge0$ in any
				compact set $K\subset\overline{\mathbb{R}_{+}^{n}}$, we apply local $C^{1}$ estimates to show
				\[
				\sup_{K}|\nabla\tilde{u}_{i}|\le C.
				\]
				This together with  $\tilde{u}_{i}(d_i/\ve_i \mathbf{e}_n)=0$ implies  $\sup_{K}|\tilde{u}_{i}|\le C$. On the other hand, by Lemma \ref{lem:finite supervolume with finite interval} we obtain
				\[
				\int_{K}e^{-(n+1)\tilde{u}_{i}}\ud \sigma_{g_{i}}=\varepsilon_{i}\int_{\Psi_{i}(K)}e^{-(n+1)u_{i}}\ud \sigma_{g}\le C\varepsilon_{i}.
				\]
				This is a contradiction.
				\vskip 4pt				
				(2) As $t\rightarrow1$, $\inf_{M}u\geq-C$ for some positive constant $C$.
				\vskip 4pt
				
				The equation in \eqref{eq:path} becomes 
				\begin{align*}
					&\sigma_{2}^{1/2}(\nabla^{2}u+\ud u\otimes\ud u-\frac{1}{2}|\nabla u|^{2}g+A_{g})\\
					=&fe^{-2u}+(1-t)(\int_{M}e^{-(n+1)u}\ud\mu_{g})^{\frac{2}{n+1}}\qquad\mathrm{~~in~~} M.
				\end{align*}
				
				Suppose not, there exists a sequence of solutions $\{u_{i}\}$ to
				(\ref{eq:path}) for some $t_i\in [0,1]$ such that $t_{i}\rightarrow1$ and $\inf_{M}u_{i}\rightarrow-\infty$.
				Notice that $\zeta(t_{i})=1$ and 	$S_{g}(t_i)=0$
				as $t_{i}\rightarrow1$.  
				By Lemma \ref{lem fudament-1-1}, we know there exist local maximum points $q_{i,a}$ of $e^{-\frac{n-2}{2}u_i}$ obtained from Lemma \ref{lem fudament-1-1} and $1\le a\le N_i$. Indeed, the global maximum point of $u_i$ belongs to $\{q_{i,a}\}_{1\le i\le N_i}$.

				Furthermore, it follows from Lemma \ref{lem:pairwize_disjoint_sigular_pts} that there exists a set of finitely many isolated simple blow-up points for $e^{-\frac{n-2}{2}u_i}$, denoted by $\mathcal{S}:=\{q_1, \cdots, q_N\}$.
				
				\begin{claim}
					Under assumption of totally non-umibilic boundary, the singular set $\mathcal{S}$ consists of only interior isolated blow-up points. 
				\end{claim}
				
				If not, without loss of generality we may assume $q_1\in \partial M$. 	
				Without loss of generality, we assume that 
				$$q_{i,1}\rightarrow q_1,$$
				where $q_{i,1}$ is a local maximum point of $e^{-\frac{n-2}{2}u_i}$ and  $\varepsilon_{i}:=e^{u_{i}(q_{i,1})}\rightarrow 0$.

				Denote $d_i:=d_g(q_{i,1},\pa M):=d_g(q_{i,1},q_{i,1}')$ for some $q_{i,1}' \in \pa M$. Clearly, $q_{i,1}' \to q_1$. By Lemma \ref{lem:boundary bounded limit}, we know that  $d_{i}/\varepsilon_{i}\le C$  for all $i$.	
				\vskip 8pt
				
				%
				%

				Let
				\[
				\Psi_{i}(x)=G_{q_{i,1}'}(\ve_i x):B_{\frac{\rho}{\ve_{i}}}^{+}\to\overline{M}
				\]
				and 
				\[
				g_{i}=\ve_{i}^{-2}\Psi_{i}^{\ast}(g).
				\]
				
				We define
				\[
				\tilde{u}_{i}(x):=u_{i}(\Psi_{i}(x))-\log\varepsilon_{i} \geq 0 \qquad\mathrm{for~~}x\in B_{\frac{\rho}{\ve_{i}}}^{+}
				\]
				and				\[
				g_{\tilde{u}_{i}}:=e^{-2\tilde{u}_{i}}g_{i}=\Psi_{i}^{\ast}(e^{-2u_{i}}g).
				\]
				
				Under Fermi coordinates around $q_{i,1}$, for any boundary point $q' $ near $q_{i,1}$ we may express
				\[
				q'=\Psi_{i}((x',0))=G_{q_{i,1}'}(\ve_i(x',0)).
				\]
				Then we conclude from conformal invariance \eqref{trace_free-2nd_ff}
				of $\mathring{L}_{g}$ that 
				\begin{equation*}
					\mathring{L}_{g_{\tilde{u}_{i}}}(x')=\mathring{L}_{g_{u_{i}}}\circ\Psi_{i}(x',0)=e^{-u_{i}\circ\Psi_{i}(x')}\mathring{L}_{g}(q')=e^{-\tilde{u}_{i}((x',0))}\ve_{i}^{-1}\mathring{L}_{g}(q').
				\end{equation*}
				In particular, there holds
				\begin{equation}\label{eq:trace_2nd_ff}
					\mathring{L}_{g_{\tilde{u}_{i}}}(0)=e^{-\tilde{u}_{i}(0)}\ve_{i}^{-1}\mathring{L}_{g}(q_{i,1}').
				\end{equation}
				
				Now $\tilde{u}_{i}$ satisfies 
				\begin{align*}
				\begin{cases}
					 \quad \sigma_{2}^{1/2}(\nabla_{g_{i}}^{2}\tilde{u}_{i}+\ud\tilde{u}_{i}\otimes\ud\tilde{u}_{i}-\frac{1}{2}|\nabla\tilde{u}_{i}|_{g_{i}}^{2}g_{i}+A_{g_{i}})&\\
					=  \varepsilon_{i}^{2}(1-t_{i})(\int_{M}e^{-(n+1)u_{i}}\ud\mu_{g})^{\frac{2}{n+1}}+f\circ \Psi_{i}e^{-2\tilde{u}_{i}} &\quad \mathrm{~~in~~} B_{\frac{\rho}{\ve_{i}}}^{+},\\
				\frac{\partial\tilde{u}_{i}}{\partial \vec n_{g_{i}}}=\zeta(t_i) c\circ \Psi_i e^{-\tilde{u}_{i}} &\quad \mathrm{~~on~~} D_{\frac{\rho}{\ve_{i}}}.
				\end{cases}
				\end{align*}
				along with the boundary condition that

				By the boundary local $C^{2}$ estimates in  Lemma \ref{lem:double_normal_derivatives}, combined with local $C^{1}$ estimates,
				for any compact $K\subset\overline{\Rn_{+}}$ we obtain 
				\begin{equation}
					\sup_{K}|\tilde{u}_{i}|+|\nabla\tilde{u}_{i}|+|\nabla^{2}\tilde{u}_{i}|\le C.\label{est: tilde_u_i}
				\end{equation}
				By \eqref{est: tilde_u_i} and Evans-Krylov theory, we know that up to a subsequence, $\tilde{u}_{i}$ converges in $C_{\mathrm{loc}}^{3}(\overline{\Rn_{+}})$
				to a smooth solution $\tilde u$ of
				\[
				\begin{cases}
					\sigma_{2}^{1/2}(\nabla^{2} \tilde u+\ud \tilde u\otimes\ud \tilde u-\frac{1}{2}|\nabla \tilde u|^{2}g_{\mathrm{E}})=f(q_1)e^{-2\tilde u} & \mathrm{~~in~~}\mathbb{R}_{+}^{n}\\
					\frac{\partial \tilde u}{\partial x_n}=c(q_1)e^{-\tilde u} & \mathrm{~~on~~}\partial\mathbb{R}_{+}^{n}.
				\end{cases}
				\]
				Moreover, $e^{-2\tilde u}|\ud x|^{2}$ is the pull-back of the round metric on a spherical cap
				and thus has umbilic boundary.
				
				On one hand, it is not hard to see that 
				\[
				\lim_{i\to\infty}\mathring{L}_{g_{\tilde{u}_{i}}}(0)=\mathring{L}_{g_{\tilde u}}(0)=0.
				\]
				On the other hand, by our assumption that $\pa M$ is totally non-umbilic with respect to $g$, we have
				\[
				|\mathring{L}_{g}(q_{i,1}')|_{g}\geq \min_{\pa M} |\mathring L|_g>0.
				\]
				This together with \eqref{est: tilde_u_i} yields 
				\[
				\lim_{i\to\infty}|e^{-\tilde{u}_{i}(0)}\ve_{i}^{-1}\mathring{L}_{g}(q_{i,1}')|_{g}=+\infty.
				\]
				
				Going back to \eqref{eq:trace_2nd_ff}, we arrive at a contradiction.
				\vskip 8pt
				
				Again it follows from Lemmas \ref{lem:boundary bounded limit} and \ref{lem fudament-1-1} that $\{q_1,\cdots, q_N\}\subset \overline M$ are isolated blow-up points. However, it was shown in Proposition \ref{prop:no blow-up pts} that interior blow-up points do not exist. 
				
				\vskip 8pt

				Until now, we have established that $\inf_{M}u\geq-C$  for all $0\le t\le 1$.		
				
				\vskip 4pt
				(3) For $t \in [0,1]$, $\inf_M u\le C$ for some positive constant $C$ independent of $t$.
				\vskip 4pt
				If not, there exists a sequence of solutions $\{u_i=u(\cdot,t_i)\}$  to \eqref{eq:path} for some $t_i$ such that $\inf_M u_i\rightarrow +\infty$ and $t_i\rightarrow t_0\in [0,1].$ Without loss of generality, we assume $\zeta(t_1)>0.$
				
				Let ${u}^*_i=u_i-u_1$ satisfy 
				\begin{align*}
					\frac{\partial{u}^*_i}{\partial\vec{n}}&=\zeta(t_1)c(e^{-u_i}-e^{-u_1})+(\zeta(t_i)-\zeta(t_1))ce^{-u_i}\\
					&=\zeta(t_1)ce^{-u_1}(e^{-{u}^*_i}-1)+(\zeta(t_i)-\zeta(t_1))ce^{-u_i}\\
					&=c\big(\zeta(t_1)e^{-u_1}(e^{-{u}^*_i}-1)+(\zeta(t_i)-\zeta(t_1))e^{-u_i}\big)\qquad \mathrm{~~on~~} \partial M.
				\end{align*}

				Notice that $e^{-{u}^*_i}-1\rightarrow -1$ due to $\inf_M {u}^*_i\rightarrow +\infty$. Let us denote $u_i^*(x_i)=\inf_M {u}^*_i$ and $u_i(x_i)\rightarrow +\infty$.
				
				For large $i$ and $x_i \in \overline M$ we have
				\[ [\zeta(t_1)e^{-u_1}(e^{-{u}^*_i}-1)+(\zeta(t_i)-\zeta(t_1))e^{-u_i}](x_i)\le 0\]
				and 	
				\[0\le \frac{\partial u^*_i}{\partial\vec{n}}(x_i)=c(x_i)\big[\zeta(t_1)e^{-u_1}(e^{-{u}^*_i}-1)+(\zeta(t_i)-\zeta(t_1))e^{-u_i}\big](x_i)\le 0, \]
				whenever  $x_i\in \partial M$. In conclusion, no matter $x_i$ is either an interior or boundary point, we both have $\nabla u_i^*(x_i)=0$ and $\nabla ^2u_i^*(x_i)\ge 0$.

				Observe that $g_{u_i}=e^{-2u_i}g=e^{-2{u}^*_i}g_{u_1}$. Then for sufficiently large $i$, at $x_i$ there holds
				\begin{align*}
					A_{e^{-2u_i}g}+S_g(t_i)=&A_{e^{-2{u}^*_i}g_{u_1}}+S_{g}(t_i)\\
					=& \nabla^2 u_i^*+A_{g_{u_1}}+S_{g}(t_i)\\
					\ge& A_{g_{u_1}}+S_{g}(t_i)\\
					\ge& A_{g_{u_1}},
				\end{align*}
				where the last inequality follows from the fact that $S_g(t_i)\ge0$ by virtue of \eqref{positive matrix}.
				The above inequality yields that 
				\begin{align*}
					0<&\sigma_2^{1/2}(A_{g_{u_1}}(x_i)) \nonumber\\
					\le&(1-t_i)\big(\int_M e^{-(n+1)u_i}\ud\mu_g\big)^{\frac{2}{n+1}}+\zeta(t_i)f(x_i)e^{-2u_i(x_i)}\\
					\le &e^{-2\inf_M u_i^*}\big(\int_M e^{-(n+1)u_1}\ud \mu_g\big)^{\frac{2}{n+1}}+f(x_i)e^{-2u_i^*(x_i)}e^{-2u_1} \to 0 \qquad \mathrm{~~as~~} i \to \infty.
				\end{align*}
				This is a contradiction. 
				\vskip 8pt
				
				We apply the gradient estimate to know $\|\nabla u\|_{C^0(\overline M)}\le C$. Thus, for any $x \in \overline M$ and $t \in [0,1]$, there holds
				$$\inf_M u-C\leq u(x)\le \inf_M u+C.$$
				This implies the uniform $C^0$ estimate of $u$ for all $t \in [0,1]$.
			\end{proof}
			As indicated in the preceding discussion, Theorem \ref{Thm:existence} immediately follows from Lemma \ref{lem:totally_nonumbilic}.
		\end{proof}

		\appendix
		\section{Appendix: Possible umbilic points of higher dimensional ellipsoids}\label{Appendix:A}
		For $n \geq 3$, consider an ellipsoid
		$$\Sigma=\left\{x \in \Rn;~~ \sum_{i=1}^{n}\frac{x_{i}^{2}}{a_{i}^{2}}=1\right\}, \quad a_i \in \R_+ \mathrm{~~for~~} 1 \leq i \leq n.$$
		Let $F(x)=\sum_{i=1}^{n}a_i^{-2}x_{i}^{2}-1$, then
 $$\nabla F=\sum_{i=1}^n \frac{2x_i}{a_i^2} \frac{\pa}{\pa x_i}, \qquad \nabla^2 F=\sum_{i=1}^n \frac{2}{a_i^2} \ud x_i \otimes \ud x_i.$$
The inward unit normal at some $p \in \Sigma$ is
 $$\vec n=-\frac{\nabla F}{|\nabla F|}.$$
 It follows that
 $$|\nabla F|^2=4 \sum_{i=1}^n \frac{x_i^2}{a_i^4} \qquad \mathrm{and} \qquad \Delta F=2 \sum_{i=1}^n \frac{1}{a_i^2}.$$
		
		\begin{proof}[Proof of Proposition \ref{prop:ellipsoid}]

For an orthonormal basis $\{e_\alpha; 1 \leq \alpha \leq n-1\}$ of $T_p \Sigma$, the second fundamental form is
 $$L_{\alpha\beta}=-\langle \nabla_{e_\alpha} \vec n,e_\beta \rangle=\frac{\nabla^2 F(e_\alpha,e_\beta)}{|\nabla F|}$$
 and the mean curvature is
 $$H=\mathrm{tr}(L)=\frac{\Delta F- \nabla^2 F(\vec n,\vec n)}{|\nabla F|}=\frac{\Delta F}{|\nabla F|}-\frac{\nabla^2 F(\nabla F,\nabla F)}{|\nabla F|^3}.$$
 This follows the exact formula of $H$.
 
 Further assume $a_1<a_2<\cdots<a_n$. Then, $p=(x_1,\cdots, x_n) \in \Sigma$ is umbilic if and only if there exists $\lambda \in \R$   such that for all 
 $(u_1,\cdots,u_n) \in T_p \Sigma$ satisfying $\sum_{i=1}^n a_i^{-2}x_iu_i=0$,
 \begin{equation}\label{eqn:umbilic}
 \sum_{i=1}^n \frac{u_i^2}{a_i^2}=\lambda \sum_{i=1}^n u_i^2.
 \end{equation}
 
 Since $p \neq 0$, without loss of generality we assume $x_n \neq 0$. We know that
 $$u_n=-\frac{a_n^2}{x_n}(\sum_{\alpha=1}^{n-1}\frac{x_\alpha}{a_\alpha^2}u_\alpha).$$
Inserting it into $\eqref{eqn:umbilic}$ to show
 \begin{align*}
  \sum_{\alpha=1}^{n-1} \frac{u_\alpha^2}{a_\alpha^2}+\frac{a_n^2}{x_n^2}\left(\sum_{\alpha=1}^{n-1}\frac{x_\alpha}{a_\alpha^2}u_\alpha\right)^2=\lambda \left[\sum_{\alpha=1}^{n-1} u_\alpha^2+\frac{a_n^4}{x_n^2}\left(\sum_{\alpha=1}^{n-1}\frac{x_\alpha}{a_\alpha^2}u_\alpha\right)^2\right]
 \end{align*}
for all $(u_1,\cdots,u_{n-1}) \in \R^{n-1}$. Then, we deduce from the above equation that
\begin{subequations}
\begin{align}
\frac{1}{a_\alpha^2}+\frac{a_n^2}{x_n^2} \frac{x_\alpha^2}{a_\alpha^4}=\lambda\left(1+\frac{a_n^4}{x_n^2} \frac{x_\alpha^2}{a_\alpha^4}\right), &\qquad 1 \leq \alpha \leq n-1;\label{eq1:umbilic_pts}\\
(1-\lambda a_n^2)x_\alpha x_\beta=0, &\qquad 1 \leq \alpha \neq \beta \leq n-1.\label{eq2:umbilic_pts}
\end{align}
\end{subequations}

If $(x_1,\cdots,x_{n-1})$ admits at least two non-zero real numbers, denoted by $x_\alpha, x_\beta \neq 0$ for some $\alpha \neq \beta$, then we derive from \eqref{eq2:umbilic_pts} that  $\lambda=\frac{1}{a_n^2}$ and then from \eqref{eq1:umbilic_pts}  that $a_\alpha^2=a_n^2$,  which violates our assumption. Thus, we assert that $(x_1,\cdots,x_{n-1})$ admits at most one $x_\alpha$ such that $x_\alpha\neq 0$.

\textit{(i)} ~For $n=3$, we can directly solve \eqref{eq1:umbilic_pts} and \eqref{eq2:umbilic_pts} with the desired solutions.

\textit{(ii)} ~For $n\geq 4$,  assume $x_\alpha=0$ for some $1 \leq \alpha \leq n-1$, it follows from \eqref{eq1:umbilic_pts} that $\lambda=\frac{1}{a_\alpha^2}$. Since $a_i$ are distinct,  we conclude again from \eqref{eq1:umbilic_pts} that $(x_1,\cdots,x_{n-1})$ admits at most one $x_\beta$ such that $x_\beta=0$.  Hence we reach a contradiction due to $n \geq 4$. So, the set of umbilic points on $\Sigma$ is empty.
\end{proof}
		
		\begin{proof}[Proof of Proposition \ref{prop:umbilic_pts_ellipsoid}]	
For an orthonormal basis $\{e_\alpha:=\sum_{i=1}^n u^i_\alpha \frac{\pa}{\pa x_i}; 1 \leq \alpha \leq n-1\}$ of $T_p \Sigma$, the second fundamental form is
\begin{align*}
 L_{\alpha\beta}=&-\langle \nabla_{e_\alpha} \vec n,e_\beta \rangle=\frac{\nabla^2 F(e_\alpha,e_\beta)}{|\nabla F|}\\
 =&\frac{1}{|\nabla F|}\sum_{i=1}^n \frac{2}{a_i^2}(\ud x_i \otimes \ud x_i)(e_\alpha,e_\beta)\\
 =&\frac{1}{|\nabla F|}\sum_{i=1}^n \frac{2}{a_i^2} u_\alpha^i u_\beta^i.
 \end{align*}
 
 We continue with computation of other components of $\nabla^2 F$:
  \begin{align*}
 \nabla^2 F(\vec n, \vec n)=\frac{\nabla^2 F(\nabla F,\nabla F)}{|\nabla F|^2}=\frac{8}{|\nabla F|^2}\sum_{i=1}^n \frac{x_i^2}{a_i^6}
 \end{align*}
 and
 \begin{align*}
 \nabla^2 F(e_\alpha, \vec n)=&\frac{1}{|\nabla F|}\sum_{i=1}^n \frac{2}{a_i^2}(\ud x_i \otimes \ud x_i)(e_\alpha,\nabla F)\\
 =&\frac{4}{|\nabla F|}\sum_{i=1}^n \frac{x_i}{a_i^4}u_\alpha^i\\
 :=&\frac{4}{|\nabla F|}\langle T, e_\alpha \rangle,
 \end{align*}
  where
 $$T=\sum_{i=1}^n \frac{x_i}{a_i^4} \frac{\pa}{\pa x_i}.$$
 These enable us to calculate
 \begin{align*}
 \sum_{\alpha}|\nabla^2 F(e_\alpha,\vec n)|^2=&\frac{16}{|\nabla F|^2} \sum_{\alpha}|\langle T,e_\alpha \rangle|^2\\
 =&\frac{16}{|\nabla F|^2} (|T|^2-|\langle T,\vec n \rangle|^2)\\
 =&\frac{16}{|\nabla F|^2}\left[\sum_{i=1}^n \frac{x_i^2}{a_i^8}-\frac{4}{|\nabla F|^2}\left(\sum_{i=1}^n \frac{x_i^2}{a_i^6}\right)^2\right]
 \end{align*}
 and
 \begin{align*}
 &\sum_{\alpha,\beta}|\nabla^2 F(e_\alpha,e_\beta)|^2\\
 =&|\nabla^2 F|^2-2\sum_{\alpha}|\nabla^2 F(e_\alpha,\vec n)|^2-|\nabla^2 F(\vec n,\vec n)|^2\\
 =&\sum_{i=1}\frac{4}{a_i^4}-\frac{32}{|\nabla F|^2}\left[\sum_{i=1}^n \frac{x_i^2}{a_i^8}-\frac{4}{|\nabla F|^2}\left(\sum_{i=1}^n \frac{x_i^2}{a_i^6}\right)^2\right]-\frac{64}{|\nabla F|^4}\left(\sum_{i=1}^n \frac{x_i^2}{a_i^6}\right)^2.
 \end{align*}

 Suppose $p$ is an umbilic point on $\Sigma$, then at $p$ we have
 \begin{equation}\label{umbilic_pt}
 0=|\mathring L|^2=|L-\frac{H}{n-1} I |^2=|L|^2-\frac{H^2}{n-1}\quad\Longleftrightarrow \quad |L|^2=\frac{H^2}{n-1}.
 \end{equation}
Therefore, inserting all above terms into \eqref{umbilic_pt}, together with the explicit formula of $H$ in Proposition \ref{prop:ellipsoid} we obtain the desired assertion.
		\end{proof}

		\section{Appendix: Proof of  Lemmas \ref{one critical point} and \ref{lem:four-dim lower bound-boundary}}\label{Appendix:B}
		
		We are now ready to give the proof of  Lemmas \ref{one critical point} and \ref{lem:four-dim lower bound-boundary}.
		
		\begin{proof}[Proof of Lemma \ref{one critical point}]
		We just provide the proof for $x_0\in \partial M$ as the  idea is similar to that for interior isolated blow-up point.

		Given $R_i \to \infty$ and $\varepsilon_{i} \to 0$, by Lemma \ref{lem fudament} we have $R_iu_i(x_i)^{-2/(n-2)} \to 0$ and
		\begin{equation}\label{key estimate}
			\left\Vert \frac{1}{u_{i}(x_{i})}u_{i}\circ G_{x_{i}'}\big(\frac{y}{u_{i}^{\frac{2}{n-2}}\left(x_{i}\right)}+d_i\mathbf{e}_n\big)-\big(\frac{2b\sqrt{\binom{n}{2}}}{f(x_0)}\big)^{\frac{n-2}{4}}(1+b|y+\tilde T_c\mathbf{e}_{n}|^{2})^{\frac{2-n}{2}}\right\Vert _{C^{2}(B^+_{3R_{i}})}<\varepsilon_{i}.
		\end{equation}
		
		For simplicity, we let 
		$$\eta_i(y)= \frac{1}{u_{i}(x_{i})}u_{i}\circ G_{x_{i}'}\left(u_{i}^{-\frac{2}{n-2}}\left(x_{i}\right)y+d_i\mathbf{e}_n\right)$$
		and 
		$$\eta(y)=\left(\frac{2b\sqrt{\binom{n}{2}}}{f(x_0)}\right)^{\frac{n-2}{4}}(1+b|y+\tilde T_c\mathbf{e}_{n}|^{2})^{\frac{2-n}{2}}.$$
		
		We first describe elementary properties involving standard bubble function. 
		
		Notice that
		\begin{align*}
			\bar{\eta}(s) & =\frac{1}{s^{n-1}}\int_{\partial^+ B_{s,\theta}^{+}}\eta(y)\ud \sigma=\int_{\partial^+ B_{1,\theta }^{+}}\eta(sz)\ud\sigma\\
			& =\bigg(\frac{2b\sqrt{\binom{n}{2}}}{f\left(x_{0}\right)}\bigg)^{\frac{n-2}{4}}\int_{\partial^+ B_{1,\theta }^{+}}\left(1+b\left(\left|sz^{\prime}\right|^{2}+\left|sz_{n}+\tilde T_c\right|^{2}\right)\right)^{\frac{2-n}{2}}\ud\sigma.
		\end{align*}
		
		A direct computation yields
		\begin{align*}
			&\bigg(\frac{2b\sqrt{\binom{n}{2}}}{f\left(x_{0}\right)}\bigg)^{\frac{2-n}{4}} \big(s^{\frac{n-2}{2}}\bar{\eta}(s)\big)'\\
			= & \frac{n-2}{2}s^{\frac{n-2}{2}-1}\int_{\partial^+ B_{1,\theta }^{+}}\bigg(1+b\left(\left|sz^{\prime}\right|^{2}+\left|sz_{n}+\tilde T_c\right|^{2}\right)\bigg)^{-\frac{n}{2}}\times\\
			& \qquad\qquad\qquad\bigg(1+b(\left|sz^{\prime}\right|^{2}+\left|sz_{n}+\tilde T_c\right|^{2})-2bs\big(s|z'|^{2}+(sz_{n}+\tilde T_c)z_{n}\big)\bigg)\ud\sigma\\
			= & \frac{n-2}{2}s^{\frac{n-2}{2}-1}\big(1-bs^{2}+b\tilde T_c^{2}\big)\int_{\partial^+ B_{1,\theta }^{+}}\big(1+b(s^{2}+2sz_{n}\tilde T_c+\tilde T_c^{2})\big)^{-\frac{n}{2}}\ud\sigma
		\end{align*}
		and
		\begin{align}
			&\bigg(\frac{2b\sqrt{\binom{n}{2}}}{f\left(x_{0}\right)}\bigg)^{\frac{2-n}{4}} \big(s^{\frac{n-2}{2}}\bar{\eta}(s)\big)''\nonumber\\
			= & \frac{n-2}{2}s^{\frac{n-2}{2}-2}\int_{\partial^+ B_{1,\theta }^{+}}\bigg(1+b\left(|sz^{\prime}|^{2}+|sz_{n}+\tilde T_c|^{2}\right)\bigg)^{-\frac{n}{2}-1}\cdot A(s,z) \ud \sigma. \label{second derivatives of standard function}
		\end{align}
		where 
		\begin{align*}
			A(s,z):=&\bigg(1+b(s^{2}+2sz_{n}\tilde T_c+\tilde T_c^{2})\bigg)\bigg((1-bs^{2}+b\tilde T_c^{2})(\frac{n}{2}-2)-2bs^{2}\bigg)\\
			&-nbs\big(s+\tilde T_cz_{n}\big)\big(1-bs^{2}+b\tilde T_c^{2}\big).
		\end{align*}
		
		Keep in mind that the sign of $1-bs^{2}+b \tilde T_c^{2}$ coincides with $\big(s^{\frac{n-2}{2}}\bar{\eta}(s)\big)'$. 
		
		Let $s_{0}:=\sqrt{\frac{1+b\tilde T_c^{2}}{b}}$, then $1-bs^{2}+b\tilde T_c^{2}\ge0$
		for $0< s\le s_{0}$. 
		
		Fix some $\alpha \in (0,1)$ to be determined later and let $\varepsilon_{0}=\alpha s_0$, then
		\begin{align*}
			1-bs^{2}+b\tilde T_c^{2} \ge bs_0^2(2\alpha-\alpha^2)>0&\quad\mathrm{for~~}\quad 0<s\le s_{0}(1-\alpha);\\
			1-bs^{2}+b\tilde T_c^{2} \leq -bs_0^2(2\alpha+\alpha^2)<0&\quad\mathrm{for~~}\quad s\ge s_{0}(1+\alpha);\\
			-bs_0^2(2\alpha+\alpha^2)< 1-bs^{2}+b\tilde T_c^{2}< bs_0^2(2\alpha-\alpha^2) &\quad\mathrm{for~~}\quad s_0(1-\alpha)<s<s_0(1+\alpha).
		\end{align*}
		
		For $n=3, 4$ and $s_0(1-\alpha)<s<s_0(1+\alpha)$, there holds 
		\[(1-bs^{2}+b\tilde T_c^{2})(\frac{n}{2}-2)-2bs^{2}\le -(\frac{n}{2}-2)bs_0^2(2\alpha+\alpha^2)-2b(1-\alpha)^2s_0^2; \]
		for $n\ge 5$ and $s_0(1-\alpha)<s<s_0(1+\alpha)$, there holds
		\[(1-bs^{2}+b\tilde T_c^{2})(\frac{n}{2}-2)-2bs^{2}\le (\frac{n}{2}-2)bs_0^2(2\alpha-\alpha^2)-2b(1-\alpha)^2s_0^2.\]
		Meanwhile, for $n\ge 3$ and $s_0(1-\alpha)<s<s_0(1+\alpha)$, 
		\[-nbs(s+\tilde T_cz_{n})(1-bs^{2}+b\tilde T_c^{2})\le n bs_0(1+\alpha)(s_0(1+\alpha)+\tilde T_c)bs_0^2(2\alpha+\alpha^2).\]
		Thus we may choose $\alpha$ small enough such that for $s_0(1-\alpha)<s<s_0(1+\alpha)$,
		\begin{align*}
			A\le & \big(1+b\tilde T_c^{2}\big)\big((2-\frac{n}{2})(2\alpha+\alpha^2)-2(1-\alpha)^2\big)bs_0^2\\
			&+nbs_0(1+\alpha)\left(s_0(1+\alpha)+\tilde T_c\right)bs_0^2(2\alpha+\alpha^2):=-c_2<0\quad 
			\mathrm{for}\quad  n=3,4;
		\end{align*}
		and 
		\begin{align*}
			A \le &\big(1+b\tilde T_c^{2}\big)\big((\frac{n}{2}-2)(2\alpha-\alpha^2)-2(1-\alpha)^2\big)bs_0^2\\
			&+nbs_0(1+\alpha)\left(s_0(1+\alpha)+\tilde T_c\right)bs_0^2(2\alpha+\alpha^2):=-c_3<0\quad \mathrm{for} \quad n\ge 5.
		\end{align*}
		Therefore, by \eqref{second derivatives of standard function}, there exits a positive constant $c_4:=\min\{c_2,c_3\}$ such that for $n\geq 3$,
		\begin{equation}\label{eq:concavity in neck}
			\big(s^{\frac{n-2}{2}}\bar{\eta}(s)\big)''\le -c_4s^{\frac{n}{2}-3}\quad \mathrm{for}\quad  s_0-\varepsilon_0<s<s_0+\varepsilon_{0}.
		\end{equation}
		
		For $s\ge (1+\alpha)s_{0}$, we have 
		\[
		s^{2}+2sz_{n}\tilde T_c+\tilde T_c^{2}\le4s^{2},\quad1-bs^{2}+b\tilde T_c^{2}=b(s_0^2-s^2)<-bs^{2}\frac{2\alpha+\alpha^2}{(1+\alpha)^2}
		\]
		and 
		\begin{align}\label{first lower derivative_2}
			\big(s^{\frac{n-2}{2}}\bar{\eta}(s)\big)'\no
			\le&  -\frac{Cs^{2}}{(1+4bs^{2})^{\frac{n}{2}}}s^{\frac{n-2}{2}-1}\no\\
			\le & -C(n,\alpha,\theta)s^{-\frac{n}{2}}\qquad\mathrm{for}\quad (1+\alpha)s_{0}\le s\le3R_{i}.
		\end{align}

		For $0<s\le s_{0}(1-\alpha)$, there exist  positive constants
		$c_{0}, c_1>0$ such that 
		\begin{equation}\label{first lower derivative_1}
			\big(s^{\frac{n-2}{2}}\bar{\eta}(s)\big)'\ge\bigg(\frac{2b\sqrt{\binom{n}{2}}}{f\left(x_{0}\right)}\bigg)^{\frac{n-2}{4}}\frac{n-2}{2}s^{\frac{n-2}{2}-1}c_{0}\ge c_1 s^{\frac{n-2}{2}-1}.
		\end{equation}
		
		Observe that 
		\begin{align*}
			\overline{\eta}_{i}(s):=\frac{1}{s^{n-1}}\int_{\partial^{+}B_{s,\theta}^{+}}\eta_{i}(y)\ud\sigma_{\widetilde G_{x_i'}^{\ast}(g)}=\int_{\partial^{+}B_{1,\theta}^{+}}\eta_{i}(sz)\ud\sigma_{g_{i}},
		\end{align*}
		where $\widetilde G_{x_i'}(y)=G_{x_i'}(u_{i}^{-2/(n-2)}(x_{i})y+d_i\mathbf{e}_n)$ and $g_{i}=s^{-2}u_i^{4/(n-2)}(x_i)\Psi_{i}^{\ast}(g) \to g_{\mathrm{E}}$ in $C_{\mathrm{loc}}^3(\overline{\Rn_+})$ with 
		$$\Psi_{i}(z):=G_{x_i'}(u_{i}^{-\frac{2}{n-2}}(x_{i})sz+d_i\mathbf{e}_n).$$
		For large $i$, 
		$$\ud\sigma_{g_i}=\ud\sigma_{\mathbb{S}^{n-1}}(1+O(s|z|u_i(x_i)^{-\frac{2}{n-2}})).$$

		For $ (1+\alpha)s_0\le s<3R_{i}$ and sufficiently large $i$, by \eqref{key estimate} and \eqref{first lower derivative_2} we have 
		\begin{align}
			& (s^{\frac{n-2}{2}}\overline{\eta}_{i}(s))'\no\label{eq:first derivative negative upper bound}\\
			= & \frac{n-2}{2}s^{\frac{n-2}{2}-1}\overline{\eta}_{i}(s)+s^{\frac{n-2}{2}}\overline{\eta}_{i}'(s)\no\\
			= & \frac{n-2}{2}s^{\frac{n-2}{2}-1}\overline{\eta}_{i}(s)+s^{\frac{n-2}{2}}\int_{\partial^{+}B_{1,\theta}^{+}}\langle(\nabla\eta_{i})(sz),z\rangle\ud\sigma_{g_{i}}\nonumber \\
			\le & \frac{n-2}{2}s^{\frac{n-2}{2}-1}\int_{\partial^{+}B_{1,\theta}^{+}}\varepsilon_{i}+\eta(sz)\ud\sigma_{g_{i}}+\frac{n-2}{2}s^{\frac{n-2}{2}-1}\overline{\eta}\no\\
			&+s^{\frac{n-2}{2}}\int_{\partial^{+}B_{1,\theta}^{+}}(\varepsilon_{i}+\langle(\nabla\eta)(sz),z\rangle)\ud\sigma_{g_{i}}\nonumber \\
			\le & C_{1}(\frac{n-2}{2}s^{\frac{n-2}{2}-1}\varepsilon_{i}+s^{\frac{n-2}{2}}\varepsilon_{i})+\big(s^{\frac{n-2}{2}}\bar{\eta}(s)\big)'+s^{-\frac{n}{2}}O(su_i(x_i)^{-\frac{2}{n-2}})\nonumber \\
			\le & (s^{\frac{n-2}{2}}\overline{\eta}(s))'+C_{1}\varepsilon_{i}s^{n-1}s^{-\frac{n}{2}}+s^{-\frac{n}{2}}O(R_iu_i(x_i)^{-\frac{2}{n-2}})\nonumber \\
			\le & -\frac{C(n,\alpha,\theta)}{2}s^{-\frac{n}{2}}.
		\end{align}
		Similarly, for $0<s\le s_0(1-\alpha)$ and sufficiently large $i$, by \eqref{first lower derivative_1} we have 
		\begin{align}
			& (s^{\frac{n-2}{2}}\overline{\eta}_{i}(s))'\no\label{eq:first derivative_2}\\
			\geq & -C_{2}(\frac{n-2}{2}s^{\frac{n-2}{2}-1}\varepsilon_{i}+s^{\frac{n-2}{2}}\varepsilon_{i})+(s^{\frac{n-2}{2}}\overline{\eta}(s))'-C_3s^{\frac{n-2}{2}}O(u_i(x_i)^{-\frac{2}{n-2}})\no\\
			\ge & c_1 s^{\frac{n-2}{2}-1}-C_{2}(\frac{n-2}{2}s^{\frac{n-2}{2}-1}\varepsilon_{i}+s^{\frac{n-2}{2}}\varepsilon_{i})-C_3s^{\frac{n-2}{2}}O(u_i(x_i)^{-\frac{2}{n-2}})\nonumber \\
			\geq & s^{\frac{n}{2}-2}\left(c_1-C_{2}(\frac{n-2}{2}\varepsilon_{i}+s\varepsilon_{i})-C_3sO(u_i(x_i)^{-\frac{2}{n-2}})\right)>0.
		\end{align}
		
		For $(1-\alpha)s_0\le s\le (1+\alpha)s_0$ and sufficiently large $i$, by \eqref{eq:concavity in neck} we have
		\begin{align}\label{eq:second derivative of standard function}
			& (s^{\frac{n-2}{2}}\overline{\eta}_{i}(s))''\no\\
			= & \frac{n-2}{2}(\frac{n-2}{2}-1)s^{\frac{n-2}{2}-2}\overline{\eta}_{i}(s)+(n-2)s^{\frac{n-2}{2}-1}\int_{\partial^{+}B_{1,\theta}^{+}}\langle(\nabla\eta_{i})(sz),z\rangle\ud\sigma_{g_{i}}\no\\
			& +s^{\frac{n-2}{2}}\int_{\partial^{+}B_{1,\theta}^{+}}\langle z(\nabla^{2}\eta_{i})(sz),z\rangle\ud\sigma_{g_{i}}\nonumber \\
			\le & \frac{n-2}{2}(\frac{n-2}{2}-1)s^{\frac{n-2}{2}-2}\int_{\partial^{+}B_{1,\theta}^{+}}(\eta(sz)+\varepsilon_{i})\ud\sigma_{g_{
					\mathbb{S}^{n-1}}}(1+O(su_i(x_i)^{-\frac{2}{n-2}}))\nonumber \\
			& +(n-2)s^{\frac{n-2}{2}-1}\int_{\partial^{+}B_{1,\theta}^{+}}(\varepsilon_{i}+\langle(\nabla\eta)(sz),z\rangle)\ud\sigma_{g_{
					\mathbb{S}^{n-1}}}(1+O(su_i(x_i)^{-\frac{2}{n-2}}))\nonumber \\
			& +s^{\frac{n-2}{2}}\int_{\partial^{+}B_{1,\theta}^{+}}(\langle z(\nabla^{2}\eta)(sz),z\rangle+\varepsilon_{i})\ud\sigma_{g_{
					\mathbb{S}^{n-1}}}(1+O(su_i(x_i)^{-\frac{2}{n-2}}))\nonumber \\
			\le &  (s^{\frac{n-2}{2}}\overline{\eta}(s))''+C_{1}\varepsilon_{i}(s^{\frac{n-2}{2}-2}+s^{\frac{n-2}{2}-1}+s^{\frac{n-2}{2}})\nonumber\\
			&+\frac{n-2}{2}(\frac{n-2}{2}-1)s^{\frac{n-2}{2}-2}\int_{\partial^{+}B_{1,\theta}^{+}}\eta(sz)\ud\sigma_{g_{
					\mathbb{S}^{n-1}}}O(su_i(x_i)^{-\frac{2}{n-2}})\nonumber \\
			&+s^{\frac{n-2}{2}}\int_{\partial^{+}B_{1,\theta}^{+}}(\langle z(\nabla^{2}\eta)(sz),z\rangle)\ud\sigma_{g_{
					\mathbb{S}^{n-1}}}(O(su_i(x_i)^{-\frac{2}{n-2}}))\nonumber\\
			&+(n-2)s^{\frac{n-2}{2}-1}\int_{\partial^{+}B_{1,\theta}^{+}}\langle(\nabla\eta)(sz),z\rangle)\ud\sigma_{g_{
					\mathbb{S}^{n-1}}}O(su_i(x_i)^{-\frac{2}{n-2}})\nonumber\\
			\leq & s^{\frac{n}{2}-3}\left(-c_4+C_{1}\varepsilon_{i}(1+s+s^{2})+u_i(x_i)^{-\frac{2}{n-2}}(s+s^2+s^3)\right)<0.
		\end{align}
		
		%
		%
		%
		%
		
		Therefore, we combine  \eqref{eq:first derivative negative upper bound}, \eqref{eq:first derivative_2} and \eqref{eq:second derivative of standard function} to conclude
		that for $n\geq3$ and sufficiently large $i$, there exists a unique point $s_{1}\in \big((1-\alpha)s_0, (1+\alpha)s_0\big)$
		such that $(s^{(n-2)/2}\overline{\eta}_{i}(s))'|_{s=s_{1}}=0$. Moreover, $(s^{(n-2)/2}\overline{\eta}_{i}(s))'>0$  for $s\in (0, s_1)$;  $(s^{(n-2)/2}\overline{\eta}_{i}(s))'<0$  for $s\in (s_1, R_i)$, this implies that $s_{1}$ is the unique critical point of $s^{(n-2)/2}\overline{\eta}_{i}(s)$
		in $(0,R_{i})$.

		Finally, with $r=u_{i}^{-2/(n-2)}(x_{i})s$ it follows from  the property
		\eqref{property_bdry_blowup_pt} that 
		\[
		(s^{\frac{n-2}{2}}\overline{\eta}_{i}(s))'=\frac{\ud}{\ud r}[r^{\frac{n-2}{2}}\overline{u}_{i}(r)]\frac{\ud r}{\ud s}=u_{i}^{-\frac{2}{n-2}}\left(x_{i}\right)\frac{\ud}{\ud r}[r^{\frac{n-2}{2}}\overline{u}_{i}(r)].
		\]
		This implies that $r^{(n-2)/2}\overline{u}_{i}(r)$ has only one critical
		point in $(0,r_{i})$.	
	\end{proof}
	
	\medskip
	
	\begin{proof}[Proof of Lemma \ref{lem:four-dim lower bound-boundary}]
		For brevity, we use $v$ and $A_\delta$ instead of $\underline{v}_{\delta}(r)$ and $A_{\underline{v}_{\delta}^{4/(n-2)}g}$, respectively. We introduce a matrix by
		\begin{align*}
			(\widetilde{L}_{ij})=\begin{pmatrix}
				L_{\alpha\beta} & 0\\
				0 & 0
			\end{pmatrix},
		\end{align*}
		whence,
		\[
		g_{ij}=\delta_{ij}-2\widetilde{L}_{ij}x_{n}+O(r^{2}),
		\]
		\[ g^{ij}=\delta_{ij}+2\widetilde{L}_{ij}x_{n}+O(r^{2}).\]
		Write $g=(\sqrt{g})^2$, then the inverse of $\sqrt{g}$ has
		\[
		\sqrt{g}^{ij}=\delta_{ij}+\tilde{L}_{ij}x_{n}+O(r^{2}).
		\]
		
		Suppose $v(x)=v(r)$ for $r=r(x)=|x|$.  A direct computation shows 
		\begin{align*}
			v_{,\alpha\beta}= & \pa_{\beta}\pa_{\alpha}v-L_{\alpha\beta}v_{,n}+O(r)|v'|,\qquad v_{,nn}=  \pa_{n}^2 v,\\
			v_{,\alpha n}= & \pa_{n}\pa_{\alpha}v+L_{\alpha\beta}v_{,\beta}+O(r)|v'|.
		\end{align*}
		
		Recall that
		\begin{align*}
			A_\delta=&-\frac{2}{n-2}\frac{\nabla^2 v}{v}+\frac{2n}{(n-2)^2} \frac{\ud v \otimes \ud v}{v^2}-\frac{2}{(n-2)^2}\frac{|\nabla v|_g^2}{v^2}g+A_g.
		\end{align*}

		Using
		\[
		|\nabla v|_{g}^{2}=g^{ij}v_{i}v_{j}=(\delta^{ij}+2\widetilde{L}_{ij}x_{n}+O(r^{2}))v_{i}v_{j},
		\]
		we have
		\begin{align*}
			\hat {D}_{ij}:= & -\frac{2}{n-2}\frac{\partial_{ij}^2 v}{v}+\frac{2n}{(n-2)^{2}}\frac{v_{i}v_{j}}{v^{2}}-\frac{2}{(n-2)^{2}}\frac{|\nabla v|_{g}^{2}}{v^{2}}g_{ij}+A_{ij}\\
			= & -\frac{2}{n-2}\frac{1}{v}\bigg(\frac{v''}{r^{2}}x_{i}x_{j}+\frac{v'}{r}\delta_{ij}-\frac{v'}{r^{3}}x_{i}x_{j}\bigg)+\frac{2n}{(n-2)^{2}}(\frac{v'}{v})^{2}\frac{x_{i}x_{j}}{r^{2}}-\frac{2}{(n-2)^{2}}(\frac{v'}{v})^{2}g_{ij}\\
			& +O(r^{2})(\frac{v'}{v})^{2}g_{ij}-\frac{4}{(n-2)^{2}}(\frac{v'}{v})^{2}L_{\alpha\beta}x_{\alpha}x_{\beta}x_{n}g_{ij}+A_{ij}\\
			=&-\frac{x_ix_j}{r^2}\left[\frac{2}{n-2}\frac{v''}{v}-\frac{2}{n-2}\frac{v'}{rv}-\frac{2n}{(n-2)^{2}}(\frac{v'}{v})^{2}\right]+\left[-\frac{2}{(n-2)^{2}}(\frac{v'}{v})^{2}-\frac{2}{n-2}\frac{v'}{rv}\right]g_{ij}\\
			&-\frac{4}{(n-2)^{2}}(\frac{v'}{v})^{2}L_{\alpha\beta}x_{\alpha}x_{\beta}x_{n}g_{ij}+\frac{2}{n-2}\frac{v'}{rv}(g_{ij}-\delta_{ij})+O(r^{2})(\frac{v'}{v})^{2}g_{ij}+A_{ij}.
		\end{align*}
		
		For clarity, we reorganize the above quantity as
		\begin{align*}
			D_{ij}=&\hat D_{ij}-\frac{2}{n-2}\frac{v'}{rv}(g_{ij}-\delta_{ij})\\
			=&-\chi_2\frac{x_i x_j}{r^2}+\chi_1 g_{ij}-\frac{4}{(n-2)^{2}}(\frac{v'}{v})^{2}L_{\alpha\beta}x_{\alpha}x_{\beta}x_{n}g_{ij}+O(r^{2})(\frac{v'}{v})^{2}g_{ij}+A_{ij}
		\end{align*}
		where
		\[
		\chi_{1}=-\frac{2}{(n-2)^{2}}(\frac{v'}{v})^{2}-\frac{2}{n-2}\frac{v'}{rv},\quad\chi_{2}=\frac{2}{n-2}\frac{v''}{v}-\frac{2}{n-2}\frac{v'}{rv}-\frac{2n}{(n-2)^{2}}(\frac{v'}{v})^{2}
		\]
		and 
		\begin{align*}
			B=&\begin{pmatrix}
				\frac{2}{n-2}\frac{v_{n}}{v}L_{\alpha\beta}+\frac{2}{n-2}\frac{v'}{rv}(g_{\alpha \beta}-\delta_{\alpha \beta}) & -\frac{2}{n-2}\frac{1}{v}L_{\alpha\beta}v_{\beta}\\
				-\frac{2}{n-2}\frac{1}{v}(L_{\alpha\beta}v_{\beta})^{\top} & 0
			\end{pmatrix}\\
			=&-\frac{2}{n-2}\frac{v'}{rv}\begin{pmatrix}
				L_{\alpha\beta}x_n+O(r^2) &L_{\alpha\beta}x_{\beta}\\
				(L_{\alpha\beta}x_{\beta})^{\top} & 0
			\end{pmatrix}=O(\big|\frac{v'}{v}\big|).
		\end{align*}
		
		This benefits us to decompose
		\begin{align*}
			\sigma_2(g^{-1}A_\delta)=& \sigma_{2}(g^{-1}({D}+B))=  \sigma_{2}\big((\sqrt{g})^{-1}({D}+B)(\sqrt{g})^{-1}\big).
		\end{align*}
		
		We choose an auxiliary function $v$ as
		$$v(r)=r^{-(n-2-\delta)}f(r),$$
		where $f(r)$ is a smooth function to be determined later.
		
		A direct computation shows 
		\begin{align*}
			\frac{v'}{v}=&-\frac{n-2-\delta}{r}+\frac{f'}{f},\\
			\frac{v''}{v}=&\frac{(n-2-\delta)(n-1-\delta)}{r^{2}}-\frac{2(n-2-\delta)}{r}\frac{f'}{f}+\frac{f''}{f}
		\end{align*}
		and
		\begin{align*}
			\chi_{1}=&\frac{2\delta(n-2-\delta)}{(n-2)^{2}r^{2}}-\frac{2}{(n-2)^{2}}(\frac{f'}{f})^{2}+\frac{2(n-2-2\delta)}{(n-2)^{2}}\frac{f'}{rf},\\
			\chi_{2}=&\frac{4\delta(n-2-\delta)}{(n-2)^{2}r^{2}}+\frac{6(n-2)-8\delta}{(n-2)^{2}r}\frac{f'}{f}-\frac{2n}{(n-2)^{2}}(\frac{f'}{f})^{2}+\frac{2}{n-2}\frac{f''}{f}.
		\end{align*}
		
		Notice that
		\begin{align*}
			& \frac{\sqrt{g}^{kp}x_{p}x_{q}\sqrt{g}^{qj}}{r^{2}}\\
			= & (\delta_{kp}+\tilde{L}_{kp}x_{n}+O(r^{2}))\frac{x_{p}x_{q}}{r^{2}}(\delta_{qj}+\tilde{L}_{qj}x_{n}+O(r^{2}))\\
			= & \frac{x_{k}x_{j}}{r^{2}}+\frac{x_{k}x_{q}x_{n}}{r^{2}}\widetilde{L}_{qj}+\widetilde{L}_{kp}\frac{x_{n}x_{p}x_{j}}{r^{2}}+O(r^{2}).
		\end{align*}		
		Then we arrange
		\begin{align*}
			& [(\sqrt{g})^{-1}({D}+B)(\sqrt{g})^{-1}]_{kj}\\
			= & -\chi_{2}\frac{x_{k}x_{j}}{r^{2}}+\chi_{1}\delta_{kj}\\
			&-\chi_{2}\bigg(\frac{x_{k}x_{q}}{r^{2}}\widetilde{L}_{qj}x_{n}+\tilde{L}_{kp}x_{n}\frac{x_{p}x_{j}}{r^{2}}\bigg)\\
			& +|\chi_{2}|O(|x|^2)-\frac{4}{(n-2)^{2}}(\frac{v'}{v})^{2}L_{\alpha\beta}x_{\alpha}x_{\beta}x_{n}\delta_{kj}+O(r^{2})(\frac{v'}{v})^{2}\delta_{kj}+\sqrt{g}^{kp}A_{pq}\sqrt{g}^{qj}\\
			&+\sqrt{g}^{kp} B_{pq}\sqrt{g}^{qj}\\
			= &(D_1)_{kj}+(D_2)_{kj}+(D_3)_{kj}+(D_4)_{kj},
		\end{align*}
		where 
		\begin{align*}
			(D_1)_{kj}=&-\chi_{2}\frac{x_{k}x_{j}}{r^{2}}+\chi_{1}\delta_{kj},\\
			(D_2)_{kj}=&-\chi_{2}\bigg(\frac{x_{k}x_{q}}{r^{2}}\widetilde{L}_{qj}x_{n}+\tilde{L}_{kp}x_{n}\frac{x_{p}x_{j}}{r^{2}}\bigg),\\
			(D_3)_{kj} =&|\chi_{2}|O(r^2)-\frac{4}{(n-2)^{2}}(\frac{v'}{v})^{2}L_{\alpha\beta}x_{\alpha}x_{\beta}x_{n}\delta_{kj}+O(r^{2})(\frac{v'}{v})^{2}\delta_{kj}+\sqrt{g}^{kp}A_{pq}\sqrt{g}^{qj},\\
			(D_4)_{kj}=&\sqrt{g}^{kp}B_{pq}\sqrt{g}^{qj}.
		\end{align*}	
		
		For $n \times n$ matrices $A,B$, we define a  mixed symmetric function by 
		$$\sigma_{k,l}(A,B)=\frac{1}{k!}\sum
		\delta\left(\begin{matrix}
			i_1&\cdots &i_k      \\
			j_1 &\cdots & j_k
		\end{matrix}
		\right)
		A_{i_1}^{j_1}\cdots A_{i_l}^{j_l}B_{i_{l+1}}^{j_{l+1}}\cdots B_{i_k}^{j_k},
$$
where $\left(\begin{matrix}
			i_1&\cdots &i_k      \\
			j_1 &\cdots & j_k
		\end{matrix}
		\right)$ is a permutation over the index set $\{1,2,\cdots,n\}$. See Reilly \cite{Reilly} or \cite[Definition 1]{Chen}.
Hence, we obtain
		\begin{align}\label{eq:summation}
			&\sigma_2(g^{-1}A_\delta)\no\\
			= & \sigma_{2}(D_{1}+D_{2}+D_{3}+D_{4})\nonumber \\
			= & \sigma_{2}(D_{1})+\sigma_{2}(D_{2}+D_{4})+\sigma_{2}(D_{3})\nonumber \\
			& +2\sigma_{2,1}(D_{1},D_{2}+D_{4})+2\sigma_{2,1}(D_{1},D_{3})+2\sigma_{2,1}(D_{2}+D_{4},D_{3})\nonumber \\
			= & \sigma_{2}(D_{1})+2\sigma_{2,1}(D_{1},D_{2}+D_{4})\nonumber \\
			& +\sigma_{2}(D_{2}+D_{4})+\sigma_{2}(D_{3})+2\sigma_{2,1}(D_{1},D_{3})+2\sigma_{2,1}(D_{2}+D_{4},D_{3}).
		\end{align}

		\vskip 4pt
		\emph{Case 1:} $n=3$. 
		\vskip 4pt
		
		It is enough to take $f(r)=e^r$.
		When $0<r<(1-\delta)\delta$ for any $0<\delta<1/2$, we know that $|\frac{v'}{v}|\le \frac{C}{r}$, $|\chi_1|\le \frac{C\delta}{r^2}$ and $|\chi_2|\le \frac{C\delta}{r^2}$ for some postive constant $C$.  Thus, $|D_1|\le \frac{C\delta}{r^2} $,  $|D_2|\le \frac{C\delta}{r}$,  $|D_3|\le C$,  $|D_4|\le \frac{C}{r} $. 
		
		Observe that
		\begin{align*}
			\sigma_{2}(D_{1}) & =(n-1)(\chi_{1}-\chi_{2})\chi_{1}+C_{n-1}^{2}\chi_{1}^{2}\\
			& =\chi_{1}(3\chi_{1}-2\chi_{2})\\
			& =\big(\frac{2\delta(1-\delta)}{r^2}-2+\frac{1}{r}(2-4\delta)\big)\big(-\frac{2\delta(1-\delta)}{r^{2}}+2+\frac{1}{r}(-6+4\delta)\big)\\
			&=-\frac{4}{r^4}\left[\delta^2(1-\delta)^2+4\delta(1-\delta)^2 r+(3-10\delta+6\delta^2)r^2-4(1-\delta)r^3+r^4\right]\\
			&=-\frac{4}{r^4}\left[\delta^2(1-\delta)^2-\delta^2 r^2+4(1-\delta)r(\delta(1-\delta) -r^2)+3(1-\delta)(1-2\delta)r^2\right]\\
			&\le -\frac{3\delta^2(1-\delta)^2}{r^4}.
		\end{align*}
		Then, for $0<r<r_1:=\min\{\frac{\delta(1-\delta)^2}{4C},(1-\delta)\delta\}$ we have
		\[\sigma_2(g^{-1}A_\delta)\leq-\frac{3\delta^2(1-\delta)^2}{r^4}+6C(\frac{\delta}{r^3}+\frac{1}{r^2}+\frac{\delta}{r^2}+\frac{1}{r})<0.\]
		
		\vskip 4pt
		\emph{Case 2:} $n=4$.	
		\vskip 4pt
		
		Notice that 
		\begin{align*}
			\sigma_{2}(D_{1}) & =(n-1)(\chi_{1}-\chi_{2})\chi_{1}+C_{n-1}^{2}\chi_{1}^{2}\\
			& =3\chi_{1}(2\chi_{1}-\chi_{2})\\
			& =3\big(-\frac{\delta(2-\delta)}{2r^{2}}+\frac{1}{2}(\frac{f'}{f})^{2}-\frac{f'}{rf}(1-\delta)\big)\big(-(\frac{f'}{f})^{2}+\frac{f''}{f}+\frac{f'}{rf}\big),
		\end{align*}
		where we have used the fact that
		\[
		2\chi_{1}-\chi_{2}=(\frac{f'}{f})^{2}-\frac{f''}{f}-\frac{f'}{rf}.
		\]
		
		Taking $f=e^{br}$ with $b \in \R_+$ to be determined later, we have
		\[
		-(\frac{f'}{f})^{2}+\frac{f''}{f}+\frac{f'}{rf}=\frac{b}{r}
		\]
		and for $0<r<\frac{\sqrt{\delta(2-\delta)}}{2b}:=r_{2}$,
		\begin{align*}
			&-\frac{\delta(2-\delta)}{2r^{2}}+\frac{1}{2}(\frac{f'}{f})^{2}-\frac{f'}{rf}(1-\delta)\\
			=&-\frac{\delta(2-\delta)}{2r^{2}}+\frac{b^2}{2}-\frac{b}{r}(1-\delta)<-\frac{\delta(2-\delta)}{4r^{2}}.
		\end{align*}
		This directly gives
		\begin{equation}\label{est:sigma_2_D_1}
		\sigma_{2}(D_{1})<-\frac{3b\delta(2-\delta)}{4r^{3}}<0.
		\end{equation}
		
		Taking $0<r<\min\{r_2,\frac{\delta(2-\delta)}{b(1-\delta)}\}:=r_{3}$,
		we obtain 
		\[
		\frac{\delta(2-\delta)}{4r^{2}}\le\chi_{1}\le\frac{3\delta(2-\delta)}{r^{2}}.
		\]
		
		Taking $0<r<\min\{r_2,\frac{\delta(2-\delta)}{b(3-2\delta)}\}:=r_{4}\leq r_3$,  we obtain
		\[
		\frac{\delta(2-\delta)}{2r^{2}}\le\chi_{2}\le\frac{3\delta(2-\delta)}{2 r^{2}}.
		\]
		
		Furthermore, taking $0<r<\min\{r_4,\frac{2-\delta}{2b}\}=r_5$, we obtain
		\begin{equation}\label{v'/v}
			-\frac{2-\delta}{r}\le \frac{v'}{v}\le -\frac{2-\delta}{2r},\quad | \frac{v'}{v}|\le \frac{2-\delta}{r}. 
		\end{equation}

		In the following, $O(f)$ means that  $O(f)\le C|f|$ for some $C \in \R_+$, which depends only on $g$ and second fundamental form  on $\pa M$, but is \emph{independent} of $b$ and $\delta$.
		
		We can obtain the following rough estimates:
		\begin{align}
			|(D_1)_{kj}|=& O(\frac{\delta}{r^{2}});\label{D1bound}\\
			|(D_2)_{kj}|=& |\chi_2|O(r);\label{D2bound}\\
			|(D_3)_{kj}|=& O(r^2|\frac{v'}{v}|^2)+O(1);\label{D3bound}\\
			|(D_4)_{kj}|=& O(|\frac{v'}{v}|).\label{D4bound}
		\end{align}

		By (\ref{eq:summation}), \eqref{D1bound}-\eqref{D4bound} and \eqref{est:sigma_2_D_1}  we have 
		\[
		\sigma_{2}(g^{-1}A_\delta)<-\frac{3b\delta(2-\delta)}{4r^{3}}+C_{1}\frac{\delta}{r^{3}}+C_{2}\frac{\delta}{r^{2}}+C_3\frac{1}{r^2}+C_4.
		\]
		Here $C_1, C_2, C_3, C_4$	are positive constants independent of $b$ and $\delta$.	
		
		Finally, we may choose $b=8(C_1+C_2+C_3+C_4)/3$, then for any $0<\delta<\frac{1}{2}$, there exists a positive constant $r_{1}=\min\{r_{5},\delta\}$ such that for all $0<r<r_1$,
		$$\sigma_{2}(g^{-1}A_\delta)<0.$$
		This means $\lambda(g^{-1}A_\delta)\in\mathbb{R}^{4}\backslash\overline{\Gamma_2^{+}}$.
		Moreover, by definition of $r_1$ we can find  two  positive constants $C_6, C_7$ depending only on $g$ and $L_{\alpha \beta}$ such that $C_6\delta \le r_1\le C_7 \sqrt{\delta}.$	
	\end{proof}
		
				\section{Appendix: Degenerate boundary $\sigma_2$-curvature equations}\label{Appendix:C}

		The purpose of this appendix is two-fold: One aims to give a delicate characterization of the limit solution to a degenerate fully nonlinear elliptic equation with Neumann boundary condition, which can be regarded as a rescaled limit equation of a  blow-up sequence of solutions to the boundary $\sigma_2$-curvature equation \eqref{PDE:sigma_2_new}. The other is to explain about an even reflection of a viscosity solution to the above degenerate fully nonlinear equation with zero Neumann boundary condition in a Euclidean unit ball. However, it is not our intention to give a general version of the convergence of viscosity solutions.  The following symbol $A_u$ is defined in \eqref{def:new_A_u}.
		\begin{proposition}
			Suppose $\Omega \subset \overline M$ is a smooth domain with boundary $\pa \Omega$ and $\partial_1\Omega \subset \pa \Omega \cap \pa M$ is a partial nonempty boundary of class $C^{3}$ and smooth metrics
			$g_{i}\rightarrow g$ in $C^{3}(\overline{\Omega})$. Let $v_{i}$
			be a smooth positive solution with $\lambda(A_{v_i})\in \Gamma_2^+$ to
			\[
			\begin{cases}
				\sigma_{2}(g_i^{-1}A_{v_i})=f_{i},~~ \lambda(g_i^{-1}A_{v_i})\in \Gamma_2^+ &\mathrm{~~in~~}\quad \Omega,\\
				\frac{\pa v_i}{\pa \vec n_{g_i}}=b_{i}v_i+q_{i} & \mathrm{~~on~~}\quad \partial_1\Omega.
			\end{cases}
			\]
			Suppose smooth functions $f_{i}\rightarrow0 $ in $C^0(\Omega\cup\partial_1\Omega)$  and $b_i\rightarrow b$, $q_i\rightarrow 0$ uniformly on $\partial_1\Omega$, and there exists a positive function $v$ such that 
			$v_{i}\rightarrow v$ in $C^0(\Omega\cup\partial_1\Omega)$. Then $v$
			is a viscosity solution to			\[
			\begin{cases}
				\sigma_{2}(A_{v})=0, ~~ \lambda(A_v)\in \pa \Gamma_2^+ & \mathrm{~~in~~}\quad\Omega,\\
				\frac{\partial v}{\partial \vec n_{g}}=bv & \mathrm{~~on~~}\quad\partial_1\Omega.
			\end{cases}
			\]
		\end{proposition}
		
		\begin{proof}
			
			\emph{Step 1:} $v$ is a viscosity supersolution. 
			\medskip
			
			Let $x_{0}\in\Omega\cup\partial_1\Omega,\varphi\in C^{2}\left(\Omega\right),(v-\varphi)\left(x_{0}\right)=0$,
			and $v-\varphi\geq0$ near $x_{0}$. 
			For simplicity, we just prove it when $x_{0}\in\partial_1\Omega$ and adopt Fermi coordinates with respect to the metric $g$ around $x_0$.  Fix a small $\delta>0$, let
			\[
			\varphi_{\delta}(x)=\varphi(x)-\delta\left|x\right|^{2}.
			\]
			Clearly, 
			\[
			\varphi_{\delta}(x)\leq\varphi(x)-\delta^{3} \qquad \mathrm{on~~} \pa^+ B_\delta^+.
			\]
			By the convergence of $v_{i}$ to $v$, for all sufficiently large $i$ we have
			\[
			v_{i}(x)\geq\varphi_{\delta}(x)+\frac{\delta^{3}}{2} \qquad \mathrm{on~~} \pa^+ B_\delta^+.
			\]
			Notice that $v_{i}\left(0\right)\rightarrow v\left(0\right)=\varphi_{\delta}\left(0\right)$,
			there exists some $\hat{x}_{i}\in B_{\delta/2}^+$ such that 
			\[
			\beta_{i}:=\left(v_{i}-\varphi_{\delta}\right)\left(\hat{x}_{i}\right)=\min_{\overline{B_\delta^+}}\left(v_{i}-\varphi_{\delta}\right)\rightarrow0\quad\mathrm{~~as~~}i\rightarrow\infty.
			\]
			Since 
			\[
			\beta_{i}=\left(v_{i}-\varphi\right)\left(\hat{x}_{i}\right)+\delta\left|\hat{x}_{i}\right|^{2}\geq\left(v_{i}-v\right)\left(\hat{x}_{i}\right)+\delta\left|\hat{x}_{i}\right|^{2},
			\]
			and $v_{i}\rightarrow v$, this forces $\hat{x}_{i}\rightarrow 0$.
			
			Let $\hat{\varphi}_{\delta}^{(i)}=\varphi_{\delta}+\beta_{i}$
			and 
			\[
			\begin{gathered}v_{i}\left(\hat{x}_{i}\right)=\hat{\varphi}_{\delta}^{(i)}\left(\hat{x}_{i}\right),\quad v_{i}(x)\geq\hat{\varphi}_{\delta}^{(i)}(x),\quad\forall\left|x-\hat{x}_{i}\right|<\frac{\delta}{2},\\
				\nabla v_{i}\left(\hat{x}_{i}\right)=\nabla\hat{\varphi}_{\delta}^{(i)}\left(\hat{x}_{i}\right),\quad\nabla^{2}v_{i}\left(\hat{x}_{i}\right)\geq\nabla^{2}\hat{\varphi}_{\delta}^{(i)}\left(\hat{x}_{i}\right)\quad\mathrm{~~if~~}\quad\hat{x}_{i}\in\Omega.
			\end{gathered}
			\]
			
			\emph{Case 1:} $\hat{x}_{i}\in\Omega$ for large $i.$ 
			
			At $\hat x_i$ we have 
			\[
			A_{v_{i}^{4/(n-2)}g_{i}}\le A_{(\hat{\varphi}_{\delta}^{(i)})^{4/(n-2)}g_{i}}
			\]
			and 
			\[
			\lambda(g_i^{-1}A_{(\hat{\varphi}_{\delta}^{(i)})^{4/(n-2)}g_{i}})\in\Gamma_{2}^{+}
			\]
			by assumption that $\lambda(g_i^{-1}A_{v_{i}^{4/(n-2)}g_{i}})\in\Gamma_{2}^{+}$. Then first taking $i\rightarrow\infty$, we obtain
			\[
			\lambda(A_{\varphi_{\delta}^{4/(n-2)}g})(0)\in \overline{\Gamma_2^+},
			\]
			and let $\delta\rightarrow0$ next to show
			\[
		\lambda(A_{\varphi^{4/(n-2)}g})(0)\in \overline{\Gamma_2^+} \quad \mathrm{or}\quad [-\frac{\partial\varphi}{\partial \vec n_g}+b\varphi](0)\ge0.
			\]
			
			\emph{Case 2:} $\hat{x}_{i}\in\partial_1\Omega$ for large $i$. 
			
			We have
			\[
			\frac{\partial}{\partial \vec n_g}(v_{i}-\hat{\varphi}_{\delta}^{(i)})\big|_{\hat{x}_{i}}\ge0
			\]
			and then as $i\rightarrow\infty$,
			\begin{align*}
				0\leftarrow q_i(\hat x_i)=[\frac{\partial v_{i}}{\partial \vec n_{g_{i}}}-b_iv_i](\hat{x}_{i})
				&\ge[\frac{\partial\hat{\varphi}_{\delta}^{(i)}}{\partial \vec n_{g_{i}}}-b_iv_i](\hat{x}_{i})\rightarrow[\frac{\partial\varphi}{\partial \vec n_g}-b\varphi](0).
			\end{align*}
			This implies that $[-\frac{\partial\varphi}{\partial \vec n_g}+b\varphi](0)\ge0$ and thus
			\[
		\lambda(A_{\varphi^{4/(n-2)}g})(0)\in \overline{\Gamma_2^+} \quad \mathrm{or}\quad [-\frac{\partial\varphi}{\partial \vec n_g}+b\varphi](0)\ge0.
		\]
			
			\emph{Step 2:} $v$ is a viscosity subsolution. 
			\medskip
			
			Let $x_{0}\in\Omega\cup\partial_1\Omega,\varphi\in C^{2}\left(\Omega\right),(v-\varphi)\left(x_{0}\right)=0$,
			and $v-\varphi\le0$ near $x_{0}$. We also simplify to prove it when $x_0 \in \pa_1 \Omega$ and use Fermi coordinates around $x_0$ as above. Fix a small $\delta>0$, consider 
			\[
			\varphi_{\delta}(x)=\varphi(x)+\delta\left|x\right|^{2}.
			\]
			Then 
			\[
			\varphi_{\delta}(x)\ge\varphi(x)+\delta^{3}\qquad \mathrm{on~~} \pa^+ B_\delta^+.
			\]
			By the convergence of $v_{i}$ to $v$, for large $i$ we have
			\[
			v_{i}(x)\le\varphi_{\delta}(x)-\frac{\delta^{3}}{2}\qquad \mathrm{on~~} \pa^+ B_\delta^+.
			\]
			Since $v_{i}\left(0\right)\rightarrow v\left(0\right)=\varphi_{\delta}\left(0\right)$,
			there exists some $\hat{x}_{i} \in B_{\delta/2}^+$ such that 
			\[
			\beta_{i}:=\left(v_{i}-\varphi_{\delta}\right)\left(\hat{x}_{i}\right)=\max_{\overline{B_\delta^+}}\left(v_{i}-\varphi_{\delta}\right)(x)\rightarrow0\quad\mathrm{~~as~~}i\rightarrow\infty.
			\]
			Since 
			\[
			\beta_{i}=\left(v_{i}-\varphi\right)\left(\hat{x}_{i}\right)-\delta\left|\hat{x}_{i}\right|^{2}\le \left(v_{i}-v\right)\left(\hat{x}_{i}\right)-\delta\left|\hat{x}_{i}\right|^{2},
			\]
			and $v_{i}\rightarrow v$, this forces $\hat{x}_{i}\rightarrow 0$.
			
			Let $\hat{\varphi}_{\delta}^{(i)}=\varphi_{\delta}+\beta_{i}$
			and 
			\[
			\begin{gathered}v_{i}\left(\hat{x}_{i}\right)=\hat{\varphi}_{\delta}^{(i)}\left(\hat{x}_{i}\right),\quad v_{i}(x)\le\hat{\varphi}_{\delta}^{(i)}(x),\quad\forall\left|x-\hat{x}_{i}\right|<\frac{\delta}{2},\\
				\nabla v_{i}\left(\hat{x}_{i}\right)=\nabla\hat{\varphi}_{\delta}^{(i)}\left(\hat{x}_{i}\right),\quad\nabla^{2}v_{i}\left(\hat{x}_{i}\right)\le\nabla^{2}\hat{\varphi}_{\delta}^{(i)}\left(\hat{x}_{i}\right)\quad\text{if}\quad\hat{x}_{i}\in\Omega
			\end{gathered}
			\]
			
			\emph{Case 1:} $\hat{x}_{i}\in\Omega$ for large $i.$ 
			
			At $\hat x_i$ we have 
			\[
			A_{v_{i}^{4/(n-2)}g_{i}}\ge A_{(\hat{\varphi}_{\delta}^{(i)})^{4/(n-2)}g_{i}},
			\]
		 and 	by assumption that $\sigma_2(A_{v_{i}^{4/(n-2)}g_{i}})=f_i\rightarrow0$ and
		 taking $i\rightarrow\infty$,
			\[
			\lambda(A_{\varphi_{\delta}^{4/(n-2)}g})(0)\in \mathbb{R}^n\backslash \Gamma_2^+.
			\]
			and let $\delta\rightarrow0$ next to show 
			\[
		\lambda(A_{\varphi})(0)\in \mathbb{R}^n\backslash \Gamma_2^+\quad \mathrm{or}\quad[-\frac{\partial\varphi}{\partial \vec n_g}+b\varphi](0)\le0.
			\]
			
			\emph{Case 2:} $\hat{x}_{i}\in\partial_1\Omega$ for large $i$.
			
			We have 
			\[
			\frac{\partial}{\partial \vec n_{g_i}}(v_{i}-\hat{\varphi}_{\delta}^{(i)})\big|_{\hat{x}_{i}}\le0
			\]
			and then for $i\rightarrow\infty$,
			\begin{align*}
				0\leftarrow q_i(\hat x_i)=[\frac{\partial v_{i}}{\partial \vec n_{g_{i}}}-b_iv_i](\hat{x}_{i})
				&\le[\frac{\partial\hat{\varphi}_{\delta}^{(i)}}{\partial \vec n_{g_{i}}}-b_iv_i](\hat{x}_{i})\rightarrow[\frac{\partial\varphi}{\partial \vec n_g}-b\varphi](0).
			\end{align*}
			This implies that $[-\frac{\partial\varphi}{\partial \vec n_g}+b\varphi](0)\le0$ and thus
				\[
		\lambda(A_{\varphi})(0)\in \mathbb{R}^n\backslash \Gamma_2^+\quad \mathrm{or}\quad[-\frac{\partial\varphi}{\partial \vec n_g}+b\varphi](0)\le0.
		\]
			
			Therefore, $v$ is a viscosity solution by definition.
		\end{proof}
		
		Let $\Omega^{+}\subset\mathbb{R}_{+}^{n}$
		be an open set and define 
		\[
		\partial^{\prime\prime}\Omega^{+}=\overline{\partial\Omega^{+}\cap\mathbb{R}_{+}^{n}},\quad\partial^{\prime}\Omega^{+}=\partial\Omega^{+}\backslash\partial^{\prime\prime}\Omega^{+}.
		\]
		To be consistent with our notations, we restate \cite[Lemma 5.3]{LiYY} as follows.

		\begin{lemma} \label{Lem:test_fcn_Li}
			
			Let $\Omega^{+} \subset \mathbb{R}_{+}^n$ be a bounded open set, and let $w \in C^2\left(\Omega^{+}\right) \cap$ $C^1\left(\Omega^{+} \cup \partial^{\prime} \Omega^{+}\right)$ satisfy
			$$
			w:=u^{-\frac{2}{n-2}} \geq c_1\quad  \mathrm{~~in~~}\quad  \Omega^{+}
			$$
			for some positive constant $c_1$, and let
			$$
			\varphi^{ \pm}(x):=e^{\delta|x|^2 \pm \delta^2 x_n} .
			$$
			Then there exist two positive constants $\delta$ depending only on $\sup \left\{|x|; x \in \Omega^{+}\right\}$ and $\bar{\ve}$ depending only on $\delta, c_1$, and $\sup \left\{|x|; x \in \Omega^{+}\right\}$ such that for any $0<\ve<\bar{\ve}$,
			\[
			\begin{array}{ll}
				(w+\ve \varphi^{ \pm})^{-\frac{4}{n-2}} A_{(w+\ve \varphi^{ \pm})^{-\frac{n-2}{2}}} \geq \left(1+\ve\varphi^{ \pm}u^{\frac{2}{n-2}}\right) u^{-\frac{4}{n-2}}A_{u}+\frac{\ve \delta}{2} \varphi^{ \pm} u^{-\frac{2}{n-2}} I & \mathrm{~~in~~} \Omega^{+}, \\
				(w-\ve \varphi^{ \pm})^{-\frac{4}{n-2}}A_{(w-\ve \varphi^{ \pm})^{-\frac{n-2}{2}}} \leq \left(1-\ve \varphi^{ \pm}u^{\frac{2}{n-2}}\right)u^{-\frac{4}{n-2}} A_u-\frac{\ve \delta}{2} \varphi^{ \pm} u^{-\frac{2}{n-2}} I & \mathrm{~~in~~} \Omega^{+}, \\
				\frac{\partial}{\partial x_n}\left(w+\ve \varphi^{ \pm}\right)=\frac{\partial w}{\partial x_n} \pm \ve \delta^2 \varphi^{\pm} & \mathrm{~~on~~} \partial^{\prime} \Omega^{+}, \\
				\frac{\partial}{\partial x_n}\left(w-\ve \varphi^{ \pm}\right)=\frac{\partial w}{\partial x_n} \mp \ve \delta^2 \varphi^{\pm} & \mathrm{~~on~~} \partial^{\prime} \Omega^{+} .
			\end{array}
			\]	
		\end{lemma}

		We follow the same lines in \cite[Proposition 5.1]{LiYY} to establish  Proposition \ref{even reflection in Appendix}.

		\noindent\begin{proof}[Proof of Proposition \ref{even reflection in Appendix}]
The proof is divided into three steps. 
			
			\emph{Step 1:} If $u^{\pm}$ is a \emph{strict} viscosity supersolution in $B_1^{\pm}$ as in Definition  \ref{vis-definition boundary}, then $u$ is a \emph{strict} viscosity supersolution in $B_1$.
			\medskip
			
			Since $u^{\pm}$ is a \emph{strict} viscosity supersolution in $B_1^{\pm}$, we know that  let $\phi\in C^{2}(\overline{B_{1}^{+}})$, $(u^{+}-\phi)(x_{0})=0$
			and $u^{+}\ge\phi$ near $x_0$. If $x_0 \in B_1^+$, then $\lambda(A_{\phi})(x_{0})\in \Gamma_2^+$. If $x_{0}\in\partial'B_{1}^{+}$,
			then \[\lambda(A_{\phi})(x_{0})\in \overline{\Gamma_2^+}\quad \mathrm{or}\quad -\frac{\partial\phi}{\partial x_n}(x_{0})>0.\]
			Similarly, let $\phi\in C^{2}(\overline{B_{1}^{-}})$, $(u^{-}-\phi)(x_{0})=0$
			and $u^{-}\ge\phi$ near $x_0$. If $x_0 \in B_1^-$, then $\lambda(A_{\phi})(x_{0})\in \Gamma_2^+$. If $x_{0}\in\partial'B_{1}^{-}$,
			then \[\lambda(A_{\phi})(x_{0})\in \overline{\Gamma_2^+}\quad \mathrm{or}\quad \frac{\partial\phi}{\partial x_n}(x_{0})>0.\]
			
			Let $\phi\in C^{2}(B_{1})$, $(u-\phi)(x_{0})=0$ and $u\ge\phi$ near $x_0$.
			
			If $x_0 \in B_1^{\pm}$, then we know that  $\lambda(A_{\phi})(x_{0})\in \overline{\Gamma_2^+}$. 
			
			If $x_{0}\in\partial'B_{1}^{+}$, then 
			$$\lambda(A_{\phi})(x_{0})\in \overline{\Gamma_2^+}\quad \mathrm{or}\quad -\frac{\partial\phi}{\partial x_n}(x_{0})>0$$
			and
			$$\lambda(A_{\phi})(x_{0})\in \overline{\Gamma_2^+}\quad \mathrm{or}\quad -\frac{\partial\phi}{\partial x_n}(x_{0})>0.$$
			This implies $\lambda(A_{\phi})(x_{0})\in \overline{\Gamma_2^+}$. Otherwise,  $-\frac{\partial\phi}{\partial x_{n}}(x_{0})>0$
			and $\frac{\partial\phi}{\partial x_{n}}(x_{0})>0$, which yields
			a contradiction. 
			
			Hence, $u$ is a \emph{strict} viscosity supersolution in $B_1$.
			
			\medskip
			\emph{Step 2:} The function 
		
			\[u_{\varepsilon}^{\pm}=\left((u^{\pm})^{-\frac{2}{n-2}}+\varepsilon\varphi^{\pm}\right)^{-\frac{n-2}{2}}\]
			is a \emph{strict} viscosity supersolution in $B_1^{\pm}$, respectively.
				\medskip
			
			To this end, for $u_\ve^+$, let $\varphi\in C^{2}(\overline{B_{1}^{+}})$ satisfy that $(u_{\varepsilon}^{+}-\varphi)(x_{0})=0$
			and $u_{\varepsilon}^{+}\ge\varphi$ near $x_0$.

			\emph{Case 1:} $x_{0}\in \partial B_{1}^{+}$. 
			
			Notice that 
			$$u^{+}=\left((u_\ve^{+})^{-\frac{2}{n-2}}-\varepsilon\varphi^{+}\right)^{-\frac{n-2}{2}}$$ and
			define
			$$\psi_\ve=\left(\varphi^{-\frac{2}{n-2}}-\varepsilon\varphi^{+}\right)^{-\frac{n-2}{2}},$$
			then $u^+-\psi_\ve \geq 0$ near $x_0$ and $(u^+-\psi_\ve)(x_0)=0$.
			Since $u^{+}$ is a viscosity supersolution, we have 
			\begin{equation}
				\lambda(A_{\psi_{\varepsilon}})(x_{0})\in \overline{\Gamma_2^+}\quad \mathrm{or}\quad -\frac{\partial\psi_{\varepsilon}}{\partial x_n}(x_{0})\ge0.\label{eq:1}
			\end{equation}
			
			At $x_0$, a direct computation yields 
			\begin{align*}
				-\frac{\partial\psi_\ve}{\partial x_{n}} =& -\frac{\partial}{\partial x_{n}}\left(\varphi^{-\frac{2}{n-2}}-\varepsilon\varphi^{+}\right)^{-\frac{n-2}{2}}\\
				=&-\frac{\pa \varphi}{\pa x_n}\left(1-\ve\varphi^+ \varphi^{\frac{2}{n-2}}\right)^{-\frac{2}{n-2}}-\varphi\frac{\pa }{\pa x_n}\left(1-\ve\varphi^+ \varphi^{\frac{2}{n-2}}\right)^{-\frac{2}{n-2}}\\
				=&-\frac{\partial\varphi}{\partial x_{n}}[1+O(\varepsilon)]-\frac{n-2}{2}\varepsilon\delta^{2}\varphi^+\varphi^{\frac{n}{n-2}}+O\left(\varepsilon^{2}\right)\\
				<&-\frac{\partial\varphi}{\partial x_{n}}[1+O(\varepsilon)].
			\end{align*}
			Also it follows from Lemma \ref{Lem:test_fcn_Li} that
			\begin{equation}\label{ineq:A_varphi}
				\psi_\ve^{-\frac{4}{n-2}}A_{\psi_\ve}<(1-\varepsilon\varphi^{+}\varphi^{\frac{2}{n-2}})\varphi^{-\frac{4}{n-2}}A_{\varphi}  \qquad \mathrm{~~in~~}\quad  B_1^+.
			\end{equation}
			These together with (\ref{eq:1}) yield for small $\ve>0$,
		
		\[\lambda(A_{\varphi})(x_{0})\in \overline{\Gamma_2^+}\quad \mathrm{or}\quad -\frac{\partial\varphi}{\partial x_n}(x_{0})>0.\]	
			This implies that $u_\ve^+$ is a \emph{strict} viscosity supersolution in $B_1^+$.
			
			\emph{Case 2:} $x_{0}\in B_{1}^{+}$, then it follows from \eqref{ineq:A_varphi} that $\lambda(A_{u_\ve^+})\in \overline{\Gamma_2^+}$.
			
			For $u_\ve^-$, repeating the above argument in place 
			of $u_\ve^+$ by 
			$$u_\ve^-=\left((u^-)^{-\frac{2}{n-2}}+\ve \varphi^-\right)^{-\frac{n-2}{2}}$$
			with minor modifications, we  can also see that $u_\ve^-$ is  a \emph{strict} viscosity supersolution  in $B_1^-$. Clearly, $u_\ve^+=u_\ve^-$ on $\pa' B_1^{+}$.
			
			\emph{Step 3:} Define
			\[
			u_{\varepsilon}\left(x^{\prime},x_{n}\right):=\begin{cases}
				u_{\varepsilon}^{+}\left(x^{\prime},x_{n}\right) & \text{ if }x_{n}\geq0,\\
				u_{\varepsilon}^{-}\left(x^{\prime},x_{n}\right) & \text{ if }x_{n}\leq0.
			\end{cases}
			\]
			
			Combining \emph{Step 1} and \emph{Step 2}, we conclude that $u_{\varepsilon}$ is a viscosity supersolution
			of $\lambda(A_{u})\in\partial\Gamma_2^+$ in $B_{1}$. Finally, letting $\varepsilon\rightarrow0$,
			by the classical convergence theory of viscosity solutions we know that $u$ is a viscosity
			supersolution of $\lambda(A_{u})\in\partial\Gamma_2^+$ in $B_{1}$.
		\end{proof}


\begin{thebibliography}{99}
			
		\bibitem{Barles} 
			G. Barles, \textit{Fully nonlinear Neumann type boundary conditions for second-order elliptic and parabolic equations}, J. Differential Equations 106 (1993), no. 1, 90-106.
			
			\bibitem{Branson-Gilkey} 
			T. Branson and P. Gilkey, \textit{The functional determinant of a four-dimensional boundary value problem}, Trans. Amer. Math. Soc. 344 (1994), no. 2, 479-531.
				
			\bibitem{Case-Wang}
			J. Case and Y. Wang, \textit{Boundary operators associated to the $\sigma_k$-curvature}, Adv. Math. 337 (2018), 83-106.
			
			
			\bibitem{CGY1}
			S.-Y. A. Chang, M. Gursky and P. Yang, \textit{An equation of Monge-Amp\`ere type in conformal geometry, and four-manifolds of positive Ricci curvature}, Ann. of Math. (2) 155 (2002), 709-787.
			
			
			
			\bibitem{CGY2} S.-Y. A. Chang, M. Gursky and P. Yang, \textit{An a
				priori estimate for a fully nonlinear equation on four-manifolds.
				Dedicated to the memory of Thomas H. Wolff}, J. Anal. Math. 87 (2002),
			151-186.
			
			\bibitem{CHY}
			S.-Y. A. Chang, Z. C. Han and P. Yang, \textit{Classification of singular radial solutions to the  $\sigma_k$  Yamabe equation on annular domains}, J. Differential Equations 216 (2005), no. 2, 482-501.
			
			
		
			\bibitem{CMW}
C.  Chen, X. N. Ma and W. Wei,
\textit{The Neumann problem of special Lagrangian equations with supercritical phase,}
J. Differential Equations 267 (2019), no.9, 5388-5409.
						
			\bibitem{CW}
		C.  Chen and W. Wei, \textit{The Neumann problem of complex Hessian quotient equations,}
Int. Math. Res. Not. IMRN (2021), no. 23, 17652-17672.



			\bibitem{Chen} S. Chen, \textit{Conformal deformation on manifolds
				with boundary}, Geom. Funct. Anal. 19 (2009), no. 4, 1029-1064.
			
			\bibitem{Chen2} S. Chen, \textit{Local estimates for some fully nonlinear
				elliptic equations}, Int. Math. Res. Not. 2005, no. 55, 3403-3425.
			
			\bibitem{Chen3} S. Chen,\textit{ Boundary value problems for some fully nonlinear elliptic equations},
			Calc. Var. Partial Differential Equations 30 (2007), no. 1, 1-15.
			
			\bibitem{Chen-Ruan-Sun}
			X. Chen, Y. Ruan and L. Sun, \textit{The Han-Li conjecture in constant scalar curvature and constant boundary mean curvature problem on compact manifolds}, Adv. Math. 358 (2019), 106854, 56 pp.
			
			\bibitem{Chen-Sun}
			X. Chen and L. Sun, \textit{Existence of conformal metrics with constant scalar curvature and constant boundary mean curvature on compact manifolds}, Commun. Contemp. Math. 21 (2019), no. 3, 1850021, 51 pp.
			
			
			\bibitem{Crandall-Ishii-Lions}
			M. G. Crandall, H. Ishii and P. L. Lions,\textit{ User's guide to viscosity solutions of second order partial differential equations}. Bull. Amer. Math. Soc. (N.S.) 27 (1992), no. 1, 1-67. 
			\bibitem{Duncan-Luc}
			J. Duncan and L. Nguyen, \textit{Differential inclusions for the Schouten tensor and nonlinear eigenvalue problems in conformal geometry}, \href{https://arxiv.org/abs/2208.00523}{arXiv:2208.00523}.
			
			\bibitem{DW}
			W. Dong and  W. Wei, 
			\textit{ The Neumann problem for a type of fully nonlinear complex equations}, J. Differential Equations 306 (2022), 525-546.
			
			\bibitem{escobar1} J. Escobar, \textit{Conformal deformation of a
				Riemannian metric to a scalar flat metric with constant mean curvature
				on the boundary}, Ann. of Math. (2) 136 (1992), no. 1, 1-50. (With
			an addendum: Ann. of Math. (2) 139 (1994), no. 3, 749-750. )
			
			
			\bibitem{escobar4} J. Escobar, \textit{The Yamabe problem on manifolds
				with boundary}, J. Differential Geom. 35 (1992), no. 1, 21-84.
				
				\bibitem{Ge-Lin-Wang}
			Y. Ge, C. S. Lin and G. Wang, \textit{On the $\sigma_2$-scalar curvature}, J. Differential Geom. 84 (2010), no. 1, 45-86.
				
			\bibitem{Ge-Wang1}
			Y. Ge and G. Wang, \textit{On a fully nonlinear Yamabe problem}, Ann. Sci. \'Ecole Norm. Sup. (4) 39 (2006), 569-598.
			
			
			\bibitem{Ge-Wang2}
			Y. Ge and G. Wang, \textit{A new conformal invariant on 3-dimensional manifolds}, Adv. Math. 249 (2013), 131-160.
			
			\bibitem{Guan-Wang} P. Guan and G. Wang, \textit{A fully nonlinear
				conformal flow on locally conformally flat manifolds}, J. Reine Angew.
			Math. 557 (2003), 219-238, \href{https://arxiv.org/abs/math/0112256}{arXiv:0112256}.
			
			\bibitem{Guan-Wang0}
			P. Guan and G. Wang, \textit{Local estimates for a class of fully nonlinear equations arising from conformal geometry}, Int. Math. Res. Not. (2003), no. 26,  1413-1432.
			
			
			
			\bibitem{Guan-Wang2}			
			P. Guan and G. Wang, \textit{ Geometric inequalities on locally conformally flat manifolds}, Duke Math. J. 124 (2004), no. 1, 177-212. 
			
			
			
			\bibitem{GVW} P. Guan, J. Viaclovsky and G. Wang, \textit{Some properties
				of the Schouten tensor and applications to conformal geometry}, Trans.
			Amer. Math. Soc. 355 (2003), no. 3, 925-933.
			
						
			
			
			\bibitem{Gursky-Viaclovsky} M. Gursky and J. Viaclovsky, \textit{Prescribing
				symmetric functions of the eigenvalues of the Ricci tensor}, Ann.
			of Math. (2) 166 (2007), no. 2, 475-531.
			
			\bibitem{GV0}
			M. Gursky and J. Viaclovsky, \textit{Volume comparison and the $\sigma_k$-Yamabe problem}, Adv. Math. 187 (2004), no.2, 447-487.
			
			
			
			\bibitem{han-li1} Z. C. Han and Y. Y. Li, \textit{The existence of
				conformal metrics with constant scalar curvature and constant boundary
				mean curvature}, Comm. Anal. Geom. 8 (2000), no. 4, 809-869.
			
			\bibitem{han-li2} Z. C. Han and Y. Y. Li, \textit{The Yamabe problem
				on manifolds with boundary: existence and compactness results}, Duke
			Math. J. 99 (1999), no. 3, 489-542.
			
			\bibitem{Han-Li-Tei}
			Z. C. Han, Y. Y. Li and E. Teixeira, \textit{Asymptotic behavior of solutions to the $\sigma_k$-Yamabe equation near isolated singularities}, Invent. Math. 182 (2010), no. 3, 635-684.
			
			\bibitem{HS} 
			Y. He and W. M. Sheng,\textit{
				On existence of the prescribing k-curvature problem on manifolds with boundary},
			Comm. Anal. Geom. 19 (2011), no. 1, 53-77.
			
			\bibitem{Ishii-Lions} 
			H. Ishii and P. L. Lions, \textit{Viscosity solutions of fully nonlinear second-order elliptic partial differential equations}, J. Differential Equations 83 (1990), no. 1, 26-78.
			
			\bibitem{JT}
			F. Jiang and N. S. Trudinger, \textit{Oblique boundary value problems for augmented Hessian equations I}, Bull. Math. Sci. 8 (2018), no. 2, 353-411.
			
			\bibitem{Jin}
			Q. Jin, 
			\textit{Local Hessian estimates for some conformally invariant fully nonlinear equations with boundary conditions,} Differential Integral Equations 20(2007), no.2, 121-132.
				
			\bibitem{Jin-Li-Li} Q. Jin, A. Li and Y. Y. Li, \textit{Estimates and existence results for a fully nonlinear Yamabe problem on manifolds with boundary}, Calc. Var. Partial Differential Equations 28 (2007),
			no. 4, 509-543.

			\bibitem{Lee-Parker}
			J. Lee and T. Parker, \textit{The Yamabe problem}, Bull. Amer. Math. Soc. (N.S.) 17 (1987), no. 1, 37-91. 
			
			\bibitem{Li-Li1} A. Li and Y. Y. Li,
			\textit{On some conformally
				invariant fully nonlinear equations}, Comm. Pure Appl. Math. 56 (2003),
			no. 10, 1416-1464.
			
			\bibitem{Li-Li2}
			A. Li and Y. Y. Li, \textit{ On some conformally invariant fully nonlinear equations. II. Liouville, Harnack and Yamabe}, Acta Math. 195 (2005), 117-154. 
			
			\bibitem{Li-Li3}				
			A. Li and Y. Y. Li, \textit{
				A fully nonlinear version of the Yamabe problem on manifolds with boundary},
			J. Eur. Math. Soc. (JEMS) 8 (2006), no. 2, 295-316.			
			
			\bibitem{LiYY0}
			Y. Y. Li, \textit{ Conformally invariant fully nonlinear elliptic equations and isolated singularities}, J. Funct. Anal. 233 (2006), no. 2, 380-425.
			
			\bibitem{LiYY}
			Y. Y. Li, \textit{Local gradient estimates of solutions to some conformally invariant fully nonlinear equations}, Comm. Pure Appl. Math. 62 (2009), no. 10, 1293-1326. 
			
			\bibitem{Li-Luc2} Y. Y. Li and L. Nguyen, \textit{A compactness theorem
				for a fully nonlinear Yamabe problem under a lower Ricci curvature
				bound}, J. Funct. Anal. 266 (2014), no. 6, 3741-3771.
				
				\bibitem{Li-Luc3}
				Y. Y. Li and L. Nguyen, \textit{A fully nonlinear version of the Yamabe problem on locally conformally flat manifolds with umbilic boundary}, Adv. Math. 251 (2014), 87-110, \href{https://arxiv.org/abs/0911.3366}{arXiv:0911.3366}.
			
			\bibitem{Li-Luc4}
			Y. Y. Li and L. Nguyen, \textit{Counterexamples to $C^2$ boundary estimates for a fully nonlinear Yamabe problem on manifolds with boundary}, Adv. Nonlinear Stud. 12 (2012), 783-797. 
			
			\bibitem{Li-Luc1} Y. Y. Li and L. Nguyen, \textit{Harnack inequalities
				and Bocher-type theorems for conformally invariant fully nonlinear degenerate elliptic equations}, Comm. Pure Appl. Math. 67 (2014), no. 11, 1843-1876.
			
			\bibitem{LTU}	
			P. L. Lions, N. S. Trudinger and J. I. Urbas, \textit{The Neumann problem for equations of Monge-Amp\`ere type},
			Comm. Pure Appl. Math. 39 (1986), no. 4, 539-563.
			
			\bibitem{MQ}
			X. N. Ma and G. H. Qiu, 
			\textit{The Neumann problem for Hessian equations},
			Comm. Math. Phys. 366 (2019), no. 1, 1-28.
			
			\bibitem{Perales} R. Perales, \textit{Volumes and limits of manifolds
				with Ricci curvature and mean curvature bounds}, Differential Geom.
			Appl. 48 (2016), 23-37.
			
			\bibitem{Reilly}
R. Reilly, \textit{On the Hessian of a function and the curvatures of its graph},
Michigan Math. J. 20 (1973), 373-383.

			
			\bibitem{Schoen}
			R. Schoen, \textit{Variational theory for the total scalar curvature functional for Riemannian metrics and related topics}, Topics in calculus of variations (Montecatini Terme, 1987), 120-154, Lecture Notes in Math., 1365, Springer, Berlin, 1989.
			
			\bibitem{Schoen-Zhang}
			R. Schoen and D. Zhang, \textit{Prescribed scalar curvature on the $n$-sphere},
			Calc. Var. Partial Differential Equations 4 (1996), no. 1, 1-25. 
			
			\bibitem{STW}
			W. M. Sheng, N.S. Trudinger and X. J. Wang, \textit{The Yamabe problem for higher order curvatures}, J. Differential Geom. 77 (2007), no.3, 515-553. 
			
			\bibitem{TW}
			N. S. Trudinger and X. J. Wang, \textit{The intermediate case of the Yamabe problem for higher order curvatures}, Int. Math. Res. Not. IMRN 2010 (2010), no.13, 2437-2458. 	
			\bibitem{Via1}
			J. Viaclovsky, \textit{Conformal geometry, contact geometry, and the calculus of variations}, Duke Math. J. 101 (2000), no.2, 283-316.
			
			
			
			\bibitem{Via2}
			J. Viaclovsky, \textit{Estimates and existence results for some fully nonlinear elliptic equations on Riemannian manifolds}, Comm. Anal. Geom. 10 (2002), 815-846. 
			
			\bibitem{Wang}
			X. J. Wang,\textit{ A priori estimates and existence for a class of fully nonlinear elliptic equations in conformal geometry}, Chinese Ann. Math. Ser. B 27 (2006), no.2, 169-178. 
			
		\end{thebibliography}
	\end{document}